\documentclass[10pt]{amsart}

\usepackage{latexsym,exscale,enumerate}
\usepackage{amsmath,amsthm,amsfonts,amssymb,amscd, stmaryrd}
\usepackage[bookmarks=false]{hyperref}
\usepackage{ulem}

\usepackage[all]{xy}
\SelectTips{cm}{}

\usepackage{graphicx}

\usepackage{tikz}
\usetikzlibrary{decorations.markings}
\usetikzlibrary{decorations.pathreplacing}
\usetikzlibrary{arrows,shapes,positioning}
\tikzstyle directed=[postaction={decorate,decoration={markings,
    mark=at position #1 with {\arrow{>}}}}]
\tikzstyle rdirected=[postaction={decorate,decoration={markings,
    mark=at position #1 with {\arrow{<}}}}]

\newcommand{\ep}{\underline{\epsilon}}
\newcommand{\onen}{{\mathbf 1}_{n}}
\newcommand{\onenn}[1]{{\mathbf 1}_{#1}}
\newcommand{\onenp}{{\mathbf 1}_{n'}}
\newcommand{\onenpp}{{\mathbf 1}_{n''}}
\newcommand{\onem}{{\mathbf 1}_{m}}
\newcommand{\onel}{{\mathbf 1}_{\lambda}}
\newcommand{\onell}[1]{{\mathbf 1}_{#1}}
\newcommand{\onea}{{\mathbf 1}_{{\bf a}}}
\newcommand{\oneaa}[1]{{\mathbf 1}_{#1}}
\newcommand{\onelp}{{\mathbf 1}_{\lambda'}}
\newcommand\rE{{\sf{E}}}
\newcommand\rF{{\sf{F}}}
\newcommand\sE{{\cal{E}}}
\newcommand\sF{{\cal{F}}}
\def\cal#1{\mathcal{#1}}%

\newcommand{\bV}{\textstyle{\bigwedge_q}}
\newcommand{\Alt}[2]{\textstyle{\bigwedge_{#1}^{#2}}}
\newcommand{\Alts}[2]{ \textstyle{\scriptstyle\bigwedge_{#1}^{#2}}}

\def\1{\mathbbm{1}}%
\newcommand\E{{\sf{E}}}
\newcommand\F{{\sf{F}}}
\newcommand\T{{\sf{T}}}
\def\l{\lambda}

\newcommand{\END}{{\rm END}}
\newcommand{\Gr}{\cat{Flag}_{N}}
\newcommand{\Grn}[1]{\cat{Flag}_{#1}}

\newcommand{\Uq}{{\bf U}_q(\mathfrak{sl}_2)}
\newcommand{\U}{\dot{{\bf U}}}
\newcommand{\Ucat}{\cal{U}}
\newcommand{\Ucatc}{\check{\cal{U}}}
\newcommand{\Ucatq}{\cal{U}_q}
\newcommand{\UcatD}{\dot{\cal{U}}}
\newcommand{\UcatDq}{\dot{\cal{U}}_q}
\newcommand{\B}{\dot{\mathbb{B}}}
\newcommand{\Bnm}{{_m\dot{\cal{B}}_n}}
\newcommand{\UA}{{_{\cal{A}}\dot{{\bf U}}}}
\newcommand{\sln}{\mf{sl}_n}
\newcommand{\slm}{\mf{sl}_m}
\newcommand{\slnn}[1]{\mf{sl}_{#1}}
\newcommand{\gln}{\mf{gl}_n}
\newcommand{\glm}{\mf{gl}_m}
\newcommand{\glnn}[1]{\mf{gl}_{#1}}

\newcommand{\und}[1]{\underline{#1}}

\newcommand{\foam}[3][N]{#2\cat{Foam}_{#3}(#1)}
\newcommand{\Bfoam}[3][N]{#2\cat{BFoam}_{#3}(#1)}
\newcommand{\nWeb}{\cat{nWeb}}

\newcommand{\qbin}[2]{
\left[
 \begin{array}{c}
 #1 \\
 #2 \\
 \end{array}
 \right]_{q^2}
}
\newcommand{\qbins}[2]{
\left[
 \begin{array}{c}
 \scs #1 \\
 \scs #2 \\
 \end{array}
 \right]
}

\newcommand{\xsum}[2]{
  \xy
  (0,.4)*{\sum};
  (0,3.7)*{\scs #2};
  (0,-2.9)*{\scs #1};
  \endxy
}

\newcommand{\refequal}[1]{\xy {\ar@{=}^{#1}
(-1,0)*{};(1,0)*{}};
\endxy}

\newcommand{\cat}[1]{\ensuremath{\mbox{\bfseries {\upshape {#1}}}}}
\newcommand{\numroman}{\renewcommand{\labelenumi}{\roman{enumi})}}
\newcommand{\numarabic}{\renewcommand{\labelenumi}{\arabic{enumi})}}
\newcommand{\numAlph}{\renewcommand{\labelenumi}{\Alph{enumi}.}}

\newcommand{\To}{\Rightarrow}
\newcommand{\TO}{\Rrightarrow}
\newcommand{\Hom}{{\rm Hom}}
\newcommand{\HOM}{{\rm HOM}}
\renewcommand{\to}{\rightarrow}
\newcommand{\maps}{\colon}
\newcommand{\op}{{\rm op}}
\newcommand{\co}{{\rm co}}
\newcommand{\iso}{\cong}
\newcommand{\id}{{\rm id}}
\newcommand{\bigb}[1]{
\begin{pspicture}(0,0)
 \rput(0,0){\psframebox[framearc=.5,fillstyle=solid]{\small $#1$}}
\end{pspicture}}
\newcommand{\del}{\partial}
\newcommand{\Res}{{\rm Res}}
\newcommand{\End}{{\rm End}}
\newcommand{\Aut}{{\rm Aut}}
\newcommand{\im}{{\rm im\ }}
\newcommand{\coim}{{\rm coim\ }}
\newcommand{\chr}{{\rm char\ }}
\newcommand{\coker}{{\rm coker\ }}
\newcommand{\spann}{{\rm span}}
\newcommand{\rk}{{\rm rk\ }}
\def\bigboxtimes{\mathop{\boxtimes}\limits}

\newcommand{\scs}{\scriptstyle}

\def\Res{{\mathrm{Res}}}
\def\Ind{{\mathrm{Ind}}}
\def\lra{{\longrightarrow}}
\def\dmod{{\mathrm{-mod}}}   
\def\fmod{{\mathrm{-fmod}}}   
\def\pmod{{\mathrm{-pmod}}}  
\def\rk{{\mathrm{rk}}}
\def\NH{{\mathrm{NH}}}
\def\pseq{{\mathrm{Seqd}}}
\def\Id{\mathrm{Id}}
\def\mc{\mathcal}
\def\mf{\mathfrak}
\def\Af{{_{\mc{A}}\mathbf{f}}}    
\def\primef{{'\mathbf{f}}}    
\def\shuffle{\,\raise 1pt\hbox{$\scriptscriptstyle\cup{\mskip
               -4mu}\cup$}\,}
\newcommand{\define}{\stackrel{\mbox{\scriptsize{def}}}{=}}

\theoremstyle{definition}
\newtheorem{thm}{Theorem}[section]
\newtheorem{cor}[thm]{Corollary}

\newtheorem{lem}[thm]{Lemma}
\newtheorem{rem}[thm]{Remark}
\newtheorem{prop}[thm]{Proposition}
\newtheorem{defn}[thm]{Definition}

\numberwithin{equation}{section}

\def\AL#1{\textcolor[rgb]{1.00,0.00,0.00}{[AL: #1]}}%
\def\HQ#1{\textcolor[rgb]{0.00,1.00,0.00}{[HQ: #1]}}%
\def\DR#1{\textcolor[rgb]{0.00,0.00,1.00}{[DR: #1]}}%
\def\MK#1{[MK: #1]}%
\def\AB#1{[AB: #1]}%
\def\b{$\blacktriangleright$}
\def\e{$\blacktriangleleft$}
\def\new#1{\b #1\e}%

\def\emph#1{{\sl #1\/}}
\def\ie{{\sl i.e.\/}}
\def\eg{{\sl e.g.\/}}
\def\Eg{{\sl E.g.\/}}
\def\etc{{\sl etc.\/}}
\def\cf{{\sl c.f.\/}}
\def\etal{\sl{et al.\/}}%
\def\vs{\sl{vs.\/}}%

\let\hat=\widehat
\let\tilde=\widetilde

\let\phi=\varphi
\let\theta=\vartheta
\let\epsilon=\varepsilon

\usepackage{bbm}
\def\C{{\mathbbm C}}
\def\N{{\mathbbm N}}
\def\R{{\mathbbm R}}
\def\Z{{\mathbbm Z}}
\def\Q{{\mathbbm Q}}
\def\H{{\mathbbm H}}
\def\P{{\mathbbm P}}

\newcommand{\Bb}{\mathbb{B}}

\def\cal#1{\mathcal{#1}}%
\def\1{\mathbbm{1}}%
\def\ev{\mathrm{ev}}%
\def\coev{\mathrm{coev}}%
\def\tr{\mathrm{tr}}%
\def\st{\mathrm{st}}%
\def\pullback#1#2#3{%
  \,\mbox{\raisebox{-.8ex}{$\scriptstyle #1$}}%
  \!\prod\!
  \mbox{\raisebox{-.8ex}{$\scriptstyle #3$}}\,}%
\def\nn{\notag}
\newcommand{\ontop}[2]{\genfrac{}{}{0pt}{2}{\scriptstyle #1}{\scriptstyle #2}}

\def\la{\langle}
\def\ra{\rangle}

\newcommand{\sdotu}[1]{\xybox{%
  (-2,0)*{};
  (2,0)*{};
  (0,0)*{}; (0,-9)*{} **\dir{-}?(.5)*{\bullet} ?(0)*\dir{<};
  (0,-11)*{\scs  #1};
}}
\newcommand{\urcurve}[1]{\xybox{
(-8,0)*{};
  (8,0)*{};
  (0,0);(0,18) **\crv{(8,9)}?(1)*\dir{>};
  (0,-2)*{ \scs #1}
}}
\newcommand{\ulcurve}[1]{\xybox{
(-8,0)*{};
  (8,)*{};
  (0,0);(0,18) **\crv{(-8,9)}?(1)*\dir{>};
  (0,-2)*{ \scs #1}
}}
\newcommand{\slineu}[1]{\xybox{%
  (-2,0)*{};
  (2,0)*{};
  (0,0)*{}; (0,-9)*{} **\dir{-}; ?(0)*\dir{<};
  (0,-11)*{\scs  #1};
}}

\newcommand{\drcurve}[1]{\xybox{
(-8,0)*{};
  (8,0)*{};
  (0,0);(0,18) **\crv{(10,9)}?(.5)*\dir{<}+(2.5,0)*{ \scs #1};
}}

\newcommand{\dlcurve}[1]{\xybox{
(-8,0)*{};
  (8,)*{};
  (0,0);(0,18) **\crv{(-10,9)}?(.5)*\dir{<}+(-2.5,0)*{ \scs #1};
}}
\newcommand{\lineu}[1]{\xybox{%
  (-2,0)*{};
  (2,0)*{};
  (0,0)*{}; (0,-18)*{} **\dir{-}; ?(.5)*\dir{<}+(1.7,-7)*{\scs #1};
}}
\newcommand{\lined}[1]{\xybox{%
  (-2,0)*{};
  (2,0)*{};
  (0,0)*{}; (0,-18)*{} **\dir{-}; ?(.5)*\dir{>}+(1.7,-7)*{\scs #1};
}}
\newcommand{\dotu}[1]{\xybox{%
  (-2,0)*{};
  (2,0)*{};
  (0,0)*{}; (0,-18)*{} **\dir{-}
  ?(.2)*{\bullet} ?(.5)*\dir{<}+(1.7,-7)*{\scs #1};
}}

\newcommand{\twod}{\xybox{%
  (-6,0)*{};
  (6,0)*{};
  (0,6)*{}="f";
  (-3,0)*{}="t1";
  (3,0)*{}="t2";
  (0,12)*{}="b";
  "t1";"f" **\crv{(-3,4)};
  "t2";"f" **\crv{(3,4)};
  "f"+(.5,0);"b"+(.5,0) **\dir{-};
  "f"+(-.5,0);"b"+(-.5,0) **\dir{-};
}}

\newcommand{\twom}{\xybox{%
  (-6,0)*{};
  (6,0)*{};
  (0,6)*{}="f";
  (-3,12)*{}="t1";
  (3,12)*{}="t2";
  (0,0)*{}="b";
  "t1";"f" **\crv{(-3,8)};
  "t2";"f" **\crv{(3,8)};
  "f"+(.5,0);"b"+(.5,0) **\dir{-};
  "f"+(-.5,0);"b"+(-.5,0) **\dir{-};
}}
\newcommand{\twoI}{\xybox{%
  (-6,0)*{};
  (6,0)*{};
  (0,6)*{}="f'";
  (0,12)*{}="f";
  (-3,18)*{}="t1";
  (3,18)*{}="t2";
    (-3,0)*{}="t1'";
  (3,0)*{}="t2'";
  (0,0)*{}="b";
  "t1";"f" **\crv{(-3,14)};
  "t2";"f" **\crv{(3,14)};
  "f"+(.5,0);"f'"+(.5,0) **\dir{-};
  "f"+(-.5,0);"f'"+(-.5,0) **\dir{-};
  "t1'";"f'" **\crv{(-3,4)};
  "t2'";"f'" **\crv{(3,4)};
}}
\newcommand{\twoIu}[1]{\xybox{%
  (-6,0)*{};
  (6,0)*{};
  (0,6)*{}="f'";
  (0,12)*{}="f";
  (-3,18)*{}="t1";
  (3,18)*{}="t2";
    (-3,0)*{}="t1'";
  (3,0)*{}="t2'";
  (0,0)*{}="b";
  "t1";"f" **\crv{(-3,14)};?(.1)*\dir{<};
  "t2";"f" **\crv{(3,14)};?(.1)*\dir{<};
  "f"+(.5,0);"f'"+(.5,0) **\dir{-};
  "f"+(-.5,0);"f'"+(-.5,0) **\dir{-};
  "t1'";"f'" **\crv{(-3,4)};
  "t2'";"f'" **\crv{(3,4)} ?(.15)*\dir{ }+(2,0)*{\scs #1};
}}

\newcommand{\lowrru}[1]{\xybox{%
  (-8,0)*{};
  (8,0)*{};
  (-6,-18)*{};(6,-9)*{} **\crv{(-6,-13) & (6,-15)} ?(1)*\dir{>};
  (6,-9)*{};(6,0)*{}  **\dir{-} ?(.3)*\dir{ }+(2,0)*{\scs {\bf j}};
}}

\newcommand{\lowllu}[1]{\xybox{%
  (-8,0)*{};
  (8,0)*{};
  (6,-18)*{};(-6,-9)*{} **\crv{(6,-13) & (-6,-15)} ?(1)*\dir{>};
  (-6,-9)*{};(-6,0)*{}  **\dir{-} ?(.3)*\dir{ }+(-2,0)*{\scs {\bf j}};
}}
\newcommand{\highrru}[1]{\xybox{%
  (-8,0)*{};
  (8,0)*{};
  (6,0)*{};(-6,-9)*{} **\crv{(6,-5) & (-6,-3)} ?(1)*\dir{<};
  (-6,-9)*{};(-6,-18)*{}  **\dir{-} ?(.3)*\dir{ }+(-2,0)*{\scs #1};
}}

\newcommand{\highllu}[1]{\xybox{%
  (-8,0)*{};
  (8,0)*{};
  (-6,0)*{};(6,-9)*{} **\crv{(-6,-5) & (6,-3)} ?(1)*\dir{<};
  (6,-9)*{};(6,-18)*{}  **\dir{-} ?(.3)*\dir{ }+(2,0)*{\scs #1};
}}

\newcommand{\bbmf}[3]{\xybox{%
  (-6,0)*{};
  (6,0)*{};
  (0,-8)*{}="f";
  (-4,0)*{}="t1";
  (4,0)*{}="t2";
  (0,-16)*{}="b";
  "t1";"f" **\crv{(-4,-4)}; ?(.35)*\dir{>}+(-2.3,0)*{\scriptstyle{#1}};
  "t2";"f" **\crv{(4,-4)}; ?(.35)*\dir{>}+(2.3,0)*{\scriptstyle{#2}};
  "f";"b" **\dir{-}; ?(.75)*\dir{>}+(3.8,0)*{\scriptstyle{#3}};;
}}

\newcommand{\bbme}[3]{\xybox{%
  (-6,0)*{};
  (6,0)*{};
  (0,-8)*{}="f";
  (-4,0)*{}="t1";
  (4,0)*{}="t2";
  (0,-16)*{}="b";
  "t1";"f" **\crv{(-4,-4)}; ?(.35)*\dir{<}+(-2.3,0)*{\scriptstyle{#1}};
  "t2";"f" **\crv{(4,-4)}; ?(.35)*\dir{<}+(2.3,0)*{\scriptstyle{#2}};
  "f";"b" **\dir{-}; ?(.6)*\dir{<}+(3.8,0)*{\scriptstyle{#3}};
}}

\newcommand{\bbdf}[3]{\xybox{%
  (-6,0)*{};
  (6,0)*{};
  (0,-8)*{}="f";
  (-4,-16)*{}="b1";
  (4,-16)*{}="b2";
  (0,-0)*{}="t";
  "f";"b1" **\crv{(-4,-12)}; ?(.75)*\dir{>}+(-2.3,0)*{\scriptstyle{#1}};
  "f";"b2" **\crv{(4,-12)}; ?(.75)*\dir{>}+(2.3,0)*{\scriptstyle{#2}};
  "t";"f" **\dir{-}; ?(.35)*\dir{>}+(3.8,0)*{\scriptstyle{#3}};
}}

\newcommand{\bbde}[3]{\xybox{%
  (-6,0)*{};
  (6,0)*{};
  (0,-8)*{}="f";
  (-4,-16)*{}="b1";
  (4,-16)*{}="b2";
  (0,-0)*{}="t";
  "f";"b1" **\crv{(-4,-12)}; ?(.6)*\dir{<}+(-2.3,0)*{\scriptstyle{#1}};
  "f";"b2" **\crv{(4,-12)}; ?(.6)*\dir{<}+(2.3,0)*{\scriptstyle{#2}};
  "t";"f" **\dir{-}; ?(.3)*\dir{<}+(3.8,0)*{\scriptstyle{#3}};
}}

\newcommand{\ccbub}[1]{
\xybox{%
 (-6,0)*{};
  (6,0)*{};
  (-4,0)*{}="t1";
  (4,0)*{}="t2";
  "t2";"t1" **\crv{(4,6) & (-4,6)};
  ?(.05)*\dir{>} ?(1)*\dir{>};
  "t2";"t1" **\crv{(4,-6) & (-4,-6)};
   ?(.25)*\dir{}+(0,0)*{\bullet}+(0,-3)*{\scs {#1}};
}}

\newcommand{\iccbub}[2]{
\xybox{%
 (-6,0)*{};
  (6,0)*{};
  (-4,0)*{}="t1";
  (4,0)*{}="t2";
  "t2";"t1" **\crv{(4,6) & (-4,6)}; ?(.7)*\dir{}+(-2,0)*{\scs #2}
  ?(.05)*\dir{>} ?(1)*\dir{>};
  "t2";"t1" **\crv{(4,-6) & (-4,-6)};
   ?(.3)*\dir{}+(0,0)*{\bullet}+(0,-3)*{\scs {#1}};
}}
\newcommand{\icbub}[2]{
\xybox{%
 (-6,0)*{};
  (6,0)*{};
  (-4,0)*{}="t1";
  (4,0)*{}="t2";
  "t2";"t1" **\crv{(4,6) & (-4,6)};?(.7)*\dir{}+(-2,0)*{\scs #2};
   ?(0)*\dir{<} ?(.95)*\dir{<};
  "t2";"t1" **\crv{(4,-6) & (-4,-6)};
   ?(.3)*\dir{}+(0,0)*{\bullet}+(0,-3)*{\scs {#1}};
}}

\newcommand{\cbub}[1]{
\xybox{%
 (-6,0)*{};
  (6,0)*{};
  (-4,0)*{}="t1";
  (4,0)*{}="t2";
  "t2";"t1" **\crv{(4,6) & (-4,6)};
   ?(0)*\dir{<} ?(.95)*\dir{<};
  "t2";"t1" **\crv{(4,-6) & (-4,-6)};
   ?(.75)*\dir{}+(0,0)*{\bullet}+(0,-3)*{\scs {#1}};
}}
\newcommand{\bbe}[1]{\xybox{%
  (-2,0)*{};
  (2,0)*{};
  (0,0);(0,-18) **\dir{-}; ?(.5)*\dir{<}+(2.3,0)*{\scriptstyle{#1}};
}}

\newcommand{\bbf}[1]{\xybox{%
  (-2,0)*{};
  (2,0)*{};
  (0,0);(0,-18) **\dir{-}; ?(.5)*\dir{>}+(2.3,0)*{\scriptstyle{#1}};
}}

\newcommand{\bbelong}[1]{\xybox{%
  (-2,0)*{};
  (2,0)*{};
  (0,0);(0,-22) **\dir{-}; ?(.5)*\dir{<}+(2.3,0)*{\scriptstyle{#1}};
}}

\newcommand{\bbflong}[1]{\xybox{%
  (-2,0)*{};
  (2,0)*{};
  (0,0);(0,-22) **\dir{-}; ?(.5)*\dir{>}+(2.3,0)*{\scriptstyle{#1}};
}}

\newcommand{\bbsid}{\xybox{%
  (-2,0)*{};
  (2,0)*{};
  (0,10);(0,4) **\dir{-};
}}
\newcommand{\bbpef}[1]{\xybox{%
  (-6,0)*{};
  (6,0)*{};
  (-4,0)*{}="t1";
  (4,0)*{}="t2";
  "t1";"t2" **\crv{(-4,-6) & (4,-6)}; ?(.15)*\dir{>} ?(.9)*\dir{>}
    ?(.5)*\dir{}+(0,-2)*{\scriptstyle{#1}};
}}
\newcommand{\bbpfe}[1]{\xybox{%
  (-6,0)*{};
  (6,0)*{};
  (-4,0)*{}="t1";
  (4,0)*{}="t2";
  "t2";"t1" **\crv{(4,-6) & (-4,-6)}; ?(.15)*\dir{>} ?(.9)*\dir{>} ?(.5)*\dir{}+(0,-2)*{\scriptstyle{#1}};
}}

\newcommand{\bbcfe}[1]{\xybox{%
  (-6,0)*{};
  (6,0)*{};
  (-4,0)*{}="t1";
  (4,0)*{}="t2";
  "t1";"t2" **\crv{(-4,6) & (4,6)}; ?(.15)*\dir{>} ?(.9)*\dir{>}
  ?(.5)*\dir{}+(0,2)*{\scriptstyle{#1}};
}}
\newcommand{\bbcef}[1]{\xybox{%
  (-6,0)*{};
  (6,0)*{};
  (-4,0)*{}="t1";
  (4,0)*{}="t2";
  "t2";"t1" **\crv{(4,6) & (-4,6)}; ?(.15)*\dir{>}
  ?(.9)*\dir{>} ?(.5)*\dir{}+(0,2)*{\scriptstyle{#1}};
}}
\newcommand{\lbbcef}[1]{\xybox{%
  (-8,0)*{};
  (8,0)*{};
  (-4,0)*{}="t1";
  (4,0)*{}="t2";
  "t2";"t1" **\crv{(4,6) & (-4,6)}; ?(.15)*\dir{>}
  ?(.9)*\dir{>} ?(.5)*\dir{}+(0,2)*{\scriptstyle{#1}};
}}

\newcommand{\sccbub}[1]{%
\xybox{%
 (-6,0)*{};
  (6,0)*{};
  (-4,0)*{}="t1";
  (4,0)*{}="t2";
  "t2";"t1" **\crv{(4,6) & (-4,6)}; ?(.05)*\dir{>} ?(1)*\dir{>};
  "t2";"t1" **\crv{(4,-6) & (-4,-6)}; ?(.3)*\dir{}+(2,-1)*{\scs #1};
}}
\newcommand{\scbub}[1]{%
\xybox{%
 (-6,0)*{};
  (6,0)*{};
  (-4,0)*{}="t1";
  (4,0)*{}="t2";
  "t2";"t1" **\crv{(4,6) & (-4,6)}; ?(0)*\dir{<} ?(.95)*\dir{<};
  "t2";"t1" **\crv{(4,-6) & (-4,-6)}; ?(.3)*\dir{}+(2,-1)*{\scs #1};
}}
\newcommand{\mcbub}[1]{%
\xybox{%
 (-12,0)*{};
  (12,0)*{};
  (-8,0)*{}="t1";
  (8,0)*{}="t2";
  "t2";"t1" **\crv{(8,10) & (-8,10)}; ?(0)*\dir{<} ?(1)*\dir{<};
  "t2";"t1" **\crv{(8,-10) & (-8,-10)}; ?(.3)*\dir{}+(2,-1)*{\scs #1};
}}
\newcommand{\mccbub}[1]{%
\xybox{%
 (-12,0)*{};
  (12,0)*{};
  (-8,0)*{}="t1";
  (8,0)*{}="t2";
  "t2";"t1" **\crv{(8,10) & (-8,10)}; ?(0)*\dir{>} ?(1)*\dir{>};
  "t2";"t1" **\crv{(8,-10) & (-8,-10)}; ?(.3)*\dir{}+(2,-1)*{\scs #1};
}}
\newcommand{\lccbub}[1]{
\xybox{%
 (-12,0)*{};
  (12,0)*{};
  (-12,0)*{}="t1";
  (12,0)*{}="t2";
  "t2";"t1" **\crv{(12,14) & (-12,14)}; ?(0)*\dir{>} ?(1)*\dir{>};
  "t2";"t1" **\crv{(12,-14) & (-12,-14)}; ?(.3)*\dir{}+(2,-1)*{\scs #1};
}}
\newcommand{\lcbub}[1]{
\xybox{%
 (-12,0)*{};
  (12,0)*{};
  (-12,0)*{}="t1";
  (12,0)*{}="t2";
  "t2";"t1" **\crv{(12,14) & (-12,14)}; ?(0)*\dir{<} ?(1)*\dir{<};
  "t2";"t1" **\crv{(12,-14) & (-12,-14)}; ?(.3)*\dir{}+(2,-1)*{\scs #1};
}}
\newcommand{\xlccbub}[1]{
\xybox{%
 (-12,0)*{};
  (12,0)*{};
  (-20,0)*{}="t1";
  (20,0)*{}="t2";
  "t2";"t1" **\crv{(20,22) & (-20,22)}; ?(0)*\dir{>} ?(1)*\dir{>};
  "t2";"t1" **\crv{(20,-22) & (-20,-22)}; ?(.3)*\dir{}+(2,-1)*{\scs #1};
}}
\newcommand{\xlcbub}[1]{
\xybox{%
 (-12,0)*{};
  (12,0)*{};
  (-20,0)*{}="t1";
  (20,0)*{}="t2";
  "t2";"t1" **\crv{(20,24) & (-20,24)}; ?(0)*\dir{<} ?(1)*\dir{<};
  "t2";"t1" **\crv{(20,-24) & (-20,-24)}; ?(.3)*\dir{}+(2,-1)*{\scs #1};
}}
\newcommand{\scap}{
\xybox{%
(-6,0)*{};
  (6,0)*{};
 (-4,0)*{};(4,0)*{}; **\crv{(4,5) & (-4,5)};
 }}
\newcommand{\mcap}{
\xybox{%
 (0,6)*{};(0,-6)*{};
 (-8,0)*{};(8,0)*{}; **\crv{(8,10) & (-8,10)};
 }}
\newcommand{\lcap}{
\xybox{%
 (-12,2)*{};(12,2)*{}; **\crv{(12,14) & (-12,14)};
 }}
\newcommand{\xlcap}{
\xybox{%
 (-20,0)*{};(20,0)*{}; **\crv{(20,22) & (-20,22)};
 }}

\newcommand{\scupfe}{\xybox{%
 (-4,0)*{};(4,0)*{}; **\crv{(4,-5) & (-4,-5)} ?(0)*\dir{<} ?(.95)*\dir{<};
 }}
\newcommand{\mcupfe}{\xybox{%
 (-8,0)*{};(8,0)*{}; **\crv{(8,-10) & (-8,-10)} ?(0)*\dir{<} ?(.95)*\dir{<};
 }}
\newcommand{\lcupfe}{\xybox{%
 (-12,0)*{};(12,0)*{}; **\crv{(12,-14) & (-12,-14)} ?(0)*\dir{<} ?(.95)*\dir{<};
 }}
\newcommand{\xlcupfe}{\xybox{%
 (-20,0)*{};(20,0)*{}; **\crv{(20,-22) & (-20,-22)} ?(0)*\dir{<} ?(.95)*\dir{<};
 }}
\newcommand{\scupef}{\xybox{%
 (-4,0)*{};(4,0)*{}; **\crv{(4,-5) & (-4,-5)} ?(0)*\dir{>} ?(1)*\dir{>};
 }}
\newcommand{\mcupef}{\xybox{%
 (-8,0)*{};(8,0)*{}; **\crv{(8,-10) & (-8,-10)} ?(0)*\dir{>} ?(1)*\dir{>};
 }}
\newcommand{\lcupef}{\xybox{%
 (12,0)*{};(-12,0)*{} **\crv{(12,-12) & (-12,-12)} ?(0)*\dir{>} ?(1)*\dir{>};
 }}
\newcommand{\xlcupef}{\xybox{%
 (-20,0)*{};(20,0)*{}; **\crv{(20,-22) & (-20,-22)} ?(0)*\dir{>} ?(1)*\dir{>};
 }}
\newcommand{\ecross}{\xybox{%
(-6,0)*{};
  (6,0)*{};
 (-4,-4)*{};(4,4)*{} **\crv{(-4,-1) & (4,1)}?(1)*\dir{>};
 (4,-4)*{};(-4,4)*{} **\crv{(4,-1) & (-4,1)}?(1)*\dir{>};
 }}
\newcommand{\fcross}{\xybox{%
 (-4,-4)*{};(4,4)*{} **\crv{(-4,-1) & (4,1)}?(1)*\dir{<};
 (4,-4)*{};(-4,4)*{} **\crv{(4,-1) & (-4,1)}?(1)*\dir{<};
 }}
\newcommand{\fecross}{\xybox{%
 (-4,-4)*{};(4,4)*{} **\crv{(-4,-1) & (4,1)}?(1)*\dir{>};
 (4,-4)*{};(-4,4)*{} **\crv{(4,-1) & (-4,1)}?(1)*\dir{<};
 }}
\newcommand{\efcross}{\xybox{%
 (-4,-4)*{};(4,4)*{} **\crv{(-4,-1) & (4,1)}?(1)*\dir{<};
 (4,-4)*{};(-4,4)*{} **\crv{(4,-1) & (-4,1)}?(1)*\dir{>};
 }}
\newcommand{\seline}{\xybox{%
 (0,-4)*{};(0,4)*{} **\dir{-} ?(1)*\dir{>};
}}
\newcommand{\sfline}{\xybox{%
 (0,-4)*{};(0,4)*{} **\dir{-} ?(1)*\dir{<};
}}
\newcommand{\meline}{\xybox{%
 (0,-8)*{};(0,8)*{} **\dir{-} ?(1)*\dir{>};
}}
\newcommand{\mfline}{\xybox{%
 (0,-8)*{};(0,8)*{} **\dir{-} ?(1)*\dir{<};
}}
\newcommand{\leline}{\xybox{%
 (0,-12)*{};(0,12)*{} **\dir{-} ?(1)*\dir{>};
}}
\newcommand{\lfline}{\xybox{%
 (0,-12)*{};(0,12)*{} **\dir{-} ?(1)*\dir{<};
}}

\newcommand{\chsheet}[1][1]{
\xy
(0,0)*{
\begin{tikzpicture} [fill opacity=0.2,decoration={markings, mark=at position 0.6 with {\arrow{>}}; }, scale=#1]
	\path [fill=red] (1,1) -- (-1,2) -- (-1,-1) -- (1,-2) -- cycle;
	\draw (-1,.58) -- (-1,2) -- (1,1) -- (1,-2) -- (-1,-1) -- (-1,-.32);
	\draw [fill=red] (0,.75) to [out=225, in=0] (-2,.75) to [out=180, in=90] (-3,0) to [out=270, in=180] (-2,-.5) to [out=0, in=135] (0,-.5);
	\filldraw [fill=white, opacity=1] (-2.5,.15) arc (-120:-60:1) -- (-1.65,.1) arc (70:110:1);
	\draw [thick, red, postaction = {decorate}] (-1,-.3) [out=0,in=270] to (-.75,.13) to [out=90,in=0] (-1,.54);
	\draw [thick, red, dashed] (-1,-.3) to [out=180,in=270] (-1.25,.13) to [out=90,in=180] (-1,.54);
	\path [fill=blue] (-1,.54) to [out=0,in=90] (-.75,.13) to [out=270,in=0] (-1,-.3) to [out=180,in=270] (-1.25,.13) to [out=90,in=180] (-1,.54);
\end{tikzpicture}};
\endxy}

\newcommand{\handlesheet}[1][1]{
\xy
(0,0)*{
\begin{tikzpicture} [fill opacity=0.2,decoration={markings, mark=at position 0.6 with {\arrow{>}}; }, scale=#1]
	\path [fill=red] (1,1) -- (-1,2) -- (-1,-1) -- (1,-2) -- cycle;
	\draw (-1,.58) -- (-1,2) -- (1,1) -- (1,-2) -- (-1,-1) -- (-1,-.32);
	\draw [fill=red] (0,.75) to [out=225, in=0] (-2,.75) to [out=180, in=90] (-3,0) to [out=270, in=180] (-2,-.5) to [out=0, in=135] (0,-.5);
	\filldraw [fill=white, opacity=1] (-2.5,.15) arc (-120:-60:1) -- (-1.65,.1) arc (70:110:1);
\end{tikzpicture}};
\endxy}

\newcommand{\diagsheet}[1][1]{
\xy
(0,0)*{
\begin{tikzpicture} [fill opacity=0.2,decoration={markings, mark=at position 0.6 with {\arrow{>}}; }, scale=#1]
	\draw [fill=red] (1,1) -- (-1,2) -- (-1,-1) -- (1,-2) -- cycle;
\end{tikzpicture}};
\endxy}

\newcommand{\dotsheet}[2][1]{
\xy
(0,0)*{
\begin{tikzpicture} [fill opacity=0.2,decoration={markings, mark=at position 0.6 with {\arrow{>}}; }, scale=#1]
	\draw [fill=red] (1,1) -- (-1,2) -- (-1,-1) -- (1,-2) -- cycle;
	\node [opacity=1] at (0,0) {$\bullet$};
	\node [opacity=1] at (-.3,0){$#2$};
\end{tikzpicture}};
\endxy}

\newcommand{\sphere}[1][1]{
\xy
(0,0)*{
\begin{tikzpicture} [fill opacity=0.2, scale=#1]
	\path [fill=red] (0,0) circle (1);
	\draw (-1,0) .. controls (-1,-.4) and (1,-.4) .. (1,0);
	\draw[dashed] (-1,0) .. controls (-1,.4) and (1,.4) .. (1,0);
	\draw[very thick] (0,0) circle (1);
\end{tikzpicture}};
\endxy}

\newcommand{\dottedsphere}[2][1]{
\xy
(0,0)*{
\begin{tikzpicture} [fill opacity=0.2, scale=#1]
	\path [fill=red] (0,0) circle (1);
	\draw (-1,0) .. controls (-1,-.4) and (1,-.4) .. (1,0);
	\draw[dashed] (-1,0) .. controls (-1,.4) and (1,.4) .. (1,0);
	\draw[very thick] (0,0) circle (1);
	\node [opacity=1]  at (0,0.6) {$\bullet$};
	\node [opacity=1] at (-0.3,0.6) {$#2$};
\end{tikzpicture}};
\endxy}

\newcommand{\twodottedsphere}[1][1]{
\xy
(0,0)*{
\begin{tikzpicture} [fill opacity=0.2, scale=#1]
	\path [fill=red] (0,0) circle (1);
	\draw (-1,0) .. controls (-1,-.4) and (1,-.4) .. (1,0);
	\draw[dashed] (-1,0) .. controls (-1,.4) and (1,.4) .. (1,0);
	\draw[very thick] (0,0) circle (1);
	\node [opacity=1]  at (0.25,0.6) {$\bullet$};
	\node [opacity=1] at (-0.25,0.6) {$\bullet$};
\end{tikzpicture}};
\endxy}

\newcommand{\thetafoam}[4][1]{
\xy
(0,0)*{
\begin{tikzpicture} [fill opacity=0.2,decoration={markings, mark=at position 0.6 with {\arrow{>}}; }, scale=#1]
	\filldraw [fill=red] (0,0) circle (1);
	\path [fill=blue] (0,0) ellipse (1 and 0.3); 	
	\draw [very thick, red, postaction={decorate}] (1,0) .. controls (1,-.4) and (-1,-.4) .. (-1,0);
	\draw[very thick, red, dashed] (-1,0) .. controls (-1,.4) and (1,.4) .. (1,0);
	\node [opacity=1]  at (0,0.7) {$\bullet$};
	\node [opacity=1]  at (0,0) {$\bullet$};
	\node [opacity=1]  at (0,-0.7) {$\bullet$};
	\node [opacity=1] at (-0.3,0.7) {$#2$};
	\node [opacity=1] at (-0.3,0) {$#3$};
	\node [opacity=1] at (-0.3,-0.7) {$#4$};
\end{tikzpicture}};
\endxy}

\newcommand{\cylinder}[1][1]{
\xy
(0,0)*{
\begin{tikzpicture} [fill opacity=0.2,  decoration={markings, mark=at position 0.6 with {\arrow{>}};  }, scale=#1]
	\path [fill=red] (0,4) ellipse (1 and 0.5);
	\path [fill=red] (0,0) ellipse (1 and 0.5);
	\path[fill=red, opacity=.3] (-1,4) .. controls (-1,3.34) and (1,3.34) .. (1,4) --
		(1,0) .. controls (1,.66) and (-1,.66) .. (-1,0) -- cycle;
	\draw [very thick] (0,4) ellipse (1 and 0.5);
	\draw [very thick] (0,0) ellipse (1 and 0.5);
	\draw[very thick] (1,4) -- (1,0);
	\draw[very thick] (-1,4) -- (-1,0);
\end{tikzpicture}};
\endxy
}

\newcommand{\bamboo}[1][1]{
\xy
(0,0)*{
\begin{tikzpicture} [fill opacity=0.2,  decoration={markings, mark=at position 0.6 with {\arrow{>}};  }, scale=#1]
	\draw [fill=red] (0,4) ellipse (1 and 0.5);
	\path [fill=red] (0,0) ellipse (1 and 0.5);
	\draw (1,0) .. controls (1,-.66) and (-1,-.66) .. (-1,0);
	\draw[dashed] (1,0) .. controls (1,.66) and (-1,.66) .. (-1,0);
	\draw (1,4) -- (1,0);
	\draw (-1,4) -- (-1,0);
	\path[fill=red, opacity=.3] (-1,4) .. controls (-1,3.34) and (1,3.34) .. (1,4) --
		(1,0) .. controls (1,.66) and (-1,.66) .. (-1,0) -- cycle;
	\draw [ultra thick, red, postaction={decorate}] (1,2) .. controls (1,1.34) and (-1,1.34) .. (-1,2);
	\draw[ultra thick, dashed, red] (1,2) .. controls (1,2.66) and (-1,2.66) .. (-1,2);
	\path [fill=blue] (0,2) ellipse (1 and 0.5);
\end{tikzpicture}};
\endxy
}

\newcommand{\bambooRHSone}[1][1]{
\xy
(0,0)*{
\begin{tikzpicture} [fill opacity=0.2,  decoration={markings, mark=at position 0.6 with {\arrow{>}};  }, scale=#1]
	\draw [fill=red] (0,4) ellipse (1 and 0.5);
	\draw (-1,4) .. controls (-1,2) and (1,2) .. (1,4);
	\path [fill=red, opacity=.3] (1,4) .. controls (1,3.34) and (-1,3.34) .. (-1,4) --
		(-1,4) .. controls (-1,2) and (1,2) .. (1,4);
	\node[opacity=1] at (0,3) {$\bullet$};
	\path [fill=red] (0,0) ellipse (1 and 0.5);
	\draw (1,0) .. controls (1,-.66) and (-1,-.66) .. (-1,0);
	\draw[dashed] (1,0) .. controls (1,.66) and (-1,.66) .. (-1,0);
	\draw (-1,0) .. controls (-1,2) and (1,2) .. (1,0);
	\path [fill=red, opacity=.3] (1,0) .. controls (1,.66) and (-1,.66) .. (-1,0) --
		(-1,0) .. controls (-1,2) and (1,2) .. (1,0);
\end{tikzpicture}};
\endxy
}

\newcommand{\bambooRHStwo}[1][1]{
\xy
(0,0)*{
\begin{tikzpicture} [fill opacity=0.2,  decoration={markings, mark=at position 0.6 with {\arrow{>}};  }, scale=#1]
	\draw [fill=red] (0,4) ellipse (1 and 0.5);
	\draw (-1,4) .. controls (-1,2) and (1,2) .. (1,4);
	\path [fill=red, opacity=.3] (1,4) .. controls (1,3.34) and (-1,3.34) .. (-1,4) --
		(-1,4) .. controls (-1,2) and (1,2) .. (1,4);
	\path [fill=red] (0,0) ellipse (1 and 0.5);
	\draw (1,0) .. controls (1,-.66) and (-1,-.66) .. (-1,0);
	\draw[dashed] (1,0) .. controls (1,.66) and (-1,.66) .. (-1,0);
	\draw (-1,0) .. controls (-1,2) and (1,2) .. (1,0);
	\path [fill=red, opacity=.3] (1,0) .. controls (1,.66) and (-1,.66) .. (-1,0) --
		(-1,0) .. controls (-1,2) and (1,2) .. (1,0);
	\node[opacity=1] at (0,1) {$\bullet$};
\end{tikzpicture}};
\endxy
}

\newcommand{\airlock}[1][1]{
\xy
(0,0)*{
\begin{tikzpicture} [fill opacity=0.2,  decoration={markings, mark=at position 0.6 with {\arrow{>}};  }, scale=#1]
	\draw [fill=red] (0,4) ellipse (1 and 0.5);
	\path [fill=red] (0,0) ellipse (1 and 0.5);
	\draw (1,0) .. controls (1,-.66) and (-1,-.66) .. (-1,0);
	\draw[dashed] (1,0) .. controls (1,.66) and (-1,.66) .. (-1,0);
	\draw (1,4) -- (1,0);
	\draw (-1,4) -- (-1,0);
	\path[fill=red, opacity=.3] (-1,4) .. controls (-1,3.34) and (1,3.34) .. (1,4) --
		(1,0) .. controls (1,.66) and (-1,.66) .. (-1,0) -- cycle;
	\draw [ultra thick, red, postaction={decorate}] (-1,1.25) .. controls (-1,.59) and (1,.59) .. (1,1.25);
	\draw[ultra thick, dashed, red] (1,1.25) .. controls (1,1.91) and (-1,1.91) .. (-1,1.25);
	\path [fill=blue] (0,1.25) ellipse (1 and 0.5);
	\draw [ultra thick, red, postaction={decorate}] (1,2.75) .. controls (1,2.09) and (-1,2.09) .. (-1,2.75);
	\draw[ultra thick, dashed, red] (1,2.75) .. controls (1,3.41) and (-1,3.41) .. (-1,2.75);
	\path [fill=blue] (0,2.75) ellipse (1 and 0.5);
\end{tikzpicture}};
\endxy
}

\newcommand{\capcup}[1][1]{
\xy
(0,0)*{
\begin{tikzpicture} [fill opacity=0.2,  decoration={markings, mark=at position 0.6 with {\arrow{>}};  }, scale=#1]
	\draw [fill=red] (0,4) ellipse (1 and 0.5);
	\draw (-1,4) .. controls (-1,2) and (1,2) .. (1,4);
	\path [fill=red, opacity=.3] (1,4) .. controls (1,3.34) and (-1,3.34) .. (-1,4) --
		(-1,4) .. controls (-1,2) and (1,2) .. (1,4);
	\path [fill=red] (0,0) ellipse (1 and 0.5);
	\draw (1,0) .. controls (1,-.66) and (-1,-.66) .. (-1,0);
	\draw[dashed] (1,0) .. controls (1,.66) and (-1,.66) .. (-1,0);
	\draw (-1,0) .. controls (-1,2) and (1,2) .. (1,0);
	\path [fill=red, opacity=.3] (1,0) .. controls (1,.66) and (-1,.66) .. (-1,0) --
		(-1,0) .. controls (-1,2) and (1,2) .. (1,0);
\end{tikzpicture}};
\endxy
}

\newcommand{\slthncone}[1][1]{
\xy
(0,0)*{
\begin{tikzpicture} [fill opacity=0.2,  decoration={markings, mark=at position 0.6 with {\arrow{>}};  }, scale=#1]
	\draw [fill=red] (0,4) ellipse (1 and 0.5);
	\draw (-1,4) .. controls (-1,2) and (1,2) .. (1,4);
	\path [fill=red, opacity=.3] (1,4) .. controls (1,3.34) and (-1,3.34) .. (-1,4) --
		(-1,4) .. controls (-1,2) and (1,2) .. (1,4);
	\node[opacity=1] at (.25,3) {$\bullet$};
	\node[opacity=1] at (-.25,3) {$\bullet$};
	\path [fill=red] (0,0) ellipse (1 and 0.5);
	\draw (1,0) .. controls (1,-.66) and (-1,-.66) .. (-1,0);
	\draw[dashed] (1,0) .. controls (1,.66) and (-1,.66) .. (-1,0);
	\draw (-1,0) .. controls (-1,2) and (1,2) .. (1,0);
	\path [fill=red, opacity=.3] (1,0) .. controls (1,.66) and (-1,.66) .. (-1,0) --
		(-1,0) .. controls (-1,2) and (1,2) .. (1,0);
\end{tikzpicture}};
\endxy
}

\newcommand{\slthnctwo}[1][1]{
\xy
(0,0)*{
\begin{tikzpicture} [fill opacity=0.2,  decoration={markings, mark=at position 0.6 with {\arrow{>}};  }, scale=#1]
	\draw [fill=red] (0,4) ellipse (1 and 0.5);
	\draw (-1,4) .. controls (-1,2) and (1,2) .. (1,4);
	\path [fill=red, opacity=.3] (1,4) .. controls (1,3.34) and (-1,3.34) .. (-1,4) --
		(-1,4) .. controls (-1,2) and (1,2) .. (1,4);
	\node[opacity=1] at (0,3) {$\bullet$};
	\path [fill=red] (0,0) ellipse (1 and 0.5);
	\draw (1,0) .. controls (1,-.66) and (-1,-.66) .. (-1,0);
	\draw[dashed] (1,0) .. controls (1,.66) and (-1,.66) .. (-1,0);
	\draw (-1,0) .. controls (-1,2) and (1,2) .. (1,0);
	\path [fill=red, opacity=.3] (1,0) .. controls (1,.66) and (-1,.66) .. (-1,0) --
		(-1,0) .. controls (-1,2) and (1,2) .. (1,0);
		\node[opacity=1] at (0,1) {$\bullet$};
\end{tikzpicture}};
\endxy
}

\newcommand{\slthncthree}[1][1]{
\xy
(0,0)*{
\begin{tikzpicture} [fill opacity=0.2,  decoration={markings, mark=at position 0.6 with {\arrow{>}};  }, scale=#1]
	\draw [fill=red] (0,4) ellipse (1 and 0.5);
	\draw (-1,4) .. controls (-1,2) and (1,2) .. (1,4);
	\path [fill=red, opacity=.3] (1,4) .. controls (1,3.34) and (-1,3.34) .. (-1,4) --
		(-1,4) .. controls (-1,2) and (1,2) .. (1,4);
	\path [fill=red] (0,0) ellipse (1 and 0.5);
	\draw (1,0) .. controls (1,-.66) and (-1,-.66) .. (-1,0);
	\draw[dashed] (1,0) .. controls (1,.66) and (-1,.66) .. (-1,0);
	\draw (-1,0) .. controls (-1,2) and (1,2) .. (1,0);
	\path [fill=red, opacity=.3] (1,0) .. controls (1,.66) and (-1,.66) .. (-1,0) --
		(-1,0) .. controls (-1,2) and (1,2) .. (1,0);
	\node[opacity=1] at (.25,1) {$\bullet$};
	\node[opacity=1] at (-.25,1) {$\bullet$};
\end{tikzpicture}};
\endxy
}

\newcommand{\slthncfour}[1][1]{
\xy
(0,0)*{
\begin{tikzpicture} [fill opacity=0.2,  decoration={markings, mark=at position 0.6 with {\arrow{>}};  }, scale=#1]
	\path [fill=red] (0,4) ellipse (1 and 0.5);
	\path [fill=red, opacity=.3] (1,4) .. controls (1,3.34) and (-1,3.34) .. (-1,4) --
		(-1,4) .. controls (-1,2) and (1,2) .. (1,4);
	\path [fill=red] (0,0) ellipse (1 and 0.5);
	\path [fill=red, opacity=.3] (1,0) .. controls (1,.66) and (-1,.66) .. (-1,0) --
		(-1,0) .. controls (-1,2) and (1,2) .. (1,0);
	\draw [very thick] (0,4) ellipse (1 and 0.5);
	\draw [very thick] (0,0) ellipse (1 and 0.5);
	\draw[very thick] (-1,0) .. controls (-1,2) and (1,2) .. (1,0);
	\draw [very thick] (-1,4) .. controls (-1,2) and (1,2) .. (1,4);
	\node[opacity=1] at (0,1) {$\bullet$};
\end{tikzpicture}};
\endxy
}

\newcommand{\slthncfive}[1][1]{
\xy
(0,0)*{
\begin{tikzpicture} [fill opacity=0.2,  decoration={markings, mark=at position 0.6 with {\arrow{>}};  }, scale=#1]
	\path [fill=red] (0,4) ellipse (1 and 0.5);
	\path [fill=red, opacity=.3] (1,4) .. controls (1,3.34) and (-1,3.34) .. (-1,4) --
		(-1,4) .. controls (-1,2) and (1,2) .. (1,4);
	\path [fill=red] (0,0) ellipse (1 and 0.5);
	\path [fill=red, opacity=.3] (1,0) .. controls (1,.66) and (-1,.66) .. (-1,0) --
		(-1,0) .. controls (-1,2) and (1,2) .. (1,0);
	\draw [very thick] (0,4) ellipse (1 and 0.5);
	\draw [very thick] (0,0) ellipse (1 and 0.5);
	\draw[very thick] (-1,0) .. controls (-1,2) and (1,2) .. (1,0);
	\draw [very thick] (-1,4) .. controls (-1,2) and (1,2) .. (1,4);
	\node[opacity=1] at (0,3) {$\bullet$};
\end{tikzpicture}};
\endxy
}

\newcommand{\tubeLHS}[1][1]{
\xy
(0,0)*{
\begin{tikzpicture} [fill opacity=0.2, scale=#1]
	\path [fill=red] (-2,4) -- (-.75,4) -- (-.75,0) -- (-2,0) -- cycle;
	\path [fill=red] (2,4) -- (.75,4) -- (.75,0) -- (2,0) -- cycle;
	\path [fill=blue] (-.75,4) .. controls (-.5,4.5) and (.5,4.5) .. (.75,4) --
		(.75,0) .. controls (.5,.5) and (-.5,.5) .. (-.75,0);
	\path [fill=blue] (-.75,4) .. controls (-.5,3.5) and (.5,3.5) .. (.75,4) --
		(.75,0) .. controls (.5,-.5) and (-.5,-.5) .. (-.75,0);
	\draw[very thick, directed=.55] (2,0) -- (.75,0);
	\draw[very thick, directed=.55] (-.75,0) -- (-2,0);
	\draw[very thick, directed=.55] (-.75,0) .. controls (-.5,-.5) and (.5,-.5) .. (.75,0);
	\draw[very thick, directed=.55] (-.75,0) .. controls (-.5,.5) and (.5,.5) .. (.75,0);
	\draw[very thick] (-2,0) -- (-2,4);
	\draw[very thick] (2,0) -- (2,4);
	\draw [very thick, red, directed=.55] (-.75,0) -- (-.75,4);
	\draw [very thick, red, directed=.55] (.75, 4) -- (.75,0);
	\draw[very thick, directed=.55] (2,4) -- (.75,4);
	\draw[very thick, directed=.55] (-.75,4) -- (-2,4);
	\draw[very thick, directed=.55] (-.75,4) .. controls (-.5,3.5) and (.5,3.5) .. (.75,4);
	\draw[very thick, directed=.55] (-.75,4) .. controls (-.5,4.5) and (.5,4.5) .. (.75,4);		
\end{tikzpicture}};
\endxy
}

\newcommand{\tubeRHStopdot}[1][1]{
\xy
(0,0)*{
\begin{tikzpicture} [fill opacity=0.2, scale=#1]
	\path [fill=red] (-2,4) -- (-.75,4) -- (-.75,0) -- (-2,0) -- cycle;
	\path [fill=red] (2,4) -- (.75,4) -- (.75,0) -- (2,0) -- cycle;
	\path [fill=red] (-.75,4) .. controls (-.75,2) and (.75,2) .. (.75,4) --
			(.75, 0) .. controls (.75,2) and (-.75,2) .. (-.75,0);
	\path [fill=blue] (-.75,4) .. controls (-.5,4.5) and (.5,4.5) .. (.75,4) --
			(.75,4) .. controls (.75,2) and (-.75,2) .. (-.75,4);
	\path [fill=blue] (-.75,4) .. controls (-.5,3.5) and (.5,3.5) .. (.75,4) --
			(.75,4) .. controls (.75,2) and (-.75,2) .. (-.75,4);
	\path [fill=blue] (-.75, 0) .. controls (-.75,2) and (.75,2) .. (.75,0) --
			(.75,0) .. controls (.5,-.5) and (-.5,-.5) .. (-.75,0);
	\path [fill=blue] (-.75, 0) .. controls (-.75,2) and (.75,2) .. (.75,0) --
			(.75,0) .. controls (.5,.5) and (-.5,.5) .. (-.75,0);
	\node [opacity=1] at (0,4) {$\bullet$};
	\draw[very thick, directed=.55] (2,0) -- (.75,0);
	\draw[very thick, directed=.55] (-.75,0) -- (-2,0);
	\draw[very thick, directed=.55] (-.75,0) .. controls (-.5,-.5) and (.5,-.5) .. (.75,0);
	\draw[very thick, directed=.55] (-.75,0) .. controls (-.5,.5) and (.5,.5) .. (.75,0);
	\draw[very thick] (-2,0) -- (-2,4);
	\draw[very thick] (2,0) -- (2,4);
	\draw [very thick, red, directed=.4] (.75,4) .. controls (.75,2) and (-.75,2) .. (-.75,4);
	\draw [very thick, red, directed=.4] (-.75, 0) .. controls (-.75,2) and (.75,2) .. (.75,0);
	\draw[very thick, directed=.55] (2,4) -- (.75,4);
	\draw[very thick, directed=.55] (-.75,4) -- (-2,4);
	\draw[very thick, directed=.55] (-.75,4) .. controls (-.5,3.5) and (.5,3.5) .. (.75,4);
	\draw[very thick, directed=.55] (-.75,4) .. controls (-.5,4.5) and (.5,4.5) .. (.75,4);
\end{tikzpicture}};
\endxy
}

\newcommand{\tubeRHSbottomdot}[1][1]{
\xy
(0,0)*{
\begin{tikzpicture} [fill opacity=0.2, scale=#1]
	\path [fill=red] (-2,4) -- (-.75,4) -- (-.75,0) -- (-2,0) -- cycle;
	\path [fill=red] (2,4) -- (.75,4) -- (.75,0) -- (2,0) -- cycle;
	\path [fill=red] (-.75,4) .. controls (-.75,2) and (.75,2) .. (.75,4) --
			(.75, 0) .. controls (.75,2) and (-.75,2) .. (-.75,0);
	\path [fill=blue] (-.75,4) .. controls (-.5,4.5) and (.5,4.5) .. (.75,4) --
			(.75,4) .. controls (.75,2) and (-.75,2) .. (-.75,4);
	\path [fill=blue] (-.75,4) .. controls (-.5,3.5) and (.5,3.5) .. (.75,4) --
			(.75,4) .. controls (.75,2) and (-.75,2) .. (-.75,4);
	\path [fill=blue] (-.75, 0) .. controls (-.75,2) and (.75,2) .. (.75,0) --
			(.75,0) .. controls (.5,-.5) and (-.5,-.5) .. (-.75,0);
	\path [fill=blue] (-.75, 0) .. controls (-.75,2) and (.75,2) .. (.75,0) --
			(.75,0) .. controls (.5,.5) and (-.5,.5) .. (-.75,0);
	\node [opacity=1] at (0,0) {$\bullet$};
	\draw[very thick, directed=.55] (2,0) -- (.75,0);
	\draw[very thick, directed=.55] (-.75,0) -- (-2,0);
	\draw[very thick, directed=.55] (-.75,0) .. controls (-.5,-.5) and (.5,-.5) .. (.75,0);
	\draw[very thick, directed=.55] (-.75,0) .. controls (-.5,.5) and (.5,.5) .. (.75,0);
	\draw[very thick] (-2,0) -- (-2,4);
	\draw[very thick] (2,0) -- (2,4);
	\draw [very thick, red, directed=.4] (.75,4) .. controls (.75,2) and (-.75,2) .. (-.75,4);
	\draw [very thick, red, directed=.4] (-.75, 0) .. controls (-.75,2) and (.75,2) .. (.75,0);
	\draw[very thick, directed=.55] (2,4) -- (.75,4);
	\draw[very thick, directed=.55] (-.75,4) -- (-2,4);
	\draw[very thick, directed=.55] (-.75,4) .. controls (-.5,3.5) and (.5,3.5) .. (.75,4);
	\draw[very thick, directed=.55] (-.75,4) .. controls (-.5,4.5) and (.5,4.5) .. (.75,4);
\end{tikzpicture}};
\endxy
}

\newcommand{\Rfoamone}[1][1]{
\xy
(0,0)*{
\begin{tikzpicture} [fill opacity=0.2,  decoration={markings,
                        mark=at position 0.5 with {\arrow{>}};    }, scale=#1]
     	\filldraw [fill=red] (2,2) rectangle (-2,-1);
      	\draw[ultra thick, red, postaction={decorate}] (-.75,-1) -- (-.75,2);
	\draw[ultra thick, red, postaction={decorate}] (.75,2) -- (.75,-1);
	\path [fill=blue] (-.75,2) -- (.25,1) -- (.25,-2) -- (-.75,-1) -- cycle;
	\path [fill=blue] (.75,2) -- (1.75,1) -- (1.75,-2) -- (.75,-1) -- cycle;
	\draw[very thick, postaction={decorate}] (-.75,-1) -- (.25,-2);
     	\draw[very thick, postaction={decorate}] (1.75,-2) --(.75,-1);
	\draw[very thick, postaction={decorate}] (-1,-2) --(.25,-2);
	\draw[very thick, postaction={decorate}] (1.75,-2) --(3,-2);
	\draw[very thick, postaction={decorate}] (1.75,-2) --(.25,-2);
	\draw[very thick, postaction={decorate}] (-.75,-1) -- (-2,-1);
	\draw[very thick, postaction={decorate}] (2,-1) -- (.75,-1);
	\draw[very thick, postaction={decorate}] (-.75,-1) -- (.75,-1);
     	\filldraw [fill=red] (-1,-2) rectangle (3,1);
     	\draw[ultra thick, red, postaction={decorate}] (.25,1) -- (.25,-2);
	\draw[ultra thick, red, postaction={decorate}] (1.75,-2) -- (1.75,1);
	\draw[very thick, postaction={decorate}] (-.75,2) -- (.25,1);
     	\draw[very thick, postaction={decorate}] (1.75,1) --(.75,2);
	\draw[very thick, postaction={decorate}] (-1,1) --(.25,1);
	\draw[very thick, postaction={decorate}] (1.75,1) --(3,1);
	\draw[very thick, postaction={decorate}] (1.75,1) --(.25,1);
	\draw[very thick, postaction={decorate}] (-.75,2) -- (-2,2);
	\draw[very thick, postaction={decorate}] (2,2) -- (.75,2);
	\draw[very thick, postaction={decorate}] (-.75,2) -- (.75,2);
\end{tikzpicture}};
\endxy
}

\newcommand{\Rfoamtwo}[1][1]{
\xy
(0,0)*{
\begin{tikzpicture} [fill opacity=0.2,  decoration={markings,
                        mark=at position 0.5 with {\arrow{>}};    }, scale=#1]
     	\filldraw [fill=red] (2,2) rectangle (-2,-1);
      	\draw[ultra thick, red, postaction={decorate}] (.25,1) arc (180:360:.75);
     	\draw[ultra thick, red, postaction={decorate}] (.75,2) arc (360:180:.75);
	\path [fill=blue] (1.75,-2) arc (0:60:0.75) -- (0.375,-0.35) arc (60:0:0.75);
      	\path [fill=blue] (0.375,-0.35) arc (60:180:0.75) -- (0.25,-2) arc (180:60:0.75);
	\draw[very thick, postaction={decorate}] (-.75,-1) -- (.25,-2);
     	\draw[very thick, postaction={decorate}] (1.75,-2) --(.75,-1);
	\draw[very thick, postaction={decorate}] (-1,-2) --(.25,-2);
	\draw[very thick, postaction={decorate}] (1.75,-2) --(3,-2);
	\draw[very thick, postaction={decorate}] (1.75,-2) --(.25,-2);
	\draw[very thick, postaction={decorate}] (-.75,-1) -- (-2,-1);
	\draw[very thick, postaction={decorate}] (2,-1) -- (.75,-1);
	\draw[very thick, postaction={decorate}] (-.75,-1) -- (.75,-1);
	\path [fill=blue] (.25,1) arc (180:240:0.75) -- (-0.375,1.35) arc (240:180:0.75);
	\path [fill=blue] (0.625,0.35) arc (240:360:0.75) -- (.75,2) arc (0:-120:0.75);
     	\filldraw [fill=red] (-1,-2) rectangle (3,1);
     	\draw[ultra thick, red, postaction={decorate}] (1.75,-2) arc (0:180:.75);
	\draw[ultra thick, red, postaction={decorate}] (-.75,-1) arc (180:0:.75);
	\draw[very thick, postaction={decorate}] (-.75,2) -- (.25,1);
     	\draw[very thick, postaction={decorate}] (1.75,1) --(.75,2);
	\draw[very thick, postaction={decorate}] (-1,1) --(.25,1);
	\draw[very thick, postaction={decorate}] (1.75,1) --(3,1);
	\draw[very thick, postaction={decorate}] (1.75,1) --(.25,1);
	\draw[very thick, postaction={decorate}] (-.75,2) -- (-2,2);
	\draw[very thick, postaction={decorate}] (2,2) -- (.75,2);
	\draw[very thick, postaction={decorate}] (-.75,2) -- (.75,2);
\end{tikzpicture}};
\endxy
}

\newcommand{\Rfoamthree}[1][1]{
\xy
(0,0)*{
\begin{tikzpicture} [fill opacity=0.2,  decoration={markings,
                        mark=at position 0.5 with {\arrow{>}};    }, scale=#1]
     	\draw [fill=red] (2,2) rectangle (.75,-1);
	\draw [fill=red] (-2,2) rectangle (-.75,-1);
	\draw[very thick, postaction={decorate}] (-.75,-1) -- (.25,-2);
     	\draw[very thick, postaction={decorate}] (1.75,-2) --(.75,-1);
	\draw[very thick, postaction={decorate}] (-1,-2) --(.25,-2);
	\draw[very thick, postaction={decorate}] (1.75,-2) --(3,-2);
	\draw[very thick, postaction={decorate}] (1.75,-2) --(.25,-2);
	\draw[very thick, postaction={decorate}] (-.75,-1) -- (-2,-1);
	\draw[very thick, postaction={decorate}] (2,-1) -- (.75,-1);
	\draw[very thick, postaction={decorate}] (-.75,-1) -- (.75,-1);
     	\draw [fill=red] (-1,-2) rectangle (.25,1);
	\draw [fill=red] (3,-2) rectangle (1.75,1);
	\path [fill=red] (-.75,2) -- (.25,1) -- (.25,-2) -- (-.75,-1) -- cycle;
	\path [fill=red] (.75,2) -- (1.75,1) -- (1.75,-2) -- (.75,-1) -- cycle;
	\draw[ultra thick, red, postaction={decorate}] (-.75,-1) .. controls (-.75,0) and (.25,0) .. (.25,-2);
	\draw[ultra thick, red, postaction={decorate}] (1.75,-2) .. controls (1.75,0) and (.75,0) .. (.75,-1);
	\draw[ultra thick, red, postaction={decorate}] (.25,1) .. controls (.25,0) and (-.75,0) .. (-.75,2);
	\draw[ultra thick, red, postaction={decorate}] (.75,2) .. controls (.75,0) and (1.75,0) .. (1.75,1);
	\path [fill=blue] (.25,-2) .. controls (.25,-.8) and (0,-.55) .. (-.25,-.35) --
			(1.25,-.35) .. controls (1.5,-.55) and (1.75,-.8) .. (1.75,-2);
	\path [fill=blue] (-.75,-1) .. controls (-.75,0) and (0,-.5) .. (-.25,-.35) --
			(1.25,-.35) .. controls (1.5,-.5) and (.75,0) .. (.75,-1);
	\path [fill=blue] (.25,1) .. controls (.25,0) and (-.5,.5) .. (-.25,.35) --
			(1.25,.35) .. controls (1,.5) and (1.75,0) .. (1.75,1);
	\path [fill=blue] (-.25,.35) .. controls (-.5,.55) and (-.75,.8) .. (-.75,2) --
			(.75,2) .. controls (.75,.8) and (1,.55) .. (1.25,.35);
	\draw[very thick, postaction={decorate}] (-.75,2) -- (.25,1);
     	\draw[very thick, postaction={decorate}] (1.75,1) --(.75,2);
	\draw[very thick, postaction={decorate}] (-1,1) --(.25,1);
	\draw[very thick, postaction={decorate}] (1.75,1) --(3,1);
	\draw[very thick, postaction={decorate}] (1.75,1) --(.25,1);
	\draw[very thick, postaction={decorate}] (-.75,2) -- (-2,2);
	\draw[very thick, postaction={decorate}] (2,2) -- (.75,2);
	\draw[very thick, postaction={decorate}] (-.75,2) -- (.75,2);
\end{tikzpicture}};
\endxy
}

\newcommand{\Efoam}[1][1]{
\xy
(0,0)*{
\begin{tikzpicture} [fill opacity=0.2,  decoration={markings,
                        mark=at position 0.6 with {\arrow{>}};    }, scale=#1]
     	\path[fill=red] (-2,0) to (2,0) to (2,3) to (-2,3);
     	\path [fill=blue] (1,-1) to (1,2) to (0,3) to (0,0);
     	\filldraw [fill=red] (-1,-1) to (3,-1) to (3,2) to (-1,2);
      	\draw[very thick] (-1,-1) -- (3,-1);
    	\draw[very thick, postaction={decorate}] (1,-1) to (0,0);
     	\draw[very thick] (-2,0) -- (2,0);
      	\draw[ultra thick, red] (0,0) to (0,3);
     	\draw[ultra thick, red] (1,-1) to (1,2);
      	\draw[very thick] (-2,3) to (2,3);
    	\draw[very thick, postaction={decorate}] (1,2) to (0,3);
     	\draw[very thick] (-1,2) to (3,2);
	\draw[very thick] (-1,-1) to (-1,2);
	\draw[very thick] (3,-1) to (3,2);
	\draw[very thick] (2,0) to (2,3);
	\draw[very thick] (-2,0) to (-2,3);
\end{tikzpicture}};
\endxy
}

\newcommand{\dotEfoam}[1][1]{
\xy
(0,0)*{
\begin{tikzpicture} [fill opacity=0.2,  decoration={markings,
                        mark=at position 0.6 with {\arrow{>}};    }, scale=#1]
     	\path[fill=red] (-2,0) to (2,0) to (2,3) to (-2,3);
     	\path [fill=blue] (1,-1) to (1,2) to (0,3) to (0,0);
     	\filldraw [fill=red] (-1,-1) to (3,-1) to (3,2) to (-1,2);
      	\draw[very thick] (-1,-1) -- (3,-1);
    	\draw[very thick, postaction={decorate}] (1,-1) to (0,0);
     	\draw[very thick] (-2,0) -- (2,0);
      	\draw[ultra thick, red] (0,0) to (0,3);
     	\draw[ultra thick, red] (1,-1) to (1,2);
      	\draw[very thick] (-2,3) to (2,3);
    	\draw[very thick, postaction={decorate}] (1,2) to (0,3);
     	\draw[very thick] (-1,2) to (3,2);
	\draw[very thick] (-1,-1) to (-1,2);
	\draw[very thick] (3,-1) to (3,2);
	\draw[very thick] (2,0) to (2,3);
	\draw[very thick] (-2,0) to (-2,3);
	\node [opacity=1]  at (0.5,1) {$\bullet$};	
\end{tikzpicture}};
\endxy
}

\newcommand{\crossingEEfoam}[1][1]{
\xy
(0,0)*{
\begin{tikzpicture} [scale=.3,fill opacity=0.2]
	\path[fill=red] (3.5,0) to (-2.5,0) to (-2.5,7.5) to (3.5,7.5);
	\path[fill=red] (2.5,1) to (-3.5,1) to (-3.5,8.5) to (2.5,8.5);
	\path[fill=blue] (.5,2.75) to (-.5,3.75) to (-.5,5.75) to (.5,4.75);
	\path[fill=blue] (-.75,7.5) to [out=270,in=130] (-.125,4.91) to (-1.125,5.91) to [out=130,in=270] (-1.75,8.5);
	\path[fill=blue] (-.125,4.91) to [out=320,in=180] (.5,4.75) to (-.5,5.75) to [out=180,in=320] (-1.125,5.91);
	\path[fill=blue] (.5,4.75) to [out=0,in=270] (1.75,7.5) to (.75,8.5) to [out=270,in=0] (-.5,5.75);
	\path[fill=blue] (-.75,0) to [out=90,in=180] (.5,2.75) to (-.5,3.75) to [out=180,in=90] (-1.75,1);
	\path[fill=blue] (.5,2.75) to [out=0,in=140] (1.125,2.59) to (.125,3.59) to [out=140,in=0] (-.5,3.75);
	\path[fill=blue] (1.125,2.59) to [out=320,in=270] (1.75,0) to (.75,1) to [out=270,in=320] (.125,3.59);
	\draw[very thick] (3.5,0) to (-2.5,0);
	\draw[very thick] (2.5,1) to (-3.5,1);
	\draw[very thick, directed=.6] (1.75,0) to (.75,1);	
	\draw[very thick, directed=.6] (-.75,0) to (-1.75,1);		
	\draw[very thick, red] (-.75,0) to [out=90,in=180] (.5,2.75);
	\draw[very thick, red] (1.75,0) to [out=90,in=0] (.5,2.75);
	\draw[very thick, red] (.5,2.75) to [out=90,in=270] (.5,4.75);
	\draw[very thick, red] (.5,4.75) to [out=0,in=270] (1.75,7.5);
	\draw[very thick, red] (.5,4.75) to [out=180,in=270] (-.75,7.5);
	\draw[very thick, red] (.5,2.75) to (-.5,3.75);
	\draw[very thick, red] (.5,4.75) to (-.5,5.75);
	\draw[very thick, red] (-1.75,1) to [out=90,in=180] (-.5,3.75);
	\draw[very thick, red] (.75,1) to [out=90,in=0] (-.5,3.75);
	\draw[very thick, red] (-.5,3.75) to [out=90,in=270] (-.5,5.75);
	\draw[very thick, red] (-.5,5.75) to [out=0,in=270] (.75,8.5);
	\draw[very thick, red] (-.5,5.75) to [out=180,in=270] (-1.75,8.5);
	\draw[very thick] (3.5,0) to (3.5,7.5);	
	\draw[very thick] (2.5,1) to (2.5,8.5);
	\draw[very thick] (-2.5,0) to (-2.5,7.5);	
	\draw[very thick] (-3.5,1) to (-3.5,8.5);
	\draw[very thick] (3.5,7.5) to (-2.5,7.5);
	\draw[very thick] (2.5,8.5) to (-3.5,8.5);
	\draw[very thick, directed=.6] (1.75,7.5) to (.75,8.5);	
	\draw[very thick, directed=.6] (-.75,7.5) to (-1.75,8.5);
	\end{tikzpicture}
};
\endxy
}

\newcommand{\crossingEtwoEonefoam}[1][1]{
\xy
(0,0)*{
\begin{tikzpicture}  [fill opacity=0.2,  decoration={markings, mark=at position 0.6 with {\arrow{>}};}, scale=#1]
	\fill [fill=red] (-3,1) rectangle (1,4);
	\path [fill=blue] (-1.75,4) .. controls (-1.75,3) and (-.25,2) .. (-.25,1) --
			(.75,0) .. controls (.75,1) and (-.75,2) .. (-.75,3);
	\fill [fill=red] (-2,0) rectangle (2,3);
	\path [fill=blue] (.75,3) .. controls (.75,2) and (-.75,1) .. (-.75,0) --
			(.25,-1) .. controls (.25,0) and (1.75,1) .. (1.75,2);
	\fill [fill=red] (-1,-1) rectangle (3,2);
	\draw[very thick, postaction={decorate}] (.75,0) -- (-.25,1);
	\draw[very thick, postaction={decorate}] (.25,-1) -- (-.75,0);
	\draw [very thick, red] (-1.75,4) .. controls (-1.75,3) and (-.25,2) .. (-.25,1);
	\draw [very thick] (-3,1) rectangle (1,4);
	\draw [very thick, red] (-.75,3) .. controls (-.75,2) and (.75,1) .. (.75,0);
	\draw [very thick, red] (.75,3) .. controls (.75,2) and (-.75,1) .. (-.75,0);
	\draw [very thick] (-2,0) rectangle (2,3);
	\draw[very thick, postaction={decorate}] (1.75,2) -- (.75,3);
	\draw[very thick, postaction={decorate}] (-.75,3) -- (-1.75,4);
	\draw[very thick, red] (.25,-1) .. controls (.25,0) and (1.75,1) .. (1.75,2);
	\draw [very thick] (-1,-1) rectangle (3,2);
\end{tikzpicture}}
\endxy
}

\newcommand{\crossingFtwoFonefoam}[1][1]{
\xy
(0,0)*{
\begin{tikzpicture}  [fill opacity=0.2, scale=#1]
	\fill [fill=red] (-3,1) rectangle (1,4);
	\path [fill=blue] (-1.75,4) .. controls (-1.75,3) and (-.25,2) .. (-.25,1) --
			(.75,0) .. controls (.75,1) and (-.75,2) .. (-.75,3);
	\fill [fill=red] (-2,0) rectangle (2,3);
	\path [fill=blue] (.75,3) .. controls (.75,2) and (-.75,1) .. (-.75,0) --
			(.25,-1) .. controls (.25,0) and (1.75,1) .. (1.75,2);
	\fill [fill=red] (-1,-1) rectangle (3,2);
	\draw[very thick, rdirected=.45] (.75,0) -- (-.25,1);
	\draw[very thick, rdirected=.45] (.25,-1) -- (-.75,0);
	\draw [very thick, red] (-1.75,4) .. controls (-1.75,3) and (-.25,2) .. (-.25,1);
	\draw [very thick] (-3,1) rectangle (1,4);
	\draw [very thick, red] (-.75,3) .. controls (-.75,2) and (.75,1) .. (.75,0);
	\draw [very thick, red] (.75,3) .. controls (.75,2) and (-.75,1) .. (-.75,0);
	\draw [very thick] (-2,0) rectangle (2,3);
	\draw[very thick, rdirected=.45] (1.75,2) -- (.75,3);
	\draw[very thick, rdirected=.45] (-.75,3) -- (-1.75,4);
	\draw[very thick, red] (.25,-1) .. controls (.25,0) and (1.75,1) .. (1.75,2);
	\draw [very thick] (-1,-1) rectangle (3,2);
\end{tikzpicture}}
\endxy
}

\newcommand{\crossingEtwoFonefoam}[1][1]{
\xy
(0,0)*{
\begin{tikzpicture}  [fill opacity=0.2, scale=#1]
	\fill [fill=red] (-3,1) rectangle (1,4);
	\path [fill=blue] (-1.75,4) .. controls (-1.75,3) and (-.25,2) .. (-.25,1) --
			(.75,0) .. controls (.75,1) and (-.75,2) .. (-.75,3);
	\fill [fill=red] (-2,0) rectangle (2,3);
	\path [fill=blue] (.75,3) .. controls (.75,2) and (-.75,1) .. (-.75,0) --
			(.25,-1) .. controls (.25,0) and (1.75,1) .. (1.75,2);
	\fill [fill=red] (-1,-1) rectangle (3,2);
	\draw[very thick, rdirected=.5] (.75,0) -- (-.25,1);
	\draw[very thick, directed=.55] (.25,-1) -- (-.75,0);
	\draw [very thick, red] (-1.75,4) .. controls (-1.75,3) and (-.25,2) .. (-.25,1);
	\draw [very thick] (-3,1) rectangle (1,4);
	\draw [very thick, red] (-.75,3) .. controls (-.75,2) and (.75,1) .. (.75,0);
	\draw [very thick, red] (.75,3) .. controls (.75,2) and (-.75,1) .. (-.75,0);
	\draw [very thick] (-2,0) rectangle (2,3);
	\draw[very thick, directed=.55] (1.75,2) -- (.75,3);
	\draw[very thick, rdirected=.5] (-.75,3) -- (-1.75,4);
	\draw[very thick, red] (.25,-1) .. controls (.25,0) and (1.75,1) .. (1.75,2);
	\draw [very thick] (-1,-1) rectangle (3,2);
\end{tikzpicture}}
\endxy
}

\newcommand{\crossingFtwoEonefoam}[1][1]{
\xy
(0,0)*{
\begin{tikzpicture}  [fill opacity=0.2, scale=#1]
	\fill [fill=red] (-3,1) rectangle (1,4);
	\path [fill=blue] (-1.75,4) .. controls (-1.75,3) and (-.25,2) .. (-.25,1) --
			(.75,0) .. controls (.75,1) and (-.75,2) .. (-.75,3);
	\fill [fill=red] (-2,0) rectangle (2,3);
	\path [fill=blue] (.75,3) .. controls (.75,2) and (-.75,1) .. (-.75,0) --
			(.25,-1) .. controls (.25,0) and (1.75,1) .. (1.75,2);
	\fill [fill=red] (-1,-1) rectangle (3,2);
	\draw[very thick, directed=.55] (.75,0) -- (-.25,1);
	\draw[very thick, rdirected=.5] (.25,-1) -- (-.75,0);
	\draw [very thick, red] (-1.75,4) .. controls (-1.75,3) and (-.25,2) .. (-.25,1);
	\draw [very thick] (-3,1) rectangle (1,4);
	\draw [very thick, red] (-.75,3) .. controls (-.75,2) and (.75,1) .. (.75,0);
	\draw [very thick, red] (.75,3) .. controls (.75,2) and (-.75,1) .. (-.75,0);
	\draw [very thick] (-2,0) rectangle (2,3);
	\draw[very thick, rdirected=.5] (1.75,2) -- (.75,3);
	\draw[very thick, directed=.55] (-.75,3) -- (-1.75,4);
	\draw[very thick, red] (.25,-1) .. controls (.25,0) and (1.75,1) .. (1.75,2);
	\draw [very thick] (-1,-1) rectangle (3,2);
\end{tikzpicture}}
\endxy
}

\newcommand{\crossingEoneEtwofoam}[1][1]{
\xy
(0,0)*{
\begin{tikzpicture}  [fill opacity=0.2,  decoration={markings, mark=at position 0.6 with {\arrow{>}};}, scale=#1]
	\fill [fill=red] (-3,1) rectangle (1,4);
	\path [fill=blue] (-.25,4) .. controls (-.25,3) and (-1.75,2) .. (-1.75,1) --
			(-.75,0) .. controls (-.75,1) and (.75,2) .. (.75,3);	
	\fill [fill=red] (-2,0) rectangle (2,3);
	\path [fill=blue] (-.75,3) .. controls (-.75,2) and (.75,1) .. (.75,0) --
			(1.75,-1) .. controls (1.75,0) and (.25,1) .. (.25,2);
	\fill [fill=red] (-1,-1) rectangle (3,2);
	\draw [very thick, red] (-.25,4) .. controls (-.25,3) and (-1.75,2) .. (-1.75,1);
	\draw [very thick] (-3,1) rectangle (1,4);
	\draw[very thick, postaction={decorate}] (-.75,0) -- (-1.75,1);
	\draw[very thick, postaction={decorate}] (1.75,-1) -- (.75,0);
	\draw [very thick, red] (-.75,3) .. controls (-.75,2) and (.75,1) .. (.75,0);
	\draw [very thick, red] (.75,3) .. controls (.75,2) and (-.75,1) .. (-.75,0);
	\draw [very thick] (-2,0) rectangle (2,3);
	\draw[very thick, postaction={decorate}] (.25,2) -- (-.75,3);
	\draw[very thick, postaction={decorate}] (.75,3) -- (-.25,4);
	\draw[very thick, red] (1.75,-1) .. controls (1.75,0) and (.25,1) .. (.25,2);
	\draw [very thick] (-1,-1) rectangle (3,2);
\end{tikzpicture}}
\endxy
}

\newcommand{\crossingFoneFtwofoam}[1][1]{
\xy
(0,0)*{
\begin{tikzpicture}  [fill opacity=0.2,  scale=#1]
	\fill [fill=red] (-3,1) rectangle (1,4);
	\path [fill=blue] (-.25,4) .. controls (-.25,3) and (-1.75,2) .. (-1.75,1) --
			(-.75,0) .. controls (-.75,1) and (.75,2) .. (.75,3);	
	\fill [fill=red] (-2,0) rectangle (2,3);
	\path [fill=blue] (-.75,3) .. controls (-.75,2) and (.75,1) .. (.75,0) --
			(1.75,-1) .. controls (1.75,0) and (.25,1) .. (.25,2);
	\fill [fill=red] (-1,-1) rectangle (3,2);
	\draw [very thick, red] (-.25,4) .. controls (-.25,3) and (-1.75,2) .. (-1.75,1);
	\draw [very thick] (-3,1) rectangle (1,4);
	\draw[very thick, rdirected=.45] (-.75,0) -- (-1.75,1);
	\draw[very thick, rdirected=.45] (1.75,-1) -- (.75,0);
	\draw [very thick, red] (-.75,3) .. controls (-.75,2) and (.75,1) .. (.75,0);
	\draw [very thick, red] (.75,3) .. controls (.75,2) and (-.75,1) .. (-.75,0);
	\draw [very thick] (-2,0) rectangle (2,3);
	\draw[very thick, rdirected=.45] (.25,2) -- (-.75,3);
	\draw[very thick, rdirected=.45] (.75,3) -- (-.25,4);
	\draw[very thick, red] (1.75,-1) .. controls (1.75,0) and (.25,1) .. (.25,2);
	\draw [very thick] (-1,-1) rectangle (3,2);
\end{tikzpicture}}
\endxy
}

\newcommand{\crossingEoneFtwofoam}[1][1]{
\xy
(0,0)*{
\begin{tikzpicture}  [fill opacity=0.2,  scale=#1]
	\fill [fill=red] (-3,1) rectangle (1,4);
	\path [fill=blue] (-.25,4) .. controls (-.25,3) and (-1.75,2) .. (-1.75,1) --
			(-.75,0) .. controls (-.75,1) and (.75,2) .. (.75,3);	
	\fill [fill=red] (-2,0) rectangle (2,3);
	\path [fill=blue] (-.75,3) .. controls (-.75,2) and (.75,1) .. (.75,0) --
			(1.75,-1) .. controls (1.75,0) and (.25,1) .. (.25,2);
	\fill [fill=red] (-1,-1) rectangle (3,2);
	\draw [very thick, red] (-.25,4) .. controls (-.25,3) and (-1.75,2) .. (-1.75,1);
	\draw [very thick] (-3,1) rectangle (1,4);
	\draw[very thick, directed=.55] (-.75,0) -- (-1.75,1);
	\draw[very thick, rdirected=.5] (1.75,-1) -- (.75,0);
	\draw [very thick, red] (-.75,3) .. controls (-.75,2) and (.75,1) .. (.75,0);
	\draw [very thick, red] (.75,3) .. controls (.75,2) and (-.75,1) .. (-.75,0);
	\draw [very thick] (-2,0) rectangle (2,3);
	\draw[very thick, rdirected=.5] (.25,2) -- (-.75,3);
	\draw[very thick, directed=.55] (.75,3) -- (-.25,4);
	\draw[very thick, red] (1.75,-1) .. controls (1.75,0) and (.25,1) .. (.25,2);
	\draw [very thick] (-1,-1) rectangle (3,2);
\end{tikzpicture}}
\endxy
}

\newcommand{\crossingFoneEtwofoam}[1][1]{
\xy
(0,0)*{
\begin{tikzpicture}  [fill opacity=0.2,  scale=#1]
	\fill [fill=red] (-3,1) rectangle (1,4);
	\path [fill=blue] (-.25,4) .. controls (-.25,3) and (-1.75,2) .. (-1.75,1) --
			(-.75,0) .. controls (-.75,1) and (.75,2) .. (.75,3);	
	\fill [fill=red] (-2,0) rectangle (2,3);
	\path [fill=blue] (-.75,3) .. controls (-.75,2) and (.75,1) .. (.75,0) --
			(1.75,-1) .. controls (1.75,0) and (.25,1) .. (.25,2);
	\fill [fill=red] (-1,-1) rectangle (3,2);
	\draw [very thick, red] (-.25,4) .. controls (-.25,3) and (-1.75,2) .. (-1.75,1);
	\draw [very thick] (-3,1) rectangle (1,4);
	\draw[very thick, rdirected=.5] (-.75,0) -- (-1.75,1);
	\draw[very thick, directed=.55] (1.75,-1) -- (.75,0);
	\draw [very thick, red] (-.75,3) .. controls (-.75,2) and (.75,1) .. (.75,0);
	\draw [very thick, red] (.75,3) .. controls (.75,2) and (-.75,1) .. (-.75,0);
	\draw [very thick] (-2,0) rectangle (2,3);
	\draw[very thick, directed=.55] (.25,2) -- (-.75,3);
	\draw[very thick, rdirected=.5] (.75,3) -- (-.25,4);
	\draw[very thick, red] (1.75,-1) .. controls (1.75,0) and (.25,1) .. (.25,2);
	\draw [very thick] (-1,-1) rectangle (3,2);
\end{tikzpicture}}
\endxy
}

\newcommand{\crossingEoneEthreefoam}[1][1]{
\xy
(0,0)*{
\begin{tikzpicture} [fill opacity=0.2,decoration={markings, mark=at position 0.6 with {\arrow{>}}; }, scale=#1]
	\fill [fill=red] (-4,2) rectangle (0,6);
	\path [fill=blue] (-1.25,6) .. controls (-1.25,5) and (-2.75,3) .. (-2.75,2) --
		 (-1.75,1) .. controls (-1.75,2) and (-.25,4) .. (-.25,5);
	\fill [fill=red] (-3,1) rectangle (1,5);
	\fill [fill=red] (-2,0) rectangle (2,4);
	\path [fill=blue] (-.75,4) .. controls (-.75,3) and (.75,1) .. (.75,0) --
		(1.75,-1) .. controls (1.75,0) and (.25,2) .. (.25,3);
	\fill [fill=red] (-1,-1) rectangle (3,3);
	\draw[very thick, red] (-1.25,6) .. controls (-1.25,5) and (-2.75,3) .. (-2.75,2);
	\draw [very thick] (-4,2) rectangle (0,6);
	\draw[very thick, postaction={decorate}] (-1.75,1) -- (-2.75,2);
	\draw[very thick, postaction={decorate}] (1.75,-1) -- (.75,0);
	\draw[very thick, red] (-1.75,1) .. controls (-1.75,2) and (-.25,4) .. (-.25,5);
	\draw [very thick] (-3,1) rectangle (1,5);
	\draw[very thick, red] (-.75,4) .. controls (-.75,3) and (.75,1) .. (.75,0);
	\draw [very thick] (-2,0) rectangle (2,4);
	\draw[very thick, postaction={decorate}] (.25,3) -- (-.75,4);
	\draw[very thick, postaction={decorate}] (-.25,5) -- (-1.25,6);
	\draw[very thick, red] (1.75,-1) .. controls (1.75,0) and (.25,2) .. (.25,3);
	\draw [very thick] (-1,-1) rectangle (3,3);
\end{tikzpicture}}
\endxy
}

\newcommand{\crossingEthreeEonefoam}[1][1]{
\xy
(0,0)*{
\begin{tikzpicture} [fill opacity=0.2,decoration={markings, mark=at position 0.6 with {\arrow{>}}; }, scale=#1]
	\fill [fill=red] (-4,2) rectangle (0,6);
	\path [fill=blue] (-2.75,6) .. controls (-2.75,5) and (-1.25,3) .. (-1.25,2) --
		 (-.25,1) .. controls (-.25,2) and (-1.75,4) .. (-1.75,5);
	\fill [fill=red] (-3,1) rectangle (1,5);
	\fill [fill=red] (-2,0) rectangle (2,4);
	\path [fill=blue] (.75,4) .. controls (.75,3) and (-.75,1) .. (-.75,0) --
		(.25,-1) .. controls (.25,0) and (1.75,2) .. (1.75,3);
	\fill [fill=red] (-1,-1) rectangle (3,3);
	\draw[very thick, red] (-2.75,6) .. controls (-2.75,5) and (-1.25,3) .. (-1.25,2);
	\draw [very thick] (-4,2) rectangle (0,6);
	\draw[very thick, postaction={decorate}] (-.25,1) -- (-1.25,2);
	\draw[very thick, postaction={decorate}] (.25,-1) -- (-.75,0);
	\draw[very thick, red] (-.25,1) .. controls (-.25,2) and (-1.75,4) .. (-1.75,5);
	\draw [very thick] (-3,1) rectangle (1,5);
	\draw[very thick, red] (.75,4) .. controls (.75,3) and (-.75,1) .. (-.75,0);
	\draw [very thick] (-2,0) rectangle (2,4);
	\draw[very thick, postaction={decorate}] (1.75,3) -- (.75,4);
	\draw[very thick, postaction={decorate}] (-1.75,5) -- (-2.75,6);
	\draw[very thick, red] (.25,-1) .. controls (.25,0) and (1.75,2) .. (1.75,3);
	\draw [very thick] (-1,-1) rectangle (3,3);
\end{tikzpicture}}
\endxy
}

\newcommand{\cupEFfoam}[3][.5]{
\xy
(0,0)*{
\begin{tikzpicture} [scale=#1,fill opacity=0.2]
\path[fill=red] (2,-2) to (-2,-2) to (-2,1) to (2,1);
\path[fill=blue] (1.5,2) to [out=270,in=0] (0,-.5) to [out=180,in=270] (-1.5,2) to
	(-.5,1) to [out=270,in=180] (0,.25) to [out=0,in=270] (.5,1);
\path[fill=red] (2.5,-1) to (-2.5,-1) to (-2.5,2) to (2.5,2);
\draw[very thick, directed=.5] (2.5,-1) to (-2.5,-1);
\draw[very thick] (2.5,-1) to (2.5,2);
\draw[very thick] (-2.5,-1) to (-2.5,2);
\draw[very thick, directed=.5] (2.5,2) to (-2.5,2);
\draw[very thick, red, rdirected=.65] (1.5,2) to [out=270,in=0] (0,-.5)
	to [out=180,in=270] (-1.5,2);
\node[red, opacity=1] at (2,1.5) {$_{#2}$};
\draw[very thick, directed=.5] (2,-2) to (-2,-2);
\draw[very thick] (2,-2) to (2,1);
\draw[very thick] (-2,-2) to (-2,1);
\draw[very thick, directed=.5] (2,1) to (-2,1);
\draw[very thick, red, rdirected=.75] (-.5,1) to [out=270,in=180] (0,.25)
	to [out=0,in=270] (.5,1);
\node[red, opacity=1] at (1.25,-1.5) {$_{#3}$};
\draw[very thick, directed=.5] (1.5,2) to (.5,1);
\draw[very thick, directed=.5] (-.5,1) to (-1.5,2);
\end{tikzpicture}};
\endxy
}

\newcommand{\cupFEfoam}[3][.5]{
\xy
(0,0)*{
\begin{tikzpicture} [scale=#1,fill opacity=0.2]
\path[fill=red] (2,-2) to (-2,-2) to (-2,1) to (2,1);
\path[fill=blue] (1.5,1) to [out=270,in=0] (0,-.25) to [out=180,in=270] (-1.5,1) to
	(-.5,2) to [out=270,in=180] (0,1.25) to [out=0,in=270] (.5,2);
\path[fill=red] (2.5,-1) to (-2.5,-1) to (-2.5,2) to (2.5,2);
\draw[very thick, directed=.5] (2.5,-1) to (-2.5,-1);
\draw[very thick] (2.5,-1) to (2.5,2);
\draw[very thick] (-2.5,-1) to (-2.5,2);
\draw[very thick, directed=.5] (2.5,2) to (-2.5,2);
\draw[very thick, red, rdirected=.5] (-.5,2) to [out=270,in=180] (0,1.25)
	to [out=0,in=270] (.5,2);
\node[red, opacity=1] at (2,1.5) {$_{#2}$};
\draw[very thick, directed=.5] (2,-2) to (-2,-2);
\draw[very thick] (2,-2) to (2,1);
\draw[very thick] (-2,-2) to (-2,1);
\draw[very thick, directed=.5] (2,1) to (-2,1);
\draw[very thick, red, rdirected=.5] (1.5,1) to [out=270,in=0] (0,-.25)
	to [out=180,in=270] (-1.5,1);
\node[red, opacity=1] at (1.25,-1.5) {$_{#3}$};
\draw[very thick, directed=.5] (1.5,1) to (.5,2);
\draw[very thick, directed=.5] (-.5,2) to (-1.5,1);
\end{tikzpicture}
}
\endxy
}

\newcommand{\capEFfoam}[3][.5]{
\xy
(0,0)*{
\begin{tikzpicture} [scale=#1,fill opacity=0.2]
\path[fill=red] (2.5,2) to (-2.5,2) to (-2.5,-1) to (2.5,-1);
\path[fill=blue] (1.5,-1) to [out=90,in=0] (0,.25) to [out=180,in=90] (-1.5,-1) to
	(-.5,-2) to [out=90,in=180] (0,-1.25) to [out=0,in=90] (.5,-2);
\path[fill=red] (2,1) to (-2,1) to (-2,-2) to (2,-2);
\draw[very thick, directed=.5] (2.5,2) to (-2.5,2);
\draw[very thick] (2.5,2) to (2.5,-1);
\draw[very thick] (-2.5,2) to (-2.5,-1);
\draw[very thick, directed=.5] (2.5,-1) to (-2.5,-1);
\draw[very thick, red, directed=.5] (1.5,-1) to [out=90,in=0] (0,.25)
	to [out=180,in=90] (-1.5,-1);
\node[red, opacity=1] at (2,1.5) {$_{#2}$};
\draw[very thick, directed=.5] (1.5,-1) to (.5,-2);
\draw[very thick, directed=.5] (-.5,-2) to (-1.5,-1);
\draw[very thick, directed=.5] (2,1) to (-2,1);
\draw[very thick] (2,1) to (2,-2);
\draw[very thick] (-2,1) to (-2,-2);
\draw[very thick, directed=.5] (2,-2) to (-2,-2);
\draw[very thick, red, directed=.5] (-.5,-2) to [out=90,in=180] (0,-1.25)
	to [out=0,in=90] (.5,-2);
\node[red, opacity=1] at (1.25,.5) {$_{#3}$};
\end{tikzpicture}};
\endxy
}

\newcommand{\capFEfoam}[3][.5]{
\xy
(0,0)*{
\begin{tikzpicture} [scale=#1,fill opacity=0.2]
\path[fill=red] (2.5,2) to (-2.5,2) to (-2.5,-1) to (2.5,-1);
\path[fill=blue] (1.5,-2) to [out=90,in=0] (0,.5) to [out=180,in=90] (-1.5,-2) to
	(-.5,-1) to [out=90,in=180] (0,-.25) to [out=0,in=90] (.5,-1);
\path[fill=red] (2,1) to (-2,1) to (-2,-2) to (2,-2);
\draw[very thick, directed=.5] (2.5,2) to (-2.5,2);
\draw[very thick] (2.5,2) to (2.5,-1);
\draw[very thick] (-2.5,2) to (-2.5,-1);
\draw[very thick, directed=.5] (2.5,-1) to (-2.5,-1);
\draw[very thick, red, directed=.75] (-.5,-1) to [out=90,in=180] (0,-.25)
	to [out=0,in=90] (.5,-1);
\node[red, opacity=1] at (1.75,1.5) {$_{#2}$};
\draw[very thick, directed=.5] (1.5,-2) to (.5,-1);
\draw[very thick, directed=.5] (-.5,-1) to (-1.5,-2);
\draw[very thick, directed=.5] (2,1) to (-2,1);
\draw[very thick] (2,1) to (2,-2);
\draw[very thick] (-2,1) to (-2,-2);
\draw[very thick, directed=.5] (2,-2) to (-2,-2);
\draw[very thick, red, directed=.65] (1.5,-2) to [out=90,in=0] (0,.5)
	to [out=180,in=90] (-1.5,-2);
\node[red, opacity=1] at (1.25,.5) {$_{#3}$};
\end{tikzpicture}};
\endxy
}

\newcommand{\bbullet}{
\begin{tikzpicture}
  \draw[fill=black] circle (0.55ex);
\end{tikzpicture}
}

\newcommand{\hackcenter}[1]{
 \xy (0,0)*{#1}; \endxy}

\usepackage{xparse}

\newcommand{\pict}[1]{
 \xy (0,0)*{ #1};
\endxy
}

\DeclareDocumentCommand{\Ustrand}{ O{black} O{semithick}}{
 \xy (0,0)*{\begin{tikzpicture}[scale=0.75] 
  \draw[ -,color=#1, #2] (0,0) -- (0,1);
\end{tikzpicture}};
(1.5,0)*{};(-1.5,0)*{}; \endxy
    }

\DeclareDocumentCommand{\sUstrand}{ O{black} O{semithick}}{
 \xy (0,0)*{\begin{tikzpicture}[scale=0.75] 
  \draw[#2,color=#1, -] (0,0) -- (0,.7);
\end{tikzpicture}};
(1.5,0)*{};(-1.5,0)*{}; \endxy
    }

\DeclareDocumentCommand{\Uup}{ O{black} O{semithick}}{
 \xy (0,0)*{\begin{tikzpicture}[scale=0.75]
  \draw[ ->,color=#1, #2] (0,0) -- (0,1);
\end{tikzpicture}};
(1.5,0)*{};(-1.5,0)*{}; \endxy
    }

\DeclareDocumentCommand{\sUup}{ O{black} O{semithick}}{
 \xy (0,0)*{\begin{tikzpicture}[scale=0.75]
  \draw[#2,color=#1, ->] (0,0) -- (0,.7);
\end{tikzpicture}};
(1.5,0)*{};(-1.5,0)*{}; \endxy
    }

\DeclareDocumentCommand{\lUup}{ O{black} O{semithick} O{1}}{
 \xy (0,0)*{\begin{tikzpicture}[scale=0.75]
  \draw[ ->,color=#1, #2] (0,0) -- (0,#3);
\end{tikzpicture}};
(1.5,0)*{};(-1.5,0)*{}; \endxy
    }

\DeclareDocumentCommand{\Udown}{ O{black} O{semithick}}{
 \xy (0,0)*{\begin{tikzpicture}[scale=0.75] 
  \draw[#2,color=#1, <-] (0,0) -- (0,1);
\end{tikzpicture}};
(1.5,0)*{};(-1.5,0)*{}; \endxy
    }

\DeclareDocumentCommand{\sUdown}{ O{black} O{semithick}}{
 \xy (0,0)*{\begin{tikzpicture}[scale=0.75] 
  \draw[#2,color=#1, <-] (0,0) -- (0,.7);
\end{tikzpicture}};
(1.5,0)*{};(-1.5,0)*{}; \endxy
    }

\DeclareDocumentCommand{\lUdown}{ O{black} O{semithick} O{1}}{
 \xy (0,0)*{\begin{tikzpicture}[scale=0.75] 
  \draw[#2,color=#1, <-] (0,0) -- (0,#3);
\end{tikzpicture}};
(1.5,0)*{};(-1.5,0)*{}; \endxy
    }

\DeclareDocumentCommand{\Uupdot}{ O{black} O{semithick}}{
 \xy (0,0)*{\begin{tikzpicture}[scale=0.75] 
  \draw[#2,color=#1, ->] (0,0) -- (0,1);
  \node at (0,.5) {$\bullet$};
\end{tikzpicture}};
(1.5,0)*{};(-2,0)*{}; \endxy
    }

\DeclareDocumentCommand{\sUupdot}{ O{black} O{semithick}}{
 \xy (0,0)*{\begin{tikzpicture}[scale=0.75]
  \draw[#2,color=#1, ->] (0,0) -- (0,.7);
  \node at (0,.25) {$ \bullet$};
\end{tikzpicture}};
(1.5,0)*{};(-2,0)*{}; \endxy
    }

\DeclareDocumentCommand{\Uupdots}{m O{black} O{semithick}}{
 \xy (0,0)*{\begin{tikzpicture}[scale=0.75] 
  \draw[#3, color=#2, ->] (0,0) -- (0,1);
  \node at (0,.5) {$\bullet$};
\node at (.3,.5) {$\scriptstyle{#1}$};
\node at (.3,-.5) {\;};
\end{tikzpicture}};
(1.5,0)*{};(-2,0)*{}; \endxy
    }

\DeclareDocumentCommand{\Udowndot}{ O{black} O{semithick}}{
 \xy (0,0)*{\begin{tikzpicture}[scale=0.75] 
  \draw[#2,color=#1, <-] (0,0) -- (0,1);
  \node at (0,.5) {$\bullet$};
\end{tikzpicture}};
(2,0)*{};(-1.5,0)*{}; \endxy
    }

\DeclareDocumentCommand{\sUdowndot}{ O{black} O{semithick}}{
 \xy (0,0)*{\begin{tikzpicture}[scale=0.75] 
  \draw[#2,color=#1, <-] (0,0) -- (0,.7);
  \node at (0,.35) {$\bullet$};
\end{tikzpicture}};
(2,0)*{};(-1.5,0)*{}; \endxy
    }

\DeclareDocumentCommand{\Ucapl}{ O{black} O{semithick}}{\;\;
\begin{tikzpicture}[scale=0.75]
  \draw[#2,color=#1, ->] (0.5,0) .. controls (0.5,0.8) and (-0.5,0.8) .. (-0.5,0);
  \node at (-.5,.8) {};
  \node at (.5,.8) {};
\end{tikzpicture}\;\; }

\DeclareDocumentCommand{\vUcapl}{ O{black} O{semithick} O{.5}}{\;\;
\begin{tikzpicture}[scale=#3]
  \draw[#2,color=#1, ->] (0.5,0) .. controls (0.5,0.8) and (-0.5,0.8) .. (-0.5,0);
  \node at (-.5,.8) {};
  \node at (.5,.8) {};
\end{tikzpicture}\;\; }

\DeclareDocumentCommand{\Ucapr}{ O{black} O{semithick}}{\;\;
\begin{tikzpicture}[scale=0.75]
  \draw[#2,color=#1, ->] (-0.5,0) .. controls (-0.5,0.8) and (0.5,0.8) .. (0.5,0)
      node[pos=0.5, shape=coordinate](X){};
  \node at (-.5,.8) {};
  \node at (.5,.8) {};
\end{tikzpicture}\;\; }

\DeclareDocumentCommand{\Ucap}{ O{black} O{semithick}}{\;\;
\begin{tikzpicture}[scale=0.75]
  \draw[#2,color=#1, -] (-0.5,0) .. controls (-0.5,0.8) and (0.5,0.8) .. (0.5,0)
      node[pos=0.5, shape=coordinate](X){};
  \node at (-.5,.8) {};
  \node at (.5,.8) {};
\end{tikzpicture}\;\; }

\DeclareDocumentCommand{\vUcapr}{ O{black} O{semithick} O{.5}}{\;\;
\begin{tikzpicture}[scale=#3]
  \draw[#2,color=#1, ->] (-0.5,0) .. controls (-0.5,0.8) and (0.5,0.8) .. (0.5,0)
      node[pos=0.5, shape=coordinate](X){};
  \node at (-.5,.8) {};
  \node at (.5,.8) {};
\end{tikzpicture}\;\; }

\DeclareDocumentCommand{\Ucupl}{ O{black} O{semithick}}{\;\;
\xy
(0,0)*{
\begin{tikzpicture}[scale=0.75]
  \draw[#2,color=#1, ->] (0.5,1) .. controls (0.5,0.2) and (-0.5,0.2) .. (-0.5,1);
  \node at (-.5,.2) {};
  \node at (.5,.2) {};
\end{tikzpicture}}; \endxy
\;\; }

\DeclareDocumentCommand{\Ucupr}{ O{black} O{semithick}}{\;\;
\xy
(0,0)*{
\begin{tikzpicture}[scale=0.75]
  \draw[#2,color=#1, ->] (-0.5,1) .. controls (-0.5,0.2) and (0.5,0.2) .. (0.5,1)
  node[pos=0.5, shape=coordinate](X){};;
  \node at (-.5,.2) {};
  \node at (.5,.2) {};
\end{tikzpicture}}; \endxy
\;\; }

\DeclareDocumentCommand{\Ucup}{ O{black} O{semithick}}{\;\;
\xy
(0,0)*{
\begin{tikzpicture}[scale=0.75]
  \draw[#2,color=#1, -] (0.5,1) .. controls (0.5,0.2) and (-0.5,0.2) .. (-0.5,1);
  \node at (-.5,.2) {};
  \node at (.5,.2) {};
\end{tikzpicture}}; \endxy
\;\; }

\DeclareDocumentCommand{\Ucross}{ O{black} O{semithick}}{\;\;
    \vcenter{\xy (0,0)*{\begin{tikzpicture}[scale=0.75]
  \draw[#2,color=#1, -] (-0.5,0) .. controls (-0.5,0.5) and (0.5,0.5) .. (0.5,1)
      node[pos=0.5, shape=coordinate](X){};
    \draw[#2,color=#1, -] (0.5,0) .. controls (0.5,0.5) and (-0.5,0.5) .. (-0.5,1);
\end{tikzpicture}}; \endxy} \;\; }

\DeclareDocumentCommand{\Ucrossr}{ O{black} O{semithick}}{\;\;
    \vcenter{\xy (0,0)*{\begin{tikzpicture}[scale=0.75]
  \draw[#2,color=#1, ->] (-0.5,0) .. controls (-0.5,0.5) and (0.5,0.5) .. (0.5,1)
      node[pos=0.5, shape=coordinate](X){};
    \draw[#2,color=#1, <-] (0.5,0) .. controls (0.5,0.5) and (-0.5,0.5) .. (-0.5,1);
\end{tikzpicture}}; \endxy} \;\; }

\DeclareDocumentCommand{\Ucrossl}{ O{black} O{semithick}}{\;\;
    \vcenter{\xy (0,0)*{\begin{tikzpicture}[scale=0.75]
  \draw[#2,color=#1, <-] (-0.5,0) .. controls (-0.5,0.5) and (0.5,0.5) .. (0.5,1)
      node[pos=0.5, shape=coordinate](X){};
  \draw[#2,color=#1, ->] (0.5,0) .. controls (0.5,0.5) and (-0.5,0.5) .. (-0.5,1);
\end{tikzpicture}}; \endxy} \;\; }

\DeclareDocumentCommand{\Ucrossu}{ O{black} O{semithick}}{\;\;
    \vcenter{\xy (0,0)*{\begin{tikzpicture}[scale=0.75]
  \draw[#2,color=#1, ->] (-0.5,0) .. controls (-0.5,0.5) and (0.5,0.5) .. (0.5,1)
      node[pos=0.5, shape=coordinate](X){};
  \draw[#2,color=#1, ->] (0.5,0) .. controls (0.5,0.5) and (-0.5,0.5) .. (-0.5,1);
\end{tikzpicture}}; \endxy} \;\; }

\DeclareDocumentCommand{\Ucrossd}{ O{black} O{semithick}}{\;\;
    \vcenter{\xy (0,0)*{\begin{tikzpicture}[scale=0.75]
  \draw[#2,color=#1, <-] (-0.5,0) .. controls (-0.5,0.5) and (0.5,0.5) .. (0.5,1)
      node[pos=0.5, shape=coordinate](X){};
  \draw[#2,color=#1, <-] (0.5,0) .. controls (0.5,0.5) and (-0.5,0.5) .. (-0.5,1);
\end{tikzpicture}}; \endxy} \;\; }

\DeclareDocumentCommand{\Umerge}{ O{black} O{semithick}}{\;\;
    \vcenter{\xy (0,0)*{\begin{tikzpicture}[scale=0.75]
  \draw[#2,color=#1, -] (-0.5,0) .. controls (-0.5,0.4) and (-.2,.4)  .. (0,0.5);
  \draw[#2,color=#1, -] (0.5,0) .. controls (0.5,0.4) and (.2,.4)  .. (0,0.5);
  \draw[ultra thick,color=#1, -] (0,.5)  -- (0,1);
\end{tikzpicture}}; \endxy} \;\; }

\DeclareDocumentCommand{\Umergeu}{ O{black} O{semithick}}{\;\;
    \vcenter{\xy (0,0)*{\begin{tikzpicture}[scale=0.75]
  \draw[#2,color=#1, -] (-0.5,0) .. controls (-0.5,0.4) and (-.2,.4)  .. (0,0.5);
  \draw[#2,color=#1, -] (0.5,0) .. controls (0.5,0.4) and (.2,.4)  .. (0,0.5);
  \draw[ultra thick,color=#1, ->] (0,.5)  -- (0,1);
\end{tikzpicture}}; \endxy} \;\; }

\DeclareDocumentCommand{\Umerged}{ O{black} O{semithick}}{\;\;
    \vcenter{\xy (0,0)*{\begin{tikzpicture}[scale=0.75]
  \draw[#2,color=#1, <-] (-0.5,0) .. controls (-0.5,0.4) and (-.2,.4)  .. (0,0.5);
  \draw[#2,color=#1, <-] (0.5,0) .. controls (0.5,0.4) and (.2,.4)  .. (0,0.5);
  \draw[ultra thick,color=#1, -] (0,.5)  -- (0,1);
\end{tikzpicture}}; \endxy} \;\; }

\DeclareDocumentCommand{\Usplit}{ O{black} O{semithick}}{\;\;
    \vcenter{\xy (0,0)*{\begin{tikzpicture}[scale=0.75]
  \draw[#2,color=#1, -] (0,0.5)  .. controls (-.2,.6) and (-.5,.6) ..  (-0.5,1);
  \draw[#2,color=#1, -] (0,0.5)  .. controls (.2,.6) and (.5,.6) ..  (0.5,1);
  \draw[ultra thick,color=#1, -] (0,0)  -- (0,.5);
\end{tikzpicture}}; \endxy} \;\; }

\DeclareDocumentCommand{\Usplitu}{ O{black} O{semithick}}{\;\;
    \vcenter{\xy (0,0)*{\begin{tikzpicture}[scale=0.75]
  \draw[#2,color=#1, ->] (0,0.5)  .. controls (-.2,.6) and (-.5,.6) ..  (-0.5,1);
  \draw[#2,color=#1, ->] (0,0.5)  .. controls (.2,.6) and (.5,.6) ..  (0.5,1);
  \draw[ultra thick,color=#1, -] (0,0)  -- (0,.5);
\end{tikzpicture}}; \endxy} \;\; }

\DeclareDocumentCommand{\vUsplitu}{ O{black} O{semithick} O{.5}}{\;\;
    \vcenter{\xy (0,0)*{\begin{tikzpicture}[scale=#3]
  \draw[#2,color=#1, ->] (0,0.5)  .. controls (-.2,.6) and (-.5,.6) ..  (-0.5,1);
  \draw[#2,color=#1, ->] (0,0.5)  .. controls (.2,.6) and (.5,.6) ..  (0.5,1);
  \draw[ultra thick,color=#1, -] (0,0)  -- (0,.5);
\end{tikzpicture}}; \endxy} \;\; }

\DeclareDocumentCommand{\Usplitd}{ O{black} O{semithick}}{\;\;
    \vcenter{\xy (0,0)*{\begin{tikzpicture}[scale=0.75]
  \draw[#2,color=#1, -] (0,0.5)  .. controls (-.2,.6) and (-.5,.6) ..  (-0.5,1);
  \draw[#2,color=#1, -] (0,0.5)  .. controls (.2,.6) and (.5,.6) ..  (0.5,1);
  \draw[ultra thick,color=#1, <-] (0,0)  -- (0,.5);
\end{tikzpicture}}; \endxy} \;\; }

\DeclareDocumentCommand{\Ucrossm}{ O{black} O{black} O{semithick}}{\;\;
    \vcenter{\xy (0,0)*{\begin{tikzpicture}[scale=0.75]
  \draw[#3,color=#1, -] (-0.5,0) .. controls (-0.5,0.5) and (0.5,0.5) .. (0.5,1)
      node[pos=0.5, shape=coordinate](X){};
    \draw[#3,color=#2, -] (0.5,0) .. controls (0.5,0.5) and (-0.5,0.5) .. (-0.5,1);
\end{tikzpicture}}; \endxy} \;\; }

\DeclareDocumentCommand{\Ucrossrm}{ O{black} O{black} O{semithick}}{\;\;
    \vcenter{\xy (0,0)*{\begin{tikzpicture}[scale=0.75]
  \draw[#3,color=#1, ->] (-0.5,0) .. controls (-0.5,0.5) and (0.5,0.5) .. (0.5,1)
      node[pos=0.5, shape=coordinate](X){};
    \draw[#3,color=#2, <-] (0.5,0) .. controls (0.5,0.5) and (-0.5,0.5) .. (-0.5,1);
\end{tikzpicture}}; \endxy} \;\; }

\DeclareDocumentCommand{\Ucrosslm}{ O{black} O{black} O{semithick}}{\;\;
    \vcenter{\xy (0,0)*{\begin{tikzpicture}[scale=0.75]
  \draw[#3,color=#1, <-] (-0.5,0) .. controls (-0.5,0.5) and (0.5,0.5) .. (0.5,1)
      node[pos=0.5, shape=coordinate](X){};
  \draw[#3,color=#2, ->] (0.5,0) .. controls (0.5,0.5) and (-0.5,0.5) .. (-0.5,1);
\end{tikzpicture}}; \endxy} \;\; }

\DeclareDocumentCommand{\Ucrossum}{ O{black} O{black} O{semithick}}{\;\;
    \vcenter{\xy (0,0)*{\begin{tikzpicture}[scale=0.75]
  \draw[#3,color=#1, ->] (-0.5,0) .. controls (-0.5,0.5) and (0.5,0.5) .. (0.5,1)
      node[pos=0.5, shape=coordinate](X){};
  \draw[#3,color=#2, ->] (0.5,0) .. controls (0.5,0.5) and (-0.5,0.5) .. (-0.5,1);
\end{tikzpicture}}; \endxy} \;\; }

\DeclareDocumentCommand{\vUcrossum}{ O{black} O{black} O{semithick} O{.5}}{\;\;
    \vcenter{\xy (0,0)*{\begin{tikzpicture}[scale=#4]
  \draw[#3,color=#1, ->] (-0.5,0) .. controls (-0.5,0.5) and (0.5,0.5) .. (0.5,1)
      node[pos=0.5, shape=coordinate](X){};
  \draw[#3,color=#2, ->] (0.5,0) .. controls (0.5,0.5) and (-0.5,0.5) .. (-0.5,1);
\end{tikzpicture}}; \endxy} \;\; }

\DeclareDocumentCommand{\Ucrossdm}{ O{black} O{black} O{semithick}}{\;\;
    \vcenter{\xy (0,0)*{\begin{tikzpicture}[scale=0.75]
  \draw[#3,color=#1, <-] (-0.5,0) .. controls (-0.5,0.5) and (0.5,0.5) .. (0.5,1)
      node[pos=0.5, shape=coordinate](X){};
  \draw[#3,color=#2, <-] (0.5,0) .. controls (0.5,0.5) and (-0.5,0.5) .. (-0.5,1);
\end{tikzpicture}}; \endxy} \;\; }

\begin{document}
%

\title[The $\mf{sl}_n$ foam $2$-category]
{The $\mf{sl}_n$ foam $2$-category: a combinatorial formulation of Khovanov-Rozansky
homology via categorical skew Howe duality}

\author{Hoel Queffelec}
\address{Univ Paris Diderot, Sorbonne Paris Cit\'e, IMJ-PRG, UMR 7586 CNRS, Sorbonne Universit\'e, UPMC Univ Paris 06, F-75013, Paris, France}
\email{hoel.queffelec@imj-prg.fr}

\author{David E. V. Rose}
\address{Department of Mathematics, University of Southern California, Los Angeles, CA 90089, USA}
\email{davidero@usc.edu}

\begin{abstract}
We give a purely combinatorial construction of colored $\sln$ link homology. 
The invariant takes values in a $2$-category where $2$-morphisms 
are given by foams, singular cobordisms between $\mf{sl}_n$ webs; applying a (TQFT-like) 
representable functor recovers (colored) Khovanov-Rozansky homology.
Novel features of the theory include the introduction of `enhanced' foam facets which fix sign issues associated with the original matrix factorization 
formulation and the use of skew Howe duality to show that (enhanced) closed foams can be evaluated in a completely combinatorial manner. 
The latter answers a question posed in \cite{MSV}.
\end{abstract}

\maketitle
\tableofcontents

%
\section{Introduction}
%

The usefulness of Khovanov homology \cite{Kh1} stems not only from its power as a topological invariant, but also 
from its relatively simple description. Indeed, one can understand Khovanov homology as assigning to a link 
in $S^3$ a (co)chain complex in which the `chain groups' are collections of disjoint circles in the plane and 
whose morphisms are given by saddle cobordisms. A simple operation assigning vector spaces to circles and 
linear maps to cobordisms (a TQFT) then gives a complex of vector spaces and the homology of this complex is the link invariant.
Despite the simplicity of this description, the link invariant is quite powerful; for example, Rasmussen uses Khovanov 
homology to give an elegant proof of the Milnor conjecture~\cite{Ras}, 
Kronheimer and Mrowka show that it detects the unknot \cite{KrMr}, and
Grigsby and Ni show that (a variant of) Khovanov homology distinguishes braids from other tangles \cite{GriNi}.

In the years since Khovanov introduced his $\mathfrak{sl}_2$ link homology, various authors have introduced 
link homology theories extending the above construction to other simple Lie algebras. Khovanov himself 
gave a construction of $\mathfrak{sl}_3$ link homology in the spirit of the invariant mentioned above \cite{Kh5}, and, in 
joint work with Rozansky, later gave a construction based on the theory of matrix factorizations \cite{KhR} which
generalizes both of the above approaches to $\mathfrak{sl}_n$ for all $n \geq 2$. Independently, Cautis and 
Kamnitzer \cite{CK01,CK02} construct $\mathfrak{sl}_n$ link homology using derived categories of coherent sheaves on 
orbits in the affine Grassmannian, Mazorchuk and Stroppel \cite{MS2} and Sussan \cite{Sussan} give 
constructions based on category $\mathcal{O}$, and Webster \cite{Webster} defines link homology associated to 
any simple Lie algebra using derived categories of diagrammatically defined algebras inspired by the 
geometry of quiver varieties.

Of the above approaches, only Khovanov's $\mathfrak{sl}_3$ link homology admits a simple, purely 
combinatorial description. This invariant assigns to a link a complex of trivalent graphs, called webs, with 
morphisms given by foams, singular cobordisms between such graphs. Again, one can pass from this complex to one 
consisting of vector spaces and linear maps and compute homology to obtain the link invariant. 

In this work, we extend Khovanov's $\mathfrak{sl}_3$ framework to give an elementary, combinatorial 
construction of (colored) $\mathfrak{sl}_n$ link homology in terms of $\sln$ webs and foams. This approach builds on earlier 
work of Mackaay-Stosic-Vaz \cite{MSV} in which a webs and foams approach was utilized; however, in that setup the 
link invariant was not shown to be purely combinatorial in nature, rather relying on the Kapustin-Li formula \cite{KapLi}, \cite{KhR3} from 
topological Landau-Ginzburg models to evaluate closed foams, the foam analogs of closed surfaces. 
We extend Blanchet's `enhanced' $\slnn{2}$ foams \cite{Blan} to the $\sln$ setting to give a simple generators and relations presentation 
which allows for an elementary description of $\mathfrak{sl}_n$ link homology and for the evaluation of any 
(enhanced) closed foam.

Our main tool is the higher representation theory of the Khovanov-Lauda categorified quantum group \cite{KL,KL2,KL3} 
(see also work of Rouquier \cite{Rou2} for an independent construction), 
in particular, the categorical skew Howe duality of Cautis-Kamnitzer-Licata \cite{CKL}. In that paper, and the
subsequent work of Cautis \cite{Cautis}, the authors show that one can obtain the braiding in the algebro-geometric 
formulation of $\mathfrak{sl}_n$ link homology as the image of Rickard complexes \cite{CR}. 
We adopt this approach for our construction of $\mathfrak{sl}_n$ link homology, also taking inspiration from diagrammatic
realization of skew Howe duality of Cautis-Kamnitzer-Morrison \cite{CKM}. In that work, the authors 
give a graphical depiction of skew Howe duality at the decategorified level, finding a generators and relations description 
for the category of representations of quantum $\mathfrak{sl}_n$. We extend this result to the $2$-categorical level, using the $2$-morphisms 
in the categorified quantum group to describe $\mathfrak{sl}_n$ foams, generalizing earlier work on the $\mathfrak{sl}_2$ and 
$\mathfrak{sl}_3$ case by the authors, joint with Lauda \cite{LQR1}. 

The following are immediate consequences of this description:
\begin{itemize}
\item A definition of colored $\sln$ link homology over the integers, which is conjecturally (properly) functorial with respect to link cobordism. 
\item A proof that (enhanced) closed foams can be evaluated combinatorially, and an algorithm for performing this evaluation. This 
can be viewed as a combinatorial analog of the Kapustin-Li formula for (enhanced) $\sln$ matrix factorizations.
\item A topological description of the (direct limit of the) categorified $q$-Schur algebra.
\end{itemize}

In Section \ref{sec:background}, we recall the requisite background information 
on link homology, skew Howe duality, and categorified 
quantum groups. Section \ref{sec:foams} contains our construction of the foam $2$-category, and in Section \ref{sec:link_hom} we define the 
link invariant and explore some consequences of our construction.

\textbf{Acknowledgements:} 
We'd like to thank Christian Blanchet, Sabin Cautis, Ben Cooper, Matt Hogancamp, Mikhail Khovanov, Tony Licata, 
Scott Morrison, Lev Rozansky, and Hao Wu for helpful discussions, 
and Daniel Tubbenhauer for his careful reading of a preliminary draft of this paper.
Both authors are indebted to Aaron Lauda for his continued support and guidance; indeed, this paper grew out of our collaboration with him. 
DR was supported by a Zumberge Fund Individual Grant Award and the John Templeton Foundation. 
%
\section{Background}\label{sec:background}
%

In this section, we discuss background information concerning link homology and (categorified) quantum groups.

\subsection{Link homology and foams} 

In \cite{Kh1}, Khovanov introduced a homology theory for links in $S^3$, now called Khovanov homology, categorifying 
the Jones polynomial. This invariant was later recast by Bar-Natan \cite{BN2} in terms of the category of tangles and 
cobordisms (see also \cite{Kh2}). In this description, the invariant assigns to a link in $S^3$, or more generally to a tangle in $\R^2 \times [0,1]$, 
a (co)chain complex in which the objects are (graded, formal direct sums of) planar tangles and whose morphisms are 
(formal matrices of linear combinations of degree zero) cobordisms, modulo local relations on such cobordisms. 
This depiction is entirely diagrammatic and combinatorial; however, the resulting invariant is extremely powerful.

In more sophisticated terms, Bar-Natan's construction can be understood in the context of the $2$-category of planar 
tangles and cobordisms. 
It will be our convention throughout that morphisms in all diagrammatically defined $2$-categories map from right to left and bottom to top. 
The objects in this $2$-category are non-negative integers, which should be viewed as specifying 
a number of boundary points. The $1$-morphisms are $\Z$-graded formal direct sums of planar tangles having appropriate 
boundary, e.g. a tangle with domain $m$ and codomain $n$ has $m$ right boundary points and $n$ left boundary points. Finally, 
the $2$-morphisms are given by matrices of linear combinations of degree zero dotted cobordisms between the relevant tangles, 
modulo the following local relations:
\begin{equation}\label{sl2closedfoam}
 \sphere[.5] \quad = \quad 0 \qquad , \qquad
\dottedsphere[.5]{}  \quad = \quad 1
\end{equation}
\begin{equation}\label{sl2neckcutting}
\cylinder[.4] \quad = \quad \slthncfour[.4] \quad + \quad \slthncfive[.4] \quad .
\end{equation}
The degree of a cobordism $C: q^{t_1}T_1 \to q^{t_2}T_2$ is given by the formula
\begin{equation}\label{sl2foamdeg}
\deg(C) = \chi(C) -2\#D - \frac{\# \partial}{2} + t_2 - t_1
\end{equation}
where $\chi$ is the Euler characteristic of the underlying surface, 
$\# D$ is the number of dots, $\# \partial$ is the number of boundary points in either $T_1$ or $T_2$, and the
powers of the formal variable $q$ denote degree.

The skein relations
\begin{equation}\label{sl2skein}
\left \llbracket
\xy
(0,0)*{
\begin{tikzpicture} [scale=.5]
\draw[very thick, directed=.99] (0,0) to (1,1);
\draw[very thick, directed=.99] (.4,.6) to (0,1);
\draw[very thick] (1,0) to (.6,.4);
\end{tikzpicture}
}
\endxy
\right \rrbracket
=
\left(
q \;
\uwave{
\xy
(0,0)*{
\begin{tikzpicture} [scale=.5]
\draw[very thick] (0,0) to [out=45,in=315] (0,1);
\draw[very thick] (1,0) to [out=135,in=225] (1,1);
\end{tikzpicture}
}
\endxy}
\stackrel{s}{\longrightarrow}
q^2 \;
\xy
(0,0)*{
\begin{tikzpicture} [scale=.5]
\draw[very thick] (0,0) to [out=45,in=135] (1,0);
\draw[very thick] (0,1) to [out=315,in=225] (1,1);
\end{tikzpicture}
}
\endxy \;
\right) 
\qquad , \qquad
\left \llbracket
\xy
(0,0)*{
\begin{tikzpicture} [scale=.5]
\draw[very thick, directed=.99] (1,0) to (0,1);
\draw[very thick, directed=.99] (.6,.6) to (1,1);
\draw[very thick] (0,0) to (.4,.4);
\end{tikzpicture}
}
\endxy
\right \rrbracket
=
\left(
q^{-2} \;
\xy
(0,0)*{
\begin{tikzpicture} [scale=.5]
\draw[very thick] (0,0) to [out=45,in=135] (1,0);
\draw[very thick] (0,1) to [out=315,in=225] (1,1);
\end{tikzpicture}
}
\endxy
\stackrel{s}{\longrightarrow}
q^{-1} \;
\uwave{
\xy
(0,0)*{
\begin{tikzpicture} [scale=.5]
\draw[very thick] (0,0) to [out=45,in=315] (0,1);
\draw[very thick] (1,0) to [out=135,in=225] (1,1);
\end{tikzpicture}
}
\endxy} \;
\right)
\end{equation}
assign complexes in the above $2$-category\footnote{To be precise, these are complexes in the $\Hom$-category
$\Hom(2,2)$.} 
to positive and negative crossings, and extend in a manner 
similar to tensor product of complexes of abelian groups to assign a complex $\llbracket T \rrbracket$ to any tangle $T$.
In these formulae, $s$ denotes a saddle cobordism and we've \uwave{underlined} the term in homological degree zero.

Given a link $L$, we can apply the representable functor $\bigoplus_{k \in \Z} \Hom(q^{-k} \varnothing,-)$ of morphisms from 
degree-shifts of the empty tangle to the complex $\llbracket L \rrbracket$ to 
obtain a complex of graded vector spaces, with degree-zero differentials. 
The homology of this complex recovers Khovanov's original invariant.

Khovanov's $\mathfrak{sl}_3$ link invariant \cite{Kh5} admits a similar description \cite{MV2},  \cite{MorrisonNieh}. 
Here, the $2$-category has sequences of $1$'s and $2$'s as objects, $1$-morphisms are $\Z$-graded formal direct sums of 
$\mathfrak{sl}_3$ webs - oriented 
trivalent graphs with boundary, in which every trivalent vertex is a `sink' or a `source' and edges are directed out from $1$'s and 
into $2$'s in the domain, and vice versa in the domain \cite{Kup}. 
The $2$-morphisms are given by matrices of linear combinations of degree-zero
foams - dotted, singular cobordisms between such graphs in which the neighborhood of a singular seam is locally modeled as the product 
of the letter `$Y$' with an interval. Again, certain local relations are placed on foams, e.g.
\begin{equation}\label{sl3tube}
\tubeLHS[.5] \quad = \quad \tubeRHStopdot[.5] \quad - \quad \tubeRHSbottomdot[.5]
\end{equation}
is a relation in the $\Hom$-category $\Hom(1,1)$; see \cite{MorrisonNieh} or \cite{LQR1} for the full collection of relations.

There are again skein relations:
\begin{align}\label{sl3skein}
\left \llbracket
\xy
(0,0)*{
\begin{tikzpicture} [scale=.5]
\draw[very thick, directed=.99] (0,0) to (1,1);
\draw[very thick, directed=.99] (.4,.6) to (0,1);
\draw[very thick] (1,0) to (.6,.4);
\end{tikzpicture}
}
\endxy
\right \rrbracket
=
\left(
q^2 \;
\uwave{
\xy
(0,0)*{
\begin{tikzpicture} [scale=.6]
\draw[very thick, directed=.99] (0,0) to [out=45,in=315] (0,1);
\draw[very thick, directed=.99] (1,0) to [out=135,in=225] (1,1);
\end{tikzpicture}
}
\endxy}
\stackrel{z}{\longrightarrow}
q^3 \;
\xy
(0,0)*{
\begin{tikzpicture} [scale=.6]
\draw[very thick, directed=.85] (.5,.75) to (1,1);
\draw[very thick, directed=.85] (.5,.75) to (0,1);
\draw[very thick, directed=.7] (.5,.75) to (.5,.25);
\draw[very thick, directed=.55] (1,0) to (.5,.25);
\draw[very thick, directed=.55] (0,0) to (.5,.25);
\end{tikzpicture}
}
\endxy \;
\right) 
\quad , \quad
\left \llbracket
\xy
(0,0)*{
\begin{tikzpicture} [scale=.5]
\draw[very thick, directed=.99] (1,0) to (0,1);
\draw[very thick, directed=.99] (.6,.6) to (1,1);
\draw[very thick] (0,0) to (.4,.4);
\end{tikzpicture}
}
\endxy
\right \rrbracket
=
\left(
q^{-3} \;
\xy
(0,0)*{
\begin{tikzpicture} [scale=.6]
\draw[very thick, directed=.85] (.5,.75) to (1,1);
\draw[very thick, directed=.85] (.5,.75) to (0,1);
\draw[very thick, directed=.7] (.5,.75) to (.5,.25);
\draw[very thick, directed=.55] (1,0) to (.5,.25);
\draw[very thick, directed=.55] (0,0) to (.5,.25);
\end{tikzpicture}
}
\endxy
\stackrel{u}{\longrightarrow}
q^{-2} \;
\uwave{
\xy
(0,0)*{
\begin{tikzpicture} [scale=.6]
\draw[very thick, directed=.99] (0,0) to [out=45,in=315] (0,1);
\draw[very thick, directed=.99] (1,0) to [out=135,in=225] (1,1);
\end{tikzpicture}
}
\endxy} \;
\right)
\end{align}
where
$
z =
\xy
(0,0)*{
\begin{tikzpicture} [scale=.37,fill opacity=0.2]
	\path [fill=red] (4.25,2) to (4.25,-.5) to [out=165,in=15] (-.5,-.5) to (-.5,2) to
		[out=0,in=225] (.75,2.5) to [out=270,in=180] (1.625,1.25) to [out=0,in=270] 
			(2.5,2.5) to [out=315,in=180] (4.25,2);
	\path [fill=red] (3.75,3) to (3.75,.5) to [out=195,in=345] (-1,.5) to (-1,3) to [out=0,in=135]
		(.75,2.5) to [out=270,in=180] (1.625,1.25) to [out=0,in=270] 
			(2.5,2.5) to [out=45,in=180] (3.75,3);
	\path[fill=blue] (2.5,2.5) to [out=270,in=0] (1.625,1.25) to [out=180,in=270] (.75,2.5);
	\draw [very thick,directed=.75] (4.25,-.5) to [out=165,in=15] (-.5,-.5);
	\draw [very thick, directed=.75] (3.75,.5) to [out=195,in=345] (-1,.5);
	\draw [very thick, red, directed=.75] (2.5,2.5) to [out=270,in=0] (1.625,1.25);
	\draw [very thick, red] (1.625,1.25) to [out=180,in=270] (.75,2.5);
	\draw [very thick] (3.75,3) to (3.75,.5);
	\draw [very thick] (4.25,2) to (4.25,-.5);
	\draw [very thick] (-1,3) to (-1,.5);
	\draw [very thick] (-.5,2) to (-.5,-.5);
	\draw [very thick,rdirected=.55] (2.5,2.5) to (.75,2.5);
	\draw [very thick,directed=.55] (.75,2.5) to [out=135,in=0] (-1,3);
	\draw [very thick,directed=.55] (.75,2.5) to [out=225,in=0] (-.5,2);
	\draw [very thick,directed=.55] (3.75,3) to [out=180,in=45] (2.5,2.5);
	\draw [very thick,directed=.55] (4.25,2) to [out=180,in=315] (2.5,2.5);		
\end{tikzpicture}
};
\endxy
$ and
$
u=
\xy
(0,0)*{
\begin{tikzpicture} [scale=.37,fill opacity=0.2]
	\path [fill=red] (4.25,-.5) to (4.25,2) to [out=165,in=15] (-.5,2) to (-.5,-.5) to 
		[out=0,in=225] (.75,0) to [out=90,in=180] (1.625,1.25) to [out=0,in=90] 
			(2.5,0) to [out=315,in=180] (4.25,-.5);
	\path [fill=red] (3.75,.5) to (3.75,3) to [out=195,in=345] (-1,3) to (-1,.5) to 
		[out=0,in=135] (.75,0) to [out=90,in=180] (1.625,1.25) to [out=0,in=90] 
			(2.5,0) to [out=45,in=180] (3.75,.5);
	\path[fill=blue] (.75,0) to [out=90,in=180] (1.625,1.25) to [out=0,in=90] (2.5,0);
	\draw [very thick,rdirected=.55] (2.5,0) to (.75,0);
	\draw [very thick,directed=.55] (.75,0) to [out=135,in=0] (-1,.5);
	\draw [very thick,directed=.55] (.75,0) to [out=225,in=0] (-.5,-.5);
	\draw [very thick,directed=.55] (3.75,.5) to [out=180,in=45] (2.5,0);
	\draw [very thick,directed=.55] (4.25,-.5) to [out=180,in=315] (2.5,0);
	\draw [very thick, red, directed=.75] (.75,0) to [out=90,in=180] (1.625,1.25);
	\draw [very thick, red] (1.625,1.25) to [out=0,in=90] (2.5,0);
	\draw [very thick] (3.75,3) to (3.75,.5);
	\draw [very thick] (4.25,2) to (4.25,-.5);
	\draw [very thick] (-1,3) to (-1,.5);
	\draw [very thick] (-.5,2) to (-.5,-.5);
	\draw [very thick,directed=.75] (4.25,2) to [out=165,in=15] (-.5,2);
	\draw [very thick, directed=.75] (3.75,3) to [out=195,in=345] (-1,3);	
\end{tikzpicture}
};
\endxy
$ ,
which assign a complex in this $2$-category to each tangle. The link homology can then be 
recovered as above by applying an appropriate representable functor to the complex $\llbracket L \rrbracket$
assigned to a link, and taking the homology of the complex of graded vector spaces.

There does not exist a purely combinatorial analog of the above constructions for $\mathfrak{sl}_n$ link homology. This invariant 
was originally described by Khovanov and Rozansky \cite{KhR} using the theory of matrix factorizations. In \cite{MSV}, the authors 
give a webs and foams interpretation of Khovanov-Rozanksy homology; however, they are unable to prove that the link invariant 
can be computed using 
their relations on foams. Rather, they appeal to the Kapustin-Li formula from topological Landau-Ginzburg models \cite{KapLi}, \cite{KhR3} to compute 
the values of `closed foams' - endomorphisms of the empty web needed to explicitly compute the differentials in the complex assigned to 
a link. Additionally, their construction does not describe the colored Khovanov-Rozanksy homology of Yonezawa \cite{Yon} and 
Wu \cite{Wu}. 
The formulation of $\mathfrak{sl}_n$ foams in this paper resolves these issues.

\subsection{Quantum groups, webs, and skew Howe duality}

The planar tangles and $\mathfrak{sl}_3$ webs in the previous section can be interpreted as graphical depictions of intertwiners between 
representations of the quantum group $U_q(\mathfrak{sl}_n)$ for $n=2,3$, see \cite{Kup} and references therein. 
In recent work of Cautis-Kamnitzer-Morrison \cite{CKM}, this was generalized to give a graphical description of the representation theory 
of quantum $\mathfrak{sl}_n$. Indeed, they show that the full subcategory of $\cat{Rep}(U_q(\mathfrak{sl}_n))$ generated (as a pivotal category) by 
the fundamental representations can be presented as follows:
\begin{itemize}
\item \textbf{Objects} are given by sequences of points labeled by the symbols $1^{\pm}, \ldots, (n-1)^{\pm}$, together with a 
zero object.

\item \textbf{Morphisms} are given by $\C(q)$-linear combinations of $\mathfrak{sl}_n$ webs - directed, labeled trivalent graphs with boundary, 
modulo planar isotopy, generated by the following:
\begin{equation}\label{webgen}
\xy
(0,0)*{
\begin{tikzpicture} [scale=.75, decoration={markings,
                        mark=at position 0.6 with {\arrow{>}};    }]
\draw[very thick, postaction={decorate}] (0,0) -- (0,1);
\node at (0,1.3) {\small$a+b$};
\draw[very thick, postaction={decorate}] (-.875,-.5) -- (0,0);
\node at (-1.2,-.5) {\small$a$};
\draw[very thick, postaction={decorate}] (.875,-.5) -- (0,0);
\node at (1.2,-.5) {\small$b$};
\end{tikzpicture}};
\endxy
\quad , \quad
\xy
(0,0)*{
\begin{tikzpicture} [scale=.75, decoration={markings,
                        mark=at position 0.6 with {\arrow{>}};    }]
\draw[very thick, postaction={decorate}] (0,1) -- (0,0);
\node at (0,1.3) {\small$a+b$};
\draw[very thick, postaction={decorate}] (0,0) -- (.875,-.5);
\node at (-1.2,-.5) {\small$a$};
\draw[very thick, postaction={decorate}] (0,0) -- (-.875,-.5);
\node at (1.2,-.5) {\small$b$};
\end{tikzpicture}};
\endxy
\quad , \quad
\xy
(0,0)*{
\begin{tikzpicture} [scale=.75, decoration={markings,
                        mark=at position 0.6 with {\arrow{>}};    }]
\draw[very thick, postaction={decorate}] (0,-.5) -- (0,.25);
\node at (0,-.8) {\small$a$};
\draw[very thick, postaction={decorate}] (0,1) -- (0,.25);
\node at (0,1.3) {\small$n-a$};
\draw[very thick] (0,.25) -- (.2,.25);
\end{tikzpicture}};
\endxy
\quad , \quad
\xy
(0,0)*{
\begin{tikzpicture} [scale=.75, decoration={markings,
                        mark=at position 0.6 with {\arrow{>}};    }]
\draw[very thick, postaction={decorate}] (0,.25) -- (0,-.5);
\node at (0,-.8) {\small$a$};
\draw[very thick, postaction={decorate}] (0,.25) -- (0,1);
\node at (0,1.3) {\small$n-a$};
\draw[very thick] (0,.25) -- (.2,.25);
\end{tikzpicture}};
\endxy
\end{equation}
modulo relations. Web edges labeled $a$ are directed into points labeled $a^+$ and out from points labeled 
$a^-$ in the codomain, and vice versa in the domain.
\end{itemize}

We'll denote this category by $\nWeb$ for the duration.
If $\C_q^n$ denotes the defining representation of $U_q(\mathfrak{sl}_n)$, the correspondence to the category of representations is 
given by associating the representation\footnote{See e.g. \cite{CKM} for the definition of the quantum exterior power.} 
$\wedge_q^a \C_q^n$ to the label $a^+$, its dual to $a^-$, 
and tensor products of such representations to a sequence of labeled points. 
The first and second graphs in \eqref{webgen} correspond to the unique (up to scalar multiple) morphisms 
in $\Hom(\wedge_q^a \C_q^n \otimes \wedge_q^b \C_q^n, \wedge_q^{a+b} \C_q^n)$ and 
$\Hom(\wedge_q^{a+b} \C_q^n, \wedge_q^a \C_q^n \otimes \wedge_q^b \C_q^n)$ 
respectively, while the latter two morphisms specify the non-canonical isomorphisms 
$\wedge_q^a \C_q^n \cong (\wedge_q^{n-a} \C_q^n)^*$. The planar tangles and $\mathfrak{sl}_3$ webs from the previous 
section are related to the above construction by choosing identifications $\C_q^2 \cong (\C_q^2)^*$ and 
$\wedge_q^2 \C_q^3 \cong (\C_q^3)^*$, eliminating the need for tag morphisms.

The relations in this category are obtained in~\cite{CKM} by adopting the geometric skew Howe duality of 
Cautis-Kamnitzer-Licata~\cite{CKL} to the diagrammatic setting. We now summarize their construction. 
The vector space $\bigwedge_q^N(\C_q^n \otimes \C_q^m)$ carries commuting actions of $U_q(\sln)$ and $U_q(\slm)$, and
the isomorphism
\begin{equation}\label{sHeq}
\textstyle \bigwedge_q^N(\C_q^n \otimes \C_q^m) \cong \bigwedge_q^N(\C_q^n \oplus \cdots \oplus \C_q^n) 
\cong \displaystyle \bigoplus_{\sum a_i = N} \wedge_q^{a_1} \C_q^n \otimes \cdots \otimes \wedge_q^{a_m} \C_q^n
\end{equation}
gives the decomposition into $\slm$ weight spaces, with $\slm$ weights $(a_1-a_2, \ldots, a_{m-1}-a_{m})$. 
The weight spaces are themselves 
tensor products of fundamental representations of $U_q(\sln)$ and the actions of $E_i$ and $F_i$ in $U_q(\slm)$ give intertwiners between 
these $\sln$ representations, which generate the full subcategory of $\cat{Rep}(U_q(\mathfrak{sl}_n))$ mentioned above. In fact, 
$U_q(\sln)$ and $U_q(\slm)$ give a Howe pair for the representation in equation \eqref{sHeq}, so all relations between these 
generating morphisms follow as consequences of the relations in $U_q(\slm)$.

To connect to the diagrammatics, recall that the Lusztig idempotent form $\dot{U}_q(\slm)$ of quantum $\slm$ can be viewed as a 
category in which objects are given by $\slm$ weights and morphisms are given by elements in $U_q(\slm)$.
Representations of $U_q(\slm)$ (with weight space decompositions) correspond to functors from $\dot{U}_q(\slm)$ to the category of 
vector spaces, so equation \eqref{sHeq} gives a functor $\phi_n^{m,N}: \dot{U}_q(\slm) \to \cat{Rep}(U_q(\mathfrak{sl}_n))$. 
In the diagrammatic description above, an $\slm$ weight $\lambda=(\lambda_1,\ldots,\lambda_{m-1})$ is sent to the sequence of labels 
$[a_1,\ldots,a_m]$ with $\sum a_i = N$ and $a_i-a_{i+1} = \lambda_i$. If no such solution exists, 
or if for some $i$ we have $a_i<0$ or $a_i>n$, then the weight is sent to the zero object. 
The functor is given by associating the following webs to the generating morphisms $1_\l, E_i 1_\l, F_i 1_\l $:
\[
1_\l 
\;\; \mapsto \;\;
\xy
(0,0)*{
\begin{tikzpicture} [scale=.75]
\draw [very thick, directed=.55] (3,0) -- (0,0);
\draw [very thick, directed=.55] (3,1) -- (0,1);
\draw [very thick, directed=.55] (3,1.5) -- (0,1.5);
\node at (3.6,0) {\small$a_m$};
\node at (3.6,1) {\small$a_2$};
\node at (3.6,1.5) {\small$a_1$};
\node at (1.5,.65) {$\vdots$};
\end{tikzpicture}};
\endxy
\]
\begin{equation}\label{eq_WsH}
E_i 1_\l 
\;\; \mapsto \;\;
\xy
(0,0)*{
\begin{tikzpicture} [scale=.75]
\draw [very thick, directed=.55] (3,0) -- (2,0);
\draw[very thick, directed=.55] (2,0) -- (0,0);
\draw [very thick, directed=.55] (3,1) -- (1,1);
\draw[very thick, directed=.55] (1,1) -- (0,1);
\draw [very thick, directed=.55] (2,0) -- (1,1);
\node at (3.6,0) {\small$a_{i+1}$};
\node at (3.6,1) {\small$a_i$};
\node at (-1,0) {\small$a_{i+1} - 1$};
\node at (-1,1) {\small$a_i+1$};
\node at (2,.5) {\tiny$1$};
\end{tikzpicture}};
\endxy
\end{equation}
\[
F_i 1_\l 
\;\; \mapsto \;\;
\xy
(0,0)*{
\begin{tikzpicture} [scale=.75]
\draw [very thick, directed=.55] (3,0) -- (1,0);
\draw[very thick, directed=.55] (1,0) -- (0,0);
\draw [very thick, directed=.55] (3,1) -- (2,1);
\draw[very thick, directed=.55] (2,1) -- (0,1);
\draw [very thick, directed=.55] (2,1) -- (1,0);
\node at (3.6,0) {\small$a_{i+1}$};
\node at (3.6,1) {\small$a_i$};
\node at (-1,0) {\small$a_{i+1} + 1$};
\node at (-1,1) {\small$a_i - 1$};
\node at (2,.5) {\tiny$1$};
\end{tikzpicture}};
\endxy
\]
where we read such diagrams\footnote{Note that our conventions for the diagrammatics of the skew Howe functor differ from 
those used in both \cite{CKM} and \cite{LQR1}. Indeed, our conventions in \cite{LQR1} were chosen to better match the diagrammatics of 
categorified quantum $\slm$, and corresponded to an unconventional choice of basis for $\C_q^m$. Our new conventions are chosen 
to match with both the categorified quantum group and the standard choice of basis for $\C_q^m$, and can be obtained from those 
in \cite{CKM} by turning their diagrams sideways and rotating about a horizontal axis.} 
from right to left. 
Edges labeled zero are erased and those labeled $n$ are truncated to the `tag' morphisms in equation \eqref{webgen}. 
If a label lies outside the range $\{0,\ldots,n\}$, then the morphism is zero, as it factors 
through the zero object. Adopting the terminology from \cite{CKM}, we'll refer to compositions of the above webs as \emph{ladder webs}.

As an example for how the web relations can be deduced from skew Howe duality, consider the $U_q(\mathfrak{sl}_2)$ relation 
$EF1_1 - FE1_1 = \id_1$. Under the functor $\phi_3^{2,3}$ this corresponds to the $\mathfrak{sl}_3$ web relation
\[
\xy
(0,0)*{
\begin{tikzpicture} [scale=.6]
\draw [very thick, directed= .15, directed=.5, directed=.85] (2.5,0) -- (-1.5,0);
\draw [very thick, directed= .1, directed=.5, directed=.95] (2.5,1) -- (-1.5,1);
\draw [very thick, directed=.6] (1.75,1) -- (1,0);
\draw[very thick, directed=.6] (0,0) to (-.75,1);
\node at (3,0) {\small$1$};
\node at (3,1) {\small$2$};
\node at (-2,0) {\small$1$};
\node at (-2,1) {\small$2$};
\node at (1.75,.5) {\tiny$1$};
\node at (-.75,.5) {\tiny$1$};
\end{tikzpicture}};
\endxy
\;\; - \;\;
\xy
(0,0)*{
\begin{tikzpicture} [scale=.6]
\draw [very thick] (2.5,0) -- (1.75,0);
\draw [very thick] (-.75,0) -- (-1.5,0);
\draw [very thick, directed=.5] (2.5,1) -- (.8,1);
\draw [very thick, directed=.5] (.2,1) -- (-1.5,1);
\draw [very thick, directed=.6] (1.75,0) -- (1,1);
\draw[very thick, directed=.6] (0,1) to (-.75,0);
\node at (3,0) {\small$1$};
\node at (3,1) {\small$2$};
\node at (-2,0) {\small$1$};
\node at (-2,1) {\small$2$};
\end{tikzpicture}};
\endxy
\;\; = \;\;
\xy
(0,0)*{
\begin{tikzpicture} [scale=.6]
\draw [very thick,directed=.55] (2.5,0) -- (-1.5,0);
\draw [very thick, directed= .55] (2.5,1) -- (-1.5,1);
\node at (3,0) {\small$1$};
\node at (3,1) {\small$2$};
\node at (-2,0) {\small$1$};
\node at (-2,1) {\small$2$};
\end{tikzpicture}};
\endxy
\]
i.e. the `square' relation from \cite{Kup}. The complete set of relations for $\sln$ webs is obtained in a similar fashion 
in \cite{CKM}.

The braidings giving quantum knot invariants can also be obtained in the above setting. Recall that the quantum Weyl group elements
\begin{equation} \label{eq_qW}
 T_i 1_{\lambda} = \sum_{s \geq 0} (-q)^s F_i^{(\lambda+s)}E_i^{(s)}1_{\lambda} \quad \text{$\lambda\geq 0$}, \qquad
 T_i 1_{\lambda} =\sum_{s\geq 0} (-q)^s  E_i^{(-\lambda+s)}F_i^{(s)}1_{\lambda} \quad \text{$\lambda\leq 0$}
\end{equation}
generate an action of the braid group on any finite-dimensional representation \cite{Lus4}. 
In \cite{CKL}, the authors show that the $\slm$ quantum Weyl group elements correspond under skew Howe duality 
to the braiding on $U_q(\sln)$ representations given by the $R$-matrix. For example, we compute that
\[
\phi_2^{2,2}(T 1_0) = 
\xy
(0,0)*{
\begin{tikzpicture} [scale=.6]
\draw [very thick,directed=.55] (1.5,0) -- (-1.5,0);
\draw [very thick, directed= .55] (1.5,1) -- (-1.5,1);
\node at (2,0) {\small$1$};
\node at (2,1) {\small$1$};
\node at (-2,0) {\small$1$};
\node at (-2,1) {\small$1$};
\end{tikzpicture}};
\endxy
\;\; - \;\; q \;\;
\xy
(0,0)*{
\begin{tikzpicture} [scale=.6]
\draw [very thick] (2.5,0) -- (1.75,0);
\draw [very thick] (-.75,0) -- (-1.5,0);
\draw [very thick, directed=.5] (2.5,1) -- (.8,1);
\draw [very thick, directed=.5] (.2,1) -- (-1.5,1);
\draw [very thick, directed=.6] (1.75,0) -- (1,1);
\draw[very thick, directed=.6] (0,1) to (-.75,0);
\node at (3,0) {\small$1$};
\node at (3,1) {\small$1$};
\node at (-2,0) {\small$1$};
\node at (-2,1) {\small$1$};
\end{tikzpicture}};
\endxy
\]
which corresponds to the skein relation for the Kauffman bracket of a positive crossing 
$\;\; \xy
(0,0)*{
\begin{tikzpicture} [scale=.5]
\draw[very thick, directed=.99] (1,0) to (0,1);
\draw[very thick] (.6,.6) to (1,1);
\draw[very thick, directed=.99] (.4,.4) to (0,0);
\end{tikzpicture}
}
\endxy \;\;$, which gives the Jones polynomial.

Our approach to $\sln$ link invariants is to categorify the above description. Playing the role of $\dot{U}_q(\slm)$ is the Khovanov-Lauda 
categorified quantum group (reviewed in the next section), a $2$-category which provides a diagrammatic description of morphisms 
between elements of $\dot{U}_q(\slm)$. 
These $2$-morphisms suggest a set of generating foams, and by insisting that a $2$-functor exists which `covers' the functor 
$\phi_n^{m,N}$, we are able to deduce a complete set of relations for $\sln$ foams. Such a $2$-functor was originally conjectured to exist by 
Khovanov and was constructed by Mackaay in the $\mathfrak{sl}_3$ case working in the restricted setting of $\Z/2\Z$ coefficients, 
before the connection to skew Howe duality was discovered. Such a $2$-functor was constructed (over arbitrary coefficients) in the 
$\mathfrak{sl}_2$ and $\mathfrak{sl}_3$ setting by the authors, joint with Lauda \cite{LQR1}, and in the $\mathfrak{sl}_3$ 
case independently by Mackaay-Pan-Tubbenhauer \cite{MPT}. 

One should note the difference between the approach taken in \cite{LQR1} and that in the present paper. 
In the former, we determined foamation functors by checking that all relations in the categorified quantum group 
map to relations in the $\slnn{2}$ and $\slnn{3}$ foam categories. 
This strategy cannot be applied in the present context since no generators and relations construction of $\sln$ foams exists.
Instead, we proceed in the same spirit as the authors did at the uncategorified level in \cite{CKM}, using relations in the categorified 
quantum group as a blueprint to give a generators and relations construction for $\sln$ foams. We then check that 
the desired foamation functor exists, and use this to study properties of foams.

In order to categorify the skein relations and obtain a categorified invariant of links, we follow \cite{CKL3}, \cite{Cautis} and utilize 
Rickard complexes \cite{CR}, complexes in the categorified quantum group categorifying the quantum Weyl group elements and which 
generate categorical braid group actions on (integrable) $2$-representations \cite{CK}. In our previous work \cite{LQR1}, 
we showed that these complexes give the categorified skein relations for the diagrammatic formulations of 
$\mathfrak{sl}_2$ and $\mathfrak{sl}_3$ link homology.

In that work, we found it beneficial to work with (a categorification of) an enhanced version of the webs described above, 
in which we retain the $n$-labeled edges, rather than truncate them to tags.
This added rigidity is necessary upon passing to the categorified level in order to define the skew Howe $2$-functor without appealing to complex 
coefficients (in the $n=2$ case).
For example, the neck cutting relation \eqref{sl2neckcutting} for $\mathfrak{sl}_2$ link homology should correspond to the nilHecke 
relation \eqref{eq_nil_dotslide} below; however, the signs don't match in the absence of edges (and foam facets) to 
keep track of the $2$-labeled data. In the $\mathfrak{sl}_2$ case, this construction is originally due to Blanchet \cite{Blan}, 
following earlier work of Clark-Morrison-Walker \cite{CMW} on the functoriality of Khovanov homology. 

As originally defined, Khovanov homology is only projectively functorial under link cobordism \cite{Kh6}, i.e. cobordisms between links induce maps between the 
homologies of those links which are well-defined only up to $\pm1$. Clark, Morrison and Walker resolved
this issue by considering the analog of Bar-Natan's construction for seamed cobordisms between oriented, tagged planar tangles, 
with coefficients in the Gaussian integers. Blanchet later introduced an enhanced version of Bar-Natan's construction, akin to $\mathfrak{sl}_3$ 
foams, in which the tags on planar tangles are extended to $2$-labeled web edges, allowing him to define a properly functorial version 
of Khovanov homology over the integers. 
In the $\slnn{3}$ case, the use of enhanced foams is not entirely necessary to define the foamation $2$-functors; however, the definition of 
these $2$-functors into non-enhanced $\slnn{3}$ foams requires mysterious sign rescalings, see \cite{MPT} and \cite{LQR1}. 
If one works in the enhanced setting, the definition of the foamation functors for $\slnn{2}$ carries over to the $\slnn{3}$ case 
mutatis mutandis. The signs giving the rescaling mentioned above can then be understood topologically, coming from a 
`forgetful' $2$-functor, which rescales enhanced $\slnn{3}$ foams according to the topology of their $3$-labeled facets.

This observation suggests that we should again utilize $n$-labeled web and foam data in the present work.
This allows us to uniformly define our foamation functors for all $n$, and should fix the problem 
of Khovanov-Rozansky homology being only projectively functorial. Indeed, the functoriality issues arise from a symmetry in the definition of the saddle morphisms, 
which is removed in the presence of $n$-labeled edges; see Remark \ref{rem:ben} for more details.

We'd like to point out that the idea of using categorical skew Howe duality to study categorifications of $\sln$ webs has been independently developed by 
Mackaay-Yonezawa \cite{MY} and Mackaay \cite{Mac2} in the setting of the homotopy category of matrix factorizations. Additionally, during the final 
preparation of this paper, the authors became aware of the preprint \cite{TubbenhauerSln}, in which Tubbenhauer uses categorical skew Howe duality to 
present an explicit basis for Mackaay's $\sln$ web algebra, by adapting the Hu-Mathas basis for the cyclotomic KLR algebra \cite{HM}. He goes on 
to use this basis to give an in-principle computable version of Khovanov-Rozansky homology using only webs in the image of the KLR algebra; 
this is a(n independently discovered) consequence of Remark \ref{rem:turningtoKLR}.
His results should translate into our framework to provide explicit bases for certain $2\Hom$ spaces of foams.

The main difference between our approach and those above is that we define and give tools to study a topological foam category, 
which allows for the use of topological intuition in analyzing morphisms in the homotopy category of matrix factorizations and in 
the (quotient of the) categorified $q$-Schur algebra through which the $2$-representation to matrix factorizations factor 
(see e.g. Lemma \ref{lem:elemIso} for an illustration of the latter). In particular, the passage from webs in image of $\Ucat_Q(\slnn{2})$ to those 
in the image of the $KLR$ algebra depicted in equation \eqref{web_trick} follows from a foam isotopy, but a less-obvious isomorphism in the 
categorified $q$-Schur algebra.

\subsection{Categorified quantum groups}

In this section we recall the Khovanov-Lauda categorification of quantum $\mathfrak{sl}_m$, the Rickard complexes, 
and the categorified $q$-Schur algebra.
Before giving the definition, we establish notation for the $\mathfrak{sl}_m$ Cartan datum.
Let $I=\{1,\ldots,m-1\}$ denote the nodes of the type $A_{m-1}$ Dynkin diagram, 
$X=\Z^{m-1}$ be the weight lattice, and let $\{\alpha_i\}_{i \in I}$ and $\{\Lambda_i\}_{i \in I}$ 
denote the collection of simple roots and fundamental weights, respectively. The $\mathfrak{sl}_m$ Cartan matrix 
\[ a_{ij} =
\left\{
\begin{array}{ll}
  2 & \text{if $i=j$}\\
  -1&  \text{if $|i-j|=1$} \\
  0 & \text{if $|i-j|>1$}
\end{array}
\right.
\]
defines a symmetric bilinear form on $X$ determined by $(\alpha_i, \alpha_j) = a_{ij}$. If $\{h_i\}_{i \in I}$ 
denote the simple coroots in $X^\vee = \Hom_{\Z}(X,\Z)$ then the canonical pairing 
$\langle -,- \rangle \maps X^{\vee} \times X\to \Z$ satisfying $\la h_i, \Lambda_i \ra = \delta_{i,j}$ is 
given on $\lambda \in X$ via $\langle h_i, \lambda \rangle = 2 \frac{(\alpha_i,\lambda)}{(\alpha_i,\alpha_i)}$.
Any weight $\lambda \in X$ can then be written as $\lambda = (\lambda_1,\lambda_2, \dots, \lambda_{m-1})$, where
$\lambda_i=\la h_i, \lambda \ra$.

\subsubsection{Categorified quantum $\slm$}
We now present the general definition of the signed version of categorified quantum $\mathfrak{sl}_m$ due to Cautis-Lauda \cite{CLau}, which extends the  
original Khovanov-Lauda construction \cite{KL,KL2,KL3}. 
Fix a field $\Bbbk$ and choice of scalars $\{t_{i,j}\}_{i,j \in I}$ satisfying 
\begin{itemize}
\item $t_{ii}=1$ for all $i \in I$ and $t_{ij} \in \Bbbk^{\times}$ for $i\neq j$,
 \item $t_{ij}=t_{ji}$ when $a_{ij}=0$.
\end{itemize}
We'll use the choice $t_{i,i+1}=-1$, $t_{i+1,i}=1$, and $t_{ij}=1$ for $|i-j|\geq 2$ in Sections \ref{sec:foams} and \ref{sec:link_hom}, 
but for now present the construction for general scalars. The $2$-categories constructed for different choices of these scalars are isomorphic.

\begin{defn} \label{defU_cat}
The 2-category $\Ucat_Q(\mathfrak{sl}_m)$ is the graded additive $\Bbbk$-linear 2-category consisting of:
\begin{itemize}
\item \textbf{Objects} $\lambda$ for $\lambda \in X$.
\item \textbf{1-morphisms} are formal direct sums of (shifts of) compositions of
$$\onel, \quad \onenn{\l+\alpha_i} \sE_i = \onenn{\l+\alpha_i} \sE_i\onel = \sE_i \onel, \quad \text{ and }\quad 
\onenn{\lambda-\alpha_i} \sF_i = \onenn{\lambda-\alpha_i} \sF_i\onel = \sF_i\onel$$
for $i \in I$ and $\l \in X$.
\item \textbf{2-morphisms} are $\Bbbk$-vector spaces spanned by compositions of colored, decorated tangle-like diagrams illustrated below.
\begin{align}
  \xy 0;/r.17pc/:
 (0,7);(0,-7); **\dir{-} ?(.75)*\dir{>};
 (0,0)*{\bullet};
 (7,3)*{ \scs \lambda};
 (-9,3)*{\scs  \lambda+\alpha_i};
 (-2.5,-6)*{\scs i};
 (-10,0)*{};(10,0)*{};
 \endxy &\maps \cal{E}_i\onel \to \cal{E}_i\onel\{ (\alpha_i,\alpha_i) \}  & \quad
 &
    \xy 0;/r.17pc/:
 (0,7);(0,-7); **\dir{-} ?(.75)*\dir{<};
 (0,0)*{\bullet};
 (7,3)*{ \scs \lambda};
 (-9,3)*{\scs  \lambda-\alpha_i};
 (-2.5,-6)*{\scs i};
 (-10,0)*{};(10,0)*{};
 \endxy\maps \cal{F}_i\onel \to \cal{F}_i\onel\{ (\alpha_i,\alpha_i) \}  \nn \\
   & & & \nn \\
   \xy 0;/r.17pc/:
  (0,0)*{\xybox{
    (-4,-4)*{};(4,4)*{} **\crv{(-4,-1) & (4,1)}?(1)*\dir{>} ;
    (4,-4)*{};(-4,4)*{} **\crv{(4,-1) & (-4,1)}?(1)*\dir{>};
    (-5.5,-3)*{\scs i};
     (5.5,-3)*{\scs j};
     (9,1)*{\scs  \lambda};
     (-10,0)*{};(10,0)*{};
     }};
  \endxy \;\;&\maps \cal{E}_i\cal{E}_j\onel  \to \cal{E}_j\cal{E}_i\onel\{ - (\alpha_i,\alpha_j) \}  &
  &
   \xy 0;/r.17pc/:
  (0,0)*{\xybox{
    (-4,4)*{};(4,-4)*{} **\crv{(-4,1) & (4,-1)}?(1)*\dir{>} ;
    (4,4)*{};(-4,-4)*{} **\crv{(4,1) & (-4,-1)}?(1)*\dir{>};
    (-6.5,-3)*{\scs i};
     (6.5,-3)*{\scs j};
     (9,1)*{\scs  \lambda};
     (-10,0)*{};(10,0)*{};
     }};
  \endxy\;\; \maps \cal{F}_i\cal{F}_j\onel  \to \cal{F}_j\cal{F}_i\onel\{ - (\alpha_i,\alpha_j) \}  \nn \\
  & & & \nn \\
     \xy 0;/r.17pc/:
    (0,-3)*{\bbpef{i}};
    (8,-5)*{\scs  \lambda};
    (-10,0)*{};(10,0)*{};
    \endxy &\maps \onel  \to \cal{F}_i\cal{E}_i\onel\{ 1 + (\l, \alpha_i) \}   &
    &
   \xy 0;/r.17pc/:
    (0,-3)*{\bbpfe{i}};
    (8,-5)*{\scs \lambda};
    (-10,0)*{};(10,0)*{};
    \endxy \maps \onel  \to\cal{E}_i\cal{F}_i\onel\{ 1 - (\l, \alpha_i) \}  \nn \\
      & & & \nn \\
  \xy 0;/r.17pc/:
    (0,0)*{\bbcef{i}};
    (8,4)*{\scs  \lambda};
    (-10,0)*{};(10,0)*{};
    \endxy & \maps \cal{F}_i\cal{E}_i\onel \to\onel\{ 1 + (\l, \alpha_i) \}  &
    &
 \xy 0;/r.17pc/:
    (0,0)*{\bbcfe{i}};
    (8,4)*{\scs  \lambda};
    (-10,0)*{};(10,0)*{};
    \endxy \maps\cal{E}_i\cal{F}_i\onel  \to\onel\{ 1 - (\l, \alpha_i) \} \nn
\end{align}
\end{itemize}
\end{defn}
Here we follow the grading conventions in \cite{CLau}
which are opposite to those from \cite{KL3} but line up nicely with the gradings on foams used later in the paper.
In order to match with the standard notions of degree in categorified quantum groups, we'll use the convention that the 
degree of the above morphisms are given by 
the degree shift of the codomain minus that of the domain, e.g. a dot has degree $2$ and an $(i,i)$ crossing has 
degree $-2$. In this $2$-category (and those throughout the paper) we
read diagrams from right to left and bottom to top.  The identity 2-morphism of the 1-morphism
$\cal{E}_i \onel$ is
represented by an upward oriented line colored by $i$ and the identity 2-morphism of $\cal{F}_i \onel$ is
represented by a downward line.

The 2-morphisms satisfy the following relations:
\begin{enumerate}
\item \label{item_cycbiadjoint} The 1-morphisms $\cal{E}_i \onel$ and $\cal{F}_i \onel$ are biadjoint (up to a specified degree shift). These conditions are expressed diagrammatically as
    \begin{equation} \label{eq_biadjoint1}
 \xy   0;/r.17pc/:
    (-8,0)*{}="1";
    (0,0)*{}="2";
    (8,0)*{}="3";
    (-8,-10);"1" **\dir{-};
    "1";"2" **\crv{(-8,8) & (0,8)} ?(0)*\dir{>} ?(1)*\dir{>};
    "2";"3" **\crv{(0,-8) & (8,-8)}?(1)*\dir{>};
    "3"; (8,10) **\dir{-};
    (12,-9)*{\lambda};
    (-6,9)*{\lambda+\alpha_i};
    \endxy
    \; =
    \;
\xy   0;/r.17pc/:
    (-8,0)*{}="1";
    (0,0)*{}="2";
    (8,0)*{}="3";
    (0,-10);(0,10)**\dir{-} ?(.5)*\dir{>};
    (5,8)*{\lambda};
    (-9,8)*{\lambda+\alpha_i};
    \endxy
\qquad \quad \xy  0;/r.17pc/:
    (8,0)*{}="1";
    (0,0)*{}="2";
    (-8,0)*{}="3";
    (8,-10);"1" **\dir{-};
    "1";"2" **\crv{(8,8) & (0,8)} ?(0)*\dir{<} ?(1)*\dir{<};
    "2";"3" **\crv{(0,-8) & (-8,-8)}?(1)*\dir{<};
    "3"; (-8,10) **\dir{-};
    (12,9)*{\lambda+\alpha_i};
    (-6,-9)*{\lambda};
    \endxy
    \; =
    \;
\xy  0;/r.17pc/:
    (8,0)*{}="1";
    (0,0)*{}="2";
    (-8,0)*{}="3";
    (0,-10);(0,10)**\dir{-} ?(.5)*\dir{<};
    (9,-8)*{\lambda+\alpha_i};
    (-6,-8)*{\lambda};
    \endxy
\end{equation}

\begin{equation}\label{eq_biadjoint2}
 \xy   0;/r.17pc/:
    (8,0)*{}="1";
    (0,0)*{}="2";
    (-8,0)*{}="3";
    (8,-10);"1" **\dir{-};
    "1";"2" **\crv{(8,8) & (0,8)} ?(0)*\dir{>} ?(1)*\dir{>};
    "2";"3" **\crv{(0,-8) & (-8,-8)}?(1)*\dir{>};
    "3"; (-8,10) **\dir{-};
    (12,9)*{\lambda};
    (-5,-9)*{\lambda+\alpha_i};
    \endxy
    \; =
    \;
    \xy 0;/r.17pc/:
    (8,0)*{}="1";
    (0,0)*{}="2";
    (-8,0)*{}="3";
    (0,-10);(0,10)**\dir{-} ?(.5)*\dir{>};
    (5,-8)*{\lambda};
    (-9,-8)*{\lambda+\alpha_i};
    \endxy
\qquad \quad \xy   0;/r.17pc/:
    (-8,0)*{}="1";
    (0,0)*{}="2";
    (8,0)*{}="3";
    (-8,-10);"1" **\dir{-};
    "1";"2" **\crv{(-8,8) & (0,8)} ?(0)*\dir{<} ?(1)*\dir{<};
    "2";"3" **\crv{(0,-8) & (8,-8)}?(1)*\dir{<};
    "3"; (8,10) **\dir{-};
    (12,-9)*{\lambda+\alpha_i};
    (-6,9)*{\lambda};
    \endxy
    \; =
    \;
\xy   0;/r.17pc/:
    (-8,0)*{}="1";
    (0,0)*{}="2";
    (8,0)*{}="3";
    (0,-10);(0,10)**\dir{-} ?(.5)*\dir{<};
   (9,8)*{\lambda+\alpha_i};
    (-6,8)*{\lambda};
    \endxy
\end{equation}

\item The 2-morphisms are $Q$-cyclic with respect to this biadjoint structure:
\begin{equation} \label{eq_cyclic_dot}
 \xy 0;/r.17pc/:
    (-8,5)*{}="1";
    (0,5)*{}="2";
    (0,-5)*{}="2'";
    (8,-5)*{}="3";
    (-8,-10);"1" **\dir{-};
    "2";"2'" **\dir{-} ?(.5)*\dir{<};
    "1";"2" **\crv{(-8,12) & (0,12)} ?(0)*\dir{<};
    "2'";"3" **\crv{(0,-12) & (8,-12)}?(1)*\dir{<};
    "3"; (8,10) **\dir{-};
    (17,-9)*{\lambda+\alpha_i};
    (-12,9)*{\lambda};
    (0,4)*{\bullet};
    (10,8)*{\scs };
    (-10,-8)*{\scs };
    \endxy
    \quad = \quad
      \xy 0;/r.17pc/:
 (0,10);(0,-10); **\dir{-} ?(.75)*\dir{<}+(2.3,0)*{\scriptstyle{}}
 ?(.1)*\dir{ }+(2,0)*{\scs };
 (0,0)*{\bullet};
 (-6,5)*{\lambda};
 (10,5)*{\lambda+\alpha_i};
 (-10,0)*{};(10,0)*{};(-2,-8)*{\scs };
 \endxy
    \quad = \quad
   \xy 0;/r.17pc/:
    (8,5)*{}="1";
    (0,5)*{}="2";
    (0,-5)*{}="2'";
    (-8,-5)*{}="3";
    (8,-10);"1" **\dir{-};
    "2";"2'" **\dir{-} ?(.5)*\dir{<};
    "1";"2" **\crv{(8,12) & (0,12)} ?(0)*\dir{<};
    "2'";"3" **\crv{(0,-12) & (-8,-12)}?(1)*\dir{<};
    "3"; (-8,10) **\dir{-};
    (17,9)*{\lambda+\alpha_i};
    (-12,-9)*{\lambda};
    (0,4)*{\bullet};
    (-10,8)*{\scs };
    (10,-8)*{\scs };
    \endxy
\end{equation}
\begin{equation} \label{eq_almost_cyclic}
   \xy 0;/r.17pc/:
  (0,0)*{\xybox{
    (-4,4)*{};(4,-4)*{} **\crv{(-4,1) & (4,-1)}?(1)*\dir{>} ;
    (4,4)*{};(-4,-4)*{} **\crv{(4,1) & (-4,-1)}?(1)*\dir{>};
    (-6.5,-3)*{\scs i};
     (6.5,-3)*{\scs j};
     (9,1)*{\scs  \lambda};
     (-10,0)*{};(10,0)*{};
     }};
  \endxy \quad = \quad
  t_{ij}^{-1}\xy 0;/r.17pc/:
  (0,0)*{\xybox{
    (4,-4)*{};(-4,4)*{} **\crv{(4,-1) & (-4,1)}?(1)*\dir{>};
    (-4,-4)*{};(4,4)*{} **\crv{(-4,-1) & (4,1)};
     (-4,4)*{};(18,4)*{} **\crv{(-4,16) & (18,16)} ?(1)*\dir{>};
     (4,-4)*{};(-18,-4)*{} **\crv{(4,-16) & (-18,-16)} ?(1)*\dir{<}?(0)*\dir{<};
     (-18,-4);(-18,12) **\dir{-};(-12,-4);(-12,12) **\dir{-};
     (18,4);(18,-12) **\dir{-};(12,4);(12,-12) **\dir{-};
     (8,1)*{ \lambda};
     (-10,0)*{};(10,0)*{};
     (-4,-4)*{};(-12,-4)*{} **\crv{(-4,-10) & (-12,-10)}?(1)*\dir{<}?(0)*\dir{<};
      (4,4)*{};(12,4)*{} **\crv{(4,10) & (12,10)}?(1)*\dir{>}?(0)*\dir{>};
      (-20,11)*{\scs j};(-10,11)*{\scs i};
      (20,-11)*{\scs j};(10,-11)*{\scs i};
     }};
  \endxy
\quad =  \quad t_{ji}^{-1}
\xy 0;/r.17pc/:
  (0,0)*{\xybox{
    (-4,-4)*{};(4,4)*{} **\crv{(-4,-1) & (4,1)}?(1)*\dir{>};
    (4,-4)*{};(-4,4)*{} **\crv{(4,-1) & (-4,1)};
     (4,4)*{};(-18,4)*{} **\crv{(4,16) & (-18,16)} ?(1)*\dir{>};
     (-4,-4)*{};(18,-4)*{} **\crv{(-4,-16) & (18,-16)} ?(1)*\dir{<}?(0)*\dir{<};
     (18,-4);(18,12) **\dir{-};(12,-4);(12,12) **\dir{-};
     (-18,4);(-18,-12) **\dir{-};(-12,4);(-12,-12) **\dir{-};
     (8,1)*{ \lambda};
     (-10,0)*{};(10,0)*{};
      (4,-4)*{};(12,-4)*{} **\crv{(4,-10) & (12,-10)}?(1)*\dir{<}?(0)*\dir{<};
      (-4,4)*{};(-12,4)*{} **\crv{(-4,10) & (-12,10)}?(1)*\dir{>}?(0)*\dir{>};
      (20,11)*{\scs i};(10,11)*{\scs j};
      (-20,-11)*{\scs i};(-10,-11)*{\scs j};
     }};
  \endxy
\end{equation}

\begin{equation} \label{eq_crossl-gen}
  \xy 0;/r.18pc/:
  (0,0)*{\xybox{
    (-4,-4)*{};(4,4)*{} **\crv{(-4,-1) & (4,1)}?(1)*\dir{>} ;
    (4,-4)*{};(-4,4)*{} **\crv{(4,-1) & (-4,1)}?(0)*\dir{<};
    (-5,-3)*{\scs j};
     (6.5,-3)*{\scs i};
     (9,2)*{ \lambda};
     (-12,0)*{};(12,0)*{};
     }};
  \endxy
\quad = \quad
 \xy 0;/r.17pc/:
  (0,0)*{\xybox{
    (4,-4)*{};(-4,4)*{} **\crv{(4,-1) & (-4,1)}?(1)*\dir{>};
    (-4,-4)*{};(4,4)*{} **\crv{(-4,-1) & (4,1)};
     (-4,4);(-4,12) **\dir{-};
     (-12,-4);(-12,12) **\dir{-};
     (4,-4);(4,-12) **\dir{-};(12,4);(12,-12) **\dir{-};
     (16,1)*{\lambda};
     (-10,0)*{};(10,0)*{};
     (-4,-4)*{};(-12,-4)*{} **\crv{(-4,-10) & (-12,-10)}?(1)*\dir{<}?(0)*\dir{<};
      (4,4)*{};(12,4)*{} **\crv{(4,10) & (12,10)}?(1)*\dir{>}?(0)*\dir{>};
      (-14,11)*{\scs i};(-2,11)*{\scs j};
      (14,-11)*{\scs i};(2,-11)*{\scs j};
     }};
  \endxy
  \quad = \quad t_{ij} \;\;
 \xy 0;/r.17pc/:
  (0,0)*{\xybox{
    (-4,-4)*{};(4,4)*{} **\crv{(-4,-1) & (4,1)}?(1)*\dir{<};
    (4,-4)*{};(-4,4)*{} **\crv{(4,-1) & (-4,1)};
     (4,4);(4,12) **\dir{-};
     (12,-4);(12,12) **\dir{-};
     (-4,-4);(-4,-12) **\dir{-};(-12,4);(-12,-12) **\dir{-};
     (16,1)*{\lambda};
     (10,0)*{};(-10,0)*{};
     (4,-4)*{};(12,-4)*{} **\crv{(4,-10) & (12,-10)}?(1)*\dir{>}?(0)*\dir{>};
      (-4,4)*{};(-12,4)*{} **\crv{(-4,10) & (-12,10)}?(1)*\dir{<}?(0)*\dir{<};
     }};
     (12,11)*{\scs j};(0,11)*{\scs i};
      (-17,-11)*{\scs j};(-5,-11)*{\scs i};
  \endxy
\end{equation}
\begin{equation} \label{eq_crossr-gen}
  \xy 0;/r.18pc/:
  (0,0)*{\xybox{
    (-4,-4)*{};(4,4)*{} **\crv{(-4,-1) & (4,1)}?(0)*\dir{<} ;
    (4,-4)*{};(-4,4)*{} **\crv{(4,-1) & (-4,1)}?(1)*\dir{>};
    (5.1,-3)*{\scs i};
     (-6.5,-3)*{\scs j};
     (9,2)*{ \lambda};
     (-12,0)*{};(12,0)*{};
     }};
  \endxy
\quad = \quad
 \xy 0;/r.17pc/:
  (0,0)*{\xybox{
    (-4,-4)*{};(4,4)*{} **\crv{(-4,-1) & (4,1)}?(1)*\dir{>};
    (4,-4)*{};(-4,4)*{} **\crv{(4,-1) & (-4,1)};
     (4,4);(4,12) **\dir{-};
     (12,-4);(12,12) **\dir{-};
     (-4,-4);(-4,-12) **\dir{-};(-12,4);(-12,-12) **\dir{-};
     (16,-6)*{\lambda};
     (10,0)*{};(-10,0)*{};
     (4,-4)*{};(12,-4)*{} **\crv{(4,-10) & (12,-10)}?(1)*\dir{<}?(0)*\dir{<};
      (-4,4)*{};(-12,4)*{} **\crv{(-4,10) & (-12,10)}?(1)*\dir{>}?(0)*\dir{>};
      (14,11)*{\scs j};(2,11)*{\scs i};
      (-14,-11)*{\scs j};(-2,-11)*{\scs i};
     }};
  \endxy
  \quad = \quad t_{ji} \;\;
  \xy 0;/r.17pc/:
  (0,0)*{\xybox{
    (4,-4)*{};(-4,4)*{} **\crv{(4,-1) & (-4,1)}?(1)*\dir{<};
    (-4,-4)*{};(4,4)*{} **\crv{(-4,-1) & (4,1)};
     (-4,4);(-4,12) **\dir{-};
     (-12,-4);(-12,12) **\dir{-};
     (4,-4);(4,-12) **\dir{-};(12,4);(12,-12) **\dir{-};
     (16,6)*{\lambda};
     (-10,0)*{};(10,0)*{};
     (-4,-4)*{};(-12,-4)*{} **\crv{(-4,-10) & (-12,-10)}?(1)*\dir{>}?(0)*\dir{>};
      (4,4)*{};(12,4)*{} **\crv{(4,10) & (12,10)}?(1)*\dir{<}?(0)*\dir{<};
      (-14,11)*{\scs i};(-2,11)*{\scs j};(14,-11)*{\scs i};(2,-11)*{\scs j};
     }};
  \endxy
\end{equation}
where the second equality in \eqref{eq_crossl-gen} and \eqref{eq_crossr-gen}
follow from \eqref{eq_almost_cyclic}.

\item The $\cal{E}$'s carry an action of the KLR algebra associated to $Q$. The KLR algebra $R=R_Q$ associated to $Q$ is defined by finite $\Bbbk$-linear combinations of 
braid--like diagrams in the plane, where each strand is colored by a vertex $i \in I$.  Strands can intersect and can carry dots but triple intersections are not allowed.  
Diagrams are considered up to planar isotopy that do not change the combinatorial type of the diagram. We recall the local relations below.
\begin{enumerate}[i)]
\item
If all strands are colored by the same $i \in I$, then the nilHecke algebra axioms hold 
 \begin{equation}
\vcenter{
\xy 0;/r.17pc/:
	(-4,-4)*{};(4,4)*{} **\crv{(-4,-1) & (4,1)}?(1)*\dir{};
	(4,-4)*{};(-4,4)*{} **\crv{(4,-1) & (-4,1)}?(1)*\dir{};
	(-4,4)*{};(4,12)*{} **\crv{(-4,7) & (4,9)}?(1)*\dir{};
	(4,4)*{};(-4,12)*{} **\crv{(4,7) & (-4,9)}?(1)*\dir{};
	(-4,12); (-4,13) **\dir{-}?(1)*\dir{>};
	(4,12); (4,13) **\dir{-}?(1)*\dir{>};
	(9,8)*{\lambda};
\endxy}
 \;\; =\;\; 0, \qquad \quad
\vcenter{\xy 0;/r.17pc/:
    (-4,-4)*{};(4,4)*{} **\crv{(-4,-1) & (4,1)}?(1)*\dir{};
    (4,-4)*{};(-4,4)*{} **\crv{(4,-1) & (-4,1)}?(1)*\dir{};
    (4,4)*{};(12,12)*{} **\crv{(4,7) & (12,9)}?(1)*\dir{};
    (12,4)*{};(4,12)*{} **\crv{(12,7) & (4,9)}?(1)*\dir{};
    (-4,12)*{};(4,20)*{} **\crv{(-4,15) & (4,17)}?(1)*\dir{};
    (4,12)*{};(-4,20)*{} **\crv{(4,15) & (-4,17)}?(1)*\dir{};
    (-4,4)*{}; (-4,12) **\dir{-};
    (12,-4)*{}; (12,4) **\dir{-};
    (12,12)*{}; (12,20) **\dir{-};
    (4,20); (4,21) **\dir{-}?(1)*\dir{>};
    (-4,20); (-4,21) **\dir{-}?(1)*\dir{>};
    (12,20); (12,21) **\dir{-}?(1)*\dir{>};
   (18,8)*{\lambda};
\endxy}
 \;\; =\;\;
\vcenter{\xy 0;/r.17pc/:
    (4,-4)*{};(-4,4)*{} **\crv{(4,-1) & (-4,1)}?(1)*\dir{};
    (-4,-4)*{};(4,4)*{} **\crv{(-4,-1) & (4,1)}?(1)*\dir{};
    (-4,4)*{};(-12,12)*{} **\crv{(-4,7) & (-12,9)}?(1)*\dir{};
    (-12,4)*{};(-4,12)*{} **\crv{(-12,7) & (-4,9)}?(1)*\dir{};
    (4,12)*{};(-4,20)*{} **\crv{(4,15) & (-4,17)}?(1)*\dir{};
    (-4,12)*{};(4,20)*{} **\crv{(-4,15) & (4,17)}?(1)*\dir{};
    (4,4)*{}; (4,12) **\dir{-};
    (-12,-4)*{}; (-12,4) **\dir{-};
    (-12,12)*{}; (-12,20) **\dir{-};
    (4,20); (4,21) **\dir{-}?(1)*\dir{>};
    (-4,20); (-4,21) **\dir{-}?(1)*\dir{>};
    (-12,20); (-12,21) **\dir{-}?(1)*\dir{>};
  (10,8)*{\lambda};
\endxy}
 \label{eq_nil_rels}
  \end{equation}

\begin{equation}
 \xy 0;/r.18pc/:
  (4,6);(4,-4) **\dir{-}?(0)*\dir{<}+(2.3,0)*{};
  (-4,6);(-4,-4) **\dir{-}?(0)*\dir{<}+(2.3,0)*{};
 \endxy
 \quad =
\xy 0;/r.18pc/:
  (0,0)*{\xybox{
    (-4,-4)*{};(4,6)*{} **\crv{(-4,-1) & (4,1)}?(1)*\dir{>}?(.25)*{\bullet};
    (4,-4)*{};(-4,6)*{} **\crv{(4,-1) & (-4,1)}?(1)*\dir{>};
     (-10,0)*{};(10,0)*{};
     }};
  \endxy
 \;\; -
\xy 0;/r.18pc/:
  (0,0)*{\xybox{
    (-4,-4)*{};(4,6)*{} **\crv{(-4,-1) & (4,1)}?(1)*\dir{>}?(.75)*{\bullet};
    (4,-4)*{};(-4,6)*{} **\crv{(4,-1) & (-4,1)}?(1)*\dir{>};
     (-10,0)*{};(10,0)*{};
     }};
  \endxy
 \;\; =
\xy 0;/r.18pc/:
  (0,0)*{\xybox{
    (-4,-4)*{};(4,6)*{} **\crv{(-4,-1) & (4,1)}?(1)*\dir{>};
    (4,-4)*{};(-4,6)*{} **\crv{(4,-1) & (-4,1)}?(1)*\dir{>}?(.75)*{\bullet};
     (-10,0)*{};(10,0)*{};
     }};
  \endxy
 \;\; -
  \xy 0;/r.18pc/:
  (0,0)*{\xybox{
    (-4,-4)*{};(4,6)*{} **\crv{(-4,-1) & (4,1)}?(1)*\dir{>} ;
    (4,-4)*{};(-4,6)*{} **\crv{(4,-1) & (-4,1)}?(1)*\dir{>}?(.25)*{\bullet};
     (-10,0)*{};(10,0)*{};
     }};
  \endxy . \label{eq_nil_dotslide}
\end{equation}

\item For $i \neq j$
\begin{equation}
 \vcenter{\xy 0;/r.17pc/:
    (-4,-4)*{};(4,4)*{} **\crv{(-4,-1) & (4,1)}?(1)*\dir{};
    (4,-4)*{};(-4,4)*{} **\crv{(4,-1) & (-4,1)}?(1)*\dir{};
    (-4,4)*{};(4,12)*{} **\crv{(-4,7) & (4,9)}?(1)*\dir{};
    (4,4)*{};(-4,12)*{} **\crv{(4,7) & (-4,9)}?(1)*\dir{};
    (8,8)*{\lambda};
    (4,12); (4,13) **\dir{-}?(1)*\dir{>};
    (-4,12); (-4,13) **\dir{-}?(1)*\dir{>};
  (-5.5,-3)*{\scs i};
     (5.5,-3)*{\scs j};
 \endxy}
 \qquad = \qquad
 \left\{
 \begin{array}{ccc}
     t_{ij}\;\xy 0;/r.17pc/:
  (3,9);(3,-9) **\dir{-}?(0)*\dir{<}+(2.3,0)*{};
  (-3,9);(-3,-9) **\dir{-}?(0)*\dir{<}+(2.3,0)*{};
  (-5,-6)*{\scs i};     (5.1,-6)*{\scs j};
 \endxy &  &  \text{if $(\alpha_i, \alpha_j)=0$,}\\ \\
 t_{ij} \vcenter{\xy 0;/r.17pc/:
  (3,9);(3,-9) **\dir{-}?(0)*\dir{<}+(2.3,0)*{};
  (-3,9);(-3,-9) **\dir{-}?(0)*\dir{<}+(2.3,0)*{};
  (-3,4)*{\bullet};(-6.5,5)*{};
  (-5,-6)*{\scs i};     (5.1,-6)*{\scs j};
 \endxy} \;\; + \;\; t_{ji}
  \vcenter{\xy 0;/r.17pc/:
  (3,9);(3,-9) **\dir{-}?(0)*\dir{<}+(2.3,0)*{};
  (-3,9);(-3,-9) **\dir{-}?(0)*\dir{<}+(2.3,0)*{};
  (3,4)*{\bullet};(7,5)*{};
  (-5,-6)*{\scs i};     (5.1,-6)*{\scs j};
 \endxy}
   &  & \text{if $(\alpha_i, \alpha_j) \neq 0$.}
 \end{array}
 \right. \label{eq_r2_ij-gen}
\end{equation}

\item For $i \neq j$ the dot sliding relations
\begin{equation} \label{eq_dot_slide_ij-gen}
\xy 0;/r.18pc/:
  (0,0)*{\xybox{
    (-4,-4)*{};(4,6)*{} **\crv{(-4,-1) & (4,1)}?(1)*\dir{>}?(.75)*{\bullet};
    (4,-4)*{};(-4,6)*{} **\crv{(4,-1) & (-4,1)}?(1)*\dir{>};
    (-5,-3)*{\scs i};
     (5.1,-3)*{\scs j};
     (-10,0)*{};(10,0)*{};
     }};
  \endxy
 \;\; =
\xy 0;/r.18pc/:
  (0,0)*{\xybox{
    (-4,-4)*{};(4,6)*{} **\crv{(-4,-1) & (4,1)}?(1)*\dir{>}?(.25)*{\bullet};
    (4,-4)*{};(-4,6)*{} **\crv{(4,-1) & (-4,1)}?(1)*\dir{>};
    (-5,-3)*{\scs i};
     (5.1,-3)*{\scs j};
     (-10,0)*{};(10,0)*{};
     }};
  \endxy
\qquad  \xy 0;/r.18pc/:
  (0,0)*{\xybox{
    (-4,-4)*{};(4,6)*{} **\crv{(-4,-1) & (4,1)}?(1)*\dir{>};
    (4,-4)*{};(-4,6)*{} **\crv{(4,-1) & (-4,1)}?(1)*\dir{>}?(.75)*{\bullet};
    (-5,-3)*{\scs i};
     (5.1,-3)*{\scs j};
     (-10,0)*{};(10,0)*{};
     }};
  \endxy
\;\;  =
  \xy 0;/r.18pc/:
  (0,0)*{\xybox{
    (-4,-4)*{};(4,6)*{} **\crv{(-4,-1) & (4,1)}?(1)*\dir{>} ;
    (4,-4)*{};(-4,6)*{} **\crv{(4,-1) & (-4,1)}?(1)*\dir{>}?(.25)*{\bullet};
    (-5,-3)*{\scs i};
     (5.1,-3)*{\scs j};
     (-10,0)*{};(12,0)*{};
     }};
  \endxy
\end{equation}
hold.

\item Unless $i = k$ and $(\alpha_i, \alpha_j) < 0$, the relation
\begin{equation}
\vcenter{\xy 0;/r.17pc/:
    (-4,-4)*{};(4,4)*{} **\crv{(-4,-1) & (4,1)}?(1)*\dir{};
    (4,-4)*{};(-4,4)*{} **\crv{(4,-1) & (-4,1)}?(1)*\dir{};
    (4,4)*{};(12,12)*{} **\crv{(4,7) & (12,9)}?(1)*\dir{};
    (12,4)*{};(4,12)*{} **\crv{(12,7) & (4,9)}?(1)*\dir{};
    (-4,12)*{};(4,20)*{} **\crv{(-4,15) & (4,17)}?(1)*\dir{};
    (4,12)*{};(-4,20)*{} **\crv{(4,15) & (-4,17)}?(1)*\dir{};
    (-4,4)*{}; (-4,12) **\dir{-};
    (12,-4)*{}; (12,4) **\dir{-};
    (12,12)*{}; (12,20) **\dir{-};
    (4,20); (4,21) **\dir{-}?(1)*\dir{>};
    (-4,20); (-4,21) **\dir{-}?(1)*\dir{>};
    (12,20); (12,21) **\dir{-}?(1)*\dir{>};
   (18,8)*{\lambda};  (-6,-3)*{\scs i};
  (6,-3)*{\scs j};
  (15,-3)*{\scs k};
\endxy}
 \;\; =\;\;
\vcenter{\xy 0;/r.17pc/:
    (4,-4)*{};(-4,4)*{} **\crv{(4,-1) & (-4,1)}?(1)*\dir{};
    (-4,-4)*{};(4,4)*{} **\crv{(-4,-1) & (4,1)}?(1)*\dir{};
    (-4,4)*{};(-12,12)*{} **\crv{(-4,7) & (-12,9)}?(1)*\dir{};
    (-12,4)*{};(-4,12)*{} **\crv{(-12,7) & (-4,9)}?(1)*\dir{};
    (4,12)*{};(-4,20)*{} **\crv{(4,15) & (-4,17)}?(1)*\dir{};
    (-4,12)*{};(4,20)*{} **\crv{(-4,15) & (4,17)}?(1)*\dir{};
    (4,4)*{}; (4,12) **\dir{-};
    (-12,-4)*{}; (-12,4) **\dir{-};
    (-12,12)*{}; (-12,20) **\dir{-};
    (4,20); (4,21) **\dir{-}?(1)*\dir{>};
    (-4,20); (-4,21) **\dir{-}?(1)*\dir{>};
    (-12,20); (-12,21) **\dir{-}?(1)*\dir{>};
  (10,8)*{\lambda};
  (-14,-3)*{\scs i};
  (-6,-3)*{\scs j};
  (6,-3)*{\scs k};
\endxy}
 \label{eq_r3_easy-gen}
\end{equation}
holds; otherwise, $(\alpha_i, \alpha_j) =-1$ and
\begin{equation}
\vcenter{\xy 0;/r.17pc/:
    (-4,-4)*{};(4,4)*{} **\crv{(-4,-1) & (4,1)}?(1)*\dir{};
    (4,-4)*{};(-4,4)*{} **\crv{(4,-1) & (-4,1)}?(1)*\dir{};
    (4,4)*{};(12,12)*{} **\crv{(4,7) & (12,9)}?(1)*\dir{};
    (12,4)*{};(4,12)*{} **\crv{(12,7) & (4,9)}?(1)*\dir{};
    (-4,12)*{};(4,20)*{} **\crv{(-4,15) & (4,17)}?(1)*\dir{};
    (4,12)*{};(-4,20)*{} **\crv{(4,15) & (-4,17)}?(1)*\dir{};
    (-4,4)*{}; (-4,12) **\dir{-};
    (12,-4)*{}; (12,4) **\dir{-};
    (12,12)*{}; (12,20) **\dir{-};
    (4,20); (4,21) **\dir{-}?(1)*\dir{>};
    (-4,20); (-4,21) **\dir{-}?(1)*\dir{>};
    (12,20); (12,21) **\dir{-}?(1)*\dir{>};
   (18,8)*{\lambda};  (-6,-3)*{\scs i};
  (6,-3)*{\scs j};
  (15,-3)*{\scs i};
\endxy}
\quad - \quad
\vcenter{\xy 0;/r.17pc/:
    (4,-4)*{};(-4,4)*{} **\crv{(4,-1) & (-4,1)}?(1)*\dir{};
    (-4,-4)*{};(4,4)*{} **\crv{(-4,-1) & (4,1)}?(1)*\dir{};
    (-4,4)*{};(-12,12)*{} **\crv{(-4,7) & (-12,9)}?(1)*\dir{};
    (-12,4)*{};(-4,12)*{} **\crv{(-12,7) & (-4,9)}?(1)*\dir{};
    (4,12)*{};(-4,20)*{} **\crv{(4,15) & (-4,17)}?(1)*\dir{};
    (-4,12)*{};(4,20)*{} **\crv{(-4,15) & (4,17)}?(1)*\dir{};
    (4,4)*{}; (4,12) **\dir{-};
    (-12,-4)*{}; (-12,4) **\dir{-};
    (-12,12)*{}; (-12,20) **\dir{-};
    (4,20); (4,21) **\dir{-}?(1)*\dir{>};
    (-4,20); (-4,21) **\dir{-}?(1)*\dir{>};
    (-12,20); (-12,21) **\dir{-}?(1)*\dir{>};
  (10,8)*{\lambda};
  (-14,-3)*{\scs i};
  (-6,-3)*{\scs j};
  (6,-3)*{\scs i};
\endxy}
 \;\; =\;\;
 t_{ij} \;\;
\xy 0;/r.17pc/:
  (4,12);(4,-12) **\dir{-}?(0)*\dir{<};
  (-4,12);(-4,-12) **\dir{-}?(0)*\dir{<}?(.25)*\dir{};
  (12,12);(12,-12) **\dir{-}?(0)*\dir{<}?(.25)*\dir{};
  (-6,-9)*{\scs i};     (6.1,-9)*{\scs j};
  (14,-9)*{\scs i};
 \endxy .
 \label{eq_r3_hard-gen}
\end{equation}
\end{enumerate}

\item When $i \ne j$ one has the mixed relations  relating $\cal{E}_i \cal{F}_j$ and $\cal{F}_j \cal{E}_i$:
\begin{equation}\label{eq_mixed_EF}
 \vcenter{   \xy 0;/r.18pc/:
    (-4,-4)*{};(4,4)*{} **\crv{(-4,-1) & (4,1)}?(1)*\dir{>};
    (4,-4)*{};(-4,4)*{} **\crv{(4,-1) & (-4,1)}?(1)*\dir{<};?(0)*\dir{<};
    (-4,4)*{};(4,12)*{} **\crv{(-4,7) & (4,9)};
    (4,4)*{};(-4,12)*{} **\crv{(4,7) & (-4,9)}?(1)*\dir{>};
  (8,8)*{\lambda};(-6,-3)*{\scs i};
     (6,-3)*{\scs j};
 \endxy}
 \;\; = \;\; t_{ji}\;\;
\xy 0;/r.18pc/:
  (3,9);(3,-9) **\dir{-}?(.55)*\dir{>}+(2.3,0)*{};
  (-3,9);(-3,-9) **\dir{-}?(.5)*\dir{<}+(2.3,0)*{};
  (8,2)*{\lambda};(-5,-6)*{\scs i};     (5.1,-6)*{\scs j};
 \endxy
\qquad \quad
    \vcenter{\xy 0;/r.18pc/:
    (-4,-4)*{};(4,4)*{} **\crv{(-4,-1) & (4,1)}?(1)*\dir{<};?(0)*\dir{<};
    (4,-4)*{};(-4,4)*{} **\crv{(4,-1) & (-4,1)}?(1)*\dir{>};
    (-4,4)*{};(4,12)*{} **\crv{(-4,7) & (4,9)}?(1)*\dir{>};
    (4,4)*{};(-4,12)*{} **\crv{(4,7) & (-4,9)};
  (8,8)*{\lambda};(-6,-3)*{\scs i};
     (6,-3)*{\scs j};
 \endxy}
 \;\;=\;\; t_{ij}\;\;
\xy 0;/r.18pc/:
  (3,9);(3,-9) **\dir{-}?(.5)*\dir{<}+(2.3,0)*{};
  (-3,9);(-3,-9) **\dir{-}?(.55)*\dir{>}+(2.3,0)*{};
  (8,2)*{\lambda};(-5,-6)*{\scs i};     (5.1,-6)*{\scs j};
 \endxy .
\end{equation}

\item \label{item_positivity} Negative degree bubbles are zero. That is, for all $m \in \Z_+$ one has
\begin{equation} \label{eq_positivity_bubbles}
\xy 0;/r.18pc/:
 (-12,0)*{\icbub{m}{i}};
 (-8,8)*{\lambda};
 \endxy
  = 0
 \qquad  \text{if $m<\lambda_i-1$,} \qquad \xy 0;/r.18pc/: (-12,0)*{\iccbub{m}{i}};
 (-8,8)*{\lambda};
 \endxy = 0\quad
  \text{if $m< -\lambda_i-1$.}
\end{equation}
On the other hand, a dotted bubble of degree zero is just  the identity 2-morphism\footnote{One can define the 2-category so that degree zero bubbles are multiplication by arbitrary scalars at the 
cost of modifying some of the other relations, see for example~\cite{Lau4,MSV2}.  However, it is shown in \cite{CLau} that the resulting 2-categories are all isomorphic.}:
\[
\xy 0;/r.18pc/:
 (0,0)*{\icbub{\lambda_i-1}{i}};
  (4,8)*{\lambda};
 \endxy
  =  \Id_{\onenn{\lambda}} \quad \text{for $\lambda_i \geq 1$,}
  \qquad \quad
  \xy 0;/r.18pc/:
 (0,0)*{\iccbub{-\lambda_i-1}{i}};
  (4,8)*{\lambda};
 \endxy  =  \Id_{\onenn{\lambda}} \quad \text{for $\lambda_i \leq -1$.}\]

\item \label{item_highersl2} For any $i \in I$ one has the extended ${\mathfrak{sl}}_2$-relations. In order to describe certain extended ${\mathfrak{sl}}_2$ relations it is convenient to use a shorthand notation from \cite{Lau1} called fake bubbles. These are diagrams for dotted bubbles where the labels of the number of dots is negative, but the total degree of the dotted bubble taken with these negative dots is still positive. They allow us to write these extended ${\mathfrak{sl}}_2$ relations more uniformly (i.e. independent on whether the weight $\lambda_i$ is positive or negative).
\begin{itemize}
 \item Degree zero fake bubbles are equal to the identity 2-morphisms
\[
 \xy 0;/r.18pc/:
    (2,0)*{\icbub{\l_i-1}{i}};
  (12,8)*{\lambda};
 \endxy
  =  \Id_{\onenn{\lambda}} \quad \text{if $\lambda_i \leq 0$,}
  \qquad \quad
\xy 0;/r.18pc/:
    (2,0)*{\iccbub{-\lambda_i-1}{i}};
  (12,8)*{\lambda};
 \endxy =  \Id_{\onenn{\lambda}} \quad  \text{if $\lambda_i \geq 0$}.\]

  \item Higher degree fake bubbles for $\lambda_i<0$ are defined inductively as
  \begin{equation} \label{eq_fake_nleqz}
  \vcenter{\xy 0;/r.18pc/:
    (2,-11)*{\icbub{\l_i-1+j}{i}};
  (12,-2)*{\l};
 \endxy} \;\; =
 \left\{
 \begin{array}{cl}
  \;\; -\;\;
\xsum{\xy (0,6)*{};  (0,1)*{\scs a+b=j}; (0,-2)*{\scs b\geq 1}; \endxy}
\;\; \vcenter{\xy 0;/r.18pc/:
    (2,0)*{\cbub{\l_i-1+a}{}};
    (20,0)*{\ccbub{-\l-1+b}{}};
  (12,8)*{\lambda};
 \endxy}  & \text{if $0 \leq j < -\l_i+1$} \\ & \\
   0 & \text{if $j < 0$. }
 \end{array}
\right.
 \end{equation}

  \item Higher degree fake bubbles for $\lambda_i>0$ are defined inductively as
   \begin{equation} \label{eq_fake_ngeqz}
  \vcenter{\xy 0;/r.18pc/:
    (2,-11)*{\iccbub{-\l_i-1+j}{i}};
  (12,-2)*{\l};
 \endxy} \;\; =
 \left\{
 \begin{array}{cl}
  \;\; -\;\;
\xsum{\xy (0,6)*{}; (0,1)*{\scs a+b=j}; (0,-2)*{\scs a\geq 1}; \endxy}
\;\; \vcenter{\xy 0;/r.18pc/:
    (2,0)*{\cbub{\l_i-1+a}{}};
    (20,0)*{\ccbub{-\l-1+b}{}};
  (12,8)*{\lambda};
 \endxy}  & \text{if $0 \leq j < \l_i+1$} \\ & \\
   0 & \text{if $j < 0$. }
 \end{array}
\right.
\end{equation}
\end{itemize}
These equations arise from the homogeneous terms in $t$ of the `infinite Grassmannian' equation
\begin{center}
\begin{eqnarray}
 \makebox[0pt]{ $
\left( \xy 0;/r.15pc/:
 (0,0)*{\iccbub{-\l_i-1}{i}};
  (4,8)*{\l};
 \endxy
 +
 \xy 0;/r.15pc/:
 (0,0)*{\iccbub{-\l_i-1+1}{i}};
  (4,8)*{\l};
 \endxy t
 + \cdots +
\xy 0;/r.15pc/:
 (0,0)*{\iccbub{-\l_i-1+\alpha}{i}};
  (4,8)*{\l};
 \endxy t^{\alpha}
 + \cdots
\right)
\left(\xy 0;/r.15pc/:
 (0,0)*{\icbub{\l_i-1}{i}};
  (4,8)*{\l};
 \endxy
 + \xy 0;/r.15pc/:
 (0,0)*{\icbub{\l_i-1+1}{i}};
  (4,8)*{\l};
 \endxy t
 +\cdots +
\xy 0;/r.15pc/:
 (0,0)*{\icbub{\l_i-1+\alpha}{i}};
 (4,8)*{\l};
 \endxy t^{\alpha}
 + \cdots
\right) = \Id_{\onel}.$ } \nn \\ \label{eq_infinite_Grass}
\end{eqnarray}
\end{center}

Now we can define the extended ${\mathfrak{sl}}_2$ relations.  
The additional curl relations given in \cite{CLau} can be derived from those below; here we provide a minimal set of relations.

If $\l_i > 0$, then we have:
\begin{equation} \label{eq_reduction-ngeqz}
  \xy 0;/r.17pc/:
  (14,8)*{\l};
  (-3,-10)*{};(3,5)*{} **\crv{(-3,-2) & (2,1)}?(1)*\dir{>};?(.15)*\dir{>};
    (3,-5)*{};(-3,10)*{} **\crv{(2,-1) & (-3,2)}?(.85)*\dir{>} ?(.1)*\dir{>};
  (3,5)*{}="t1";  (9,5)*{}="t2";
  (3,-5)*{}="t1'";  (9,-5)*{}="t2'";
   "t1";"t2" **\crv{(4,8) & (9, 8)};
   "t1'";"t2'" **\crv{(4,-8) & (9, -8)};
   "t2'";"t2" **\crv{(10,0)} ;
   (-6,-8)*{\scs i};
 \endxy\;\; = \;\; 0
   \qquad \quad
 \vcenter{\xy 0;/r.17pc/:
  (-8,0)*{};(-6,-8)*{\scs i};(6,-8)*{\scs i};
  (8,0)*{};
  (-4,10)*{}="t1";
  (4,10)*{}="t2";
  (-4,-10)*{}="b1";
  (4,-10)*{}="b2";
  "t1";"b1" **\dir{-} ?(.5)*\dir{>};
  "t2";"b2" **\dir{-} ?(.5)*\dir{<};
  (10,2)*{\l};
  \endxy}
\;\; = \;\; -\;\;
   \vcenter{\xy 0;/r.17pc/:
    (-4,-4)*{};(4,4)*{} **\crv{(-4,-1) & (4,1)}?(1)*\dir{<};?(0)*\dir{<};
    (4,-4)*{};(-4,4)*{} **\crv{(4,-1) & (-4,1)}?(1)*\dir{>};
    (-4,4)*{};(4,12)*{} **\crv{(-4,7) & (4,9)}?(1)*\dir{>};
    (4,4)*{};(-4,12)*{} **\crv{(4,7) & (-4,9)};
  (8,8)*{\l};(-6.5,-3)*{\scs i};  (6,-3)*{\scs i};
 \endxy}
\end{equation}
\begin{equation}
 \vcenter{\xy 0;/r.17pc/:
  (-8,0)*{};
  (8,0)*{};
  (-4,10)*{}="t1";
  (4,10)*{}="t2";
  (-4,-10)*{}="b1";
  (4,-10)*{}="b2";(-6,-8)*{\scs i};(6,-8)*{\scs i};
  "t1";"b1" **\dir{-} ?(.5)*\dir{<};
  "t2";"b2" **\dir{-} ?(.5)*\dir{>};
  (10,2)*{\l};
  \endxy}
\;\; = \;\; -\;\;
 \vcenter{   \xy 0;/r.17pc/:
    (-4,-4)*{};(4,4)*{} **\crv{(-4,-1) & (4,1)}?(1)*\dir{>};
    (4,-4)*{};(-4,4)*{} **\crv{(4,-1) & (-4,1)}?(1)*\dir{<};?(0)*\dir{<};
    (-4,4)*{};(4,12)*{} **\crv{(-4,7) & (4,9)};
    (4,4)*{};(-4,12)*{} **\crv{(4,7) & (-4,9)}?(1)*\dir{>};
  (8,8)*{\l};
     (-6,-3)*{\scs i};
     (6.5,-3)*{\scs i};
 \endxy}
  \;\; + \;\;
   \sum_{ \xy  (0,3)*{\scs f_1+f_2+f_3}; (0,0)*{\scs =\lambda_i-1};\endxy}
    \vcenter{\xy 0;/r.17pc/:
    (-12,10)*{\l};
    (-8,0)*{};
  (8,0)*{};
  (-4,-15)*{}="b1";
  (4,-15)*{}="b2";
  "b2";"b1" **\crv{(5,-8) & (-5,-8)}; ?(.05)*\dir{<} ?(.93)*\dir{<}
  ?(.8)*\dir{}+(0,-.1)*{\bullet}+(-3,2)*{\scs f_3};
  (-4,15)*{}="t1";
  (4,15)*{}="t2";
  "t2";"t1" **\crv{(5,8) & (-5,8)}; ?(.15)*\dir{>} ?(.95)*\dir{>}
  ?(.4)*\dir{}+(0,-.2)*{\bullet}+(3,-2)*{\scs \; f_1};
  (0,0)*{\iccbub{\scs \quad -\l_i-1+f_2}{i}};
  (7,-13)*{\scs i};
  (-7,13)*{\scs i};
  \endxy} \label{eq_ident_decomp-ngeqz}
\end{equation}

If $\lambda_i < 0$, then we have:
\begin{equation} \label{eq_reduction-nleqz}
  \xy 0;/r.17pc/:
  (-14,8)*{\l};
  (3,-10)*{};(-3,5)*{} **\crv{(3,-2) & (-2,1)}?(1)*\dir{>};?(.15)*\dir{>};
    (-3,-5)*{};(3,10)*{} **\crv{(-2,-1) & (3,2)}?(.85)*\dir{>} ?(.1)*\dir{>};
  (-3,5)*{}="t1";  (-9,5)*{}="t2";
  (-3,-5)*{}="t1'";  (-9,-5)*{}="t2'";
   "t1";"t2" **\crv{(-4,8) & (-9, 8)};
   "t1'";"t2'" **\crv{(-4,-8) & (-9, -8)};
   "t2'";"t2" **\crv{(-10,0)} ;
   (6,-8)*{\scs i};
 \endxy\;\; = \;\;
0
\qquad \qquad
\vcenter{\xy 0;/r.17pc/:
  (-8,0)*{};
  (8,0)*{};
  (-4,10)*{}="t1";
  (4,10)*{}="t2";
  (-4,-10)*{}="b1";
  (4,-10)*{}="b2";(-6,-8)*{\scs i};(6,-8)*{\scs i};
  "t1";"b1" **\dir{-} ?(.5)*\dir{<};
  "t2";"b2" **\dir{-} ?(.5)*\dir{>};
  (10,2)*{\l};
  \endxy}
\;\; = \;\; -\;\;
\vcenter{   \xy 0;/r.17pc/:
    (-4,-4)*{};(4,4)*{} **\crv{(-4,-1) & (4,1)}?(1)*\dir{>};
    (4,-4)*{};(-4,4)*{} **\crv{(4,-1) & (-4,1)}?(1)*\dir{<};?(0)*\dir{<};
    (-4,4)*{};(4,12)*{} **\crv{(-4,7) & (4,9)};
    (4,4)*{};(-4,12)*{} **\crv{(4,7) & (-4,9)}?(1)*\dir{>};
  (8,8)*{\l};(-6,-3)*{\scs i};
     (6.5,-3)*{\scs i};
 \endxy}
\end{equation}
\begin{equation}
 \vcenter{\xy 0;/r.17pc/:
  (-8,0)*{};(-6,-8)*{\scs i};(6,-8)*{\scs i};
  (8,0)*{};
  (-4,10)*{}="t1";
  (4,10)*{}="t2";
  (-4,-10)*{}="b1";
  (4,-10)*{}="b2";
  "t1";"b1" **\dir{-} ?(.5)*\dir{>};
  "t2";"b2" **\dir{-} ?(.5)*\dir{<};
  (10,2)*{\l};
  (-10,2)*{\l};
  \endxy}
\;\; = \;\;
  -\;\;\vcenter{\xy 0;/r.17pc/:
    (-4,-4)*{};(4,4)*{} **\crv{(-4,-1) & (4,1)}?(1)*\dir{<};?(0)*\dir{<};
    (4,-4)*{};(-4,4)*{} **\crv{(4,-1) & (-4,1)}?(1)*\dir{>};
    (-4,4)*{};(4,12)*{} **\crv{(-4,7) & (4,9)}?(1)*\dir{>};
    (4,4)*{};(-4,12)*{} **\crv{(4,7) & (-4,9)};
  (8,8)*{\l};(-6.5,-3)*{\scs i};  (6,-3)*{\scs i};
 \endxy}
  \;\; + \;\;
    \sum_{ \xy  (0,3)*{\scs g_1+g_2+g_3}; (0,0)*{\scs =-\l_i-1};\endxy}
    \vcenter{\xy 0;/r.17pc/:
    (-8,0)*{};
  (8,0)*{};
  (-4,-15)*{}="b1";
  (4,-15)*{}="b2";
  "b2";"b1" **\crv{(5,-8) & (-5,-8)}; ?(.1)*\dir{>} ?(.95)*\dir{>}
  ?(.8)*\dir{}+(0,-.1)*{\bullet}+(-3,2)*{\scs g_3};
  (-4,15)*{}="t1";
  (4,15)*{}="t2";
  "t2";"t1" **\crv{(5,8) & (-5,8)}; ?(.15)*\dir{<} ?(.9)*\dir{<}
  ?(.4)*\dir{}+(0,-.2)*{\bullet}+(3,-2)*{\scs g_1};
  (0,0)*{\icbub{\scs \quad\; \l_i-1 + g_2}{i}};
    (7,-13)*{\scs i};
  (-7,13)*{\scs i};
  (-10,10)*{\l};
  \endxy} \label{eq_ident_decomp-nleqz}
\end{equation}

If $\lambda_i =0$, then we have:
\begin{equation}\label{eq_reduction-neqz}
  \xy 0;/r.17pc/:
  (14,8)*{\l};
  (-3,-10)*{};(3,5)*{} **\crv{(-3,-2) & (2,1)}?(1)*\dir{>};?(.15)*\dir{>};
    (3,-5)*{};(-3,10)*{} **\crv{(2,-1) & (-3,2)}?(.85)*\dir{>} ?(.1)*\dir{>};
  (3,5)*{}="t1";  (9,5)*{}="t2";
  (3,-5)*{}="t1'";  (9,-5)*{}="t2'";
   "t1";"t2" **\crv{(4,8) & (9, 8)};
   "t1'";"t2'" **\crv{(4,-8) & (9, -8)};
   "t2'";"t2" **\crv{(10,0)} ;
   (-6,-8)*{\scs i};
 \endxy\;\; = \;\;-
   \xy 0;/r.17pc/:
  (-8,8)*{\l};
  (0,0)*{\bbe{}};
  (-3,-8)*{\scs i};
 \endxy
\qquad \qquad
  \xy 0;/r.17pc/:
  (-14,8)*{\l};
  (3,-10)*{};(-3,5)*{} **\crv{(3,-2) & (-2,1)}?(1)*\dir{>};?(.15)*\dir{>};
    (-3,-5)*{};(3,10)*{} **\crv{(-2,-1) & (3,2)}?(.85)*\dir{>} ?(.1)*\dir{>};
  (-3,5)*{}="t1";  (-9,5)*{}="t2";
  (-3,-5)*{}="t1'";  (-9,-5)*{}="t2'";
   "t1";"t2" **\crv{(-4,8) & (-9, 8)};
   "t1'";"t2'" **\crv{(-4,-8) & (-9, -8)};
   "t2'";"t2" **\crv{(-10,0)} ;
   (6,-8)*{\scs i};
 \endxy \;\; = \;\;
  \xy 0;/r.17pc/:
  (-8,8)*{\l};
  (0,0)*{\bbe{}};
  (-3,-8)*{\scs i};
 \endxy
\end{equation}
\begin{equation}\label{eq_reduction-neqz_2}
 \vcenter{\xy 0;/r.17pc/:
  (-8,0)*{};
  (8,0)*{};
  (-4,10)*{}="t1";
  (4,10)*{}="t2";
  (-4,-10)*{}="b1";
  (4,-10)*{}="b2";(-6,-8)*{\scs i};(6,-8)*{\scs i};
  "t1";"b1" **\dir{-} ?(.5)*\dir{<};
  "t2";"b2" **\dir{-} ?(.5)*\dir{>};
  (10,2)*{\l};
  \endxy}
\;\; = \;\;
 \;\; - \;\;
 \vcenter{   \xy 0;/r.17pc/:
    (-4,-4)*{};(4,4)*{} **\crv{(-4,-1) & (4,1)}?(1)*\dir{>};
    (4,-4)*{};(-4,4)*{} **\crv{(4,-1) & (-4,1)}?(1)*\dir{<};?(0)*\dir{<};
    (-4,4)*{};(4,12)*{} **\crv{(-4,7) & (4,9)};
    (4,4)*{};(-4,12)*{} **\crv{(4,7) & (-4,9)}?(1)*\dir{>};
  (8,8)*{l};(-6,-3)*{\scs i};
     (6.5,-3)*{\scs i};
 \endxy}
   \qquad \quad
 \vcenter{\xy 0;/r.17pc/:
  (-8,0)*{};(-6,-8)*{\scs i};(6,-8)*{\scs i};
  (8,0)*{};
  (-4,10)*{}="t1";
  (4,10)*{}="t2";
  (-4,-10)*{}="b1";
  (4,-10)*{}="b2";
  "t1";"b1" **\dir{-} ?(.5)*\dir{>};
  "t2";"b2" **\dir{-} ?(.5)*\dir{<};
  (10,2)*{\l};
  \endxy}
\;\; = \;\;
 \;\; - \;\;
   \vcenter{\xy 0;/r.17pc/:
    (-4,-4)*{};(4,4)*{} **\crv{(-4,-1) & (4,1)}?(1)*\dir{<};?(0)*\dir{<};
    (4,-4)*{};(-4,4)*{} **\crv{(4,-1) & (-4,1)}?(1)*\dir{>};
    (-4,4)*{};(4,12)*{} **\crv{(-4,7) & (4,9)}?(1)*\dir{>};
    (4,4)*{};(-4,12)*{} **\crv{(4,7) & (-4,9)};
  (8,8)*{\l};(-6,-3)*{\scs i};  (6,-3)*{\scs i};
 \endxy}
\end{equation}
\end{enumerate}

Next, recall that an idempotent $e \maps b\to b$ in a category $\cal{C}$ is a morphism such that
$e^2 = e$.  The idempotent is said to split if there exist morphisms
$
 \xymatrix{ b \ar[r]^g & b' \ar[r]^h &b}
$
such that $e=h  g$ and $gh = \id_{b'}$.    
The Karoubi envelope $Kar(\cal{C})$ of a category $\cal{C}$ is a minimal 
enlargement of $\cal{C}$ in which all idempotents split.  Explicitly, in the category $Kar(\cal{C})$:
\begin{itemize}
  \item \textbf{Objects} are pairs $(b,e)$ where $e \maps b \to b$ is an idempotent of $\cal{C}$.
\item \textbf{Morphisms} are triples $(e,f,e') \maps (b,e) \to (b',e')$
where $f \maps b \to b'$ in $\cal{C}$ making the diagram
\begin{equation} \label{eq_Kar_morph}
 \xymatrix{
 b \ar[r]^f \ar[d]_e \ar[dr]^{f} & b' \ar[d]^{e'} \\ b \ar[r]_f & b'
 }
\end{equation}
commute, i.e. $fe=f=e'f$.
\end{itemize}
The identity morphism of $(b,e)$ is given by the idempotent $(e,e,e) \maps (b,e) \to (b,e)$.
When $\cal{C}$ is an additive category, we write $(b,e)\in Kar(\cal{C})$ as $\im e$ and we have 
$b \cong \im e \oplus \im (1-e)$ in $Kar(\cal{C})$.

The Karoubi envelope $\UcatD_Q(\slm) := Kar(\Ucat_Q(\mf{sl}_m))$ of the 2-category $\Ucat_Q(\mf{sl}_m)$ is the 2-category with the 
same objects as $\Ucat_Q(\mf{sl}_m)$ whose $\Hom$-categories are given by
\[
\UcatD_Q(\onel,\onelp) := Kar\big(\Ucat_Q(\onel,\onelp)\big).
\]
In particular, all idempotent 2-morphisms split in $\UcatD_Q(\onel,\onelp)$. The main result of \cite{KL3} is that there is an isomorphism of
$\Z[q,q^{-1}]$-algebras
\begin{equation} \label{eq_gamma}
  \gamma \maps K_0(\UcatD_Q(\mf{sl}_m)) \longrightarrow \; _{\Z[q^{\pm}]}\U(\mf{sl}_m)
\end{equation}
between the split Grothendieck ring $K_0(\UcatD_Q(\mf{sl}_m))$ and the integral form 
$_{\Z[q^{\pm}]}\U(\mf{sl}_m)$ of the Lusztig idempotent quantum enveloping algebra. 
Furthermore, Webster ~\cite{Web4} has shown that the images of the indecomposable 1-morphisms in $\UcatD_Q(\mf{sl}_m)$ in 
$K_0(\Ucat_Q(\mf{sl}_m))$ agree with the Lusztig canonical basis in $_{\Z[q^{\pm}]}\U(\mf{sl}_m)$.

\subsubsection{Thick calculus in $\UcatD_Q(\mf{sl}_m)$}

Typically the passage from a diagrammatically defined category to its Karoubi envelope results in the loss of a completely diagrammatic 
description of the resulting category.  However, the Karoubi envelope $\UcatD_Q(\mf{sl}_2)$ of the 2-category $\Ucat_Q(\mf{sl}_2)$ 
still admits a completely diagrammatic description~\cite{KLMS}.  In this case, one defines idempotent 2-morphisms 
$e_a \maps \cal{E}^a\onel \to \cal{E}^a\onel$ given by the composite of any reduced presentation of the longest braid word on $a$ strands 
together with a specific pattern of dots starting with $a-1$ dots on the top left-most strand, $a-2$ on the next strand, 
and ending with no dots on the right-most of the $a$ strands. An example is shown below for $a=4$:
\[
\xy 0;/r.15pc/:
 (-12,-20)*{}; (12,20) **\crv{(-12,-8) & (12,8)}?(1)*\dir{>};
 (-4,-20)*{}; (4,20) **\crv{(-4,-13) & (12,2) & (12,8)&(4,13)}?(1)*\dir{>};?(.88)*\dir{}+(0.1,0)*{\bullet};
 (4,-20)*{}; (-4,20) **\crv{(4,-13) & (12,-8) & (12,-2)&(-4,13)}?(1)*\dir{>}?(.86)*\dir{}+(0.1,0)*{\bullet};
 ?(.92)*\dir{}+(0.1,0)*{\bullet};
 (12,-20)*{}; (-12,20) **\crv{(12,-8) & (-12,8)}?(1)*\dir{>}?(.70)*\dir{}+(0.1,0)*{\bullet};
 ?(.90)*\dir{}+(0.1,0)*{\bullet};?(.80)*\dir{}+(0.1,0)*{\bullet};
 \endxy
 \qquad =: \qquad   \xy
 (0,0)*{\includegraphics[scale=0.4]{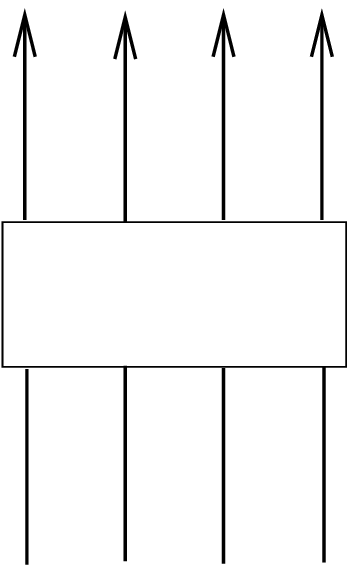}};
 (0,-0.5)*{e_a};
  \endxy
\]
in which we've used the box notation from \cite{KLMS} for this composite 2-morphism. Note that when considering 
$\Ucat_Q(\mf{sl}_2)$, the string diagrams need not carry a color $i$, since there is only one Dynkin node.

The divided power $\cal{E}^{(a)}\onel$ is defined in the Karoubi envelope $\UcatD_Q(\mf{sl}_2)$ as the pair
\[
\cal{E}^{(a)}\onel:= (\cal{E}^a\onel \{ \frac{a(a-1)}{2}\} , e_a)
\]
where the grading shift is necessary to get an isomorphism $\cal{E}^a\onel \cong \oplus_{[a]!}\cal{E}^{(a)}\onel$.
Here, and throughout, we set $\oplus_{f(q)} C = C^{\oplus a_k} \{k\} \oplus \cdots \oplus C^{\oplus a_l} \{l\}$ 
for $f(q) = a_{k}q^k + \cdots + a_{l}q^l \in \Z[q,q^{-1}]$ when $C$ is an object of a graded, additive category.
The divided power $\onel\cal{F}^{(a)}$ is then defined as the adjoint of $\cal{E}^{(a)}\onel$.  It was shown in \cite{KLMS} that splitting the idempotents $e_a$ by adding $\cal{E}^{(a)}\onel$ and $\cal{F}^{(b)}\onel$  gives rise 
to explicit decompositions of arbitrary 1-morphisms into indecomposable 1-morphisms using only the relations from $\Ucat_Q(\mf{sl}_2)$.  This allows for a strengthening of the categorification result to the case when 
we define $\Ucat_Q(\mf{sl}_2)$ by taking $\Z$-linear combinations of 2-morphisms.

It is possible to represent the 1-morphisms $\cal{E}^{(a)}\onel$ in $\UcatD_Q(\mf{sl}_2)$ by introducing an augmented graphical 
calculus of colored, thickened strands.  For example, the identity 2-morphism for $\cal{E}^{(a)}\onel$ is given
by the triple
\begin{equation}
 (e_a,e_a,e_a)
 \quad
 =
 \quad
 \left( e_a,\;
 \xy
 (0,0)*{\includegraphics[scale=0.5]{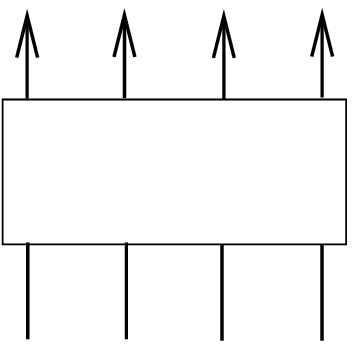}};
 (0,0)*{e_a}; (11,-5)*{\lambda};
  \endxy
 , e_a \right)
  \;\; =: \;\;     
\xy
(0,0)*{
\begin{tikzpicture} [scale=.5]
\draw[ultra thick, black!30!green, directed=.95] (0,0) to (0,4); 
\node[black!30!green,opacity=1] at (0,-.5) {$a$};
\node at (2,3) {$\lambda$};
\end{tikzpicture}
}
\endxy
\end{equation}
where we think of the colored label $a$ on the right as describing the thickness of the strand.
A downward oriented line of thickness $a$ similarly describes the identity
2-morphism for $\cal{F}^{(a)}\onel$ in $\UcatD_Q(\mf{sl}_2)$.

One can introduce further notation to describe natural 2-morphisms in $\UcatD_Q(\mf{sl}_2)$.  For example, using the shorthands
\[ 
  \xy
 (0,0)*{\includegraphics[scale=0.5]{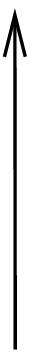}};
 (-3,-3)*{a};
  \endxy
  \quad : = \quad
  \xy
 (0,0)*{\includegraphics[scale=0.5]{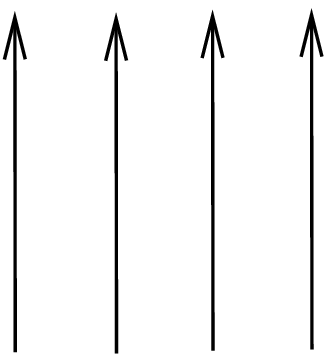}};
 (0,-11)*{\underbrace{\hspace{0.7in}}};  (0,-14)*{a};
  \endxy
 \qquad
 \qquad
   \xy
 (0,0)*{\includegraphics[scale=0.5]{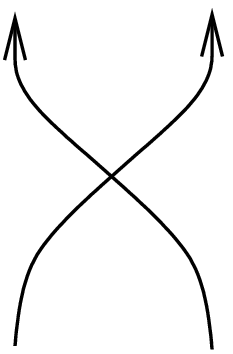}};
 (-7,-6)*{a};(7,-6)*{b};
  \endxy
  \quad := \quad
  \xy
 (0,0)*{\includegraphics[scale=0.5]{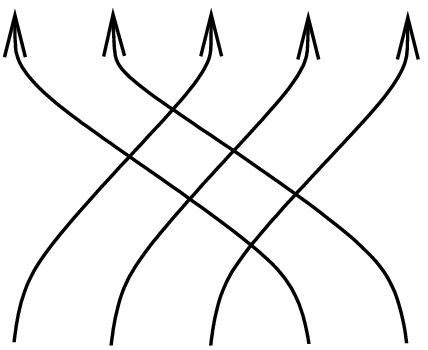}};
 (-5.5,-11)*{\underbrace{\hspace{0.45in}}};  (-5.5,-14)*{a};
 (7.5,-11)*{\underbrace{\hspace{0.25in}}};  (7.5,-14)*{b};
  \endxy
\]
there are 2-morphisms in $\UcatD_Q(\mf{sl}_2)$ given by
\begin{align}
\xy
(0,0)*{
\begin{tikzpicture} [scale=.6,fill opacity=0.2]
\draw[ultra thick, black!30!green] (0,-1.5) to (0,0);
\draw[ultra thick, black!30!green, directed=1] (0,0) to [out=30,in=270] (1,1.5);
\draw[ultra thick, black!30!green, directed=1] (0,0) to [out=150,in=270] (-1,1.5);
\node[black!30!green, opacity=1] at (0,-2) {\tiny$_{a+b}$};
\node[black!30!green, opacity=1] at (-1,2) {\tiny$a$};
\node[black!30!green, opacity=1] at (1,2) {\tiny$b$};
\end{tikzpicture}
}
\endxy
& :=
     \left(e_{a+b}, \;\xy
 (0,0)*{\includegraphics[scale=0.5]{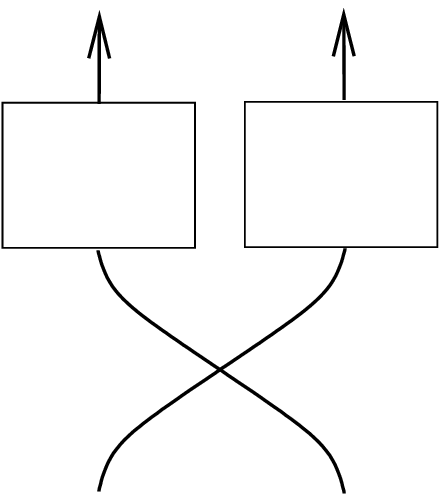}};
 (-7,-10)*{b};(7,-10)*{a};(-6,4)*{e_a};(6,4)*{e_b}; (14,-2)*{\l};
  \endxy \;, e_ae_b \right) \maps \cal{E}^{(a+b)}\onel \{t\} \to
  \cal{E}^{(a)}\cal{E}^{(b)}\onel\{t-ab\} \nn
  \\
\xy
(0,0)*{
\begin{tikzpicture} [scale=.6,fill opacity=0.2]
\draw[ultra thick, black!30!green, directed=1] (0,0) to (0,1.5);
\draw[ultra thick, black!30!green] (1,-1.5) to [out=90,in=330] (0,0);
\draw[ultra thick, black!30!green] (-1,-1.5) to [out=90,in=210] (0,0);
\node[black!30!green, opacity=1] at (0,2) {\tiny$_{a+b}$};
\node[black!30!green, opacity=1] at (-1,-2) {\tiny$a$};
\node[black!30!green, opacity=1] at (1,-2) {\tiny$b$};
\end{tikzpicture}
}
\endxy
  &:= 
  \left( e_a e_b,\;
     \xy
 (0,0)*{\includegraphics[scale=0.5,angle=180]{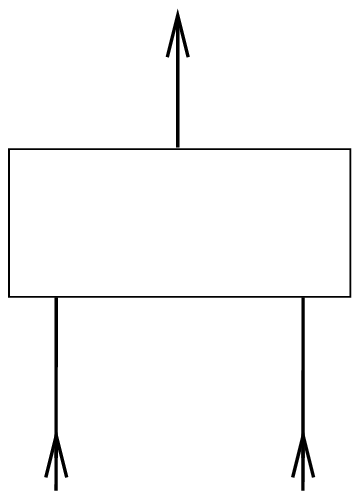}};
 (-8,-7)*{a};(8,-7)*{b};(0,1)*{e_{a+b}};(-5,10)*{a+b}; (11,6)*{\l};
  \endxy\;, e_{a+b} \right) \maps \cal{E}^{(a)}\cal{E}^{(b)}\onel \{t\} \to
  \cal{E}^{(a+b)}\onel\{t-ab\}. \nn
\end{align}
In computing the degree of the above diagrams, one must account for the shift in the definition of divided powers, e.g.
in the first diagram the degree shift in the divided power for $\cal{E}^{(a+b)}\onel$ is
$\frac{(a+b)(a+b-1)}{2}$, while the degree shift in the composite
$\cal{E}^{(a)}\cal{E}^{(b)}\onel$ is $\frac{a(a-1)}{2}+\frac{b(b-1)}{2}$, so that the net difference is
$\frac{2ab}{2} = ab$.  Both of the above diagrams in the thick calculus have degree $-ab$. 

The dot $2$-morphisms in $\Ucat_Q(\slnn{2})$ generalize to $2$-morphisms denoted by decorating a 
strand of thickness $a$ by an element of the algebra of symmetric polynomials in $a$ variables, e.g. 
$
\xy
(0,0)*{
\begin{tikzpicture} [scale=.3]
\draw[ultra thick, black!30!green, directed=.95] (0,0) to (0,4); 
\node[black!30!green,opacity=1] at (0,-.5) {$a$};
\node at (2,3) {$\lambda$};
\node at (0.35,2) {$\bullet^f$};
\end{tikzpicture}
}
\endxy
$
is a $2$-morphism in $\UcatD_Q(\mf{sl}_2)$ for any $f \in \Z[t_1,\ldots,t_a]^{S_a}$, the ring of symmetric polynomials in the variables $t_1, \dots, t_a$.
Recall that both the elementary symmetric polynomials $\{e_i\}_{i=1}^a$ and the complete symmetric polynomials
$\{h_i\}_{i=1}^a$ give generating sets for this algebra, and the Schur polynomials $\pi_\alpha$ give a linear basis,
where $\alpha=(\alpha_1,\ldots,\alpha_a)$ is a partition with at most $a$ parts.
The set $P(a,b)$ consists of the set of partitions $\alpha$ with at most $a$ parts and so that $\alpha_1 \leq b$. Given 
a partition $\alpha$, we denote the dual of the complementary partition by $\hat{\alpha}$; if the former is an element of 
$P(a,b)$ then the latter lies in $P(b,a)$. We refer the reader to \cite{KLMS} for the requisite facts about symmetric functions, 
and for complete details of the graphical calculus of $\UcatD_Q(\mf{sl}_2)$.

For general $m$ there does not yet exist a 
completely diagrammatic description of the Karoubi envelope of $\Ucat_Q(\slm)$.  In this case one currently lacks a 
set of diagrammatic relations needed to decompose arbitrary 1-morphisms into indecomposables, though explicit isomorphisms giving 
higher Serre relations were defined by Sto{\v{s}}i{\'c}~\cite{Stosic}. It will nevertheless be convenient to introduce a version of the 2-category 
$\Ucat_Q(\slm)$ where we have split the idempotents needed to define divided powers $\cal{E}_i^{(a)}\onel$ and $\cal{F}_i^{(a)}\onel$, 
but where we have not passed to the full Karoubi envelope.  
\begin{defn}
Let $\Ucatc_Q(\slm)$ denote the full sub-2-category of  $\UcatD_Q(\slm)$ with the same objects $\lambda \in X$ as $\UcatD_Q(\slm)$ and 
with 1-morphisms generated as a graded additive $\Bbbk$-linear category by  the 1-morphisms 
$\cal{E}_i\onel:=(\cal{E}_i \onel,\id_{\cal{E}_i \onel})$ and $\cal{E}^{(a)}_i\onel:= (\cal{E}_i^a\onel \{\frac{a(a-1)}{2}\} , e_a)$ 
and their adjoints.
\end{defn}
Diagrammatically, this 2-category can be described using thick strands carrying two labels, one indicating the thickness of the strand, and 
the other indicating the color $i \in I$ of the strand. Since this notation can become confusing, we'll typically use differing colors to denote 
the Dynkin node coloring a strand. We'll generally use the following `chromatic' ordering:
\[
\xy
(0,0)*{
\begin{tikzpicture} [scale=1]
\draw (-1.625,0) to (2.625,0);
\node[red] at (-1,0) {\Large$\bullet$};
\node[black!30!green] at (0,0) {\Large$\bullet$};
\node[blue] at (1,0) {\Large$\bullet$};
\node[black] at (2,0) {\Large$\bullet$};
\node at (-2,0) {$\cdots$};
\node at (3,0) {$\cdots$};
\node[red] at (-1,.5) {$i-1$};
\node[black!30!green] at (0,.5) {$i$};
\node[blue] at (1,.5) {$i+1$};
\node[black] at (2,.5) {$i+2$};
\end{tikzpicture}
}
\endxy
\]
of the Dynkin diagram in our thick calculus computations.
Since the thick strands are defined in terms of idempotents in thin strands, all the 2-morphisms can be studied using only the 
relations from $\Ucat_Q(\slm)$, together with the `pitchfork relation' from~\cite{Stosic}.

We'll study our foam $2$-categories using $2$-functors from $\Ucat_Q(\slm)$, $\Ucatc_Q(\slm)$, and other related $2$-categories.
\begin{defn}\label{def:2rep}
Let $\Ucat_Q$ denote one of the aforementioned $2$-categories. A \emph{$2$-representation} of $\Ucat_Q$ is a graded additive $\Bbbk$-linear 2-functor 
$\Ucat_Q \to \cal{K}$ for some graded, additive 2-category $\cal{K}$.
\end{defn}
To aid in our study of such $2$-functors, we collect here useful relations which follow from the defining relations of $\Ucat_Q(\slm)$ and the thick calculus of \cite{KLMS}; 
many of these relations are signed versions of those appearing in \cite{Stosic}. Throughout, we use the color scheme adopted above.

The following relations hold in $\Ucatc_Q(\slm)$:
\[
\xy
(0,0)*{

};
\endxy
\]

The KLR relation \eqref{eq_r2_ij-gen} for $(i,i+1)$ crossings generalizes to the thick calculus relations:
\begin{equation} \label{eq_r2_ij-gen-thick}
\xy
(0,0)*{
 %
};
\endxy
\end{equation}
and their downward counterparts:
\[
\xy
(0,0)*{
 %
};
\endxy \quad .
\]
Similarly, the mixed crossing relations \eqref{eq_mixed_EF} give:
\begin{align}
 \xy
(0,0)*{
 %
};
\endxy
. \nonumber
\end{align}
We will also extensively use the following special cases of the ``higher reduction to bubbles" relations from \cite[Formula 5.10]{KLMS} 
for the duration:
\begin{equation} \label{KLMSHigherBubbleRed}
 \xy
(0,0)*{
%
};
\endxy \quad
\end{equation}
Finally, the following signed isotopies of elementary diagrams can easily be deduced;  
below, the green line has thickness $k$ while the blue one has thickness $l$:
\begin{eqnarray*}
\xy
(0,0)*{
%
};
\endxy .
\end{eqnarray*}

\subsubsection{Useful relations in the categorified $q$-Schur algebra}

We will later see that our $2$-functors from $\Ucatc_Q(\slm)$ to the $\sln$ foam $2$-categories actually factor through (a version of) 
the categorified q-Schur algebra from \cite{MSV2}. We next recall the definition of this $2$-category and describe some basic diagrammatic 
relations which will be useful later.

Let $\Ucat_Q(\glm) = \displaystyle \bigoplus_{N \in \Z} \Ucat_Q(\slm)$, i.e. the direct sum of an infinite number of 
copies of $\Ucat_Q(\slm)$, each indexed by an integer $N$. The objects in this $2$-category correspond exactly to $\glm$ weights, 
denoted by $m$-tuples $\mathbf{a} = [a_1,\ldots,a_m]$ of 
integers, corresponding to the $\slm$ weight $\lambda$ with $\lambda_i=a_{i}-a_{i+1}$ in the copy of $\Ucat_Q(\slm)$ indexed by 
$N= \sum_{i=1}^m a_i$. This $2$-category admits an (almost) identical diagrammatic description to $\Ucat_Q(\slm)$, with the $\slm$ weights 
labeling the regions replaced by the appropriate $\glm$ weights.

\begin{defn}[\cite{MSV2,CKM}]\label{def:qSchur}
Let $\Ucat_Q(\glm)^{\geq 0}$ be the quotient of $\Ucat_Q(\glm)$ by 
(the ideal generated by the identity $2$-morphisms of the identity $1$-morphisms for) 
weights $\mathbf{a} = [a_1,\ldots,a_m]$ such that there exists $i$ with $a_i<0$. 
Let $\Ucatc_Q(\glm)^{(n)}$ be the $n$-bounded quotient of $\Ucat_Q(\glm)^{\geq 0}$, 
defined by additionally taking the quotient with respect to $\glm$ weights $\mathbf{a}$ such that there exists $i$ such that $a_i>n$.
\end{defn}

These $2$-categories admit analogous idempotent completed versions to $\Ucat_Q(\slm)$, e.g. the $2$-category $\Ucatc_Q(\glm)^{(n)}$ 
is obtained by taking the subcategory of the idempotent completion of $\Ucat_Q(\glm)$ generated by divided powers, then taking the 
quotient which kills the relevant weights.
The $2$-category $\Ucat_Q(\glm)^{\geq 0}$ is precisely the direct sum $\displaystyle \bigoplus_{d \in \N} \cal{S}(n,d)$ of the 
categorified q-Schur algebras from \cite{MSV2}. The definition of the $2$-category $\Ucatc_Q(\glm)^{(n)}$ is motivated by the 
fact that the $\sln$ foamation functors defined in Section \ref{sec:foams} factor through this $2$-category.

Our next result presents some formulas that hold in $\Ucatc_Q(\glm)^{\geq 0}$, hence in all categories $\Ucatc_Q(\glm)^{(n)}$. 
These are often reformulations or simplifications of formulas that appear in \cite{KLMS} and \cite{Stosic}.

\begin{prop} \label{prop:relBoundedQuotient}
In $\Ucatc_q(\glm)^{\geq 0}$ we have:
\begin{equation} \label{eq:EFqSchur}
\xy
(0,0)*{

};
\endxy
\end{equation}
and the generalizations\footnote{
Recall that a bubble decorated with $\pi_{\alpha}^{\spadesuit}$ has degree equal to $\deg(\pi_\alpha)$; 
see \cite[Section 4]{KLMS} for the precise definition of these decorated bubbles.}:
\begin{align}\label{eq:EFqSchur:gen}
(-1)^{\frac{a(a-1)}{2}+\frac{k(k-1)}{2}}
\xy
(0,0)*{
 %
};
\endxy
\quad 
\end{align}
where $\pi$ decorating the red bubble stands for 
$\pi_{\hat{\alpha}+\hat{\beta}}$. 
We also have:
\[
\xy
(0,0)*{
 %
};
\endxy
\]
\end{prop}

\begin{proof}
We'll begin by proving the first relation in \eqref{eq:EFqSchur}. 
To do so, we first deduce the relation
\begin{equation}\label{qS:1rel}
\xy
(0,0)*{
 %
};
\endxy
\end{equation}
in  $\Ucat_Q(\glm)$, and observe that the rightmost term is zero in the quotient $\Ucat_Q(\glm)^{\geq 0}$.
To verify this relation, first `explode' the thick strand
\[
\xy
(0,0)*{
 %
};
\endxy
\]
and recall that if $\lambda\geq 0$, we have:
$
\xy
(0,0)*{
 %
};
\endxy
$. We can perform this operation provided we have a non-negative weight to the right of the strand we wish to cross over to obtain:
\[
\xy
(0,0)*{
 %
};
\endxy
\]
where in the above picture, the dots come with multiplicity, decreasing from $a-1$ on the leftmost strand to zero on the rightmost strand. 
The number of strands on the left of the downward one in the middle part of the diagram is $\lfloor \frac{a}{2}\rfloor$, so that the number of strands on the right is $\lceil \frac{a}{2} \rceil$ and 
the $\slnn{2}$ weight right next to the downward strand is thus either $0$ or $1$. 
To continue moving the downward strand to the right, we use \eqref{eq_ident_decomp-nleqz} and thus at the end of this procedure we
obtain:
\[
\xy
(0,0)*{
 %
};
\endxy
\]
where in the sum the dots may come with multiplicity. 
Focusing on the bottom part of these morphisms, we have:
\[
\xy
(0,0)*{
 %
};
\endxy
\]
with $\lambda=a-2k-2$ and $r\leq 2k-a+1$. Using equations \eqref{eq_nil_rels} and \eqref{eq_nil_dotslide}, we find that such morphisms are zero unless $r \geq k$, hence 
$k = a-1 = r$. It follows that the top portion of the only surviving $2$-morphism in the sum is
\[
\xy
(0,0)*{
 %
};
\endxy
\]
which pairs with the bottom to give the desired term in \eqref{qS:1rel}.

We hence compute in  $\Ucat_Q(\glm)^{\geq 0}$:
\[
\xy
(0,0)*{
 %
};
\endxy
\]
and then one-by-one pull the downward strands to the right, as above until they merge with an upward one. 
This again produces a sum of terms, with a typical one taking the form:
\[
\xy
(0,0)*{
 %
};
\endxy
\]
where dots and other bubbles may be added. Each time a downward strand doesn't fuse with an upward one, 
it pulls all the way through to the far right, and hence the corresponding element is zero in $\Ucat_Q(\glm)^{\geq 0}$.

In the end, we find that
\[
\xy
(0,0)*{
 %
};
\endxy
\;=\;
\sum (\text{coeff})\;
\left(
\text{Dotted versions of}\;\xy
(0,0)*{
 %
};
\endxy
\right)
\]
and for degree reasons, the only admissible term is the one with no dots. 
Since $\Ucat_Q(\glm)^{\geq 0}$ admits $2$-representations with the left hand term mapping to a non-zero $2$-morphism, we know that 
the coefficient on the remaining term on the right-hand side is non-zero, and the coefficient in equation \eqref{eq:EFqSchur} can be deduced 
by composing both sides with thick caps and cups. The second formula in this equation can be deduced using similar methods.

We now turn our attention to equation \eqref{eq:EFqSchur:gen}. We compute
\[
\xy
(0,0)*{
 %
};
\endxy
\]
while on the other hand
\begin{align*}
\xy
(0,0)*{
 %
};
\endxy \quad .
\end{align*}
We deduce the first part of the equation by closing the red strand as a degree-zero bubble and observing that $t_{i,i+1}^{|\hat{\alpha}|}t_{i+1,i}^{|\hat{\alpha}|}=(-1)^{|\hat{\alpha}|}$. 
For the second part, we close the red strand in the middle of the previous computation and slide it back between both green strands.

We conclude by noting that the final relations in the Proposition follow from \eqref{eq:EFqSchur} using the higher reduction to bubbles equation \eqref{KLMSHigherBubbleRed}.
\end{proof}

Equation \eqref{eq:EFqSchur}, together with the value of degree zero thick bubbles, provides an isomorphism $[0,a] \cong [a,0]$ in $\Ucat_Q(\glm)^{\geq 0}$.  It will later correspond to a foam isotopy allowing us to 
`push away' a facet. Up to increasing the value of $m$, that is, up to the $2$-functor $\Ucat_Q(\glm)^{\geq 0} \rightarrow \Ucat_Q(\glnn{m+k})^{\geq 0}$ induced by $[a_1,\dots,a_m] \rightarrow [a_1,\dots,a_m,0,\dots,0]$, 
we can always assume that there is enough room to push facets away. We'll use this idea to assume that certain entries of $\glm$ weights are zero when proving certain relations in $\Ucat_Q(\glm)^{\geq 0}$, 
e.g. see Lemma \ref{lem:elemIso}.

\subsubsection{Rickard complexes}

For $\mf{sl}_2$, the quantum Weyl group action in \eqref{eq_qW} gives a reflection isomorphism between the $\lambda$ weight space of an $U_q(\mf{sl}_2)$-representation 
$V$ and the $-\lambda$ weight space. This isomorphism was categorified by Chuang and Rouquier in their elegant solution to the Brou\'{e} conjecture~\cite{CR}. 
Their work is closely related to a variant of the 2-category $\Ucat_Q(\mf{sl}_2)$ where the nilHecke algebra is replaced by the affine Hecke algebra and there is no grading. 
Cautis, Kamnitzer and Licata~\cite{CKL3} later developed analogous complexes in the context of $\Ucat_Q(\mf{sl}_2)$, generalizing Chuang and Rouquier's construction.

To categorify the elements $T_i1_{\lambda}$, it is clear from \eqref{eq_qW} that we must work in $\Ucatc_Q(\slm)$ so that we have lifts of divided powers. 
Also, the minus signs in the definition of the braid group 
generators suggest that we will have to pass to the $2$-category $Kom(\Ucatc_Q(\slm))$ of complexes in the $\Hom$-categories of $\Ucatc_Q(\slm)$. Indeed, Cautis 
presents\footnote{We take the opposite of Cautis's quantum and homological gradings in the 
definition of these complexes, to better fit with previous diagrammatic constructions of link homology. 
Sadly, this convention for homological degree is opposite to that used in (colored) Khovanov-Rozansky homology.
} 
the Rickard complexes in this context~\cite{Cautis}.
The braid group generators $T_i1_{\lambda}$ are categorified by complexes $\cal{T}_i\onel$ in $Kom(\Ucatc_Q(\slm))$ of the form:
\begin{equation} \label{Rickardn}
\cal{T}_i \onel =
\xymatrix{ \uwave{\cal{E}_i^{(-\l_i)} \onel} \ar[r]^-{d_1} & \cal{E}_i^{(-\l_i+1)} \cal{F}_i \onel \{1\} \ar[r]^-{d_2} & \cdots
\ar[r]^-{d_s} & \cal{E}_i^{(-\l_i+s)} \cal{F}_i^{(s)} \onel \{s\} \ar[r]^-{d_{s+1}} &\cdots}
\end{equation}
when $\l_i \leq 0$ and
\begin{equation} \label{Rickardp}
\cal{T}_i \onel =
\xymatrix{ \uwave{\cal{F}_i^{(\l_i)} \onel} \ar[r]^-{d_1} & \cal{F}_i^{(\l_i+1)} \cal{E}_i \onel \{1\} \ar[r]^-{d_2} & \cdots
\ar[r]^-{d_s} & \cal{F}_i^{(\l_i+s)} \cal{E}_i^{(s)} \onel \{s\} \ar[r]^-{d_{s+1}} &\cdots}
\end{equation}
when $\l_i \geq 0$.
The differential $d_k$ that appears in the first complex is conveniently expressed in thick calculus as
\[
d_k=
\xy
(0,0)*{
\begin{tikzpicture}[scale=1]
\draw [ultra thick, black!30!green, ->] (0,0) -- (0,2);
\draw [ultra thick, black!30!green, ->] (1,2) -- (1,0);
\draw [semithick, black!30!green, directed= .55] (1,1.2) .. controls (1,.5) and (0,.5) .. (0,1.2);
\node [black!30!green] at (0,-.5) {\footnotesize $-\lambda+k$}; 
\node [black!30!green] at (.5,1) {\footnotesize $1$}; 
\node [black!30!green] at (1,-.5) {\footnotesize $k-1$}; 
\node [black!30!green] at (1,2.5) {\footnotesize $k$}; 
\node at (1.7,1) {$\lambda$};
\end{tikzpicture}
};
\endxy
\]
where all strands are colored by the index $i \in I$. The differential in the second complex is defined similarly and the equality
$d^2=0$ follows easily from a computation in thick calculus.
Results of Cautis-Kamnitzer show the images of the complexes $\cal{T}_i\onel$ under any integrable 2-representation $\Ucatc_Q(\slm) \to \cal{K}$ 
satisfy braid relations up to homotopy in $Kom(\cal{K})$~\cite[Section 6]{CK}.

The images of the complexes $\cal{T}_i\onel$ under such a $2$-representation are invertible, up to homotopy, with inverses given by taking the left adjoint of the 
complex $\cal{T}_i\onel$ in the 2-category $Kom(\Ucatc_Q(\slm)$. More explicitly, they are given by
\[
\onel \cal{T}_i^{-1} =
\xymatrix{ \cdots \ar[r]^-{d_{s+1}^*} & \onel \cal{E}_i^{(s)} \cal{F}_i^{(-\l_i+s)} \{-s\} \ar[r]^-{d_s^*}
& \cdots \ar[r]^-{d_2^*}
& \onel \cal{E}_i \cal{F}_i^{(-\l_i+1)} \{-1\} \ar[r]^-{d_1^*} & \uwave{\onel \cal{F}_i^{(-\l_i)}} }
\]
when $\l_i \leq 0$ and
\[
\onel \cal{T}_i^{-1} =
\xymatrix{ \cdots \ar[r]^-{d_{s+1}^*} & \onel \cal{F}_i^{(s)} \cal{E}_i^{(\l_i+s)} \{-s\} \ar[r]^-{d_s^*}
& \cdots \ar[r]^-{d_2^*}
& \onel \cal{F}_i \cal{E}_i^{(\l_i+1)} \{-1\} \ar[r]^-{d_1^*} &\uwave{ \onel \cal{E}_i^{(\l_i)}} }
\]
when $\l_i \geq 0$.

Given a 2-representation $\cal{F} \maps \Ucatc_Q(\slm) \to \cal{K}$, the equivalences given by these complexes are highly 
nontrivial and have been applied to a variety of contexts ranging from the representation theory of the symmetric group~\cite{CR} to coherent sheaves 
on cotangent bundles~\cite{CKL2,CKL3}. We'll later see that the skein relations for $\sln$ link homology also follow from these complexes.

%
\section{Foams and foamation} \label{sec:foams}
%

In this section we define the $\mathfrak{sl}_n$ foam 2-categories and show that there exists a 2-representations from 
categorified quantum $\slm$ to this $2$-category. Using properties of these foamation functors, we'll show that 
analogs of closed foams can be evaluated in a purely combinatorial manner, and that the foam 2-category is 
non-degenerate.

\subsection{The $\mathfrak{sl}_n$ foam $2$-category}
To define the $2$-category of (enhanced) $\mathfrak{sl}_n$ foams, we'll begin by defining 
a pre-additive $2$-category $\tilde{\foam{n}{m}}$. After showing that the $2$-morphisms in this $2$-category
are graded, we'll then define $\foam{n}{m}$ to be the ``free additive $2$-category'' of the 
degree-zero part of $\tilde{\foam{n}{m}}$.

\begin{defn} \label{defn:foams}
$\tilde{\foam{n}{m}}$ is the $2$-category defined as follows:
\begin{itemize}
\item \textbf{Objects} are sequences $(a_1,\ldots,a_m)$ labeling $m$ points in the interval $[0,1]$
with $a_i \in \{0,1, \ldots, n\}$ and $N = \sum_{i=1}^m a_i$, together with a zero object. 
\item \textbf{1-morphisms} are left-directed enhanced $\sln$ webs, 
where $a$-labeled edges are directed out from points labeled by $a$ in the
domain and directed into such points in the codomain. No web edges are attached to points labeled by $0$.
\item \textbf{2-morphisms} are $\Bbbk$-linear combinations of enhanced $\sln$ foams - labeled, decorated 
singular surfaces with oriented seams whose generic slices are left-directed webs, 
generated by the following basic foams:
\begin{equation} 
\xy
(0,0)*{
\begin{tikzpicture} [scale=.6,fill opacity=0.2]
	\path[fill=blue] (2.25,3) to (.75,3) to (.75,0) to (2.25,0);
	\path[fill=red] (.75,3) to [out=225,in=0] (-.5,2.5) to (-.5,-.5) to [out=0,in=225] (.75,0);
	\path[fill=red] (.75,3) to [out=135,in=0] (-1,3.5) to (-1,.5) to [out=0,in=135] (.75,0);	
	\draw [very thick,directed=.55] (2.25,0) to (.75,0);
	\draw [very thick,directed=.55] (.75,0) to [out=135,in=0] (-1,.5);
	\draw [very thick,directed=.55] (.75,0) to [out=225,in=0] (-.5,-.5);
	\draw[very thick, red, directed=.55] (.75,0) to (.75,3);
	\draw [very thick] (2.25,3) to (2.25,0);
	\draw [very thick] (-1,3.5) to (-1,.5);
	\draw [very thick] (-.5,2.5) to (-.5,-.5);
	\draw [very thick,directed=.55] (2.25,3) to (.75,3);
	\draw [very thick,directed=.55] (.75,3) to [out=135,in=0] (-1,3.5);
	\draw [very thick,directed=.55] (.75,3) to [out=225,in=0] (-.5,2.5);
	\node [blue, opacity=1]  at (1.5,2.5) {\tiny{$_{a+b}$}};
	\node[red, opacity=1] at (-.75,3.25) {\tiny{$b$}};
	\node[red, opacity=1] at (-.25,2.25) {\tiny{$a$}};		
\end{tikzpicture}
};
\endxy
\quad \quad , \quad \quad
\xy
(0,0)*{
\begin{tikzpicture} [scale=.6,fill opacity=0.2]
	\path[fill=blue] (-2.25,3) to (-.75,3) to (-.75,0) to (-2.25,0);
	\path[fill=red] (-.75,3) to [out=45,in=180] (.5,3.5) to (.5,.5) to [out=180,in=45] (-.75,0);
	\path[fill=red] (-.75,3) to [out=315,in=180] (1,2.5) to (1,-.5) to [out=180,in=315] (-.75,0);	
	\draw [very thick,rdirected=.55] (-2.25,0) to (-.75,0);
	\draw [very thick,rdirected=.55] (-.75,0) to [out=315,in=180] (1,-.5);
	\draw [very thick,rdirected=.55] (-.75,0) to [out=45,in=180] (.5,.5);
	\draw[very thick, red, rdirected=.55] (-.75,0) to (-.75,3);
	\draw [very thick] (-2.25,3) to (-2.25,0);
	\draw [very thick] (1,2.5) to (1,-.5);
	\draw [very thick] (.5,3.5) to (.5,.5);
	\draw [very thick,rdirected=.55] (-2.25,3) to (-.75,3);
	\draw [very thick,rdirected=.55] (-.75,3) to [out=315,in=180] (1,2.5);
	\draw [very thick,rdirected=.55] (-.75,3) to [out=45,in=180] (.5,3.5);
	\node [blue, opacity=1]  at (-1.5,2.5) {\tiny{$_{a+b}$}};
	\node[red, opacity=1] at (.25,3.25) {\tiny{$b$}};
	\node[red, opacity=1] at (.75,2.25) {\tiny{$a$}};		
\end{tikzpicture}
};
\endxy
\end{equation}
\begin{equation}\label{ccgens} 
\xy
(0,0)*{
\begin{tikzpicture} [scale=.6,fill opacity=0.2]
	\path[fill=blue] (-.75,4) to [out=270,in=180] (0,2.5) to [out=0,in=270] (.75,4) .. controls (.5,4.5) and (-.5,4.5) .. (-.75,4);
	\path[fill=red] (-.75,4) to [out=270,in=180] (0,2.5) to [out=0,in=270] (.75,4) -- (2,4) -- (2,1) -- (-2,1) -- (-2,4) -- (-.75,4);
	\path[fill=blue] (-.75,4) to [out=270,in=180] (0,2.5) to [out=0,in=270] (.75,4) .. controls (.5,3.5) and (-.5,3.5) .. (-.75,4);
	\draw[very thick, directed=.55] (2,1) -- (-2,1);
	\path (.75,1) .. controls (.5,.5) and (-.5,.5) .. (-.75,1); 
	\draw [very thick, red, directed=.65] (-.75,4) to [out=270,in=180] (0,2.5) to [out=0,in=270] (.75,4);
	\draw[very thick] (2,4) -- (2,1);
	\draw[very thick] (-2,4) -- (-2,1);
	\draw[very thick,directed=.55] (2,4) -- (.75,4);
	\draw[very thick,directed=.55] (-.75,4) -- (-2,4);
	\draw[very thick,directed=.55] (.75,4) .. controls (.5,3.5) and (-.5,3.5) .. (-.75,4);
	\draw[very thick,directed=.55] (.75,4) .. controls (.5,4.5) and (-.5,4.5) .. (-.75,4);
	\node [red, opacity=1]  at (1.5,3.5) {\tiny{$_{a+b}$}};
	\node[blue, opacity=1] at (.25,3.4) {\tiny{$a$}};
	\node[blue, opacity=1] at (-.25,4.1) {\tiny{$b$}};	
\end{tikzpicture}
};
\endxy
\quad \quad , \quad \quad
\xy
(0,0)*{
\begin{tikzpicture} [scale=.6,fill opacity=0.2]
	\path[fill=blue] (-.75,-4) to [out=90,in=180] (0,-2.5) to [out=0,in=90] (.75,-4) .. controls (.5,-4.5) and (-.5,-4.5) .. (-.75,-4);
	\path[fill=red] (-.75,-4) to [out=90,in=180] (0,-2.5) to [out=0,in=90] (.75,-4) -- (2,-4) -- (2,-1) -- (-2,-1) -- (-2,-4) -- (-.75,-4);
	\path[fill=blue] (-.75,-4) to [out=90,in=180] (0,-2.5) to [out=0,in=90] (.75,-4) .. controls (.5,-3.5) and (-.5,-3.5) .. (-.75,-4);
	\draw[very thick, directed=.55] (2,-1) -- (-2,-1);
	\path (.75,-1) .. controls (.5,-.5) and (-.5,-.5) .. (-.75,-1); 
	\draw [very thick, red, directed=.65] (.75,-4) to [out=90,in=0] (0,-2.5) to [out=180,in=90] (-.75,-4);
	\draw[very thick] (2,-4) -- (2,-1);
	\draw[very thick] (-2,-4) -- (-2,-1);
	\draw[very thick,directed=.55] (2,-4) -- (.75,-4);
	\draw[very thick,directed=.55] (-.75,-4) -- (-2,-4);
	\draw[very thick,directed=.55] (.75,-4) .. controls (.5,-3.5) and (-.5,-3.5) .. (-.75,-4);
	\draw[very thick,directed=.55] (.75,-4) .. controls (.5,-4.5) and (-.5,-4.5) .. (-.75,-4);
	\node [red, opacity=1]  at (1.25,-1.25) {\tiny{$_{a+b}$}};
	\node[blue, opacity=1] at (-.25,-3.4) {\tiny{$b$}};
	\node[blue, opacity=1] at (.25,-4.1) {\tiny{$a$}};
\end{tikzpicture}
};
\endxy
\end{equation}
\begin{equation} \label{uzgens}
\xy
(0,0)*{
\begin{tikzpicture} [scale=.6,fill opacity=0.2]
	\path [fill=red] (4.25,-.5) to (4.25,2) to [out=165,in=15] (-.5,2) to (-.5,-.5) to 
		[out=0,in=225] (.75,0) to [out=90,in=180] (1.625,1.25) to [out=0,in=90] 
			(2.5,0) to [out=315,in=180] (4.25,-.5);
	\path [fill=red] (3.75,.5) to (3.75,3) to [out=195,in=345] (-1,3) to (-1,.5) to 
		[out=0,in=135] (.75,0) to [out=90,in=180] (1.625,1.25) to [out=0,in=90] 
			(2.5,0) to [out=45,in=180] (3.75,.5);
	\path[fill=blue] (.75,0) to [out=90,in=180] (1.625,1.25) to [out=0,in=90] (2.5,0);
	\draw [very thick,directed=.55] (2.5,0) to (.75,0);
	\draw [very thick,directed=.55] (.75,0) to [out=135,in=0] (-1,.5);
	\draw [very thick,directed=.55] (.75,0) to [out=225,in=0] (-.5,-.5);
	\draw [very thick,directed=.55] (3.75,.5) to [out=180,in=45] (2.5,0);
	\draw [very thick,directed=.55] (4.25,-.5) to [out=180,in=315] (2.5,0);
	\draw [very thick, red, directed=.75] (.75,0) to [out=90,in=180] (1.625,1.25);
	\draw [very thick, red] (1.625,1.25) to [out=0,in=90] (2.5,0);
	\draw [very thick] (3.75,3) to (3.75,.5);
	\draw [very thick] (4.25,2) to (4.25,-.5);
	\draw [very thick] (-1,3) to (-1,.5);
	\draw [very thick] (-.5,2) to (-.5,-.5);
	\draw [very thick,directed=.55] (4.25,2) to [out=165,in=15] (-.5,2);
	\draw [very thick, directed=.55] (3.75,3) to [out=195,in=345] (-1,3);
	\node [blue, opacity=1]  at (1.625,.5) {\tiny{$_{a+b}$}};
	\node[red, opacity=1] at (3.5,2.65) {\tiny{$b$}};
	\node[red, opacity=1] at (4,1.85) {\tiny{$a$}};		
\end{tikzpicture}
};
\endxy
\quad \quad , \quad \quad
\xy
(0,0)*{
\begin{tikzpicture} [scale=.6,fill opacity=0.2]
	\path [fill=red] (4.25,2) to (4.25,-.5) to [out=165,in=15] (-.5,-.5) to (-.5,2) to
		[out=0,in=225] (.75,2.5) to [out=270,in=180] (1.625,1.25) to [out=0,in=270] 
			(2.5,2.5) to [out=315,in=180] (4.25,2);
	\path [fill=red] (3.75,3) to (3.75,.5) to [out=195,in=345] (-1,.5) to (-1,3) to [out=0,in=135]
		(.75,2.5) to [out=270,in=180] (1.625,1.25) to [out=0,in=270] 
			(2.5,2.5) to [out=45,in=180] (3.75,3);
	\path[fill=blue] (2.5,2.5) to [out=270,in=0] (1.625,1.25) to [out=180,in=270] (.75,2.5);
	\draw [very thick,directed=.55] (4.25,-.5) to [out=165,in=15] (-.5,-.5);
	\draw [very thick, directed=.55] (3.75,.5) to [out=195,in=345] (-1,.5);
	\draw [very thick, red, directed=.75] (2.5,2.5) to [out=270,in=0] (1.625,1.25);
	\draw [very thick, red] (1.625,1.25) to [out=180,in=270] (.75,2.5);
	\draw [very thick] (3.75,3) to (3.75,.5);
	\draw [very thick] (4.25,2) to (4.25,-.5);
	\draw [very thick] (-1,3) to (-1,.5);
	\draw [very thick] (-.5,2) to (-.5,-.5);
	\draw [very thick,directed=.55] (2.5,2.5) to (.75,2.5);
	\draw [very thick,directed=.55] (.75,2.5) to [out=135,in=0] (-1,3);
	\draw [very thick,directed=.55] (.75,2.5) to [out=225,in=0] (-.5,2);
	\draw [very thick,directed=.55] (3.75,3) to [out=180,in=45] (2.5,2.5);
	\draw [very thick,directed=.55] (4.25,2) to [out=180,in=315] (2.5,2.5);
	\node [blue, opacity=1]  at (1.625,2) {\tiny{$_{a+b}$}};
	\node[red, opacity=1] at (3.5,2.65) {\tiny{$b$}};
	\node[red, opacity=1] at (4,1.85) {\tiny{$a$}};		
\end{tikzpicture}
};
\endxy
\end{equation}
\begin{equation}\label{MVgensup}
\xy
(0,0)*{
\begin{tikzpicture} [scale=.6,fill opacity=0.2]
	\path[fill=red] (-2.5,4) to [out=0,in=135] (-.75,3.5) to [out=270,in=90] (.75,.25)
		to [out=135,in=0] (-2.5,1);
	\path[fill=blue] (-.75,3.5) to [out=270,in=125] (.29,1.5) to [out=55,in=270] (.75,2.75) 
		to [out=135,in=0] (-.75,3.5);
	\path[fill=blue] (-.75,-.5) to [out=90,in=235] (.29,1.5) to [out=315,in=90] (.75,.25) 
		to [out=225,in=0] (-.75,-.5);
	\path[fill=red] (-2,3) to [out=0,in=225] (-.75,3.5) to [out=270,in=125] (.29,1.5)
		to [out=235,in=90] (-.75,-.5) to [out=135,in=0] (-2,0);
	\path[fill=red] (-1.5,2) to [out=0,in=225] (.75,2.75) to [out=270,in=90] (-.75,-.5)
		to [out=225,in=0] (-1.5,-1);
	\path[fill=red] (2,3) to [out=180,in=0] (.75,2.75) to [out=270,in=55] (.29,1.5)
		to [out=305,in=90] (.75,.25) to [out=0,in=180] (2,0);
	\draw[very thick, directed=.55] (2,0) to [out=180,in=0] (.75,.25);
	\draw[very thick, directed=.55] (.75,.25) to [out=225,in=0] (-.75,-.5);
	\draw[very thick, directed=.55] (.75,.25) to [out=135,in=0] (-2.5,1);
	\draw[very thick, directed=.55] (-.75,-.5) to [out=135,in=0] (-2,0);
	\draw[very thick, directed=.55] (-.75,-.5) to [out=225,in=0] (-1.5,-1);
	\draw[very thick, red, rdirected=.85] (-.75,3.5) to [out=270,in=90] (.75,.25);
	\draw[very thick, red, rdirected=.75] (.75,2.75) to [out=270,in=90] (-.75,-.5);	
	\draw[very thick] (-1.5,-1) -- (-1.5,2);	
	\draw[very thick] (-2,0) -- (-2,3);
	\draw[very thick] (-2.5,1) -- (-2.5,4);	
	\draw[very thick] (2,3) -- (2,0);
	\draw[very thick, directed=.55] (2,3) to [out=180,in=0] (.75,2.75);
	\draw[very thick, directed=.55] (.75,2.75) to [out=135,in=0] (-.75,3.5);
	\draw[very thick, directed=.65] (.75,2.75) to [out=225,in=0] (-1.5,2);
	\draw[very thick, directed=.55]  (-.75,3.5) to [out=225,in=0] (-2,3);
	\draw[very thick, directed=.55]  (-.75,3.5) to [out=135,in=0] (-2.5,4);
	\node[red, opacity=1] at (-2.25,3.375) {\tiny$c$};
	\node[red, opacity=1] at (-1.75,2.75) {\tiny$b$};	
	\node[red, opacity=1] at (-1.25,1.75) {\tiny$a$};
	\node[blue, opacity=1] at (0,2.75) {\tiny$_{b+c}$};
	\node[blue, opacity=1] at (0,.25) {\tiny$_{a+b}$};
	\node[red, opacity=1] at (1.37,2.5) {\tiny$_{a+b+c}$};	
\end{tikzpicture}
};
\endxy
\quad \quad , \quad \quad
\xy
(0,0)*{
\begin{tikzpicture} [scale=.6,fill opacity=0.2]
	\path[fill=red] (-2.5,4) to [out=0,in=135] (.75,3.25) to [out=270,in=90] (-.75,.5)
		 to [out=135,in=0] (-2.5,1);
	\path[fill=blue] (-.75,2.5) to [out=270,in=125] (-.35,1.5) to [out=45,in=270] (.75,3.25) 
		to [out=225,in=0] (-.75,2.5);
	\path[fill=blue] (-.75,.5) to [out=90,in=235] (-.35,1.5) to [out=315,in=90] (.75,-.25) 
		to [out=135,in=0] (-.75,.5);	
	\path[fill=red] (-2,3) to [out=0,in=135] (-.75,2.5) to [out=270,in=125] (-.35,1.5) 
		to [out=235,in=90] (-.75,.5) to [out=225,in=0] (-2,0);
	\path[fill=red] (-1.5,2) to [out=0,in=225] (-.75,2.5) to [out=270,in=90] (.75,-.25)
		to [out=225,in=0] (-1.5,-1);
	\path[fill=red] (2,3) to [out=180,in=0] (.75,3.25) to [out=270,in=45] (-.35,1.5) 
		to [out=315,in=90] (.75,-.25) to [out=0,in=180] (2,0);				
	\draw[very thick, directed=.55] (2,0) to [out=180,in=0] (.75,-.25);
	\draw[very thick, directed=.55] (.75,-.25) to [out=135,in=0] (-.75,.5);
	\draw[very thick, directed=.55] (.75,-.25) to [out=225,in=0] (-1.5,-1);
	\draw[very thick, directed=.45]  (-.75,.5) to [out=225,in=0] (-2,0);
	\draw[very thick, directed=.35]  (-.75,.5) to [out=135,in=0] (-2.5,1);	
	\draw[very thick, red, rdirected=.75] (-.75,2.5) to [out=270,in=90] (.75,-.25);
	\draw[very thick, red, rdirected=.85] (.75,3.25) to [out=270,in=90] (-.75,.5);
	\draw[very thick] (-1.5,-1) -- (-1.5,2);	
	\draw[very thick] (-2,0) -- (-2,3);
	\draw[very thick] (-2.5,1) -- (-2.5,4);	
	\draw[very thick] (2,3) -- (2,0);
	\draw[very thick, directed=.55] (2,3) to [out=180,in=0] (.75,3.25);
	\draw[very thick, directed=.55] (.75,3.25) to [out=225,in=0] (-.75,2.5);
	\draw[very thick, directed=.55] (.75,3.25) to [out=135,in=0] (-2.5,4);
	\draw[very thick, directed=.55] (-.75,2.5) to [out=135,in=0] (-2,3);
	\draw[very thick, directed=.55] (-.75,2.5) to [out=225,in=0] (-1.5,2);
	\node[red, opacity=1] at (-2.25,3.75) {\tiny$c$};
	\node[red, opacity=1] at (-1.75,2.75) {\tiny$b$};	
	\node[red, opacity=1] at (-1.25,1.75) {\tiny$a$};
	\node[blue, opacity=1] at (-.125,2.25) {\tiny$_{a+b}$};
	\node[blue, opacity=1] at (-.125,.75) {\tiny$_{b+c}$};
	\node[red, opacity=1] at (1.35,2.75) {\tiny$_{a+b+c}$};	
\end{tikzpicture}
};
\endxy
\end{equation}
modulo isotopies (which preserve the left-directed generic slice condition) and local relations. 
Components of the complement of the singular set, which we'll refer to as \emph{facets}, 
carry labels in $\{1, \ldots, n\}$ and each $a$-labeled component is decorated by an element of the ring 
of symmetric functions in $a$ variables.
If a component is undecorated it is implicitly decorated by the unit element and, for a $1$-labeled facet, the decoration $\bullet^p$ 
denotes the polynomial $t^p$ in the ring $\C[t]$ of symmetric polynomials in one variable 
(if the dot doesn't carry any label, then it corresponds to $t$). 
Vertical and horizontal composition is given by glueing 
along the relevant boundaries (and multiplying the symmetric functions decorating the same facet).
\end{itemize}

Note that the foams
\begin{equation}\label{MVgensdown} 
\xy
(0,0)*{

};
\endxy
\end{equation}
can be constructed via isotopy from the generators. The local foam relations are as follows, where 
any foam carrying a facet label outside the set $\{0,\ldots,n\}$ is zero and $0$-labeled edges/facets are 
understood to be erased: \\

\noindent \textbf{The Matveev-Piergalini\footnote{We adopt the terminology for this relation from \cite{MSV}.}
relations:}
\begin{equation}\label{MVrel}
\xy
(0,-8.8)*{
%
};
\endxy
\end{equation}

\noindent \textbf{The dot migration relation:} Let $c_{\alpha, \beta}^\gamma$ be the Littlewood-Richardson coefficients \cite{Macdonald}, then
\begin{equation}\label{dotmigrel}
\xy
(0,0)*{
%
};
\endxy
\end{equation}

\noindent \textbf{The blister relations:} Let $\alpha=(\alpha_1, \ldots,\alpha_{k+1})$ be a partition, $\alpha' = (\alpha_2, \ldots,\alpha_{k+1})$, 
and let $f$ and $g$ be (any) symmetric functions, then
\begin{equation}\label{blisterrel}
\xy
(0,0)*{
%
};
\endxy
\end{equation}
\noindent \textbf{The KLR relations:} For $a,c \geq 1$,
\begin{equation}\label{R2KLR1rel}
\xy
(0,0)*{
%
};
\endxy
\end{equation}
\noindent \textbf{The pitchfork relations:}
\begin{equation}\label{PFrel1}
\xy
(0,0)*{
%
};
\endxy
\end{equation}
\end{defn}

\begin{rem}\label{rem:rel}
The above equations do not constitute a minimal collection of relations; for example, we could have only imposed 
\begin{equation*}
\xy
(0,0)*{
%
};
\endxy
\end{equation*}
in equations \eqref{dotmigrel} and \eqref{blisterrel}. However, each of the foam relation equations contains a necessary relation (for certain values of the 
labelings and decorations), and we've decided to include the general relations of each form.
Additionally, one can check that each of these foam relations is a consequence of a relation in $\Ucatc_Q(\slm)$ under the $2$-representation 
given in Theorem \ref{thm:2functor}.
\end{rem}

The $2$-morphisms in this category admit a grading. In order to define this, we first introduce a 
weighted Euler characteristics $\chi_n$ for $\sln$ webs and foams. Note that both webs and foams are in particular CW-complexes;
placing an arbitrary CW-complex structure on a web or foam, we assign each cell $c$ a number $\gamma(c)$ 
depending on its labeling, and for $0$- and $1$-cells depending also on its neighborhood in the CW-complex:
\begin{itemize}
 \item $2$-cells: \[ c=\;\;\xy
(0,0)*{
\begin{tikzpicture} [scale=.6,fill opacity=0.2,decoration={markings, mark=at position 0.5 with {\arrow{>}}; }]
\filldraw [dashed, fill=red] (0,0) rectangle (2,2);
\node[red,opacity=1] at (1,1) {\tiny $k$};
\end{tikzpicture}};
\endxy \quad \mapsto \quad \gamma(c)=k(n-k).\]

\item $1$-cells: \[ c=\;\;
\xy
(0,0)*{
\begin{tikzpicture} [scale=.6,fill opacity=0.2,decoration={markings, mark=at position 0.5 with {\arrow{>}}; }]
\filldraw [dashed, fill=red] (0,0) rectangle (2,2);
\draw [very thick] (1,0) -- (1,2);
\node[red,opacity=1] at (.5,1) {\tiny $k$};
\node[red,opacity=1] at (1.5,1) {\tiny $k$};
\end{tikzpicture}};
\endxy 
\quad, \;\;
\xy
(0,0)*{
\begin{tikzpicture} [scale=.6,fill opacity=0.2,decoration={markings, mark=at position 0.5 with {\arrow{>}}; }]
\filldraw [dashed, fill=red] (1,0) rectangle (2,2);
\draw [very thick] (1,0) -- (1,2);
\node[red,opacity=1] at (1.5,1) {\tiny $k$};
\end{tikzpicture}};
\endxy 
\quad, \;\;
\text{or} \quad
\xy
(0,0)*{
\begin{tikzpicture} [scale=.6,fill opacity=0.2,decoration={markings, mark=at position 0.5 with {\arrow{>}}; }]
\draw [very thick] (1,0) -- (1,2);
\node[opacity=1] at (1.5,1) {\tiny $k$};
\end{tikzpicture}};
\endxy 
\quad \mapsto \quad \gamma(c)=k(n-k),\] 
and
\[ c=\;\;\xy
(0,0)*{
\begin{tikzpicture} [scale=.6,fill opacity=0.2,decoration={markings, mark=at position 0.5 with {\arrow{>}}; }]
\filldraw [dashed, fill=blue] (1,0) to (2.2,.5) -- (2.2,2.5) to (1,2) -- (1,0);
\filldraw [dashed, fill=blue] (1,0) to (1.8,-.5) -- (1.8,1.5) to (1,2) -- (1,0);
\filldraw [dashed, fill=red] (-.25,0) rectangle (1,2);
\draw [very thick] (1,0) -- (1,2);
\node[red,opacity=1] at (.35,1) {\tiny $k+l$};
\node[blue,opacity=1] at (1.4,.9) {\tiny $k$};
\node[blue,opacity=1] at (2,1.2) {\tiny $l$};
\end{tikzpicture}};
\endxy \quad \mapsto \quad \gamma(c)=(k+l)(n-k-l)+kl.\]

\item $0$-cells: \[ c=\;\;
\xy
(0,0)*{
\begin{tikzpicture} [scale=.6,fill opacity=0.2,decoration={markings, mark=at position 0.5 with {\arrow{>}}; }]
\filldraw [dashed, fill=red] (0,0) rectangle (2,2);
\draw [very thick] (1,0) -- (1,2);
\node[red,opacity=1] at (.5,1) {\tiny $k$};
\node[red,opacity=1] at (1.5,1) {\tiny $k$};
\node[opacity=1] at (1,1) {$\bullet$};
\end{tikzpicture}};
\endxy \quad, \;\;
\xy
(0,0)*{
\begin{tikzpicture} [scale=.6,fill opacity=0.2,decoration={markings, mark=at position 0.5 with {\arrow{>}}; }]
\filldraw [dashed, fill=red] (1,0) rectangle (2,2);
\draw [very thick] (1,0) -- (1,2);
\node[red,opacity=1] at (1.5,1) {\tiny $k$};
\node[opacity=1] at (1,1) {$\bullet$};
\end{tikzpicture}};
\endxy \quad, \;\;
\xy
(0,0)*{
\begin{tikzpicture} [scale=.6,fill opacity=0.2,decoration={markings, mark=at position 0.5 with {\arrow{>}}; }]
\draw [very thick] (1,0) -- (1,2);
\node[opacity=1] at (1.5,1.5) {\tiny $k$};
\node[opacity=1] at (1.5,.5) {\tiny $k$};
\node[opacity=1] at (1,1) {$\bullet$};
\end{tikzpicture}};
\endxy \;\; , \;\; \text{or} \quad
\xy
(0,0)*{
\begin{tikzpicture} [scale=.6,fill opacity=0.2,decoration={markings, mark=at position 0.5 with {\arrow{>}}; }]
\draw [very thick] (1,.25) -- (1,1.75);
\node[opacity=1] at (1.5,1) {\tiny $k$};
\node[opacity=1] at (1,1.75) {$\bullet$};
\end{tikzpicture}};
\endxy \quad
\mapsto \quad \gamma(c)=k(n-k),\]
\[ c=\;\;
\xy
(0,0)*{
\begin{tikzpicture} [scale=.6,fill opacity=0.2,decoration={markings, mark=at position 0.5 with {\arrow{>}}; }]
\filldraw [dashed, fill=blue] (1,0) to (2.2,.5) -- (2.2,2.5) to (1,2) -- (1,0);
\filldraw [dashed, fill=blue] (1,0) to (1.8,-.5) -- (1.8,1.5) to (1,2) -- (1,0);
\filldraw [dashed, fill=red] (-.25,0) rectangle (1,2);
\draw [very thick] (1,0) -- (1,2);
\node[red,opacity=1] at (.35,1) {\tiny $k+l$};
\node[blue,opacity=1] at (1.4,.9) {\tiny $k$};
\node[blue,opacity=1] at (2,1.2) {\tiny $l$};
\node[opacity=1] at (1,1) {$\bullet$};
\end{tikzpicture}};
\endxy \quad, \;\; 
\xy
(0,0)*{
\begin{tikzpicture} [scale=.6,fill opacity=0.2,decoration={markings, mark=at position 0.5 with {\arrow{>}}; }]
\filldraw [dashed, fill=blue] (1,0) to (2.2,.5) -- (2.2,2.5) to (1,2) -- (1,0);
\filldraw [dashed, fill=blue] (1,0) to (1.8,-.5) -- (1.8,1.5) to (1,2) -- (1,0);
\filldraw [dashed, fill=red] (-.25,0) rectangle (1,2);
\draw [dashed] (1,0) -- (1,2);
\draw[very thick] (1,2) to (-.25,2);
\draw[very thick] (1,2) to (1.8,1.5);
\draw[very thick] (1,2) to (2.2,2.5);
\node[red,opacity=1] at (.35,1) {\tiny $k+l$};
\node[blue,opacity=1] at (1.4,.9) {\tiny $k$};
\node[blue,opacity=1] at (2,1.2) {\tiny $l$};
\node[opacity=1] at (1,2) {$\bullet$};
\end{tikzpicture}};
\endxy \quad, \;\; \text{or} \quad
\xy
(0,0)*{
\begin{tikzpicture} [scale=.6,fill opacity=0.2,decoration={markings, mark=at position 0.5 with {\arrow{>}}; }]
\draw[very thick] (0,1) to (0,2);
\draw[very thick] (.865,.5) to (0,1); 
\draw[very thick] (-.865,.5) to (0,1);
\node[opacity=1] at (.5,2) {\tiny $k+l$};
\node[opacity=1] at (-.865,.875) {\tiny $k$};
\node[opacity=1] at (.865,.875) {\tiny $l$};
\node[opacity=1] at (0,1) {$\bullet$};
\end{tikzpicture}};
\endxy 
\mapsto \quad \gamma(c)=(k+l)(n-k-l)+kl,\] and
\[ c=\;\;\xy
(0,0)*{
\begin{tikzpicture} [scale=.6,fill opacity=0.2,decoration={markings, mark=at position 0.5 with {\arrow{>}}; }]
\filldraw [dashed, fill=red] (0,1) -- (2.2,1) to (2.5,2.3) -- (.3,2.3) to (0,1);
\filldraw [dashed, fill=blue] (1,0) -- (2.2,0) -- (2.2,2) -- (1,2) -- (1,0);
\draw [very thick] (0,1) -- (2.2,1);
\filldraw [dashed, fill=blue] (1,0) to (1.8,-.7) -- (1.8,1.3) to (1,2) -- (1,0);
\filldraw [dashed, fill=red] (0,0) rectangle (1,2);
\draw [very thick] (1,0) -- (1,2);
\node[blue,opacity=1] at (1.6,-.2) {\tiny $k$};
\node[blue,opacity=1] at (2,1.7) {\tiny $l$};
\node[red,opacity=1] at (2.2,2.1) {\tiny $m$};
\node[blue,opacity=1] at (1.7,.7) {\tiny $l+m$};
\node[red,opacity=1] at (.45,1.5) {\tiny $k+l$};
\node[red,opacity=1] at (.45,.5) {\tiny $k+l$};
\node[red,opacity=1] at (.4,.2) {\tiny $+m$};
\node[opacity=1] at (1,1) {$\bullet$};
\end{tikzpicture}};
\endxy \quad \mapsto \quad \gamma(c)=(k+l+m)(n-k-l-m)+kl+km+lm.\] 
\end{itemize}
Given a web or a foam $C$ with decomposition into cells $c_{ij}$, where $j=0,1,2$ denotes the cell dimension, 
we define the weighted Euler characteristic of $C$ to be:
\begin{equation} \label{Euler_char}
\chi_n(C)=\sum_{j=0}^2 (-1)^j \left(\sum_i \gamma(c_{ij})\right).
\end{equation}

\begin{defn}
Let $F$ be an undecorated foam between webs $W_{Bottom}$ and $W_{Top}$. Define the degree of $F$ as:
\begin{equation} \label{degree}
\deg(F)= - \chi_n(F)+\frac{1}{2} (\chi_n(W_{Top})+\chi_n(W_{Bottom})).
\end{equation}
\end{defn}

The next lemma ensures that the degree defined above is well-behaved with respect to the 2-categorical structure 
of $\tilde{\foam{n}{m}}$. In particular, identity foams are in degree zero and degree is additive with respect to both 
horizontal and vertical composition.

\begin{lem}
Let $W$ be a web and let $F_i$ be foams which are composable as in the schematic below, then
\begin{align*}
 \deg(id_W)&=0 \nonumber \\
\deg \left(\; \xy
(0,0)*{
\begin{tikzpicture} [scale=0.8]
  \draw (0,0) to (1,0) to (1,1);
  \draw (0,1) rectangle (1,2);
  \draw (-1,0) rectangle (0,1);
 \draw (-1,0) to (-1.25,.25) to (-1.25,1.25) to (-.25,1.25) to (-.25,2.25) to (.75,2.25) to (1,2);
 \draw (-1.25,1.25) to (-1,1);
 \draw (-.25,1.25) to (0,1);
 \draw (-.25,2.25) to (0,2);
  \node at (0.5,0.5) {$F_1$};
   \node at (0.5,1.5) {$F_2$};
  \node at (-0.5,0.5) {$F_3$};
\end{tikzpicture}};
\endxy\;\right) &= \deg(F_1)+\deg(F_2)+\deg(F_3). 
\end{align*}
\end{lem}
\begin{proof}
The foam $id_W$ has same web $W$ on the top and the bottom, which is a deformation retract of the foam. The formulae 
for $\gamma(c)$ above then imply that 
\[
\deg(id_{W})= - \chi_n(id_{W})+\frac{1}{2}(\chi_n(W)+\chi_n(W)) =\chi_n(W) - \chi_n(W) =0.
\]
For the second equality, we denote the top and bottom webs of $F_i$ by $(W_i)_{Top}$ and $(W_i)_{Bottom}$ and compute
\begin{align*}
\deg \left(\; \xy
(0,0)*{
\begin{tikzpicture} [scale=0.8]
  \draw (0,1) to (0,0) to (1,0) to (1,1);
  \draw (0,1) rectangle (1,2);
\draw (0,0) to (-.25,.25) to (-.25,2.25) to (.75,2.25) to (1,2);
\draw (0,1) to (-.25,1.25);
\draw (0,2) to (-.25,2.25);
  \node at (0.5,0.5) {$F_1$};
   \node at (0.5,1.5) {$F_2$};
\end{tikzpicture}};
\endxy \;
\right)&= - \chi_n \left(\; \xy
(0,0)*{
\begin{tikzpicture} [scale=0.8]
  \draw (0,1) to (0,0) to (1,0) to (1,1);
  \draw (0,1) rectangle (1,2);
\draw (0,0) to (-.25,.25) to (-.25,2.25) to (.75,2.25) to (1,2);
\draw (0,1) to (-.25,1.25);
\draw (0,2) to (-.25,2.25);
  \node at (0.5,0.5) {$F_1$};
   \node at (0.5,1.5) {$F_2$};
\end{tikzpicture}};
\endxy \;
\right) + \frac{1}{2}(\chi_n((W_2)_{Top})+\chi_n((W_1)_{Bottom})) \\
&= - \chi_n(F_1) - \chi_n(F_2) + \chi_n((W_1)_{Top}) + \frac{1}{2}\chi_n((W_2)_{Top}) + \frac{1}{2}\chi_n((W_1)_{Bottom}) \\
& = \deg(F_1)+ \deg(F_2) .
\end{align*}
The result for horizontal composition with $F_3$ follows similarly.
\end{proof}

Extending the above definition of degree to decorated foams by shifting the grading by the sum of the degrees of the decorations 
(in the respective rings of symmetric functions, where the generators are of degree $2$), we see that the relations in $\tilde{\foam{n}{m}}$ are 
degree-homogeneous, so this $2$-category is graded. We can now define the $\mathfrak{sl}_n$ foam $2$-category by introducing a 
formal grading on $1$-morphisms and considering only degree-preserving $2$-morphisms.

\begin{defn}
Let $\foam{n}{m}$ be the $2$-category defined as follows:
\begin{itemize}
\item \textbf{Objects} are the same as in $\tilde{\foam{n}{m}}$.
\item \textbf{$1$-morphisms} are (formally) $\Z$-graded formal direct sums of $1$-morphisms $W$ in $\tilde{\foam{n}{m}}$.
\item \textbf{$2$-morphisms} are matrices whose entries are $2$-morphisms 
$q^{d_1}W_1 \stackrel{\sum \alpha_i F_i}{\longrightarrow} q^{d_2}W_2$ in $\tilde{\foam{n}{m}}$
with $\deg{F_i} = d_2 - d_1$ for all $i$.
\end{itemize}
\end{defn}
Here, the powers of $q$ denote the formal grading on webs.

\begin{rem}
Using arguments similar to the proof of \cite[Theorem 5.1]{KLMS},  equations \eqref{dotmigrel}, \eqref{blisterrel}, and \eqref{nH3rel} imply that 
\begin{equation}\label{digon_web_iso}
\xy
(0,0)*{
\begin{tikzpicture}
\node[rotate=-45] at (0,0){$\xy
(0,0)*{
\begin{tikzpicture} [scale=.6]
	\draw[very thick, directed=.55] (2,-2) -- (.75,-2);
	\draw[very thick, directed=.55] (-.75,-2) -- (-2,-2);
	\draw[very thick, directed=.55] (.75,-2) .. controls (.5,-2.5) and (-.5,-2.5) .. (-.75,-2);
	\draw[very thick, directed=.55] (.75,-2) .. controls (.5,-1.5) and (-.5,-1.5) .. (-.75,-2);
\end{tikzpicture}
};
\endxy$};
\node at (-.25,-.5) {\small$a$};
\node at (.25,.5) {\small$b$};
\end{tikzpicture}
};
\endxy
\cong \quad
\bigoplus_{a+b \brack a}
\xy
(0,0)*{
\begin{tikzpicture}
\node[rotate=-45] at (0,0){$\xy
(0,0)*{
\begin{tikzpicture} [scale=.4]
	\draw[very thick, directed=.55] (2,-2) -- (-2,-2);
\end{tikzpicture}
};
\endxy$};
\node at (.25,.5) {\small$a+b$};
\end{tikzpicture}
};
\endxy 
\end{equation}
in $\foam[a+b]{n}{1}$.
\end{rem}

In the following subsections, we prove that $\foam{n}{m}$ is non-degenerate, using methods from higher representation theory. The reader only 
interested in the construction of the link homology theory can skip ahead to Section \ref{sec:link_hom}.

\subsection{The foamation $2$-functors}

We now aim to show that for each $m,n \geq 2$ and $N  \geq 1$ there exists a $2$-representation
$\Phi_n  : \Ucat_Q(\mf{sl}_m)  \to  \foam{n}{m}$. This $2$-functor extends the diagrammatic skew-Howe functor 
of Cautis-Kamnitzer-Morrison \cite{CKM}.
We'll establish this result via a sequence of lemmata; recall that for the duration we take $t_{i,i+1}=-1$, $t_{i+1,i}=1$, 
and $t_{i,j}=1$ when $|i-j| \geq 2$.
Our first result says that there is a $2$-functor from the `KLR' part 
of $\Ucat_Q(\slm)$, i.e. the full $2$-subcategory $\Ucat_Q^+(\mf{sl}_m)$ generated by the $\cal{E}_i$'s, to $\foam{n}{m}$.

\begin{lem} \label{lemKLR}
For each $m,n\geq 2$ and $N\geq 1$, there is a $2$-functor
$\Phi_n  \maps \Ucat_Q^+(\mf{sl}_m)  \to  \foam{n}{m}$ defined on objects and $1$-morphisms 
analogously\footnote{The labeling of webs in \eqref{eq_WsH} also determines the labeling of facets and boundary webs in the equations below.}
to equation \eqref{eq_WsH}, on single strand $2$-morphisms by:
\[
\Phi_n \left(
 \xy 0;/r.17pc/:
 (0,7);(0,-7); **\dir{-} ?(.75)*\dir{>};
 (7,3)*{ \scs \lambda};
 (-2.5,-6)*{\scs i};
 (-10,0)*{};(10,0)*{};
 \endxy
\right) = \quad \Efoam[.5] \quad , \quad
\Phi_n \left(
 \xy 0;/r.17pc/:
 (0,7);(0,-7); **\dir{-} ?(.75)*\dir{>};
 (0,0)*{\bullet};
 (7,3)*{ \scs \lambda};
 (-2.5,-6)*{\scs i};
 (-10,0)*{};(10,0)*{};
 \endxy
\right) = \quad \dotEfoam[.5]
\]
and on crossings by:
\[
\qquad
\Phi_n \left(
\xy 0;/r.20pc/:
	(0,0)*{\xybox{
    	(-4,-4)*{};(4,4)*{} **\crv{(-4,-1) & (4,1)}?(1)*\dir{>} ;
    	(4,-4)*{};(-4,4)*{} **\crv{(4,-1) & (-4,1)}?(1)*\dir{>};
    	(-5.5,-3)*{\scs i};
     	(5.5,-3)*{\scs i};
     	(9,1)*{\scs  \lambda};
     	(-10,0)*{};(10,0)*{};}};
\endxy
\right) = \quad \crossingEEfoam[.5]
\]
\[
\Phi_n \left(
\xy 0;/r.20pc/:
	(0,0)*{\xybox{
    	(-4,-4)*{};(4,4)*{} **\crv{(-4,-1) & (4,1)}?(1)*\dir{>} ;
    	(4,-4)*{};(-4,4)*{} **\crv{(4,-1) & (-4,1)}?(1)*\dir{>};
    	(-5.5,-3)*{\scs i};
     	(5.5,-3)*{\scs \ \ \ i+1};
     	(9,1)*{\scs  \lambda};
     	(-10,0)*{};(10,0)*{};}};
\endxy
\right) = \quad \crossingEoneEtwofoam[.5] \quad , \qquad
\Phi_n \left(
\xy 0;/r.20pc/:
	(0,0)*{\xybox{
    	(-4,-4)*{};(4,4)*{} **\crv{(-4,-1) & (4,1)}?(1)*\dir{>} ;
    	(4,-4)*{};(-4,4)*{} **\crv{(4,-1) & (-4,1)}?(1)*\dir{>};
    	(-5.5,-3)*{\scs i+1 \ \ \ };
     	(5.5,-3)*{\scs i};
     	(9,1)*{\scs  \lambda};
     	(-10,0)*{};(10,0)*{};}};
\endxy
\right) = \quad \crossingEtwoEonefoam[.5]
\]

\[
\Phi_n \left(
\xy 0;/r.20pc/:
	(0,0)*{\xybox{
    	(-4,-4)*{};(4,4)*{} **\crv{(-4,-1) & (4,1)}?(1)*\dir{>} ;
    	(4,-4)*{};(-4,4)*{} **\crv{(4,-1) & (-4,1)}?(1)*\dir{>};
    	(-5.5,-3)*{\scs j};
     	(5.5,-3)*{\scs i};
     	(9,1)*{\scs  \lambda};
     	(-10,0)*{};(10,0)*{};}};
\endxy
\right) = \quad \crossingEthreeEonefoam[.5] \quad , \qquad
\Phi_n \left(
\xy 0;/r.20pc/:
	(0,0)*{\xybox{
    	(-4,-4)*{};(4,4)*{} **\crv{(-4,-1) & (4,1)}?(1)*\dir{>} ;
    	(4,-4)*{};(-4,4)*{} **\crv{(4,-1) & (-4,1)}?(1)*\dir{>};
    	(-5.5,-3)*{\scs i};
     	(5.5,-3)*{\scs j};
     	(9,1)*{\scs  \lambda};
     	(-10,0)*{};(10,0)*{};}};
\endxy
\right) = \quad \crossingEoneEthreefoam[.5]
\]
where in the above diagrams $j-i > 1$ and the $i^{th}$ sheet is always in the back.
\end{lem}

Note that the diagram for the image of the $(i,i)$-crossing is shorthand for 
the composition of cap and cup generators from equation \eqref{ccgens} with 
four seam-crossing generators from \eqref{MVgensup}. 
We'll occasionally use the notation $\Phi_n^{m,N}$ when we wish to specify the target $2$-category, and 
occasionally omit these sub- and superscripts when it's clear from the context.

Before proving Lemma \ref{lemKLR}, we first deduce the following foam relations, which are consequences 
of the relations in $\foam{n}{m}$:

\begin{align}
\label{mergefoamrel}
\xy
(0,0)*{

};
\endxy
\end{align}

To this end, consider the diagram:
\[
\xymatrix{ 
& \xy
(0,0)*{
%
};
\endxy  
\ar@{=>}[u]_-{\delta_r} \ar@{=>}[ur]^-{\gamma_r}
& 
}
\]
in which the left- and right-hand sides of \eqref{mergefoamrel} are given by $\gamma_l \circ \beta_l \circ \alpha_l$ and 
$\gamma_r \circ \beta_r \circ \alpha_r$, respectively. The relation follows since 
$\delta_l \beta_l \alpha_l = \delta_r \beta_r \alpha_r$ using equation \eqref{nH2rel} and 
$\gamma_l \delta_l^{-1} = \gamma_r  \delta_r^{-1}$ since they are isotopic. Equation \eqref{splitfoamrel} follows 
similarly.

Next, note that for any $p,q$ we have the isomorphism
\[
\xy
(0,0)*{

};
\endxy
\]
so rather than prove equation \eqref{rightsquarefoamrel}, it suffices to check the corresponding 
relation in 
\[
\Hom \left( 
\xy
(0,0)*{
%
};
\endxy \; \; \right). 
\]
Using equation \eqref{nH2rel}, it then suffices to check the foam relation:
\[
\xy
(0,0)*{
%
};
\endxy \quad .
\]
Applying equation \eqref{nH3rel}, isotopy, and equation \eqref{MVrel}, the latter foam equals:
\begin{equation*}
\sum_{\alpha \in P(a+b,c)} (-1)^{| \hat{\alpha} |} \quad
\xy
(0,0)*{
%
};
\endxy 
\quad 
\end{equation*}
as desired, where $\hat{\alpha}' = (b,\ldots,b)$ and $\gamma' = (c,\ldots,c)$. 
Equation \eqref{leftsquarefoamrel} follows from a similar computation.

We now proceed with the proof of Lemma \ref{lemKLR}.
\begin{proof}
It suffices to check that the images of the KLR relations hold in $\foam{n}{m}$. 

The first nilHecke relation in \eqref{eq_nil_rels} holds since applying relation \eqref{MVrel} 
yields a blister in a $2$-labeled sheet with no dots on its facets. 
This foam is zero by equation \eqref{blisterrel} with $a=1=b$ and $\alpha=\emptyset=\beta$.
The second relation in \eqref{eq_nil_rels} follows from equations \eqref{mergefoamrel}, \eqref{splitfoamrel}, 
\eqref{rightsquarefoamrel}, and \eqref{leftsquarefoamrel} with $a=b=c=1$.
The final nilHecke relation \eqref{eq_nil_dotslide} holds since, using
\eqref{MVrel}, this is equivalent to \eqref{nH3rel} with $a=1=b$.

The $j=i+1$ case of relation \eqref{eq_r2_ij-gen} follows from equation \eqref{R2KLR1rel} and the $j=i-1$ case 
follows from \eqref{R2KLR2rel}; the $(\alpha_i,\alpha_j)=0$ case follows via isotopy, as does 
relation \eqref{eq_dot_slide_ij-gen}.

Relation \eqref{eq_r3_easy-gen} holds via isotopy in the case that $(\alpha_i,\alpha_k)=0$, so we first consider the 
situation where $k=i \pm 1$. The relation again holds by isotopy if $(\alpha_i,\alpha_j)=0$ or $(\alpha_j,\alpha_k)=0$, 
so we can assume that $j=i$ or $j=k$. The case when $i=j=k+1$ follows from equation \eqref{PFrel1}. Indeed, the 
foam relation:
\begin{equation}\label{PFexplicitrel1}
\xy
(0,0)*{
\begin{tikzpicture} [scale=.5,fill opacity=0.2]
	\path[fill=blue] (-.5,5.5) to [out=270,in=90] (-2.25,3.5) to [out=270,in=135] 
		(-1.9,2.6) to (-.9,1.6) to [out=135,in=270] (-1.25,2.5) to [out=90,in=270] (.5,4.5);
	\path[fill=blue] (-1.5,2.5) to [out=180,in=315] (-1.9,2.6) to (-.9,1.6) to [out=315,in=180] (-.5,1.5);
	\path[fill=blue] (1,5.5) to [out=270,in=90] (-.75,3.5) to [out=270,in=0] (-1.5,2.5) to (-.5,1.5) to 
		[out=0,in=270] (.25,2.5) to [out=90,in=270] (2,4.5);
	\path[fill=blue] (-.5,0) to (-.5,1.5) to (-1.5,2.5) to (-1.5,1);
	\path[fill=blue] (-3,6.5) to [out=270,in=135] (-2.5,5.37) to (-1.5,4.37) to [out=135,in=270] (-2,5.5);
	\path[fill=blue] (-2.5,5.37) to [out=325, in=135] (-.25,3.85) to (.75,2.85) to [out=135, in=325] (-1.5,4.37);
	\path[fill=blue] (0,2) to [out=90,in=295] (-.25,3.85) to (.75,2.85) to [out=295,in=90] (1,1);
	\path[fill=red] (1.5,2) to (1.5,6.5) to (-4.5,6.5) to (-4.5,2);
	\path[fill=red] (2.5,1) to (2.5,5.5) to (-3.5,5.5) to (-3.5,1);
	\path[fill=red] (3.5,0) to (3.5,4.5) to (-2.5,4.5) to (-2.5,0);
	\draw[very thick] (1.5,2) to (1.5,6.5);	
	\draw[very thick] (2.5,1) to (2.5,5.5);
	\draw[very thick] (3.5,0) to (3.5,4.5);
	\draw[very thick] (-4.5,2) to (-4.5,6.5);	
	\draw[very thick] (-3.5,1) to (-3.5,5.5);
	\draw[very thick] (-2.5,0) to (-2.5,4.5);
	\draw[very thick,directed=.15,directed=.8] (1.5,2) to (-4.5,2);	
	\draw[very thick,directed=.5] (1,1) to (0,2);
	\draw[very thick,directed=.15,directed=.45,directed=.8] (2.5,1) to (-3.5,1);
	\draw[very thick,directed=.5] (-.5,0) to (-1.5,1);
	\draw[very thick,directed=.4,directed=.9] (3.5,0) to (-2.5,0);
	\draw[very thick, red] (-.5,0) to (-.5,1.5);
	\draw[very thick, red] (-.5,1.5) to [out=0,in=270] (.25,2.5) to [out=90,in=270] (2,4.5);
	\draw[very thick, red] (-.5,1.5) to [out=180,in=270] (-1.25,2.5) to [out=90,in=270] (.5,4.5);
	\draw[very thick, red] (-.5,1.5) to (-1.5,2.5);
	\draw[very thick, red] (-1.5,1) to (-1.5,2.5);
	\draw[very thick, red] (-1.5,2.5) to [out=0,in=270] (-.75,3.5) to [out=90,in=270] (1,5.5);
	\draw[very thick, red] (-1.5,2.5) to [out=180,in=270] (-2.25,3.5) to [out=90,in=270] (-.5,5.5);
	\draw[very thick, red] (1,1) to (1,2) to [out=90,in=270] (-2,5.5);
	\draw[very thick, red] (0,2) to (0,3) to [out=90,in=270] (-3,6.5);
	\draw[very thick,directed=.5,directed=.9] (1.5,6.5) to (-4.5,6.5);
	\draw[very thick,directed=.5] (-2,5.5) to (-3,6.5);
	\draw[very thick,directed=.15,directed=.4,directed=.6,directed=.9] (2.5,5.5) to (-3.5,5.5);
	\draw[very thick,directed=.5] (2,4.5) to (1,5.5);
	\draw[very thick,directed=.5] (.5,4.5) to (-.5,5.5);
	\draw[very thick,directed=.15,directed=.4,directed=.8] (3.5,4.5) to (-2.5,4.5);
	\node[blue,opacity=1] at (1.125,4.125) {\tiny$b$};
	\node[blue,opacity=1] at (-1.5,3.5) {\tiny$a$};
	\node[blue,opacity=1] at (.5,2.5) {\tiny$d$};
	\node[blue,opacity=1] at (-1,1.25) {\tiny$_{a+b}$};
	\end{tikzpicture}
};
\endxy
\quad = \quad
\xy
(0,0)*{
\begin{tikzpicture} [scale=.5,fill opacity=0.2]
	\path[fill=blue] (.5,4.5) to [out=270,in=120] (.77,3.15) to (-.23,4.15) to [out=120,in=270] (-.5,5.5);
	\path[fill=blue] (1.25,3) to [out=180,in=315] (.77,3.15) to (-.23,4.15) to [out=315,in=180] (.25,4);
	\path[fill=blue] (2,4.5) to [out=270,in=0] (1.25,3) to (.25,4) to [out=0,in=270] (1,5.5);
	\path[fill=blue] (1.25,3) to [out=270,in=90] (-.5,0) to  (-1.5,1) to [out=90,in=270] (.25,4);
	\path[fill=blue] (-3,6.5) to [out=270,in=135] (-2.6,4.5) to (-1.6,3.5) to [out=135,in=270] (-2,5.5);
	\path[fill=blue] (-.27,2.85) to [out=135,in=315] (-2.6,4.5) to (-1.6,3.5) to [out=315,in=135] (.73,1.85);
	\path[fill=blue] (0,2) to [out=90,in=315] (-.27,2.85) to (.73,1.85) to [out=315,in=90] (1,1);
	\path[fill=red] (1.5,2) to (1.5,6.5) to (-4.5,6.5) to (-4.5,2);
	\path[fill=red] (2.5,1) to (2.5,5.5) to (-3.5,5.5) to (-3.5,1);
	\path[fill=red] (3.5,0) to (3.5,4.5) to (-2.5,4.5) to (-2.5,0);
	\draw[very thick] (1.5,2) to (1.5,6.5);	
	\draw[very thick] (2.5,1) to (2.5,5.5);
	\draw[very thick] (-4.5,2) to (-4.5,6.5);	
	\draw[very thick] (-3.5,1) to (-3.5,5.5);
	\draw[very thick,directed=.15,directed=.8] (1.5,2) to (-4.5,2);	
	\draw[very thick,directed=.5] (1,1) to (0,2);
	\draw[very thick,directed=.15,directed=.45,directed=.8] (2.5,1) to (-3.5,1);
	\draw[very thick,directed=.5] (-.5,0) to (-1.5,1);
	\draw[very thick,directed=.4,directed=.9] (3.5,0) to (-2.5,0);
	\draw[very thick, red] (1.25,3) to [out=0,in=270] (2,4.5);
	\draw[very thick, red] (1.25,3) to [out=180,in=270] (.5,4.5);
	\draw[very thick, red] (-.5,0) to [out=90,in=270] (1.25,3);
	\draw[very thick, red] (1.25,3) to (.25,4);
	\draw[very thick, red] (.25,4) to [out=0,in=270] (1,5.5);
	\draw[very thick, red] (.25,4) to [out=180,in=270] (-.5,5.5);
	\draw[very thick, red] (-1.5,1) to [out=90,in=270] (.25,4);
	\draw[very thick, red] (1,1) to [out=90,in=270] (-2,4.5) to (-2,5.5);
	\draw[very thick, red] (0,2) to [out=90,in=270] (-3,5.5) to (-3,6.5);
	\draw[very thick,directed=.5,directed=.9] (1.5,6.5) to (-4.5,6.5);
	\draw[very thick,directed=.5] (-2,5.5) to (-3,6.5);
	\draw[very thick,directed=.15,directed=.4,directed=.6,directed=.9] (2.5,5.5) to (-3.5,5.5);
	\draw[very thick,directed=.5] (2,4.5) to (1,5.5);
	\draw[very thick,directed=.5] (.5,4.5) to (-.5,5.5);
	\draw[very thick,directed=.15,directed=.4,directed=.8] (3.5,4.5) to (-2.5,4.5);
	\draw[very thick] (3.5,0) to (3.5,4.5);
	\draw[very thick] (-2.5,0) to (-2.5,4.5);
	\node[blue,opacity=1] at (1.125,4.125) {\tiny$b$};
	\node[blue,opacity=1] at (-.125,4.75) {\tiny$a$};
	\node[blue,opacity=1] at (-2.5,5) {\tiny$d$};
	\node[blue,opacity=1] at (-.75,1.5) {\tiny$_{a+b}$};
	\end{tikzpicture}
};
\endxy
\end{equation}
is a direct consequence of \eqref{PFrel1} and the foam relation
\begin{equation}\label{PFexplicitrel2}
\xy
(0,0)*{
\begin{tikzpicture} [scale=.5,fill opacity=0.2]
	\path[fill=blue] (-3,6.5) to [out=270,in=135] (-2.5,5.37) to (-1.5,4.37) to [out=135,in=270] (-2,5.5);
	\path[fill=blue] (-2.5,5.37) to [out=325, in=135] (-.25,3.85) to (.75,2.85) to [out=135, in=325] (-1.5,4.37);
	\path[fill=blue] (0,2) to [out=90,in=295] (-.25,3.85) to (.75,2.85) to [out=295,in=90] (1,1);
	\path[fill=blue] (-.25,1.5) to [out=0,in=135] (.3,1.25) to (-.7,2.25) to [out=135,in=0] (-1.25,2.5);
	\path[fill=blue] (.5,0) to [out=90,in=300] (.3,1.25) to (-.7,2.25) to [out=300,in=90] (-.5,1);
	\path[fill=blue] (-1,0) to [out=90,in=180] (-.25,1.5) to (-1.25,2.5) to [out=180,in=90] (-2,1);
	\path[fill=blue] (-.25,1.5) to [out=90,in=270] (1,4.5) to (0,5.5) to [out=270,in=90] (-1.25,2.5);
	\path[fill=red] (1.5,2) to (1.5,6.5) to (-4.5,6.5) to (-4.5,2);
	\path[fill=red] (2.5,1) to (2.5,5.5) to (-3.5,5.5) to (-3.5,1);
	\path[fill=red] (3.5,0) to (3.5,4.5) to (-2.5,4.5) to (-2.5,0);
	\draw[very thick] (1.5,2) to (1.5,6.5);	
	\draw[very thick] (2.5,1) to (2.5,5.5);
	\draw[very thick] (3.5,0) to (3.5,4.5);
	\draw[very thick] (-4.5,2) to (-4.5,6.5);	
	\draw[very thick] (-3.5,1) to (-3.5,5.5);
	\draw[very thick] (-2.5,0) to (-2.5,4.5);
	\draw[very thick,directed=.15,directed=.8] (1.5,2) to (-4.5,2);
	\draw[very thick, directed=.5] (1,1) to (0,2);
	\draw[very thick,directed=.15,directed=.4,directed=.6,directed=.9] (2.5,1) to (-3.5,1);
	\draw[very thick, directed=.5] (.5,0) to (-.5,1);
	\draw[very thick, directed=.5] (-1,0) to (-2,1);
	\draw[very thick,directed=.3,directed=.7,directed=.9] (3.5,0) to (-2.5,0);
	\draw[very thick, red] (-1,0) to [out=90,in=180] (-.25,1.5);
	\draw[very thick, red] (.5,0) to [out=90,in=0] (-.25,1.5);
	\draw[very thick, red] (-.25,1.5) to [out=90,in=270] (1,4.5);
	\draw[very thick, red] (-.25,1.5) to (-1.25,2.5);
	\draw[very thick, red] (-2,1) to [out=90,in=180] (-1.25,2.5);
	\draw[very thick, red] (-.5,1) to [out=90,in=0] (-1.25,2.5);
	\draw[very thick, red] (-1.25,2.5) to [out=90,in=270] (0,5.5);
	\draw[very thick, red] (1,1) to (1,2) to [out=90,in=270] (-2,5.5);
	\draw[very thick, red] (0,2) to (0,3) to [out=90,in=270] (-3,6.5);
	\draw[very thick,directed=.5,directed=.9] (1.5,6.5) to (-4.5,6.5);
	\draw[very thick,directed=.5] (-2,5.5) to (-3,6.5);
	\draw[very thick,directed=.25,directed=.6,directed=.9] (2.5,5.5) to (-3.5,5.5);
	\draw[very thick, directed=.5] (1,4.5) to (0,5.5);
	\draw[very thick,directed=.25,directed=.7] (3.5,4.5) to (-2.5,4.5);
	\node[blue,opacity=1] at (.125,.75) {\tiny$b$};
	\node[blue,opacity=1] at (-1.5,1.5) {\tiny$a$};
	\node[blue,opacity=1] at (.5,2.5) {\tiny$d$};
	\node[blue,opacity=1] at (.25,4.125) {\tiny$_{a+b}$};
	\end{tikzpicture}
};
\endxy
\quad = \quad
\xy
(0,0)*{
\begin{tikzpicture} [scale=.5,fill opacity=0.2]
	\path[fill=blue] (-1,0) to [out=90,in=270] (.25,2) to [out=90,in=180] (1,3) to (0,4) to [out=180,in=90]
		(-.75,3) to [out=270,in=90] (-2,1);
	\path[fill=blue] (1.5,2.8) to [out=135,in=0] (1,3) to (0,4) to [out=0,in=135] (.5,3.8);
	\path[fill=blue] (.5,0) to [out=90,in=270] (1.75,2) to [out=90,in=300] (1.5,2.8) to (.5,3.8) to [out=300,in=90]
		(.75,3) to [out=270,in=90] (-.5,1);
	\path[fill=blue] (1,4.5) to (0,5.5) to (0,4) to (1,3);
	\path[fill=blue] (-3,6.5) to [out=270,in=135] (-2.6,4.5) to (-1.6,3.5) to [out=135,in=270] (-2,5.5);
	\path[fill=blue] (-.27,2.85) to [out=135,in=315] (-2.6,4.5) to (-1.6,3.5) to [out=315,in=135] (.73,1.85);
	\path[fill=blue] (0,2) to [out=90,in=315] (-.27,2.85) to (.73,1.85) to [out=315,in=90] (1,1);
	\path[fill=red] (1.5,2) to (1.5,6.5) to (-4.5,6.5) to (-4.5,2);
	\path[fill=red] (2.5,1) to (2.5,5.5) to (-3.5,5.5) to (-3.5,1);
	\path[fill=red] (3.5,0) to (3.5,4.5) to (-2.5,4.5) to (-2.5,0);
	\draw[very thick] (1.5,2) to (1.5,6.5);	
	\draw[very thick] (2.5,1) to (2.5,5.5);
	\draw[very thick] (-4.5,2) to (-4.5,6.5);	
	\draw[very thick] (-3.5,1) to (-3.5,5.5);
	\draw[very thick,directed=.15,directed=.8] (1.5,2) to (-4.5,2);
	\draw[very thick, directed=.5] (1,1) to (0,2);
	\draw[very thick,directed=.15,directed=.4,directed=.6,directed=.9] (2.5,1) to (-3.5,1);
	\draw[very thick, directed=.5] (.5,0) to (-.5,1);
	\draw[very thick, directed=.5] (-1,0) to (-2,1);
	\draw[very thick,directed=.3,directed=.7,directed=.9] (3.5,0) to (-2.5,0);
	\draw[very thick, red] (.5,0) to [out=90,in=270] (1.75,2) to [out=90,in=0] (1,3);
	\draw[very thick, red] (-1,0) to [out=90,in=270] (.25,2) to [out=90,in=180] (1,3);
	\draw[very thick, red] (1,3) to (1,4.5);	
	\draw[very thick, red] (1,3) to (0,4);
	\draw[very thick, red] (-.5,1) to [out=90,in=270] (.75,3) to [out=90,in=0] (0,4);
	\draw[very thick, red] (-2,1) to [out=90,in=270] (-.75,3) to [out=90,in=180] (0,4);
	\draw[very thick, red] (0,4) to (0,5.5);
	\draw[very thick, red] (1,1) to [out=90,in=270] (-2,4.5) to (-2,5.5);
	\draw[very thick, red] (0,2) to [out=90,in=270] (-3,5.5) to (-3,6.5);
	\draw[very thick,directed=.5,directed=.9] (1.5,6.5) to (-4.5,6.5);
	\draw[very thick,directed=.5] (-2,5.5) to (-3,6.5);
	\draw[very thick,directed=.25,directed=.6,directed=.9] (2.5,5.5) to (-3.5,5.5);
	\draw[very thick, directed=.5] (1,4.5) to (0,5.5);
	\draw[very thick,directed=.25,directed=.7] (3.5,4.5) to (-2.5,4.5);
	\draw[very thick] (3.5,0) to (3.5,4.5);
	\draw[very thick] (-2.5,0) to (-2.5,4.5);
	\node[blue,opacity=1] at (.125,.75) {\tiny$b$};
	\node[blue,opacity=1] at (-1,1.5) {\tiny$a$};
	\node[blue,opacity=1] at (-2.5,5) {\tiny$d$};
	\node[blue,opacity=1] at (.5,4.25) {\tiny$_{a+b}$};
	\end{tikzpicture}
};
\endxy
\end{equation}
follows by composing \eqref{PFrel1} with foams of the form \eqref{MVgensdown} and applying \eqref{MVrel}. 
The $a=b=d=1$ case of these relations gives the aforementioned case of relation \eqref{eq_r3_easy-gen}. 
In a similar manner, the $j=k=i+1$ case of \eqref{eq_r3_easy-gen} follows from \eqref{PFrel2}, the $j=k=i-1$ case from 
\eqref{PFrel3}, and the $i=j=k-1$ case from \eqref{PFrel4}.

Finally, relation \eqref{eq_r3_hard-gen} follows from \eqref{R3KLRrel1} when $j=i+1$ and from \eqref{R3KLRrel2} when $j=i-1$.
\end{proof}

The above $2$-functor actually induces one from $\Ucatc_Q^+(\mf{sl}_m)$, the full $2$-subcategory of 
$\UcatD_Q^+(\mf{sl}_m)$ with $1$-morphisms given by compositions of divided powers 
$\cal{E}^{(a)}\onel\{t\}$.

\begin{lem}\label{lemKLRc}
The $2$-functor in Lemma \ref{lemKLR} extends to a $2$-functor 
$\Phi_n  \maps \Ucatc_Q^+(\mf{sl}_m)  \to  \foam{n}{m}$
\end{lem}
\begin{proof}
The $2$-functors $\Phi_n$ induce $2$-functors $Kar(\Ucat_Q(\mf{sl}_m))  \to  Kar(\foam{n}{m})$.
Since any category $\mathcal{C}$ embeds fully faithfully into its Karoubi envelope by sending an 
object $c$ to the pair $(c,\id_c)$, it suffices to check that the images of 
$\cal{E}^{(a)}\onel$ in $Kar(\foam{n}{m})$ are isomorphic to $1$-morphisms of this form.

It follows from equations \eqref{MVrel}, \eqref{blisterrel}, \eqref{rightsquarefoamrel}, and \eqref{leftsquarefoamrel} that
\[
\Phi_n(\cal{E}^{(a)}\onel) \cong 
\left( 
\xy
(0,0)*{
};
\endxy
\right) \quad ,
\]
which we'll demonstrate in the case that $a=3$ (the case for general $a$ proceeds analogously, 
however, drawing the relevant foams becomes cumbersome).

Indeed, equations \eqref{rightsquarefoamrel} and \eqref{leftsquarefoamrel} imply that
\[
\Phi_n(\cal{E}^{(a)}\onel) =
\left(
\xy
(0,0)*{
%
};
\endxy
\;\;
\right)
\]
give the isomorphism in $Kar(\foam{n}{m})$ by equations \eqref{MVrel} and \eqref{blisterrel}.
\end{proof}

One can check that, unsurprisingly, the $2$-functor $\Phi_n$ is given on splitter and merger thick-calculus 
$2$-morphisms by
\[
\xy
(0,0)*{
%
};
\endxy
\]

Next, we aim to show that the above $2$-functors extend to the entire categorified quantum group. 
In \cite{LQR1}, we were able to deduce this result from the main theorem of \cite{CLau}; however, in 
the present case we can't guarantee a priori that the hypotheses of that theorem hold. Instead, we'll 
employ the following trick: note that we have web isotopies: 
\begin{align}
\xy
(0,0)*{
%
}
\endxy \nn
\end{align}
which give a correspondence between foams mapping between the webs on the left-hand side 
and foams between the webs on the right. We'll use this correspondence to show that all of the 
$\mathfrak{sl}_2$ relations are satisfied in the foam $2$-category.

\begin{thm} \label{thm:2functor}
Extending $\Phi_n$ to cap and cup $1$-morphisms via
\[
\Phi_n \left(
\xy 0;/r.20pc/:
	(0,0)*{\bbcef{i}};
	(8,4)*{\scs  \lambda};
	(-10,0)*{};(10,0)*{};
\endxy
\right) =  \quad \capFEfoam{}{} \quad , \qquad
\Phi_n \left(
\xy 0;/r.20pc/:
	(0,0)*{\bbcfe{i}};
	(8,4)*{\scs  \lambda};
	(-10,0)*{};(10,0)*{};
\endxy
\right) = (-1)^{a_{i+1}} \quad \capEFfoam{}{}
\]
\[
\Phi_n \left(
\xy 0;/r.20pc/:
	(0,0)*{\bbpef{i}};
	(8,-4)*{\scs \lambda};
	(-10,0)*{};(10,0)*{};
\endxy
\right) = (-1)^{a_{i+1} + 1} \quad \cupFEfoam{}{} \quad , \qquad
\Phi_n \left(
\xy 0;/r.20pc/:
	(0,0)*{\bbpfe{i}};
	(8,-4)*{\scs \lambda};
	(-10,0)*{};(10,0)*{};
\endxy
\right) = \quad \cupEFfoam{}{}
\]
defines a $2$-representation $\Phi_n  \maps \Ucatc_Q(\mf{sl}_m)  \to  \foam{n}{m}$.
\end{thm}

One can then immediately check that $\Phi_n$ acts on thick-calculus elements as follows.
\begin{cor}
\[
\Phi_n \left(
\;
\xy 
(0,0)*{
\begin{tikzpicture}[scale=.75]
\draw[ultra thick, black!30!green, ->] (1,0) .. controls (1,.8) and (0,.8) .. (0,0);
\node at (1,.7) {$\scs \lambda$};
\node [black!30!green] at (1,-.2) {$\scs k$};
\end{tikzpicture}
};
\endxy
\;
\right) = \quad \capFEfoam{}{} \quad , \qquad
\Phi_n \left(
\;
\xy 
(0,0)*{
\begin{tikzpicture}[scale=.75]
\draw[ultra thick, black!30!green, <-] (1,0) .. controls (1,.8) and (0,.8) .. (0,0);
\node at (1,.7) {$\scs \lambda$};
\node [black!30!green] at (0,-.2) {$\scs k$};
\end{tikzpicture}
};
\endxy
\;
\right) = (-1)^{\frac{k(k-1)}{2}+ka_{i+1}} \quad \capEFfoam{}{}
\]

\[
\Phi_n \left(
\;
\xy 
(0,0)*{
\begin{tikzpicture}[scale=.75]
\draw[ultra thick, black!30!green, <-] (1,1) .. controls (1,.2) and (0,.2) .. (0,1);
\node at (1,.3) {$\scs \lambda$};
\node [black!30!green] at (0,1.2) {$\scs k$};
\end{tikzpicture}
};
\endxy
\;
\right) = (-1)^{\frac{k(k-1)}{2}+k(a_{i+1} + 1)} \quad \cupFEfoam{}{} \quad , \qquad
\Phi_n \left(
\;
\xy 
(0,0)*{
\begin{tikzpicture}[scale=.75]
\draw[ultra thick, black!30!green, ->] (1,1) .. controls (1,.2) and (0,.2) .. (0,1);
\node at (1,.3) {$\scs \lambda$};
\node [black!30!green] at (1,1.2) {$\scs k$};
\end{tikzpicture}
};
\endxy
\;
\right) = \quad \cupEFfoam{}{}
\]

\[
\Phi_n \left(
\;
\xy 
(0,0)*{
\begin{tikzpicture}[scale=1]
\draw[ultra thick, black!30!green, ->] (1,1) .. controls (1,.6) and (0,.4) .. (0,0);
\draw[ultra thick, blue, ->] (0,1) .. controls (0,.6) and (1,.4) .. (1,0);
\node [black!30!green] at (1,1.2) {$\scs k$};
\node [blue] at (0,1.2) {$\scs l$};
\end{tikzpicture}
};
\endxy
\;
\right) =  \quad 
\xy
(0,0)*{
\begin{tikzpicture} [scale=.35,fill opacity=0.2]
	\path[fill=red] (3,6) rectangle (-4,1);
	\path[fill=red] (4,5) rectangle (-3,0);
	\path[fill=red] (5,4) rectangle (-2,-1);
	\path[fill=blue] (0,1) to [out=90,in=270] (2,6) to (1,5) to [out=270,in=90] (-1,0);
	\path[fill=blue] (2,0) to [out=90,in=270] (0,5) to (-1,4) to [out=270,in=90] (1,-1);
	\draw[very thick] (3,1) to (-4,1);
	\draw[very thick] (4,0) to (-3,0);
	\draw[very thick] (5,-1) to (-2,-1);
	\draw[very thick, directed=.55] (2,0) to (1,-1);
	\draw[very thick, directed=.55] (0,1) to (-1,0);
	\draw[very thick] (-2,-1) to (-2,4);
	\draw[very thick] (-3,0) to (-3,5);
	\draw[very thick] (-4,1) to (-4,6);
	\draw[very thick] (3,1) to (3,6);
	\draw[very thick] (4,0) to (4,5);
	\draw[very thick] (5,-1) to (5,4);
	\draw[very thick, red] (1,-1) to [out=90,in=270] (-1,4);
	\draw[very thick, red] (2,0) to [out=90,in=270] (0,5);	
	\draw[very thick, red] (0,1) to [out=90,in=270] (2,6);		
	\draw[very thick, red] (-1,0) to [out=90,in=270] (1,5);
	\draw[very thick] (3,6) to (-4,6);
	\draw[very thick] (4,5) to (-3,5);
	\draw[very thick] (5,4) to (-2,4);
	\draw[very thick, directed=.55] (2,6) to (1,5);
	\draw[very thick, directed=.55] (0,5) to (-1,4);	
\end{tikzpicture}};
\endxy
\quad , \qquad
\Phi_n \left(
\;
\xy 
(0,0)*{
\begin{tikzpicture}[scale=1]
\draw[ultra thick, black!30!green, ->] (0,1) .. controls (0,.6) and (1,.4) .. (1,0);
\draw[ultra thick, blue, ->] (1,1) .. controls (1,.6) and (0,.4) .. (0,0);
\node [black!30!green] at (0,1.2) {$\scs k$};
\node [blue] at (1,1.2) {$\scs l$};
\end{tikzpicture}
};
\endxy
\;
\right) = (-1)^{kl}\quad 
\xy
(0,0)*{
\begin{tikzpicture} [scale=.35,fill opacity=0.2]
	\path[fill=red] (3,6) rectangle (-4,1);
	\path[fill=red] (4,5) rectangle (-3,0);
	\path[fill=red] (5,4) rectangle (-2,-1);
	\path[fill=blue] (0,0) to [out=90,in=270] (2,5) to (1,4) to [out=270,in=90] (-1,-1);
	\path[fill=blue] (2,1) to [out=90,in=270] (0,6) to (-1,5) to [out=270,in=90] (1,0);
	\draw[very thick] (3,1) to (-4,1);
	\draw[very thick] (4,0) to (-3,0);
	\draw[very thick] (5,-1) to (-2,-1);
	\draw[very thick, directed=.55] (2,1) to (1,0);
	\draw[very thick, directed=.55] (0,0) to (-1,-1);
	\draw[very thick] (-2,-1) to (-2,4);
	\draw[very thick] (-3,0) to (-3,5);
	\draw[very thick] (-4,1) to (-4,6);
	\draw[very thick] (3,1) to (3,6);
	\draw[very thick] (4,0) to (4,5);
	\draw[very thick] (5,-1) to (5,4);
	\draw[very thick, red] (1,0) to [out=90,in=270] (-1,5);
	\draw[very thick, red] (2,1) to [out=90,in=270] (0,6);	
	\draw[very thick, red] (0,0) to [out=90,in=270] (2,5);		
	\draw[very thick, red] (-1,-1) to [out=90,in=270] (1,4);
	\draw[very thick] (3,6) to (-4,6);
	\draw[very thick] (4,5) to (-3,5);
	\draw[very thick] (5,4) to (-2,4);
	\draw[very thick, directed=.55] (2,5) to (1,4);
	\draw[very thick, directed=.55] (0,6) to (-1,5);	
\end{tikzpicture}};
\endxy
\]

\[
\Phi_n \left(
\;
\xy 
(0,0)*{
\begin{tikzpicture}[scale=1]
\draw[ultra thick, black!30!green, ->] (0,1) .. controls (0,.6) and (1,.4) .. (1,0);
\draw[ultra thick, blue, ->] (0,0) .. controls (0,.4) and (1,.6) .. (1,1);
\node [black!30!green] at (0,1.2) {$\scs k$};
\node [blue] at (0,-.2) {$\scs l$};
\end{tikzpicture}
};
\endxy
\;
\right) = (-1)^{kl} \quad 
\xy
(0,0)*{
\begin{tikzpicture} [scale=.35,fill opacity=0.2]
	\path[fill=red] (3,6) rectangle (-4,1);
	\path[fill=red] (4,5) rectangle (-3,0);
	\path[fill=red] (5,4) rectangle (-2,-1);
	\path[fill=blue] (0,-1) to [out=90,in=270] (2,4) to (1,5) to [out=270,in=90] (-1,0);
	\path[fill=blue] (2,1) to [out=90,in=270] (0,6) to (-1,5) to [out=270,in=90] (1,0);
	\draw[very thick] (3,1) to (-4,1);
	\draw[very thick] (4,0) to (-3,0);
	\draw[very thick] (5,-1) to (-2,-1);
	\draw[very thick, directed=.55] (2,1) to (1,0);
	\draw[very thick, directed=.55] (0,-1) to (-1,0);
	\draw[very thick] (-2,-1) to (-2,4);
	\draw[very thick] (-3,0) to (-3,5);
	\draw[very thick] (-4,1) to (-4,6);
	\draw[very thick] (3,1) to (3,6);
	\draw[very thick] (4,0) to (4,5);
	\draw[very thick] (5,-1) to (5,4);
	\draw[very thick, red] (1,0) to [out=90,in=270] (-1,5);
	\draw[very thick, red] (2,1) to [out=90,in=270] (0,6);	
	\draw[very thick, red] (0,-1) to [out=90,in=270] (2,4);		
	\draw[very thick, red] (-1,0) to [out=90,in=270] (1,5);
	\draw[very thick] (3,6) to (-4,6);
	\draw[very thick] (4,5) to (-3,5);
	\draw[very thick] (5,4) to (-2,4);
	\draw[very thick, directed=.55] (2,4) to (1,5);
	\draw[very thick, directed=.55] (0,6) to (-1,5);	
\end{tikzpicture}};
\endxy
\quad , \qquad
\Phi_n \left(
\;
\xy 
(0,0)*{
\begin{tikzpicture}[scale=1]
\draw[ultra thick, black!30!green, ->] (1,0) .. controls (1,.4) and (0,.6) .. (0,1);
\draw[ultra thick, blue, ->] (1,1) .. controls (1,.6) and (0,.4) .. (0,0);
\node [black!30!green] at (1,-.2) {$\scs k$};
\node [blue] at (1,1.2) {$\scs l$};
\end{tikzpicture}
};
\endxy
\;
\right) = \quad 
\xy
(0,0)*{
\begin{tikzpicture} [scale=.35,fill opacity=0.2]
	\path[fill=red] (2.5,6) rectangle (-4.5,1);
	\path[fill=red] (3.5,5) rectangle (-3.5,0);
	\path[fill=red] (4.5,4) rectangle (-2.5,-1);
	\path[fill=blue] (0,0) to [out=90,in=270] (2,5) to (1,4) to [out=270,in=90] (-1,-1);
	\path[fill=blue] (2,0) to [out=90,in=270] (0,5) to (-1,6) to [out=270,in=90] (1,1);
	\draw[very thick] (2.5,1) to (-4.5,1);
	\draw[very thick] (3.5,0) to (-3.5,0);
	\draw[very thick] (4.5,-1) to (-2.5,-1);
	\draw[very thick, directed=.55] (2,0) to (1,1);
	\draw[very thick, directed=.55] (0,0) to (-1,-1);
	\draw[very thick] (-2.5,-1) to (-2.5,4);
	\draw[very thick] (-3.5,0) to (-3.5,5);
	\draw[very thick] (-4.5,1) to (-4.5,6);
	\draw[very thick] (2.5,1) to (2.5,6);
	\draw[very thick] (3.5,0) to (3.5,5);
	\draw[very thick] (4.5,-1) to (4.5,4);
	\draw[very thick, red] (1,1) to [out=90,in=270] (-1,6);
	\draw[very thick, red] (0,0) to [out=90,in=270] (2,5);		
	\draw[very thick, red] (2,0) to [out=90,in=270] (0,5);	
	\draw[very thick, red] (-1,-1) to [out=90,in=270] (1,4);
	\draw[very thick] (2.5,6) to (-4.5,6);
	\draw[very thick] (3.5,5) to (-3.5,5);
	\draw[very thick] (4.5,4) to (-2.5,4);
	\draw[very thick, directed=.55] (2,5) to (1,4);
	\draw[very thick, directed=.55] (0,5) to (-1,6);	
\end{tikzpicture}};
\endxy
\]

\[
\Phi_n \left(
\;
\xy 
(0,0)*{
\begin{tikzpicture}[scale=1]
\draw[ultra thick, black!30!green, ->] (0,0) .. controls (0,.4) and (1,.6) .. (1,1);
\draw[ultra thick, blue, ->] (0,1) .. controls (0,.6) and (1,.4) .. (1,0);
\node [black!30!green] at (0,-.2) {$\scs k$};
\node [blue] at (0,1.2) {$\scs l$};
\end{tikzpicture}
};
\endxy
\;
\right) = \quad 
\xy
(0,0)*{
\begin{tikzpicture} [scale=.35,fill opacity=0.2]
	\path[fill=red] (3,6) rectangle (-4,1);
	\path[fill=red] (4,5) rectangle (-3,0);
	\path[fill=red] (5,4) rectangle (-2,-1);
	\path[fill=blue] (0,0) to [out=90,in=270] (2,5) to (1,6) to [out=270,in=90] (-1,1);
	\path[fill=blue] (2,0) to [out=90,in=270] (0,5) to (-1,4) to [out=270,in=90] (1,-1);
	\draw[very thick] (3,1) to (-4,1);
	\draw[very thick] (4,0) to (-3,0);
	\draw[very thick] (5,-1) to (-2,-1);
	\draw[very thick, directed=.55] (2,0) to (1,-1);
	\draw[very thick, directed=.55] (0,0) to (-1,1);
	\draw[very thick] (-2,-1) to (-2,4);
	\draw[very thick] (-3,0) to (-3,5);
	\draw[very thick] (-4,1) to (-4,6);
	\draw[very thick] (3,1) to (3,6);
	\draw[very thick] (4,0) to (4,5);
	\draw[very thick] (5,-1) to (5,4);
	\draw[very thick, red] (1,-1) to [out=90,in=270] (-1,4);
	\draw[very thick, red] (2,0) to [out=90,in=270] (0,5);	
	\draw[very thick, red] (0,0) to [out=90,in=270] (2,5);		
	\draw[very thick, red] (-1,1) to [out=90,in=270] (1,6);
	\draw[very thick] (3,6) to (-4,6);
	\draw[very thick] (4,5) to (-3,5);
	\draw[very thick] (5,4) to (-2,4);
	\draw[very thick, directed=.55] (2,5) to (1,6);
	\draw[very thick, directed=.55] (0,5) to (-1,4);	
\end{tikzpicture}};
\endxy
\quad , \qquad
\Phi_n \left(
\;
\xy 
(0,0)*{
\begin{tikzpicture}[scale=1]
\draw[ultra thick, black!30!green, ->] (1,1) .. controls (1,.6) and (0,.4) .. (0,0);
\draw[ultra thick, blue, ->] (1,0) .. controls (1,.4) and (0,.6) .. (0,1);
\node [black!30!green] at (1,1.2) {$\scs k$};
\node [blue] at (1,-.2) {$\scs l$};
\end{tikzpicture}
};
\endxy
\;
\right) = \quad 
\xy
(0,0)*{
\begin{tikzpicture} [scale=.35,fill opacity=0.2]
	\path[fill=red] (3,6) rectangle (-4,1);
	\path[fill=red] (4,5) rectangle (-3,0);
	\path[fill=red] (5,4) rectangle (-2,-1);
	\path[fill=blue] (0,1) to [out=90,in=270] (2,6) to (1,5) to [out=270,in=90] (-1,0);
	\path[fill=blue] (2,-1) to [out=90,in=270] (0,4) to (-1,5) to [out=270,in=90] (1,0);
	\draw[very thick] (3,1) to (-4,1);
	\draw[very thick] (4,0) to (-3,0);
	\draw[very thick] (5,-1) to (-2,-1);
	\draw[very thick, directed=.55] (2,-1) to (1,0);
	\draw[very thick, directed=.55] (0,1) to (-1,0);
	\draw[very thick] (-2,-1) to (-2,4);
	\draw[very thick] (-3,0) to (-3,5);
	\draw[very thick] (-4,1) to (-4,6);
	\draw[very thick] (3,1) to (3,6);
	\draw[very thick] (4,0) to (4,5);
	\draw[very thick] (5,-1) to (5,4);
	\draw[very thick, red] (1,0) to [out=90,in=270] (-1,5);
	\draw[very thick, red] (0,1) to [out=90,in=270] (2,6);		
	\draw[very thick, red] (2,-1) to [out=90,in=270] (0,4);	
	\draw[very thick, red] (-1,0) to [out=90,in=270] (1,5);
	\draw[very thick] (3,6) to (-4,6);
	\draw[very thick] (4,5) to (-3,5);
	\draw[very thick] (5,4) to (-2,4);
	\draw[very thick, directed=.55] (2,6) to (1,5);
	\draw[very thick, directed=.55] (0,4) to (-1,5);	
\end{tikzpicture}};
\endxy
\]
\end{cor}

Adopting the terminology from the decategorified setting, we'll refer to foams in the image of $\Phi_n$ 
as \emph{ladder foams}.

\begin{proof}[Proof (of Theorem \ref{thm:2functor})]
The biadjointness relations \eqref{eq_biadjoint1}, \eqref{eq_biadjoint2} and the $Q$-cyclic relation
for dots \eqref{eq_cyclic_dot} are satisfied since the corresponding foams are isotopic and the factors of $-1$ 
in the definition of $\Phi_n$ coming from the right-oriented caps and cups cancel.
The $Q$-cyclic relation for crossings \eqref{eq_almost_cyclic} again holds via isotopy; the factors of $-1$ coming from the 
scalings on the right-oriented caps and cups cancel with each other provided 
$j \neq i \pm1$ and with the factor of $-1$ coming from $t_{ij}$ or $t_{ji}$ when 
$j=i \pm 1$. The mixed $\cal{E}_i \cal{F}_j$ relation \eqref{eq_mixed_EF} follows similarly from the cap/cup 
rescalings and equation \eqref{MVrel}.

This $2$-functor is graded, as can be checked on a case-by-case basis. Indeed,
\[
\deg \left(\;
\xy
(0,0)*{
\begin{tikzpicture}[scale=0.75]
  \draw[ultra thick,color=black!30!green, ->] (-0.5,1) .. controls (-0.5,0.2) and (0.5,0.2) .. (0.5,1)
  node[pos=0.5, shape=coordinate](X){};;
  \node at (.5,.2) {$\scs \lambda$};
  \node at (-.2,.7) {\textcolor{black!30!green}{$\scs k$}};
\end{tikzpicture}};
\endxy
\;\right) = k(k+\lambda_i) = k(k+a_{i}-a_{i+1})
\]
and the corresponding foam $F$ is (up to sign):
\[
\xy
(0,0)*{
\begin{tikzpicture} [scale=.5,fill opacity=0.2]
\path[fill=red] (2,-2) to (-2,-2) to (-2,1) to (2,1);
\path[fill=blue] (1.5,1) to [out=270,in=0] (0,-.25) to [out=180,in=270] (-1.5,1) to
	(-.5,2) to [out=270,in=180] (0,1.25) to [out=0,in=270] (.5,2);
\path[fill=red] (2.5,-1) to (-2.5,-1) to (-2.5,2) to (2.5,2);
\draw[very thick, directed=.5] (2.5,-1) to (-2.5,-1);
\draw[very thick] (2.5,-1) to (2.5,2);
\draw[very thick] (-2.5,-1) to (-2.5,2);
\draw[very thick, directed=.5] (2.5,2) to (-2.5,2);
\draw[very thick, red, rdirected=.5] (-.5,2) to [out=270,in=180] (0,1.25)
	to [out=0,in=270] (.5,2);
\node[red, opacity=1] at (2,1.5) {$_{a_i}$};
\draw[very thick, directed=.5] (2,-2) to (-2,-2);
\draw[very thick] (2,-2) to (2,1);
\draw[very thick] (-2,-2) to (-2,1);
\draw[very thick, directed=.5] (2,1) to (-2,1);
\draw[very thick, red, rdirected=.5] (1.5,1) to [out=270,in=0] (0,-.25)
	to [out=180,in=270] (-1.5,1);
\node[red, opacity=1] at (1.25,-1.5) {$_{a_{i+1}}$};
\node[blue,opacity=1] at (0,.5) {$_k$};
\draw[very thick, directed=.5] (1.5,1) to (.5,2);
\draw[very thick, directed=.5] (-.5,2) to (-1.5,1);
\end{tikzpicture}
}
\endxy
\]
so we compute 
\begin{align*}
 \chi^{Top}&= a_i(n-a_i)+a_{i+1}(n-a_{i+1}) + k(n-a_i)+k(n-a_{i+1})\\ &\quad + 2k(a_{i+1}-k)-2k(n-k) -a_{i+1}k+a_ik \\
&= a_i(n-a_i)+a_{i+1}(n-a_{i+1})  \\
\chi^{Bot}&= a_i(n-a_i)+a_{i+1}(n-a_{i+1})  \\
\chi&= a_i (n-a_i) + a_{i+1}(n-a_{i+1}) -k(k+a_i-a_{i+1})
\end{align*}
which gives that $\deg(F) = k(k + a_i - a_{i+1}).$ 
The result for all other caps and cups following similarly.

For splitters we have 
\[
\deg \left(
\xy
(0,0)*{
\begin{tikzpicture} [scale=.4,fill opacity=0.2]
\draw[ultra thick, black!30!green] (0,-1.5) to (0,0);
\draw[ultra thick, black!30!green, directed=1] (0,0) to [out=30,in=270] (1,1.5);
\draw[ultra thick, black!30!green, directed=1] (0,0) to [out=150,in=270] (-1,1.5);
\node[black!30!green, opacity=1] at (0,-2) {\tiny$_{k+l}$};
\node[black!30!green, opacity=1] at (-1,2) {\tiny$k$};
\node[black!30!green, opacity=1] at (1,2) {\tiny$l$};
\end{tikzpicture}
}
\endxy
\right) = -kl.
\]
and the corresponding foam is:
\[
\xy
(0,0)*{
\begin{tikzpicture} [scale=.4,fill opacity=0.2]
	\path[fill=red] (3.5,2.5) to (-2.5,2.5) to (-2.5,7.5) to (3.5,7.5);
	\path[fill=red] (2.5,3.5) to (-3.5,3.5) to (-3.5,8.5) to (2.5,8.5);
	\path[fill=blue] (.5,2.5) to (-.5,3.5) to (-.5,5.75) to (.5,4.75);
	\path[fill=blue] (-.75,7.5) to [out=270,in=130] (-.125,4.91) to (-1.125,5.91) to [out=130,in=270] (-1.75,8.5);
	\path[fill=blue] (-.125,4.91) to [out=320,in=180] (.5,4.75) to (-.5,5.75) to [out=180,in=320] (-1.125,5.91);
	\path[fill=blue] (.5,4.75) to [out=0,in=270] (1.75,7.5) to (.75,8.5) to [out=270,in=0] (-.5,5.75);
	\draw[very thick, directed=.25, directed=.75] (3.5,2.5) to (-2.5,2.5);
	\draw[very thick, directed=.25, directed=.75] (2.5,3.5) to (-3.5,3.5);
	\draw[very thick, directed=.55] (.5,2.5) to (-.5,3.5);	
	\draw[very thick, red] (.5,2.5) to [out=90,in=270] (.5,4.75);
	\draw[very thick, red] (.5,4.75) to [out=0,in=270] (1.75,7.5);
	\draw[very thick, red] (.5,4.75) to [out=180,in=270] (-.75,7.5);
	\draw[very thick, red] (.5,4.75) to (-.5,5.75);
	\draw[very thick, red] (-.5,3.5) to [out=90,in=270] (-.5,5.75);
	\draw[very thick, red] (-.5,5.75) to [out=0,in=270] (.75,8.5);
	\draw[very thick, red] (-.5,5.75) to [out=180,in=270] (-1.75,8.5);
	\draw[very thick] (3.5,2.5) to (3.5,7.5);	
	\draw[very thick] (2.5,3.5) to (2.5,8.5);
	\draw[very thick] (-2.5,2.5) to (-2.5,7.5);	
	\draw[very thick] (-3.5,3.5) to (-3.5,8.5);
	\draw[very thick, directed=.5, directed=.9, directed=.1] (3.5,7.5) to (-2.5,7.5);
	\draw[very thick, directed=.5, directed=.9, directed=.1] (2.5,8.5) to (-3.5,8.5);
	\draw[very thick, directed=.5] (1.75,7.5) to (.75,8.5);	
	\draw[very thick, directed=.5] (-.75,7.5) to (-1.75,8.5);
	\node[blue, opacity=1] at (-1.2,7) {\tiny$k$};
	\node[blue, opacity=1] at (1.25,7) {\tiny$l$};
	\node[blue, opacity=1] at (0,4.25) {\tiny$_{k+l}$};
	\node[red, opacity=1] at (2,8) {\tiny$a_i$};
	\node[red, opacity=1] at (2.75,3) {\tiny$a_{i+1}$};	
	\end{tikzpicture}
};
\endxy
\]
We again compute:
\begin{align*}
 \chi^{Top}&= (a_{i+1}+k+l)(n-a_{i+1}-k-l) + (a_{i}+k+l)(n-a_{i}-k-l)\\ &\quad +a_{i+1}k +a_{i+1}l +2kl+a_{i}l+ ka_{i}+ k^2+l^2 -kn -ln \\
\chi^{Bot}&= (a_{i+1}+k+l)(n-a_{i+1}-k-l) + (a_{i}+k+l) (n-a_{i}-k-l)\\ &\quad +a_{i+1}k+a_{i+1}l+a_{i}k+a_{i}l+k^2+l^2+2kl-kn-ln \\
\chi&= (a_{i+1}+k+l)(n-a_{i+1}-k-l)+a_{i+1}k +a_{i+1}l+a_{i}l\\ &\quad +(a_{i}+k+l)(n-a_{i}-k-l)+a_{i}k+kl-(k+l)(n-k-l)
\end{align*}
from which we deduce that the degree of the foam is $-kl$.

We finally check that
\[
\deg \left(
\xy
(0,0)*{
\begin{tikzpicture} [scale=.4,fill opacity=0.2]
\draw[ultra thick, blue, directed=.99] (0,0) to [out=90,in=270] (2,3);
\draw[ultra thick, black!30!green, rdirected=.05] (0,3) to [out=270,in=90] (2,0);
\node[blue, opacity=1] at (0,-.5) {\tiny$k$};
\node[black!30!green, opacity=1] at (2,-.5) {\tiny$l$};
\node[opacity=1] at (2.75,1.5) {\tiny$\lambda$};
\end{tikzpicture}
}
\endxy \;
\right) = -kl,
\]
and that the degree of the foam
\[
\xy
(0,0)*{
\begin{tikzpicture}  [fill opacity=0.2,  decoration={markings, mark=at position 0.6 with {\arrow{>}};}, scale=.7]
	\fill [fill=red] (-3,1) rectangle (1,4);
	\path [fill=blue] (-1.75,4) .. controls (-1.75,3) and (-.25,2) .. (-.25,1) --
			(.75,0) .. controls (.75,1) and (-.75,2) .. (-.75,3);
	\fill [fill=red] (-2,0) rectangle (2,3);
	\path [fill=blue] (.75,3) .. controls (.75,2) and (-.75,1) .. (-.75,0) --
			(.25,-1) .. controls (.25,0) and (1.75,1) .. (1.75,2);
	\fill [fill=red] (-1,-1) rectangle (3,2);
	\draw[very thick, postaction={decorate}] (.75,0) -- (-.25,1);
	\draw[very thick, postaction={decorate}] (.25,-1) -- (-.75,0);
	\draw [very thick, red] (-1.75,4) .. controls (-1.75,3) and (-.25,2) .. (-.25,1);
	\draw [very thick] (-3,1) rectangle (1,4);
	\draw [very thick, red] (-.75,3) .. controls (-.75,2) and (.75,1) .. (.75,0);
	\draw [very thick, red] (.75,3) .. controls (.75,2) and (-.75,1) .. (-.75,0);
	\draw [very thick] (-2,0) rectangle (2,3);
	\draw[very thick, postaction={decorate}] (1.75,2) -- (.75,3);
	\draw[very thick, postaction={decorate}] (-.75,3) -- (-1.75,4);
	\draw[very thick, red] (.25,-1) .. controls (.25,0) and (1.75,1) .. (1.75,2);
	\draw [very thick] (-1,-1) rectangle (3,2);
	\node[red,opacity=1] at (2.5,1.75) {\footnotesize $a_{i+2}$};
	\node[red,opacity=1] at (1.5,.25) {\footnotesize $a_{i+1}$};
	\node[red,opacity=1] at (.75,3.75) {\footnotesize $a_i$};
	\node[blue,opacity=1] at (1.25,1.5) {\footnotesize $k$};
	\node[blue,opacity=1] at (-1,2.75) {\footnotesize $l$};
\end{tikzpicture}}
\endxy
\]
is $-kl$, which follows from the computations
\begin{align*}
 \chi^{Top}&= (a_{i}+l)(n-a_{i}-l)+a_{i}l+(a_{i+1}+k)(n-a_{i+1}-k)+l(a_{i+1}+k-l) +ka_{i+1}\\ &\quad +a_{i+2}(n-a_{i+2}) + k(a_{i+2}-k) - k(n-k) -l(n-l) \\
\chi^{Bot}&= a_{i+2}(n-a_{i+2})+k(a_{i+2}-k)+(a_{i+1}+k-l)(n-a_{i+1}-k+l)+k(a_{i+1}-l)+a_{i+1}(n-a_{i+1})\\ 
	&\quad+l(a_{i+1}-l)+(a_{i}+l)(n-a_{i}-l) +a_{i}l-k(n-k)-(a_{i+1}-l)(n-a_{i+1}+l)-l(n-l) \\
\chi&= (a_{i+1}+k)(n-a_{i+1}-k)+ka_{i+1}+a_{i+1}l-l^2+a_{i+2}(n-a_{i+2})+k(a_{i+2}-k)\\ &\quad +(a_{i}+l)(n-a_{i}-l)+a_{i}l-k(n-k)-l(n-l) \;\; .
\end{align*}

It remains to show that the $\mathfrak{sl}_2$ relations and bubble relations 
hold in $\foam{n}{m}$. Let $\Phi_n^N$ denote either of the $2$-functors $\Phi_n^{2,N}$ or $\Phi_n^{3,N}$ for the duration.
For $\lambda > 0$, we have that
\[
\xy
(0,0)*{

}
\endxy
\]
which is isotopic (under the above correspondence) 
to the image of the $\Ucatc_Q^+(\mathfrak{sl}_3)$ $2$-morphism
\[
(-1)^{a-\lambda} \quad
\xy
(0,0)*{
%
}
\endxy 
\]
have isotopic images. It follows that the image of the dotted bubble
\[
\xy
(0,0)*{
%
}
\endxy 
\]
is isotopic to the image of
\begin{equation}\label{cwKLRbubble}
(-1)^{a-\lambda} \quad
\xy
(0,0)*{
%
}
\endxy 
\quad = 0
\]
unless $\hat{\alpha}+l \geq a-1$, the inequality $\hat{\alpha} \leq a-\lambda$ gives that the image of an
$l$-dotted bubble is zero unless $l \geq \lambda-1$, i.e. the image of a negative degree 
clockwise bubble is zero.

Equation \eqref{cwKLRbubble} gives that the image of a non-negative degree bubble 
\[
\xy
(0,0)*{
%
}
\endxy 
\]
is isotopic to the image of
\begin{equation}\label{cwKLR+bubble}
\sum_{i=0}^c (-1)^i \quad
\xy
(0,0)*{
%
}
\endxy 
\end{equation}
under $\Phi_n^{2a-\lambda}$, which shows that a degree-zero counterclockwise bubble
maps to the identity. This also implies that the image of the fake bubble
\[
\xy
(0,0)*{
%
}
\endxy 
\]
is equal (up to isotopy) to the image of the $2$-morphism
\begin{equation}\label{ccwKLRfakebubble}
\sum_{i=0}^c (-1)^i \quad
\xy
(0,0)*{
%
}
\endxy
\end{equation}
for $\lambda>0$.
Indeed, this follows from \eqref{eq_infinite_Grass} and the computation:
\begin{align*}
\sum_{c=0}^{r} \left(
\sum_{i=0}^c (-1)^i \quad
\xy
(0,0)*{
%
}
\endxy = 0
\end{align*}
for $0<r<\l-1$, where we've used the standard fact that $\sum_{j=0}^s (-1)^j h_{s-j}e_j = 0$ if $s>0$.

Noting that the $2$-morphisms
\[
\xy
(0,0)*{
%
}
\endxy 
\]
have isotopic images, we find that the $\mathfrak{sl}_2$ equation \eqref{eq_ident_decomp-ngeqz} holds in $\foam{n}{m}$ if and 
only if the image of the $2$-morphism
\[
\xy
(0,0)*{
%
}
\endxy 
\]
is the identity foam. Since
\begin{align*}
\sum_{\xy  (0,3)*{\scs b+c+d}; (0,0)*{\scs =\l-1};\endxy} \ \ \sum_{i=0}^{c} (-1)^i \quad
\xy
(0,0)*{
%
}
\endxy
\end{align*}
this follows from the thick calculus relation:
\begin{equation} \label{KLRbGB}
\xy
(0,0)*{
%
}
\endxy 
\end{equation}
which can be extracted from the computations in \cite{Stosic} with signs conveniently adjusted.

Note that the second part of equation \eqref{eq_reduction-nleqz} and the first part of equation \eqref{eq_reduction-neqz_2} 
also follow from equation \eqref{KLRbGB} since the last term is zero when $\lambda \leq 0$.

Similarly, if $\lambda < 0$ then
\[
\xy
(0,0)*{
%
}
\endxy
\]
which is isotopic to the image of the $2$-morphism
\[
\xy
(0,0)*{
%
}
\endxy 
\]
as are the images of
\[
\xy
(0,0)*{
%
}
\endxy
\]
under $\Phi_n^{2a+\l}$. We again deduce that the image of the dotted bubble:
\[
\xy
(0,0)*{
%
}
\endxy 
\]
is equal (up to isotopy) to the image of:
\begin{align*}\label{ccwKLRbubble}
(-1)^{a+1} \quad
\xy
(0,0)*{
%
}
\endxy 
\end{align*}
which gives that negative-degree bubbles are zero and that the image 
of a non-negative degree bubble 
\[
\xy
(0,0)*{
%
}
\endxy 
\]
is isotopic to the image of
\begin{equation}\label{ccwKLR+bubble}
\sum_{i=0}^c (-1)^i \quad
\xy
(0,0)*{
%
}
\endxy .
\end{equation}
As above, this implies that the fake bubble
\[
\xy
(0,0)*{
%
}
\endxy 
\]
is isotopic to the image of
\begin{equation}\label{cwKLRfakebubble}
\sum_{i=0}^c (-1)^i \quad
\xy
(0,0)*{
%
}
\endxy 
\end{equation}
under $\Phi_n^{2a+\lambda}$. We again find that \eqref{eq_ident_decomp-nleqz} holds in $\foam{n}{m}$ if and 
only if the image of the $2$-morphism
\[
\xy
(0,0)*{
%
}
\endxy 
\]
is the identity foam; this in turn follows from the thick calculus relation
\begin{equation}\label{KLRGBg}
\xy
(0,0)*{
%
}
\endxy .
\end{equation}
Equation \eqref{KLRGBg} also shows that the second part of equation \eqref{eq_reduction-ngeqz} and the second part of 
equation \eqref{eq_reduction-neqz_2} since the last term is zero when $\lambda \geq 0$.

Finally, we analyze the images of the left and right curls. The images of the $2$-morphisms
\[
  \xy 0;/r.17pc/:
  (14,8)*{\l};
  (-3,-10)*{};(3,5)*{} **\crv{(-3,-2) & (2,1)}?(1)*\dir{>};?(.15)*\dir{>};
    (3,-5)*{};(-3,10)*{} **\crv{(2,-1) & (-3,2)}?(.85)*\dir{>} ?(.1)*\dir{>};
  (3,5)*{}="t1";  (9,5)*{}="t2";
  (3,-5)*{}="t1'";  (9,-5)*{}="t2'";
   "t1";"t2" **\crv{(4,8) & (9, 8)};
   "t1'";"t2'" **\crv{(4,-8) & (9, -8)};
   "t2'";"t2" **\crv{(10,0)} ;
 \endxy
\quad \quad \text{and} \quad \quad
(-1)^{a-\l} \quad
\xy
(0,0)*{
\begin{tikzpicture} [scale=.5]
\draw[ultra thick, blue, directed=.95] (0,-4) to (0,-3) to [out=90,in=270] (1,0)
	to [out=90,in=270] (0,3) to (0,4);
\draw[black!30!green, directed=.7] (-1,-4) to (-1,-1) to [out=90,in=270] (0,1) to [out=90,in=210] (1,3);
\draw[black!30!green, directed=.95] (1,-3) to [out=150,in=270] (0,-1) to [out=90,in=270] (-1,1) to (-1,4);
\draw[ultra thick, black!30!green, directed=.95] (1,-3) to [out=30,in=270] (2,0);
\draw[ultra thick, black!30!green] (2,0) to [out=90,in=330] (1,3);
\draw[ultra thick, black!30!green, directed=.5] (1,-4) to (1,-3);
\draw[ultra thick, black!30!green, directed=.5] (1,3) to (1,4);
\node[opacity=1] at (4,2) {$\mu = (\textcolor{black!30!green}{-a} ,\textcolor{blue}{\lambda})$};
\node[opacity=1, black!30!green] at (-1,-4.5) {\tiny$1$};
\node[opacity=1, blue] at (0,-4.5) {\tiny$a-\l$};
\node[opacity=1, black!30!green] at (-1,4.5) {\tiny$1$};
\node[opacity=1, black!30!green] at (1,4.5) {\tiny$a$};
\draw[ultra thick, blue,directed=.95] (-2,4) to (-2,-4);
\draw[ultra thick, black!30!green,directed=.95] (-3,4) to (-3,-4);
\node[blue,opacity=1] at (-2,-4.5) {\tiny$a-\lambda$};
\node[blue,opacity=1] at (-2,-5) {\tiny$-1$};
\node[black!30!green,opacity=1] at (-3,4.5) {\tiny$a+1$};
\end{tikzpicture}
}
\endxy 
\]
under $\Phi_n^{2a-\l}$ are isotopic; the latter $2$-morphism is equal to
\[
(-1)^{a-\l} \quad
\xy
(0,0)*{
\begin{tikzpicture} [scale=.5]
\draw[ultra thick, blue, directed=.95] (0,-4) to (0,-3) to [out=90,in=270] (-1.5,0)
	to [out=90,in=270] (0,3) to (0,4);
\draw[black!30!green, directed=.7] (-1,-4) to (-1,-1) to [out=90,in=270] (0,1) to [out=90,in=210] (1,3);
\draw[black!30!green, directed=.95] (1,-3) to [out=150,in=270] (0,-1) to [out=90,in=270] (-1,1) to (-1,4);
\draw[ultra thick, black!30!green, directed=.95] (1,-3) to [out=30,in=270] (2,0);
\draw[ultra thick, black!30!green] (2,0) to [out=90,in=330] (1,3);
\draw[ultra thick, black!30!green, directed=.5] (1,-4) to (1,-3);
\draw[ultra thick, black!30!green, directed=.5] (1,3) to (1,4);
\node[opacity=1, black!30!green] at (-1,-4.5) {\tiny$1$};
\node[opacity=1, blue] at (0,-4.5) {\tiny$a-\l$};
\node[opacity=1, black!30!green] at (-1,4.5) {\tiny$1$};
\node[opacity=1, black!30!green] at (1,4.5) {\tiny$a$};
\end{tikzpicture}
}
\endxy 
\quad + \quad (-1)^{a-\l} 
\sum_{\alpha,\beta,\gamma \in P(1,a-\l-1)}
(-1)^{\alpha+\beta} c_{\alpha \beta \gamma}^{(a-\l-1)} \quad
\xy
(0,0)*{
\begin{tikzpicture} [scale=.5]
\draw[ultra thick, blue, directed=.95] (-1,-4) to (-1,4);
\draw[black!30!green, directed=.7] (-2,-4) to (-2,4);
\draw[ultra thick, black!30!green, directed=.5] (1,-2) to [out=30,in=330] (1,2);
\draw[black!30!green, directed=.5] (1,-2) to [out=150,in=210] (1,2);
\draw[ultra thick, black!30!green, directed=.5] (1,-4) to (1,-2);
\draw[ultra thick, black!30!green, directed=.5] (1,2) to (1,4);
\node[opacity=1, black!30!green] at (0,-1.5) {\tiny$1$};
\node[opacity=1, blue] at (-1,-4.5) {\tiny$a-\l$};
\node[opacity=1, black!30!green] at (-2,4.5) {\tiny$1$};
\node[opacity=1, black!30!green] at (1,4.5) {\tiny$a$};
\node[opacity=1] at (-2,0) {$\bullet$};
\node[opacity=1] at (-1,2) {$\bullet$};
\node[opacity=1] at (0,0) {$\bullet$};
\node[opacity=1] at (-1.5,0) {$\alpha$};
\node[opacity=1] at (-.5,2.25) {$\pi_{\bar{\gamma}}$};
\node[opacity=1] at (.5,0) {$\beta$};
\end{tikzpicture}
}
\endxy \quad .
\]
The relevant relations are then satisfied since the image of the first term is always zero (for weight reasons), 
the summation is zero when $\l>0$, and equals negative the identity $2$-morphism when $\l=0$.

The $2$-morphisms
\[
  \xy 0;/r.17pc/:
  (-14,8)*{\l};
  (3,-10)*{};(-3,5)*{} **\crv{(3,-2) & (-2,1)}?(1)*\dir{>};?(.15)*\dir{>};
    (-3,-5)*{};(3,10)*{} **\crv{(-2,-1) & (3,2)}?(.85)*\dir{>} ?(.1)*\dir{>};
  (-3,5)*{}="t1";  (-9,5)*{}="t2";
  (-3,-5)*{}="t1'";  (-9,-5)*{}="t2'";
   "t1";"t2" **\crv{(-4,8) & (-9, 8)};
   "t1'";"t2'" **\crv{(-4,-8) & (-9, -8)};
   "t2'";"t2" **\crv{(-10,0)} ;
 \endxy
\quad \quad \text{and} \quad \quad
(-1)^{a+1} \quad
\xy
(0,0)*{
\begin{tikzpicture} [scale=.5]
\draw[ultra thick, black!30!green, directed=.95] (0,-4) to (0,-3) to [out=90,in=270] (-1,0)
	to [out=90,in=270] (0,3) to (0,4);
\draw[blue, directed=.7] (1,-4) to (1,-1) to [out=90,in=270] (0,1) to [out=90,in=330] (-1,3);
\draw[blue, directed=.95] (-1,-3) to [out=30,in=270] (0,-1) to [out=90,in=270] (1,1) to (1,4);
\draw[ultra thick, blue, directed=.95] (-1,-3) to [out=150,in=270] (-2,0);
\draw[ultra thick, blue] (-2,0) to [out=90,in=210] (-1,3);
\draw[ultra thick, blue, directed=.5] (-1,-4) to (-1,-3);
\draw[ultra thick, blue, directed=.5] (-1,3) to (-1,4);
\node[opacity=1] at (4.75,2) {$\mu = (\textcolor{black!30!green}{-a-\lambda+1} ,\textcolor{blue}{\lambda-2})$};
\node[opacity=1, blue] at (1,-4.5) {\tiny$1$};
\node[opacity=1, black!30!green] at (0,-4.5) {\tiny$a+\l$};
\node[opacity=1, blue] at (-1,4.5) {\tiny$a$};
\node[opacity=1, blue] at (1,4.5) {\tiny$1$};
\draw[ultra thick, blue,directed=.95] (-3,4) to (-3,-4);
\draw[ultra thick, black!30!green,directed=.95] (-4,4) to (-4,-4);
\node[blue,opacity=1] at (-3,-4.5) {\tiny$a$};
\node[black!30!green,opacity=1] at (-4,4.5) {\tiny$a+\lambda$};
\end{tikzpicture}
}
\endxy 
\]
have isotopic images under $\Phi_n^{2a+\l}$. We find that the latter $2$-morphism equals
\[
(-1)^{a+1} \quad
\xy
(0,0)*{
\begin{tikzpicture} [scale=.5]
\draw[ultra thick, black!30!green, directed=.95] (0,-4) to (0,-3) to [out=90,in=270] (1.5,0)
	to [out=90,in=270] (0,3) to (0,4);
\draw[blue, directed=.7] (1,-4) to (1,-1) to [out=90,in=270] (0,1) to [out=90,in=330] (-1,3);
\draw[blue, directed=.95] (-1,-3) to [out=30,in=270] (0,-1) to [out=90,in=270] (1,1) to (1,4);
\draw[ultra thick, blue, directed=.95] (-1,-3) to [out=150,in=270] (-2,0);
\draw[ultra thick, blue] (-2,0) to [out=90,in=210] (-1,3);
\draw[ultra thick, blue, directed=.5] (-1,-4) to (-1,-3);
\draw[ultra thick, blue, directed=.5] (-1,3) to (-1,4);
\node[opacity=1, blue] at (1,-4.5) {\tiny$1$};
\node[opacity=1, black!30!green] at (0,-4.5) {\tiny$a+\l$};
\node[opacity=1, blue] at (-1,4.5) {\tiny$a$};
\node[opacity=1, blue] at (1,4.5) {\tiny$1$};
\draw[ultra thick, blue,directed=.95] (-3,4) to (-3,-4);
\draw[ultra thick, black!30!green,directed=.95] (-4,4) to (-4,-4);
\node[blue,opacity=1] at (-3,-4.5) {\tiny$a$};
\node[black!30!green,opacity=1] at (-4,4.5) {\tiny$a+\lambda$};
\end{tikzpicture}
}
\endxy 
\quad + \quad (-1)^{a+1} 
\sum_{\alpha,\beta,\gamma \in P(1,a+\l-1)}
(-1)^{\gamma} c_{\alpha \beta \gamma}^{(a+\l-1)} \quad
\xy
(0,0)*{
\begin{tikzpicture} [scale=.5]
\draw[blue, directed=.95] (1,-4) to (1,4);
\draw[ultra thick, black!30!green, directed=.95] (0,-4) to (0,4); 
\draw[ultra thick, blue, directed=.5] (-2,-4) to (-2,-2);
\draw[ultra thick, blue, directed=.5] (-2,2) to (-2,4);
\draw[blue, directed=.5] (-2,-2) to [out=30,in=330] (-2,2);
\draw[ultra thick, blue, directed=.5] (-2,-2) to [out=150,in=210] (-2,2);
\node[opacity=1] at (-1,0) {$\bullet$};
\node[opacity=1] at (0,2) {$\bullet$};
\node[opacity=1] at (1,0) {$\bullet$};
\node[opacity=1, blue] at (-1,-1.5) {\tiny$1$};
\node[opacity=1, black!30!green] at (0,-4.5) {\tiny$a+\l$};
\node[opacity=1, blue] at (-2,4.5) {\tiny$a$};
\node[opacity=1, blue] at (1,4.5) {\tiny$1$};
\node[opacity=1] at (-1.5,0) {$\alpha$};
\node[opacity=1] at (-.5,2.25) {$\pi_{\bar{\gamma}}$};
\node[opacity=1] at (1.5,0) {$\beta$};
\draw[ultra thick, blue,directed=.95] (-4,4) to (-4,-4);
\draw[ultra thick, black!30!green,directed=.95] (-5,4) to (-5,-4);
\node[blue,opacity=1] at (-4,-4.5) {\tiny$a$};
\node[black!30!green,opacity=1] at (-5,4.5) {\tiny$a+\lambda$};
\end{tikzpicture}
}
\endxy
\]
which again confirms the $\mathfrak{sl}_2$ relations since the image of the first $2$-morphism is 
zero, the image of the summation is zero if $\l <0$, and equals the identity for $\l=0$.
\end{proof}

\begin{rem} \label{rem:turningtoKLR}
Remark \ref{rem:easyisotopy} below shows that the web isotopies in equation \eqref{web_trick} actually follow from 
isomorphisms in $\Ucatc(\glm)^{\geq0}$. Repeated use of this allows one to identify $1$-morphisms in $\Ucatc(\glm)^{\geq0}$ with those 
in the KLR part $\Ucatc^+(\glnn{m'})^{\geq0}$ for $m'$ sufficiently\footnote{The specific value of $m'$ is determined by the combinatorial data of the web.}
large. Under closer inspection, the above proof then actually 
shows that the $\slnn{2}$ relations in $\Ucatc(\glm)^{\geq0}$ follow from the relations in $\Ucatc^+(\glnn{m'})^{\geq0}$. A version of this result was 
independently discovered in \cite{TubbenhauerSln} after passing to matrix factorizations, where Tubbenhauer uses only the KLR algebra to describe the complexes 
giving Khovanov-Rozansky homology. 
The naturality of this identification explains the `$F$-braiding' used in that work\footnote{Note that $F$'s in \cite{TubbenhauerSln} correspond to $E$'s here, by our 
diagrammatic conventions.}; it is equivalent to the Rickard complex in $\Ucatc(\glm)^{\geq0}$ 
via the isomorphisms presented in Lemma \ref{easyisotopy}.
\end{rem}

\begin{rem} \label{rem:SquareFormulas}
The existence of $2$-representations $\Ucatc_Q(\mathfrak{sl}_m) \to \foam{n}{m}$ implies that isomorphic $1$-morphisms in the former are sent to 
isomorphic $1$-morphisms in the latter. In particular, when $a-b\geq l - k$ we have isomorphisms 
\begin{equation}\label{EF_to_cb}
\xy
(0,0)*{
\begin{tikzpicture} [scale=.6]
\draw[very thick, directed=.15, directed=.5, directed=.9] (2.25,1) to (-2.25,1);
\draw[very thick, directed=.2, directed=.5, directed=.8] (2.25,0) to (-2.25,0);
\draw[very thick, directed=.55] (1.25,1) to (.5,0);
\draw[very thick, directed=.55] (-.5,0) to (-1.25,1);
\node at (2.5,1) {\small$a$};
\node at (2.5,0) {\small$b$};
\node at (-1.25,.5) {\small$k$};
\node at (1.25,.5) {\small$l$};
\end{tikzpicture}
}
\endxy
\quad \cong
\bigoplus_{j=0}^{\min(k,l)}
\bigoplus_{\alpha \in P(j,a+k-b-l-j)}
q^{j(a+k-b-l)-2|\alpha|} \quad
\xy
(0,0)*{
\begin{tikzpicture} [scale=.6]
\draw[very thick, directed=.2, directed=.5, directed=.8] (2.25,1) to (-2.25,1);
\draw[very thick, directed=.15, directed=.5, directed=.9] (2.25,0) to (-2.25,0);
\draw[very thick, directed=.55] (1.25,0) to (.5,1);
\draw[very thick, directed=.55] (-.5,1) to (-1.25,0);
\node at (2.5,1) {\small$a$};
\node at (2.5,0) {\small$b$};
\node at (-1.5,.5) {\small$_{l-j}$};
\node at (1.5,.5) {\small$_{k-j}$};
\end{tikzpicture}
}
\endxy
\end{equation}
and when $a-b \leq l-k$ we have isomorphisms
\begin{equation}\label{FE_to_cb}
\xy
(0,0)*{
\begin{tikzpicture} [scale=.6]
\draw[very thick, directed=.2, directed=.5, directed=.8] (2.25,1) to (-2.25,1);
\draw[very thick, directed=.15, directed=.5, directed=.9] (2.25,0) to (-2.25,0);
\draw[very thick, directed=.55] (1.25,0) to (.5,1);
\draw[very thick, directed=.55] (-.5,1) to (-1.25,0);
\node at (2.5,1) {\small$a$};
\node at (2.5,0) {\small$b$};
\node at (-1.5,.5) {\small$l$};
\node at (1.5,.5) {\small$k$};
\end{tikzpicture}
}
\endxy
\quad \cong
\bigoplus_{j=0}^{\min(k,l)}
\bigoplus_{\alpha \in P(j,b+l-a-k-j)}
q^{j(b+l-a-k)-2|\alpha|} \quad
\xy
(0,0)*{
\begin{tikzpicture} [scale=.6]
\draw[very thick, directed=.15, directed=.5, directed=.9] (2.25,1) to (-2.25,1);
\draw[very thick, directed=.2, directed=.5, directed=.8] (2.25,0) to (-2.25,0);
\draw[very thick, directed=.55] (1.25,1) to (.5,0);
\draw[very thick, directed=.55] (-.5,0) to (-1.25,1);
\node at (2.5,1) {\small$a$};
\node at (2.5,0) {\small$b$};
\node at (-1.5,.5) {\small$_{k-j}$};
\node at (1.5,.5) {\small$_{l-j}$};
\end{tikzpicture}
}
\endxy
\end{equation}
which follow from \cite[Theorem 5.9]{KLMS}. 
Moreover, one can obtain the explicit foams giving 
these isomorphisms by applying $\Phi_n$ to the relevant $2$-morphisms in that paper.
Note that if any label in the above webs lies outside 
the range $\{0,\ldots,n\}$, then that web is zero in $\foam{n}{m}$.
\end{rem}

%
\subsection{Properties of the foamation functors}
%

In this section, we prove certain fullness and faithfulness properties of the foamation functors. These will 
later serve to show that our link homology theory can be computed entirely combinatorially, and that 
the $2$-categories $\foam{n}{m}$ are sufficiently non-degenerate to give a non-trivial invariant. 
For the duration of the discussion, we'll fix the index $n$ and let $\Phi^{m,N}$ denote the $2$-functor
$\Ucatc_Q(\mathfrak{sl}_m) \to \foam{n}{m}$. By definition, this $2$-functor factors through the 
$2$-categories $\Ucatc_Q(\glm)^{\geq0}$ and $\Ucatc_Q(\glm)^{(n)}$ from Definition \eqref{def:qSchur}, 
facts which we'll use to analyze properties of the foamation $2$-functors.

\begin{lem}\label{easyisotopy}
Let $W_1$ and $W_2$ be two ladder webs with equal boundary and which are planar isotopic, 
then there exists $m$ and a foam $F:W_1 \to W_2$ in the image of $\Phi^{m,N}$ realizing the web isotopy.
\end{lem}
\begin{proof}
Any two isotopic webs in ladder form are related via a sequence of the following moves: \\
\begin{equation} \label{EdgeBend}
\xymatrix{
\xy
(0,0)*{
\begin{tikzpicture} [scale=.45]
\draw[very thick, directed=.5] (2.5,0) to (1.25,0) to (.5,-1) to (-.5,-1) to (-1.25,0) to (-2.5,0);
\draw[dotted] (2.5,1) to (-2.5,1);
\draw[dotted] (2.5,-1) to (-2.5,-1);
\node at (3,0) {\small $a$};
\end{tikzpicture}};
\endxy
\ar@/^/[rr]^{ (-1)^{\gamma}
\xy 
(0,0)*{
\begin{tikzpicture}[scale=.6]
\draw[ultra thick, black!30!green, <-] (1,0) .. controls (1,.8) and (0,.8) .. (0,0);
\end{tikzpicture}
};
\endxy} 
\quad & & \quad
\xy
(0,0)*{
\begin{tikzpicture} [scale=.45]
\draw[very thick, directed=.55] (2.5,0) to (-2.5,0);
\draw[dotted] (2.5,1) to (-2.5,1);
\draw[dotted] (2.5,-1) to (-2.5,-1);
\node at (3,0) {\small $a$};
\end{tikzpicture}};
\endxy
\ar@/^/[ll]^{
\xy 
(0,0)*{
\begin{tikzpicture}[scale=.6]
\draw[ultra thick, black!30!green, ->] (1,1) .. controls (1,.2) and (0,.2) .. (0,1);
\end{tikzpicture}
};
\endxy}
\ar@/^/[rr]^{ (-1)^\gamma
\xy 
(0,0)*{
\begin{tikzpicture}[scale=.6]
\draw[ultra thick, black!30!green, <-] (1,1) .. controls (1,.2) and (0,.2) .. (0,1);
\node [black!30!green] at (0,1.3) {\small $a$};
\end{tikzpicture}
};
\endxy} 
\quad & & \quad
\xy
(0,0)*{
\begin{tikzpicture} [scale=.45]
\draw[very thick, directed=.5] (2.5,0) to (1.25,0) to (.5,1) to (-.5,1) to (-1.25,0) to (-2.5,0);
\draw[dotted] (2.5,1) to (-2.5,1);
\draw[dotted] (2.5,-1) to (-2.5,-1);
\node at (3,0) {\small $a$};
\end{tikzpicture}};
\endxy
\ar@/^/[ll]^{
\xy 
(0,0)*{
\begin{tikzpicture}[scale=.6]
\draw[ultra thick, black!30!green, ->] (1,0) .. controls (1,.8) and (0,.8) .. (0,0);
\end{tikzpicture}
};
\endxy}
}
\end{equation} \\
\begin{equation} \label{CornerChange}
\xymatrix{
\xy
(0,0)*{
\begin{tikzpicture} [scale=.45]
\draw[very thick, directed=.4, directed=.85] (2.5,0) to (1.25,0) to (.5,-1) to (-2.5,-1);
\draw[very thick, directed=.3] (-.5,-1) to (-1.25,0) to (-2.5,0);
\draw[dotted] (2.5,0) to (-2.5,0);
\draw[dotted] (2.5,-1) to (-2.5,-1);
\node at (-3,0) {\small $a$};
\node at (-3,-1) {\small $b$};
\end{tikzpicture}};
\endxy
\ar@/^/[rr]^{ (-1)^\gamma
\xy 
(0,0)*{
\begin{tikzpicture}[scale=.4]
\draw[ultra thick, black!30!green, directed=.35] (0,.25) to [out=90,in=180] (.375,1.5) to [out=0,in=180] (1,1);
\draw[ultra thick, black!30!green, directed=.55] (1.5,2.25) to [out=270,in=0] (1,1);
\draw[ultra thick, black!30!green, directed=.85] (1,1) to (1,.25);
\end{tikzpicture}
};
\endxy}
\quad & & \quad
\xy
(0,0)*{
\begin{tikzpicture} [scale=.45,fill opacity=.2]
\draw[very thick, directed=.2, directed=.8] (2,0) to (-2,0);
\draw[very thick, directed=.25] (.375,0) to (-.375,-1) to (-2,-1);
\draw[dotted] (2,0) to (-2,0);
\draw[dotted] (2,-1) to (-2,-1);
\end{tikzpicture}};
\endxy 
\ar@/^/[ll]^{
\xy 
(0,0)*{
\begin{tikzpicture}[scale=.4]
\draw[ultra thick, black!30!green, rdirected=.35] (0,-.25) to [out=270,in=180] (.375,-1.5) to [out=0,in=180] (1,-1);
\draw[ultra thick, black!30!green, rdirected=.55] (1.5,-2.25) to [out=90,in=0] (1,-1);
\draw[ultra thick, black!30!green, rdirected=.85] (1,-1) to (1,-.25);
\end{tikzpicture}
};
\endxy}
}
,
\xymatrix{
\xy
(0,0)*{
\begin{tikzpicture} [scale=.45]
\draw[very thick, directed=.35, directed=.7] (2.5,0) to (1.25,0) to (.5,-1) to (-.5,-1) to (-1.25,0) to (-2.5,0);
\draw[very thick, directed=.5] (2.5,-1) to (.5,-1);
\draw[dotted] (2.5,0) to (-2.5,0);
\draw[dotted] (2.5,-1) to (-2.5,-1);
\end{tikzpicture}};
\endxy
\quad 
\ar@/^/[rr]^{ (-1)^{\gamma + ab} \; \;
\xy 
(0,0)*{
\begin{tikzpicture}[scale=.4]
\draw[ultra thick, black!30!green, rdirected=.35] (0,.25) to [out=90,in=0] (-.375,1.5) to [out=180,in=0] (-1,1);
\draw[ultra thick, black!30!green, rdirected=.55] (-1.5,2.25) to [out=270,in=180] (-1,1);
\draw[ultra thick, black!30!green, rdirected=.85] (-1,1) to (-1,.25);
\end{tikzpicture}
};
\endxy}
& & \quad
\xy
(0,0)*{
\begin{tikzpicture} [scale=.45]
\draw[very thick, directed=.2, directed=.8] (2,0) to (-2,0);
\draw[very thick, directed=.85] (2,-1) to (.375,-1) to (-.375,0);
\draw[dotted] (2,0) to (-2,0);
\draw[dotted] (2,-1) to (-2,-1);
\node at (2.5,0) {\small $a$};
\node at (2.5,-1) {\small $b$};
\end{tikzpicture}};
\endxy
\ar@/^/[ll]^{
\xy 
(0,0)*{
\begin{tikzpicture}[scale=.4]
\draw[ultra thick, black!30!green, directed=.35] (0,-.25) to [out=270,in=0] (-.375,-1.5) to [out=180,in=0] (-1,-1);
\draw[ultra thick, black!30!green, directed=.55] (-1.5,-2.25) to [out=90,in=180] (-1,-1);
\draw[ultra thick, black!30!green, directed=.85] (-1,-1) to (-1,-.25);
\end{tikzpicture}
};
\endxy}
}
\end{equation} \\
\begin{equation*}
\xymatrix{
\xy
(0,0)*{
\begin{tikzpicture} [scale=.45]
\draw[very thick, directed=.2, directed=.6] (2.5,0) to (0,0) to (-.75,1) to (-2.25,1);
\draw[very thick, directed=.25] (.75,0) to (0,-1) to (-2.25,-1);
\draw[dotted] (2.5,1) to (-2.25,1);
\draw[dotted] (2.5,0) to (-2.25,0);
\draw[dotted] (2.5,-1) to (-2.25,-1);
\node at (-2.75,1) {\small $a$};
\node at (-2.75,-1) {\small $b$};
\end{tikzpicture}};
\endxy
\quad 
\ar@/^/[rr]^{ 
\xy
(0,0)*{
\begin{tikzpicture} [scale=.25]
\draw[ultra thick, black!30!green, directed=.99] (0,0) to [out=90,in=270] (2,3);
\draw[ultra thick, blue, directed=.99] (0,3) to [out=270,in=90] (2,0);
\end{tikzpicture}
};
\endxy}
& & \quad
\xy
(0,0)*{
\begin{tikzpicture} [scale=.45]
\draw[very thick, directed=.2, directed=.6] (2.5,0) to (0,0) to (-.75,-1) to (-2.25,-1);
\draw[very thick, directed=.25] (.75,0) to (0,1) to (-2.25,1);
\draw[dotted] (2.5,1) to (-2.25,1);
\draw[dotted] (2.5,0) to (-2.25,0);
\draw[dotted] (2.5,-1) to (-2.25,-1);
\end{tikzpicture}};
\endxy
\ar@/^/[ll]^{
\xy
(0,0)*{
\begin{tikzpicture} [scale=.25]
\draw[ultra thick, black!30!green, directed=.99] (2,0) to [out=90,in=270] (0,3);
\draw[ultra thick, blue, directed=.99] (2,3) to [out=270,in=90] (0,0);
\end{tikzpicture}
};
\endxy}
}
, 
\xymatrix{
\xy
(0,0)*{
\begin{tikzpicture} [scale=.45]
\draw[very thick, directed=.85, directed=.45] (2.25,1) to (.75,1) to (0,0) to (-2.5,0);
\draw[very thick, directed=.85] (2.25,-1) to (0,-1) to  (-.75,0);
\draw[dotted] (2.25,1) to (-2.5,1);
\draw[dotted] (2.25,0) to (-2.5,0);
\draw[dotted] (2.25,-1) to (-2.5,-1);
\end{tikzpicture}};
\endxy
\quad
\ar@/^/[rr]^{ (-1)^{ab}
\xy
(0,0)*{
\begin{tikzpicture} [scale=.25]
\draw[ultra thick, blue, directed=.99] (0,0) to [out=90,in=270] (2,3);
\draw[ultra thick, black!30!green, directed=.99] (0,3) to [out=270,in=90] (2,0);
\end{tikzpicture}
};
\endxy}
& & \quad
\xy
(0,0)*{
\begin{tikzpicture} [scale=.45]
\draw[very thick, directed=.85, directed=.45] (2.25,-1) to (.75,-1) to (0,0) to (-2.5,0);
\draw[very thick, directed=.85] (2.25,1) to (0,1) to  (-.75,0);
\draw[dotted] (2.25,1) to (-2.5,1);
\draw[dotted] (2.25,0) to (-2.5,0);
\draw[dotted] (2.25,-1) to (-2.5,-1);
\node at (2.75,1) {\small $a$};
\node at (2.75,-1) {\small $b$};
\end{tikzpicture}};
\endxy
\ar@/^/[ll]^{
\xy
(0,0)*{
\begin{tikzpicture} [scale=.25]
\draw[ultra thick, blue, directed=.99] (2,0) to [out=90,in=270] (0,3);
\draw[ultra thick, black!30!green, directed=.99] (2,3) to [out=270,in=90] (0,0);
\end{tikzpicture}
};
\endxy}
}
\end{equation*} \\
all of which are realized by foams given by the image of the indicated thick calculus morphisms under $\Phi^{m,N}$. 
In these formulae, the thickness of the green and blue strands can be deduced from the webs, and $\gamma = \frac{a(a-1)}{2}$.
Note that we need to take $m$ sufficiently large in order to have enough $0$-labeled edges to perform the requisite moves.
\end{proof}

\begin{rem}\label{rem:easyisotopy}
One can check that in fact the thick calculus morphisms indicated above are isomorphisms in $\Ucatc_Q(\glm)^{\geq0}$,
hence each web isotopy is actually induced by a thick calculus isomorphism.
In Lemma \ref{lem:elemIso} below, we'll prove certain coherence relations between these web isotopies and foam generators, which will 
allow us to deduce that foam isotopy is always induced by relations in $\Ucatc_Q(\glm)^{\geq0}$. 
\end{rem}

For now, we'll use Lemma \ref{easyisotopy} to give an ``eventual fullness" result for the family of $2$-functors $\{\Phi^{m,N}\}_{m \geq 2}$. Note 
that for $m' >m$ we have a natural inclusion $2$-functor $\foam{n}{m} \stackrel{\iota}{\lhook\joinrel\longrightarrow} \foam{n}{m'}$ given on objects by 
$[a_1,\ldots,a_m] \mapsto [a_1,\ldots,a_m,0,\ldots,0]$.
\begin{prop}\label{prop:fullness}
Let $F:q^{k_1}W_1 \to q^{k_2} W_2$ be a foam in $\foam{n}{m}$ between ladder webs, then there exists $m' \geq m$ and 
a $2$-morphism $\mathcal{D}$ in $\Ucatc_Q(\mathfrak{sl}_{m'})$ so that $\Phi^{m',N}(\mathcal{D}) = \iota(F)$.
\end{prop}

\begin{proof}
We can locally realize the foam generators from equations \eqref{ccgens}, \eqref{uzgens}, and \eqref{MVgensup} 
in the image of $\Phi$ as follows:
\begin{align*}
\Phi \left(
\xy
(0,0)*{
\begin{tikzpicture} [scale=.4,fill opacity=0.2]
\draw[ultra thick, black!30!green] (0,-1.5) to (0,0);
\draw[ultra thick, black!30!green, directed=1] (0,0) to [out=30,in=270] (1,1.5);
\draw[ultra thick, black!30!green, directed=1] (0,0) to [out=150,in=270] (-1,1.5);
\node[black!30!green, opacity=1] at (0,-2) {\tiny$_{a+b}$};
\node[black!30!green, opacity=1] at (-1,2) {\tiny$a$};
\node[black!30!green, opacity=1] at (1,2) {\tiny$b$};
\node[opacity=1] at (1.75,-.5) {\tiny$[0,a+b]$};
\end{tikzpicture}
}
\endxy
\right)
=
\xy
(0,0)*{
\begin{tikzpicture} [scale=.4,fill opacity=0.2]
	\path[fill=red] (3.5,7.5) to (1.75,7.5) to [out=270,in=0] (0,5) to [out=180,in=270] (-1.75,8.5) to (-3.5,8.5) to 
		(-3.5,3.5) to (-.5,3.5) to (.5,2.5) to (3.5,2.5);
	\path[fill=blue] (1.75,7.5) to [out=270,in=0] (0,5) to [out=180,in=270] (-1.75,8.5) to (-.75,7.5);
	\path[fill=blue] (1.75,7.5) to [out=270,in=0] (0,5) to [out=180,in=270] (-1.75,8.5) to (.75,8.5);
	\draw[very thick] (3.5,2.5) to (.5,2.5);
	\draw[very thick] (-.5,3.5) to (-3.5,3.5);
	\draw[very thick, directed=.55] (.5,2.5) to (-.5,3.5);	
	\draw[very thick, red, rdirected=.5] (1.75,7.5) to [out=270,in=0] (0,5) to [out=180,in=270] (-1.75,8.5);
	\draw[very thick] (3.5,2.5) to (3.5,7.5);	
	\draw[very thick] (-3.5,3.5) to (-3.5,8.5);
	\draw[very thick, directed=.2] (3.5,7.5) to (-.75,7.5);
	\draw[very thick, directed=.8] (.75,8.5) to (-3.5,8.5);
	\draw[very thick, directed=.5] (1.75,7.5) to (.75,8.5);	
	\draw[very thick, directed=.5] (-.75,7.5) to (-1.75,8.5);
	\node[blue, opacity=1] at (0,6.5) {\tiny$a$};
	\node[blue, opacity=1] at (0,8) {\tiny$b$};
	\node[red, opacity=1] at (2.75,7) {\tiny$_{a+b}$};
\end{tikzpicture}
};
\endxy 
\quad &, \quad
\Phi \left(
\xy
(0,0)*{
\begin{tikzpicture} [scale=.4,fill opacity=0.2]
\draw[ultra thick, black!30!green, directed=1] (0,0) to (0,1.5);
\draw[ultra thick, black!30!green] (1,-1.5) to [out=90,in=330] (0,0);
\draw[ultra thick, black!30!green] (-1,-1.5) to [out=90,in=210] (0,0);
\node[black!30!green, opacity=1] at (0,2) {\tiny$_{a+b}$};
\node[black!30!green, opacity=1] at (-1,-2) {\tiny$a$};
\node[black!30!green, opacity=1] at (1,-2) {\tiny$b$};
\node[opacity=1] at (1.75,.5) {\tiny$[0,a+b]$};
\end{tikzpicture}
}
\endxy
\right)
=
\xy
(0,0)*{
\begin{tikzpicture} [scale=.4,fill opacity=0.2]
	\path[fill=red] (3.5,0) to (1.75,0) to [out=90,in=0] (0,3.5) to [out=180,in=90] (-1.75,1) to (-3.5,1) to (-3.5,6) to
		(-.5,6) to (.5,5) to (3.5,5);
	\path[fill=blue] (1.75,0) to [out=90,in=0] (0,3.5) to [out=180,in=90] (-1.75,1) to (-.75,0);
	\path[fill=blue] (1.75,0) to [out=90,in=0] (0,3.5) to [out=180,in=90] (-1.75,1) to (.75,1);
	\draw[very thick, directed=.2] (3.5,0) to (-.75,0);
	\draw[very thick, directed=.8] (.75,1) to (-3.5,1);
	\draw[very thick, directed=.5] (1.75,0) to (.75,1);	
	\draw[very thick, directed=.5] (-.75,0) to (-1.75,1);		
	\draw[very thick, red, directed=.55] (1.75,0) to [out=90,in=0] (0,3.5) to [out=180,in=90] (-1.75,1);
	\draw[very thick] (3.5,0) to (3.5,5);	
	\draw[very thick] (-3.5,1) to (-3.5,6);
	\draw[very thick] (3.5,5) to (.5,5);
	\draw[very thick] (-.5,6) to (-3.5,6);
	\draw[very thick, directed=.55] (.5,5) to (-.5,6);
	\node[blue, opacity=1] at (0,.5) {\tiny$a$};
	\node[red, opacity=1] at (2.75,3) {\tiny$_{a+b}$};
	\node[blue, opacity=1] at (0,2.5) {\tiny$b$};
\end{tikzpicture}
};
\endxy \\ 
\Phi \left(
\xy
(0,0)*{
\begin{tikzpicture} [scale=.4,fill opacity=0.2]
\draw[ultra thick, black!30!green, directed=.99] (0,0) to [out=90,in=180] (1,1.5) to [out=0,in=90] (2,0);
\node[black!30!green, opacity=1] at (0,-.5) {\tiny$b$};
\node[opacity=1] at (1.5,2) {\tiny$[b,a]$};
\end{tikzpicture}
}
\endxy \;
\right)
=
\xy
(0,0)*{
\begin{tikzpicture} [scale=.4,fill opacity=0.2]
	\path[fill=blue] (-1.5,0) to [out=90,in=180] (-.5,2) to [out=0,in=90] (.5,0);
	\path[fill=red] (3.5,4) to (-3.5,4) to (-3.5,0) to (-1.5,0) to [out=90,in=180] (-.5,2) to [out=0,in=90] (.5,0) to (3.5,0);
	\path[fill=red] (2.5,5) to (-4.5,5) to (-4.5,1) to (-2.5,1) to (-1.5,0) to [out=90,in=180] (-.5,2) to [out=0,in=90] (.5,0)
		to (1.5,1) to (2.5,1);
	\draw[very thick, directed=.2, directed=.6, directed=.9] (3.5,0) to (-3.5,0);
	\draw[very thick, directed=.65] (2.5,1) to (1.5,1) to (.5,0);
	\draw[very thick, directed=.35] (-1.5,0) to (-2.5,1) to (-4.5,1);
	\draw[very thick] (3.5,0) to (3.5,4);
	\draw[very thick] (2.5,1) to (2.5,5);
	\draw[very thick] (-3.5,0) to (-3.5,4);
	\draw[very thick] (-4.5,1) to (-4.5,5);
	\draw (-2.5,1) to [out=90,in=180] (-.5,3) to [out=0,in=90] (1.5,1);
	\draw[very thick, red, directed=.3] (-1.5,0) to [out=90,in=180] (-.5,2) to [out=0,in=90] (.5,0);
	\draw[very thick, directed=.55] (3.5,4) to (-3.5,4);
	\draw[very thick, directed=.55] (2.5,5) to (-4.5,5);
	\node[red, opacity=1] at (-3,3.5) {\tiny$a$};
	\node[red, opacity=1] at (-4,4.5) {\tiny$b$};
\end{tikzpicture}
};
\endxy
\quad &, \quad
\Phi \left(
\xy
(0,0)*{
\begin{tikzpicture} [scale=.4,fill opacity=0.2]
\draw[ultra thick, black!30!green, directed=.99] (2,3) to [out=270,in=0] (1,1.5) to [out=180,in=270] (0,3);
\node[black!30!green, opacity=1] at (2,3.5) {\tiny$b$};
\node[opacity=1] at (1.5,1) {\tiny$[b,a]$};
\end{tikzpicture}
}
\endxy \;
\right)
=
\xy
(0,0)*{
\begin{tikzpicture} [scale=.4,fill opacity=0.2]
	\path[fill=blue] (.5,4) to [out=270,in=0] (-.5,2.5) to [out=180,in=270] (-1.5,4);
	\path[fill=red] (3.5,4) to (.5,4) to [out=270,in=0] (-.5,2.5) to [out=180,in=270] (-1.5,4) to (-3.5,4) to (-3.5,0) to (3.5,0);
	\path[fill=red] (2.5,5) to (1.5,5) to (.5,4) to [out=270,in=0] (-.5,2.5) to [out=180,in=270] (-1.5,4) to (-2.5,5) to (-4.5,5)
		to (-4.5,1) to (2.5,1);
	\draw[very thick, directed=.55] (3.5,0) to (-3.5,0);
	\draw[very thick, directed=.55] (2.5,1) to (-4.5,1);
	\draw[very thick] (3.5,0) to (3.5,4);
	\draw[very thick] (2.5,1) to (2.5,5);
	\draw[very thick] (-3.5,0) to (-3.5,4);
	\draw[very thick] (-4.5,1) to (-4.5,5);
	\draw (-2.5,5) to [out=270,in=180] (-.5,1.75) to [out=0,in=270] (1.5,5);
	\draw[very thick, red, directed=.3] (.5,4) to [out=270,in=0] (-.5,2.5) to [out=180,in=270] (-1.5,4);
	\draw[very thick, directed=.2, directed=.6, directed=.9] (3.5,4) to (-3.5,4);
	\draw[very thick, directed=.65] (2.5,5) to (1.5,5) to (.5,4);
	\draw[very thick, directed=.35] (-1.5,4) to (-2.5,5) to (-4.5,5);
	\node[red, opacity=1] at (-3,3.5) {\tiny$a$};
	\node[red, opacity=1] at (-4,4.5) {\tiny$b$};
\end{tikzpicture}
};
\endxy \\
\Phi \left(
\xy
(0,0)*{
\begin{tikzpicture} [scale=.4,fill opacity=0.2]
\draw[ultra thick, black!30!green, directed=.99] (0,0) to [out=90,in=270] (2,3);
\draw[ultra thick, blue, directed=.99] (0,3) to [out=270,in=90] (2,0);
\node[black!30!green, opacity=1] at (0,-.5) {\tiny$c$};
\node[blue, opacity=1] at (0,3.5) {\tiny$a$};
\node[opacity=1] at (2.75,1.5) {\tiny$[0,d,0]$};
\end{tikzpicture}
}
\endxy \;
\right)
=
\xy
(0,0)*{
\begin{tikzpicture} [scale=.4,fill opacity=0.2]
	\path[fill=red] (3.5,5) to (2,5) to [out=270,in=60] (1,2.5) to [out=300,in=90] (2,0) to (3.5,0);
	\path[fill=blue] (2,5) to [out=270,in=60] (1,2.5) to [out=120,in=270] (0,5);
	\path[fill=blue] (2,0) to [out=90,in=300] (1,2.5) to [out=240,in=90] (0,0);
	\path[fill=red] (-3.5,0) to (0,0) to [out=90,in=240] (1,2.5) to [out=120,in=270] (0,5) to (-3.5,5);
	\path[fill=red] (0,0) to [out=90,in=270] (2,5) to (1,6) to (-4.5,6) to (-4.5,1) to (-1,1);
	\path[fill=red] (2,0) to [out=90,in=270] (0,5) to (-1,4) to (-2.5,4) to (-2.5,-1) to (1,-1);
	\draw[very thick, directed=.1, directed=.4,directed=.8] (3.5,0) to (-3.5,0);
	\draw[very thick, directed=.55] (2,0) to (1,-1) to (-2.5,-1);
	\draw[very thick, directed=.55] (0,0) to (-1,1) to (-4.5,1);
	\draw[very thick] (-2.5,-1) to (-2.5,4);
	\draw[very thick] (-3.5,0) to (-3.5,5);
	\draw[very thick] (-4.5,1) to (-4.5,6);
	\draw[very thick] (3.5,0) to (3.5,5);
	\draw (-1,1) to [out=90,in=270] (1,6);
	\draw[very thick, red, directed=.7] (0,0) to [out=90,in=270] (2,5);		
	\draw[very thick, red, directed=.7] (2,0) to [out=90,in=270] (0,5);	
	\draw (1,-1) to [out=90,in=270] (-1,4);
	\draw[very thick, directed=.1,directed=.4,directed=.8] (3.5,5) to (-3.5,5);
	\draw[very thick, directed=.55] (2,5) to (1,6) to (-4.5,6);
	\draw[very thick, directed=.3] (0,5) to (-1,4) to (-2.5,4);
	\node[red, opacity=1] at (-2,3.5) {\tiny$a$};
	\node[red, opacity=1] at (-3,4.5) {\tiny$b$};
	\node[red, opacity=1] at (-4,5.5) {\tiny$c$};
\end{tikzpicture}
};
\endxy
\quad &, \quad
\Phi \left(
\xy
(0,0)*{
\begin{tikzpicture} [scale=.4,fill opacity=0.2]
\draw[ultra thick, black!30!green, directed=.99] (2,0) to [out=90,in=270] (0,3);
\draw[ultra thick, blue, directed=.99] (2,3) to [out=270,in=90] (0,0);
\node[black!30!green, opacity=1] at (2,-.5) {\tiny$c$};
\node[blue, opacity=1] at (2,3.5) {\tiny$a$};
\node[opacity=1] at (2.75,1.5) {\tiny$[0,d,0]$};
\end{tikzpicture}
}
\endxy \;
\right)
=
\xy
(0,0)*{
\begin{tikzpicture} [scale=.4,fill opacity=0.2]
	\path[fill=red] (3.5,5) to (2,5) to [out=270,in=60] (1,2.5) to [out=300,in=90] (2,0) to (3.5,0);
	\path[fill=blue] (2,5) to [out=270,in=60] (1,2.5) to [out=120,in=270] (0,5);
	\path[fill=blue] (2,0) to [out=90,in=300] (1,2.5) to [out=240,in=90] (0,0);
	\path[fill=red] (-3.5,0) to (0,0) to [out=90,in=240] (1,2.5) to [out=120,in=270] (0,5) to (-3.5,5);
	\path[fill=red] (0,0) to [out=90,in=270] (2,5) to (1,4) to (-2.5,4) to (-2.5,-1) to (-1,-1);
	\path[fill=red] (2,0) to [out=90,in=270] (0,5) to (-1,6) to (-4.5,6) to (-4.5,1) to (1,1);
	\draw[very thick, directed=.1,directed=.4,directed=.8] (3.5,0) to (-3.5,0);
	\draw[very thick, directed=.55] (2,0) to (1,1) to (-4.5,1);
	\draw[very thick, directed=.3] (0,0) to (-1,-1) to (-2.5,-1);
	\draw[very thick] (-2.5,-1) to (-2.5,4);
	\draw[very thick] (-3.5,0) to (-3.5,5);
	\draw[very thick] (-4.5,1) to (-4.5,6);
	\draw[very thick] (3.5,0) to (3.5,5);
	\draw (1,1) to [out=90,in=270] (-1,6);
	\draw[very thick, red, directed=.7] (0,0) to [out=90,in=270] (2,5);		
	\draw[very thick, red, directed=.7] (2,0) to [out=90,in=270] (0,5);	
	\draw (-1,-1) to [out=90,in=270] (1,4);
	\draw[very thick, directed=.1, directed=.4,directed=.8] (3.5,5) to (-3.5,5);
	\draw[very thick, directed=.55] (2,5) to (1,4) to (-2.5,4);
	\draw[very thick, directed=.55] (0,5) to (-1,6) to (-4.5,6);	
	\node[red, opacity=1] at (-2,3.5) {\tiny$a$};
	\node[red, opacity=1] at (-3,4.5) {\tiny$b$};
	\node[red, opacity=1] at (-4,5.5) {\tiny$c$};
\end{tikzpicture}
};
\endxy
\end{align*}
where $d=a+b+c$ in the last equation. The result then follows from Lemma \ref{easyisotopy} since we can 
apply web isotopies in the image of $\Phi$ to obtain a web to which one of the above generators 
can be applied.
\end{proof}

We now aim to show a similar ``eventual faithfulness" result. This will aid in later showing that 
the foam $2$-categories are sufficiently non-degenerate to construct link homology. The main idea is to observe that one can realize several foam isotopies 
as the images under $\Phi_n$ of relations holding in $\Ucatc_Q(\glm)^{\geq 0}$ and then use these to deduce that all foam isotopies arise as consequences of 
relations in $\Ucatc_Q(\glm)^{\geq 0}$. It is fairly easy to check that elementary isotopies between ladder foams are induced by relations in 
$\Ucatc_Q(\glm)^{\geq 0}$ \emph{up to sign}, so the main difficulty lies in checking that the signs coming from relations in the categorified 
quantum group compensate with the rescalings on caps and cups in Theorem \ref{thm:2functor}.

Before trying to prove as many relations as possible, it is beneficial to consider symmetries of $\Ucatc_Q(\slm)$ (or $\Ucatc_Q(\glm)$) 
which correspond to foam symmetries. In addition to the duality given by biadjointness and $Q$-cyclicity, 
the following three symmetries of $\Ucatc_Q(\slm)$ are described in \cite{KL3}:
\begin{enumerate}
\item reflect across a horizontal axis and reverse orientation. \label{sym1}
\item rescale, reverse orientation, and send $\lambda\mapsto -\lambda$, \label{sym2}
\item rescale, reflect across a vertical axis, and send $\lambda\mapsto -\lambda$, and \label{sym3}
\end{enumerate}
Our next result describes the first symmetry, adapted from the non-signed version of $\Ucatc_Q(\slm)$ to the signed setting, 
where an additional rescaling is required.
\begin{prop} \label{prop:sym1}
The signed version of symmetry \eqref{sym1}, given by reflecting across a horizontal axis, reversing orientation, and rescaling by:
\begin{align*}
\xy
(0,0)*{
\begin{tikzpicture}[scale=1]
\draw[ultra thick, black!30!green, ->] (0,1) .. controls (0,.2) and (1,.2) .. (1,1);
\node at (1.3,.5) {\tiny $[a,b]$};
\end{tikzpicture}};
\endxy
\quad \mapsto \quad
(-1)^{b+1} \quad
\xy
(0,0)*{
\begin{tikzpicture}[scale=1]
\draw[ultra thick, black!30!green, ->] (1,0) .. controls (1,.8) and (0,.8) .. (0,0);
\node at (1.3,.5) {\tiny $[a,b]$};
\end{tikzpicture}};
\endxy
\quad&,\quad
\xy
(0,0)*{
\begin{tikzpicture}[scale=1]
\draw[ultra thick, black!30!green, ->] (1,0) .. controls (1,.8) and (0,.8) .. (0,0);
\node at (1.3,.5) {\tiny $[a,b]$};
\end{tikzpicture}};
\endxy
\quad \mapsto \quad
(-1)^{b} \quad
\xy
(0,0)*{
\begin{tikzpicture}[scale=1]
\draw[ultra thick, black!30!green, ->] (0,1) .. controls (0,.2) and (1,.2) .. (1,1);
\node at (1.3,.5) {\tiny $[a,b]$};
\end{tikzpicture}};
\endxy
\nonumber \\
\xy
(0,0)*{
\begin{tikzpicture}[scale=1]
\draw[ultra thick, black!30!green, <-] (0,1) .. controls (0,.2) and (1,.2) .. (1,1);
\node at (1.3,.5) {\tiny $[a,b]$};
\end{tikzpicture}};
\endxy
\quad \mapsto \quad
(-1)^{b} \quad
\xy
(0,0)*{
\begin{tikzpicture}[scale=1]
\draw[ultra thick, black!30!green, <-] (1,0) .. controls (1,.8) and (0,.8) .. (0,0);
\node at (1.3,.5) {\tiny $[a,b]$};
\end{tikzpicture}};
\endxy
\quad&,\quad
\xy
(0,0)*{
\begin{tikzpicture}[scale=1]
\draw[ultra thick, black!30!green, <-] (1,0) .. controls (1,.8) and (0,.8) .. (0,0);
\node at (1.3,.5) {\tiny $[a,b]$};
\end{tikzpicture}};
\endxy
\quad \mapsto \quad
(-1)^{b+1} \quad
\xy
(0,0)*{
\begin{tikzpicture}[scale=1]
\draw[ultra thick, black!30!green, <-] (0,1) .. controls (0,.2) and (1,.2) .. (1,1);
\node at (1.3,.5) {\tiny $[a,b]$};
\end{tikzpicture}};
\endxy
\\
\xy
(0,0)*{
\begin{tikzpicture}[scale=.8]
\draw [ultra thick, black!30!green, ->] (0,1) to [out=270,in=90] (1,0);
\draw [ultra thick, blue,->] (0,0) to [out=90,in=270] (1,1);
\end{tikzpicture}
};
\endxy
\quad \mapsto \quad (-1) \quad
\xy
(0,0)*{
\begin{tikzpicture}[scale=.8]
\draw [ultra thick, black!30!green, <-] (0,0) to [out=90,in=270] (1,1);
\draw [ultra thick, blue, <-] (0,1) to [out=270,in=90] (1,0);
\end{tikzpicture}
};
\endxy
\quad&, \quad
\xy
(0,0)*{
\begin{tikzpicture}[scale=.8]
\draw [ultra thick, black!30!green, <-] (0,0) to [out=90,in=270] (1,1);
\draw [ultra thick, blue, <-] (0,1) to [out=270,in=90] (1,0);
\end{tikzpicture}
};
\endxy
\quad \mapsto \quad (-1) \quad
\xy
(0,0)*{
\begin{tikzpicture}[scale=.8]
\draw [ultra thick, black!30!green, ->] (0,1) to [out=270,in=90] (1,0);
\draw [ultra thick, blue,->] (0,0) to [out=90,in=270] (1,1);
\end{tikzpicture}
};
\endxy
\nonumber \\
\xy
(0,0)*{
\begin{tikzpicture}[scale=.8]
\draw [ultra thick, black!30!green, ->] (0,1) to [out=270,in=90] (1,0);
\draw [ultra thick, blue,<-] (0,0) to [out=90,in=270] (1,1);
\end{tikzpicture}
};
\endxy
\quad \mapsto \quad (-1) \quad
\xy
(0,0)*{
\begin{tikzpicture}[scale=.8]
\draw [ultra thick, blue, ->] (0,1) to [out=270,in=90] (1,0);
\draw [ultra thick, black!30!green,<-] (0,0) to [out=90,in=270] (1,1);
\end{tikzpicture}
};
\endxy
\quad&, \quad
\xy
(0,0)*{
\begin{tikzpicture}[scale=.8]
\draw [ultra thick, blue, ->] (0,1) to [out=270,in=90] (1,0);
\draw [ultra thick, black!30!green,<-] (0,0) to [out=90,in=270] (1,1);
\end{tikzpicture}
};
\endxy
\quad \mapsto \quad (-1) \quad
\xy
(0,0)*{
\begin{tikzpicture}[scale=.8]
\draw [ultra thick, black!30!green, ->] (0,1) to [out=270,in=90] (1,0);
\draw [ultra thick, blue,<-] (0,0) to [out=90,in=270] (1,1);
\end{tikzpicture}
};
\endxy
\end{align*}
is an automorphism of $\Ucatc_Q(\glm)$. 
Under the foamation $2$-functors, this corresponds to reflecting a ladder foam across a plane parallel to the boundary webs and reversing seam orientation.
\end{prop}

The change from $\lambda$ to $-\lambda$ in symmetries \eqref{sym2} and \eqref{sym3} 
corresponds to replacing $a=[a_1,\dots,a_m]$ by $[n-a_1,\dots,n-a_m]$ 
(which changes the value of $N$ to $nm-N$) and encodes the duality between $\bigwedge_q^k(\C^n)$ and 
$\bigwedge_q^{n-k}(\C^n)$. The following composition of these two symmetries corresponds to a horizontal foam symmetry. 

\begin{prop} \label{prop:sym2}
The map given by reflecting across a vertical axis, reversing orientation, and rescaling by:
\begin{align*}
\xy
(0,0)*{
\begin{tikzpicture}[scale=.8]
\draw [ultra thick, black!30!green, ->] (0,0) to [out=90,in=270] (1,1);
\draw [ultra thick, blue,->] (1,0) to [out=90,in=270] (0,1);
\end{tikzpicture}
};
\endxy
\quad \mapsto \quad (-1) \quad
\xy
(0,0)*{
\begin{tikzpicture}[scale=.8]
\draw [ultra thick, black!30!green, ->] (0,1) to [out=270,in=90] (1,0);
\draw [ultra thick, blue,<-] (0,0) to [out=90,in=270] (1,1);
\end{tikzpicture}
};
\endxy
\quad&, \quad
\xy
(0,0)*{
\begin{tikzpicture}[scale=.8]
\draw [ultra thick, black!30!green, ->] (0,1) to [out=270,in=90] (1,0);
\draw [ultra thick, blue,<-] (0,0) to [out=90,in=270] (1,1);
\end{tikzpicture}
};
\endxy
\quad \mapsto \quad (-1) \quad
\xy
(0,0)*{
\begin{tikzpicture}[scale=.8]
\draw [ultra thick, black!30!green, ->] (0,0) to [out=90,in=270] (1,1);
\draw [ultra thick, blue,->] (1,0) to [out=90,in=270] (0,1);
\end{tikzpicture}
};
\endxy \\
\xy
(0,0)*{
\begin{tikzpicture}[scale=.8]
\draw [ultra thick, black!30!green, ->] (0,0) to [out=90,in=270] (1,1);
\draw [ultra thick, blue, ->] (0,1) to [out=270,in=90] (1,0);
\end{tikzpicture}
};
\endxy
\quad \mapsto \quad (-1) \quad
\xy
(0,0)*{
\begin{tikzpicture}[scale=.8]
\draw [ultra thick, black!30!green, ->] (0,1) to [out=270,in=90] (1,0);
\draw [ultra thick, blue,->] (0,0) to [out=90,in=270] (1,1);
\end{tikzpicture}
};
\endxy
\quad&, \quad
\xy
(0,0)*{
\begin{tikzpicture}[scale=.8]
\draw [ultra thick, blue , ->] (0,0) to [out=90,in=270] (1,1);
\draw [ultra thick, black!30!green, ->] (0,1) to [out=270,in=90] (1,0);
\end{tikzpicture}
};
\endxy
\quad \mapsto \quad (-1) \quad
\xy
(0,0)*{
\begin{tikzpicture}[scale=.8]
\draw [ultra thick, blue, ->] (0,1) to [out=270,in=90] (1,0);
\draw [ultra thick, black!30!green,->] (0,0) to [out=90,in=270] (1,1);
\end{tikzpicture}
};
\endxy
\nonumber
\end{align*}
is an automorphism of $\Ucatc_Q(\glm)$. Under the foamation functor, this corresponds to reflecting a ladder foam across a plane
parallel to the left and right boundary planes and reversing seams.
\end{prop}

\begin{rem}
We have given ways to realize symmetries along two elementary planes (horizontal and vertical, orthogonal to the sheets). One may wish to have also a 
symmetry across a vertical plane, parallel to the sheets. However, it is harder to make precise sense of the corresponding symmetry in the quantum group. 
One reason for this is that the symmetry reverses the local ordering of sheets along a vertical seam from back-to-front (which is our usual convention) to front-to-back 
(which \emph{a priori} doesn't lie in the image of the $2$-functor). 
As we do not need this extra symmetry for our considerations, we omit it, rather than making all necessary rescalings and adjustments for it to work. 
It is interesting to notice that this symmetry would map $\glm$ weights $[a_1,\ldots,a_m] \mapsto [a_m,\ldots,a_1]$, 
which does not correspond to one of Khovanov-Lauda's symmetries on $\Ucat_Q(\slm)$.
\end{rem}

\begin{lem} \label{lem:elemIso}
The following equalities in $\Ucatc_Q(\glm)^{\geq 0}$ induce isotopies between ladder foams. 
\newline
\noindent \textbf{Rung bend:} For $\sigma = \frac{1}{2}[a(a-1)+k(k-1)+(b-k)(b-k-1)]+ak$, 
\begin{equation} \label{RungBend1}
  \xy
 (0,15)*{
};
{\ar@{=>} (0,-7)*{}; (0,7)*{}}; ?(.5)*\dir{}+(3,0)*{\scriptstyle{id}}; 
 \endxy
\end{equation}
and the isotopy whose domain and codomain are given by the middle terms on the 
left-hand side holds as well, with the same sign. \newline

\noindent \textbf{Point-singularity turn:}
For $\rho=a(k+l)+\frac{(a)(a-1)}{2}$,
\begin{equation} \label{Isotopy_6ptTurn1}
  \xy
 (0,27)*{
 %
};
{\ar@{=>} (0,-7)*{}; (0,7)*{}}; ?(.5)*\dir{}+(3,0)*{\scriptstyle{id}}; 
 \endxy
\end{equation}
and a similar relation holds with a merger instead of a splitter. Similarly, for
$\rho'=\frac{(a)(a-1)}{2}$,
\begin{equation} \label{Isotopy_6ptTurn2}
  \xy
 (0,27)*{
 %
};
{\ar@{=>} (0,-19.5)*{}; (0,-16.5)*{}}; 
{\ar@{=>} (0,-1.5)*{}; (0,1.5)*{}}; 
{\ar@{=>} (0,16.5)*{}; (0,19.5)*{}}; 
 \endxy
\end{equation}
and similar relations hold for cups and the other cap.

\noindent \textbf{Cap split:}
For $\zeta= (-1)^{\frac{a(a-1)}{2}+\frac{k(k-1)}{2}}$ and $\zeta'=(-1)^{\frac{k(k-1)}{2}+k(a+1)}$,

\begin{equation} \label{CapSplit}
  \xy
 (0,27)*{
 %
};
{\ar@{=>} (0,-16.5)*{}; (0,-13.5)*{}}; 
{\ar@{=>} (0,-1.5)*{}; (0,1.5)*{}}; 
{\ar@{=>} (0,13.5)*{}; (0,16.5)*{}}; 
 \endxy
\end{equation}
and similarly for splitters. \newline

\noindent \textbf{Merger/splitter change:}
\begin{equation} \label{SplitCornerChange}
  \xy
 (0,15)*{
 %
};
{\ar@{=>} (0,-18)*{}; (0,-15)*{}}; 
{\ar@{=>} (0,-1.5)*{}; (0,1.5)*{}}; 
{\ar@{=>} (0,15)*{}; (0,18)*{}}; 
 \endxy
\end{equation}
and similarly for mergers. 

\noindent \textbf{Cap to Merge}\footnote{This formula is simply a rewriting of Relation \eqref{KLMSHigherBubbleRed} extracted from \cite{KLMS}.
}
For $\xi= (-1)^{\frac{a(a-1)}{2}+\frac{k(k-1)}{2}}$,

\begin{equation} \label{CapToMerge}
  \xy
 (0,36)*{
};
{\ar@{=>} (0,-7)*{}; (0,7)*{}}; ?(.5)*\dir{}+(3,0)*{\scriptstyle{id}}; 
 \endxy
\end{equation}

\end{lem}

\begin{proof}
The proof follows from repeated use of the relations in $\Ucatc_Q(\glm)$ and of Proposition \ref{prop:relBoundedQuotient}. 
We'll prove Relation \eqref{Isotopy_6ptbend} and leave the rest to the reader. We compute:
\[
\xy
(0,0)*{
%
};
\endxy
\]
The encircled blue cycle can be evaluated, and we obtain:
\[
(-1)^{\frac{(a-k)(a-k-1)}{2}}\;
\xy
(0,0)*{
%
};
\endxy
\]
In addition to simplifying the highlighted area, we also perform a ``Reidemeister 3''-like relation to pull off the blue strand, and obtain:
\[
=\quad (-1)^{\frac{(a-k)(a-k-1)}{2}}\sum_{\alpha\in P(k,b-l)}(-1)^{k(b-l) - |\alpha|}\;
\xy
(0,0)*{
%
};
\endxy
\]
The highlighted diagram can then be reduced as follows:
\[
\xy
(0,0)*{
%
};
\endxy\]
where the last equality holds true only if $\alpha=(0,\dots,0)$ (so that $\hat{\alpha}=(k,\dots,k)$), and otherwise it is zero. Our diagram becomes:
\[
(-1)^{\frac{(a-k)(a-k-1)}{2}}\;
\xy
(0,0)*{
%
};
\endxy
\]
and the latter equals $(-1)^{\frac{(a-k)(a-k-1)}{2}+\frac{(b+k-l)(b+k-l-1)}{2}}$ times the following:
\[
\sum_{\alpha\in P(l,c)} (-1)^{lc-|\alpha|}\;
\xy
(0,0)*{
%
};
\endxy
\]
The red part can be reduced as follows:
\[
\xy
(0,0)*{
%
};
\endxy
\;=
(-1)^{lc+\frac{(c+l)(c+l-1)}{2}}.
\]
where again the middle equality holds only if $\alpha=0$, and otherwise this term is zero. 
The last equality comes from the fact that this bubble is degree zero.
Combining all steps, we obtain that the term we started from equals
\[
(-1)^{\frac{(a-k)(a-k-1)}{2}+\frac{(b+k-l)(b+k-l-1)}{2}+\frac{(c+l)(c+l-1)}{2}}\;
\xy
(0,0)*{
%
};
\endxy
\]
\end{proof}

Using the duality given by cap/cup biadjointness and $Q$-cyclicity, together with the symmetries described in 
Propositions \ref{prop:sym1} and \ref{prop:sym2}, we can deduce analogs of the above relations for all possible symmetries of the webs. 

\begin{rem}\label{rem:nIsotopies}
In most of the above isotopies, and in particular in equations \eqref{RungBend1}, \eqref{RungBend2}, \eqref{Isotopy_6ptbend}, \eqref{PushCap}, and 
\eqref{PushSplitter}, one could replace the $0$-labeled `edges' by $n$-labeled ones, and adapt the orientation and labeling 
of the web parts accordingly. Up to sign adjustments, the relations still hold, but this time in $\Ucatc_Q(\glm)^{(n)}$, and the proofs are very similar 
(because $[n-k,n]$ and $[0,k]$ correspond to the same $\sln$ weight, and the $a_i\geq 0$ condition in the above proofs is replaced by the $n$-bounded condition). 
This shows that $n$-labeled sheets only encode signs, and don't prevent isotopies of the underlying non-$n$-labeled foam. 
This also shows that we can always `pull in' an $n$-labeled sheet from the front of the foam to any place inside it. This is a 
categorified version of a process we use in defining the link invariant in Section \ref{sec:link_hom}. We'll utilize these facts later in comparing 
our construction to that given in \cite{MSV}.
\end{rem}

Having shown that the above local foam isotopies are consequences of relations in $\Ucatc_Q(\glm)^{\geq 0}$, we can now prove 
the following result, which shows that, up to increasing the value of $m$, all isotopies between ladder foams are induced by 
equalities in $\Ucatc_Q(\glm)^{\geq 0}$.

\begin{prop} \label{prop:foam_isotopy}
Let $F_1 = \Phi(\cal{D}_1)$ and $F_2 = \Phi(\cal{D}_2)$ be two isotopic ladder foams in $\foam{n}{m}$ mapping between (the same) ladder webs, 
then there exists $m' > m$ so that $\iota(\cal{D}_1) = \iota(\cal{D}_2)$, where $\iota$ is the $2$-functor
$\Ucatc_Q(\glm)^{\geq 0} \stackrel{\iota}{\longrightarrow} \Ucatc_Q(\glnn{m'})^{\geq 0}$ given on objects by $[a_1,\ldots,a_m] \mapsto [a_1,\ldots,a_m,0,\ldots,0]$.
\end{prop}

\begin{proof}
The result follows by observing that, using Lemmata \ref{easyisotopy} and \ref{lem:elemIso}, all foam singularities can be pushed or pulled between the 
vertical planes giving a ladder foam its ladder structure, denoted for the duration by `slices.'
By applying these isotopies, which are the images under $\Phi$ of equalities in $\Ucatc_Q(\glm)^{\geq 0}$, 
we'll construct ladder foams $F'_1$ and $F'_2$ so that $F_1 \sim F'_1$, $F_2 \sim F'_2$, and so that 
the preimages of $F'_1$ and $F'_2$ are obviously equal.

We'll begin with the case where the ladder foams are undecorated, and then extend the proof to the decorated case. 
Throughout, we'll implicitly apply the foam isotopies given in Lemma \ref{easyisotopy} and assume that $m'$ is sufficiently large, allowing us 
enough room (i.e. enough `$0$-labeled' edges/facets) to perform the relevant ladder foam isotopies.

First, note that we can apply isotopies in the image of $\Phi$ to ensure that all point singularities, i.e. those lying in the intersection of two singular seams and locally modeled by \eqref{MVgensup}, 
are given by the images of splitter/merger morphisms in $\Ucatc_Q(\glm)^{\geq 0}$. If a point singularity is given by a `mixed crossing' of $i$ and $i+1$ colored (thick) strands, 
we can apply equation\footnote{In this proof, when we reference an equation
giving a foam isotopy as a consequence of an equality in $\Ucatc_Q(\glm)^{\geq 0}$, we implicitly refer to all relations obtained from the indicated 
one by symmetry and duality.} \eqref{Isotopy_6ptbend} and the isotopy
\begin{equation} \label{Isotopy_Mixed6pt}
  \xy
 (0,36)*{
 \begin{tikzpicture} [scale=.35]
\draw [dotted] (9,0) -- (0,0);
\draw [very thick, ->] (9,1) -- (0,1);
\draw [very thick, ->] (9,2) -- (0,2);
\draw [very thick, ->] (9,3) -- (0,3);
\draw [dotted] (9,4) -- (0,4);
\draw [very thick, directed= .6] (5,1) -- (5,2);
\draw [very thick, directed= .6] (4,2) -- (4,3);
\node [opacity=1] at (9.3,0) {\tiny $0$};
\node [opacity=1] at (9.3,1) {\tiny $a$};
\node [opacity=1] at (9.3,2) {\tiny $b$};
\node [opacity=1] at (9.3,3) {\tiny $c$};
\node [opacity=1] at (9.3,4) {\tiny $0$};
\node [opacity=1] at (3.5,2.5) {\tiny $l$};
\node [opacity=1] at (5.5,1.5) {\tiny $k$};
\end{tikzpicture}};
 (0,13.5)*{
 \begin{tikzpicture} [scale=.35]
\draw [dotted] (9,0) -- (8,0);
\draw [very thick] (8,0) -- (1,0);
\draw [dotted] (1,0) -- (0,0);
\draw [very thick] (9,1) -- (8,1);
\draw [dotted] (8,1) -- (6,1);
\draw [very thick] (6,1) -- (5,1);
\draw [dotted] (5,1) -- (1,1);
\draw [very thick, ->] (1,1) -- (0,1);
\draw [very thick, ->] (9,2) -- (0,2);
\draw [very thick] (9,3) -- (7,3);
\draw [dotted] (7,3) -- (4,3);
\draw [very thick] (4,3) -- (3,3);
\draw [dotted] (3,3) -- (2,3);
\draw [very thick, ->] (2,3) -- (0,3);
\draw [dotted] (9,4) -- (7,4);
\draw [very thick] (7,4) -- (2,4);
\draw [dotted] (2,4) -- (0,4);
\draw [very thick, directed=.6] (1,0) -- (1,1);
\draw [very thick, rdirected=.6] (2,3) -- (2,4);
\draw [very thick, directed=.6] (3,3) -- (3,4);
\draw [very thick, directed= .6] (4,2) -- (4,3);
\draw [very thick, directed= .6] (5,1) -- (5,2);
\draw [very thick, directed=.6] (6,0) -- (6,1);
\draw [very thick, directed=.6] (7,3) -- (7,4);
\draw [very thick, rdirected=.6] (8,0) -- (8,1);
\node [opacity=1] at (9.3,0) {\tiny $0$};
\node [opacity=1] at (9.3,1) {\tiny $a$};
\node [opacity=1] at (9.3,2) {\tiny $b$};
\node [opacity=1] at (9.3,3) {\tiny $c$};
\node [opacity=1] at (9.3,4) {\tiny $0$};
\node [opacity=1] at (3.5,2.5) {\tiny $l$};
\node [opacity=1] at (5.5,1.5) {\tiny $k$};
\end{tikzpicture}};
 (0,-13.5)*{
 \begin{tikzpicture} [scale=.35]
\draw [dotted] (9,0) -- (8,0);
\draw [very thick] (8,0) -- (1,0);
\draw [dotted] (1,0) -- (0,0);
\draw [very thick] (9,1) -- (8,1);
\draw [dotted] (8,1) -- (6,1);
\draw [very thick] (6,1) -- (4,1);
\draw [dotted] (4,1) -- (1,1);
\draw [very thick, ->] (1,1) -- (0,1);
\draw [very thick, ->] (9,2) -- (0,2);
\draw [very thick] (9,3) -- (7,3);
\draw [dotted] (7,3) -- (5,3);
\draw [very thick] (5,3) -- (3,3);
\draw [dotted] (3,3) -- (2,3);
\draw [very thick, ->] (2,3) -- (0,3);
\draw [dotted] (9,4) -- (7,4);
\draw [very thick] (7,4) -- (2,4);
\draw [dotted] (2,4) -- (0,4);
\draw [very thick, directed=.6] (1,0) -- (1,1);
\draw [very thick, rdirected=.6] (2,3) -- (2,4);
\draw [very thick, directed=.6] (3,3) -- (3,4);
\draw [very thick, directed= .6] (5,2) -- (5,3);
\draw [very thick, directed= .6] (4,1) -- (4,2);
\draw [very thick, directed=.6] (6,0) -- (6,1);
\draw [very thick, directed=.6] (7,3) -- (7,4);
\draw [very thick, rdirected=.6] (8,0) -- (8,1);
\node [opacity=1] at (9.3,0) {\tiny $0$};
\node [opacity=1] at (9.3,1) {\tiny $a$};
\node [opacity=1] at (9.3,2) {\tiny $b$};
\node [opacity=1] at (9.3,3) {\tiny $c$};
\node [opacity=1] at (9.3,4) {\tiny $0$};
\node [opacity=1] at (3.5,1.5) {\tiny $k$};
\node [opacity=1] at (5.5,2.5) {\tiny $l$};
\end{tikzpicture}};
 (0,-36)*{
 \begin{tikzpicture} [scale=.35]
\draw [dotted] (9,0) -- (0,0);
\draw [very thick, ->] (9,1) -- (0,1);
\draw [very thick, ->] (9,2) -- (0,2);
\draw [very thick, ->] (9,3) -- (0,3);
\draw [dotted] (9,4) -- (0,4);
\draw [very thick, directed= .6] (4,1) -- (4,2);
\draw [very thick, directed= .6] (5,2) -- (5,3);
\node [opacity=1] at (9.3,0) {\tiny $0$};
\node [opacity=1] at (9.3,1) {\tiny $a$};
\node [opacity=1] at (9.3,2) {\tiny $b$};
\node [opacity=1] at (9.3,3) {\tiny $c$};
\node [opacity=1] at (9.3,4) {\tiny $0$};
\node [opacity=1] at (5.5,2.5) {\tiny $l$};
\node [opacity=1] at (3.5,1.5) {\tiny $k$};
\end{tikzpicture}};
{\ar@{=>} (0,-27)*{}; (0,-23)*{}}; 
{\ar@{=>} (0,-2)*{}; (0,2)*{}}; 
{\ar@{=>} (0,23)*{}; (0,27)*{}}; 
 \endxy
\; \longmapsfrom \;
(-1)^\vartheta
\;\;
\xy
(0,0)*{
\begin{tikzpicture}[scale=.6]
\draw [ultra thick, black!30!green, ->] (1,-4) .. controls (1,0) and (-1,0) ..  (-1,4);
\draw [ultra thick, blue, ->] (-1,-4) .. controls (-1,0) and (1,0) .. (1,4);
\draw [ultra thick, red, directed=.55] (-2,.5) -- (-2,-.5);
\draw [ultra thick, red, directed=.55] (-2,-.5) .. controls (-1.9,-.6) and (-1.7,-.7) .. (-1.7,-1) .. controls (-1.7,-1.5) and (-1.2,-1.5) .. (-1.2,-1) -- (-1.2,1) .. controls (-1.2,1.5) and (-1.7,1.5) .. (-1.7,1) .. controls (-1.7,.7) and (-1.9,.6) .. (-2,.5);
\draw [ultra thick, red, directed= .52] (-2,-.5) .. controls (-2.1,-.6) and (-2.3,-.7) .. (-2.3,-1) .. controls (-2.3,-2) and (2,-3.5) .. (2,-2.5) -- (2,2.5) .. controls (2,3.5) and (-2.3,2) .. (-2.3,1) .. controls (-2.3,.7) and (-2.1,.6) .. (-2,.5);
\draw [ultra thick, black, directed=.55] (2.5,.5) -- (2.5,-.5);
\draw [ultra thick, black, directed=.55] (2.5,-.5) .. controls (2.4,-.6) and (2.3,-.7) .. (2.3,-1) .. controls (2.3,-1.5) and (1.7,-1.5) .. (1.7,-1) -- (1.7,1) .. controls (1.7,1.5) and (2.2,1.5) .. (2.2,1) .. controls (2.2,.7) and (2.4,.6) .. (2.5,.5);
\draw [ultra thick, black, directed= .52] (2.5,-.5) .. controls (2.6,-.6) and (2.8,-.7) .. (2.8,-1) -- (2.8,-2.5) .. controls (2.8,-4) and (-2.5,-3) .. (-2.5,-2) -- (-2.5,2) .. controls (-2.5,3) and (2.8,4) .. (2.8,2.5) -- (2.8,1) .. controls (2.8,.7) and (2.6,.6) .. (2.5,.5);
\node [blue] at (1,4.3) {\tiny $k$};
\node [black!30!green] at (-1,4.3) {\tiny $l$};
\node at (2.8,0) {\tiny $a$};
\node [red] at (2.2,0) {\tiny $c$};
\node [red] at (-1.7,0) {\tiny $l$};
\node at (1.4,0) {\tiny $k$};
\node [rotate=90] at (-2.8,0) {\tiny $a-k$};
\node at (3,3.7) {\tiny $[0,c,b,a,0]$};
\end{tikzpicture}
};
\endxy
\; = \;
\xy
(0,0)*{
\begin{tikzpicture}[scale=.6]
\draw [ultra thick, black!30!green, ->] (1,-4) .. controls (1,0) and (-1,0) ..  (-1,4);
\draw [ultra thick, blue, ->] (-1,-4) .. controls (-1,0) and (1,0) .. (1,4);
\node [blue] at (1,4.3) {\tiny $k$};
\node [black!30!green] at (-1,4.3) {\tiny $l$};
\end{tikzpicture}
};
\endxy
\; \longmapsto \;
 \xy
 (0,19)*{
 \begin{tikzpicture} [scale=.35]
\draw [dotted] (9,0) -- (0,0);
\draw [very thick, ->] (9,1) -- (0,1);
\draw [very thick, ->] (9,2) -- (0,2);
\draw [very thick, ->] (9,3) -- (0,3);
\draw [dotted] (9,4) -- (0,4);
\draw [very thick, directed= .6] (5,1) -- (5,2);
\draw [very thick, directed= .6] (4,2) -- (4,3);
\node [opacity=1] at (9.3,0) {\tiny $0$};
\node [opacity=1] at (9.3,1) {\tiny $a$};
\node [opacity=1] at (9.3,2) {\tiny $b$};
\node [opacity=1] at (9.3,3) {\tiny $c$};
\node [opacity=1] at (9.3,4) {\tiny $0$};
\node [opacity=1] at (3.5,2.5) {\tiny $l$};
\node [opacity=1] at (5.5,1.5) {\tiny $k$};
\end{tikzpicture}};
 (0,-19)*{
 \begin{tikzpicture} [scale=.35]
\draw [dotted] (9,0) -- (0,0);
\draw [very thick, ->] (9,1) -- (0,1);
\draw [very thick, ->] (9,2) -- (0,2);
\draw [very thick, ->] (9,3) -- (0,3);
\draw [dotted] (9,4) -- (0,4);
\draw [very thick, directed= .6] (4,1) -- (4,2);
\draw [very thick, directed= .6] (5,2) -- (5,3);
\node [opacity=1] at (9.3,0) {\tiny $0$};
\node [opacity=1] at (9.3,1) {\tiny $a$};
\node [opacity=1] at (9.3,2) {\tiny $b$};
\node [opacity=1] at (9.3,3) {\tiny $c$};
\node [opacity=1] at (9.3,4) {\tiny $0$};
\node [opacity=1] at (5.5,2.5) {\tiny $l$};
\node [opacity=1] at (3.5,1.5) {\tiny $k$};
\end{tikzpicture}};
{\ar@{=>} (0,-7)*{}; (0,7)*{}}?(.5)*\dir{}+(3,0)*{\scriptstyle{id}}; 
 \endxy
\end{equation}
where $\vartheta=\frac{(a-k)(a-k-1)+k(k-1)+c(c-1)+l(l-1)}{2}+k(a-k)+lc$, 
to ensure that locally the remainder of the foam doesn't occupy the slices adjacent to the point singularity. 
Equations \eqref{Isotopy_6ptTurn1} and \eqref{Isotopy_6ptTurn2} then show that the mixed crossing and splitter/merger 
realizations of point singularities are equivalent.
Moreover, the isotopy 
\begin{equation} \label{Isotopy_Split6pt}
  \xy
 (0,40)*{
 \begin{tikzpicture} [scale=.35]
\draw [very thick, ->] (7,0) -- (0,0);
\draw [very thick, ->] (7,1) -- (0,1);
\draw [dotted] (7,2) -- (0,2);
\draw [very thick, directed= .6] (4,0) -- (4,1);
\draw [very thick, directed= .6] (3,0) -- (3,1);
\node [opacity=1] at (7.3,0) {\tiny $a$};
\node [opacity=1] at (7.3,1) {\tiny $b$};
\node [opacity=1] at (7.3,2) {\tiny $0$};
\node [opacity=1] at (2.5,.5) {\tiny $l$};
\node [opacity=1] at (4.5,.5) {\tiny $k$};
\end{tikzpicture}};
 (0,20)*{
 \begin{tikzpicture} [scale=.35]
\draw [very thick, ->] (7,0) -- (0,0);
\draw [very thick] (7,1) -- (6,1);
\draw [dotted] (6,1) -- (5,1);
\draw [ultra thick] (5,1) -- (4,1);
\draw [dotted] (4,1) -- (3,1);
\draw [ultra thick] (3,1) -- (2,1);
\draw [dotted] (2,1) -- (1,1);
\draw [very thick, ->] (1,1) -- (0,1);
\draw [dotted] (7,2) -- (6,2);
\draw [very thick] (6,2) -- (1,2);
\draw [dotted] (1,2) -- (0,2);
\draw [very thick, directed= .6] (6,1) -- (6,2);
\draw [very thick, directed= .6] (5,0) -- (5,1);
\draw [very thick, directed= .6] (4,1) -- (4,2);
\draw [very thick, directed= .6] (3,0) -- (3,1);
\draw [very thick, directed= .6] (2,1) -- (2,2);
\draw [very thick, directed= .6] (1,2) -- (1,1);
\node [opacity=1] at (7.3,0) {\tiny $a$};
\node [opacity=1] at (7.3,1) {\tiny $b$};
\node [opacity=1] at (7.3,2) {\tiny $0$};
\node [opacity=1] at (2.5,.5) {\tiny $l$};
\node [opacity=1] at (5.5,.5) {\tiny $k$};
\end{tikzpicture}};
 (0,0)*{
 \begin{tikzpicture} [scale=.35]
\draw [very thick, ->] (7,0) -- (0,0);
\draw [very thick] (7,1) -- (6,1);
\draw [dotted] (6,1) -- (5,1);
\draw [ultra thick] (5,1) -- (2,1);
\draw [dotted] (2,1) -- (1,1);
\draw [very thick, ->] (1,1) -- (0,1);
\draw [dotted] (7,2) -- (6,2);
\draw [very thick] (6,2) -- (1,2);
\draw [dotted] (1,2) -- (0,2);
\draw [very thick, directed= .6] (6,1) -- (6,2);
\draw [very thick, directed= .6] (5,0) -- (5,1);
\draw [very thick, directed= .6] (4,0) -- (4,1);
\draw [very thick, directed= .6] (3,1) -- (3,2);
\draw [very thick, directed= .6] (2,1) -- (2,2);
\draw [very thick, directed= .6] (1,2) -- (1,1);
\node [opacity=1] at (7.3,0) {\tiny $a$};
\node [opacity=1] at (7.3,1) {\tiny $b$};
\node [opacity=1] at (7.3,2) {\tiny $0$};
\node [opacity=1] at (3.5,.5) {\tiny $l$};
\node [opacity=1] at (5.5,.5) {\tiny $k$};
\end{tikzpicture}};
 (0,-20)*{
 \begin{tikzpicture} [scale=.35]
\draw [very thick, ->] (7,0) -- (0,0);
\draw [very thick] (7,1) -- (6,1);
\draw [dotted] (6,1) -- (4.5,1);
\draw [ultra thick] (4.5,1) -- (2.5,1);
\draw [dotted] (2.5,1) -- (1,1);
\draw [very thick, ->] (1,1) -- (0,1);
\draw [dotted] (7,2) -- (6,2);
\draw [very thick] (6,2) -- (1,2);
\draw [dotted] (1,2) -- (0,2);
\draw [very thick, directed= .6] (6,1) -- (6,2);
\draw [very thick, directed= .6] (4.5,0) -- (4.5,1);
\draw [very thick, directed= .6] (2.5,1) -- (2.5,2);
\draw [very thick, directed= .6] (1,2) -- (1,1);
\node [opacity=1] at (7.3,0) {\tiny $a$};
\node [opacity=1] at (7.3,1) {\tiny $b$};
\node [opacity=1] at (7.3,2) {\tiny $0$};
\node [opacity=1] at (3.5,.5) {\tiny $k+l$};
\end{tikzpicture}};
 (0,-40)*{
 \begin{tikzpicture} [scale=.35]
\draw [very thick, ->] (7,0) -- (0,0);
\draw [very thick, ->] (7,1) -- (0,1);
\draw [dotted] (7,2) -- (0,2);
\draw [very thick, directed= .6] (3.5,0) -- (3.5,1);
\node [opacity=1] at (7.3,0) {\tiny $a$};
\node [opacity=1] at (7.3,1) {\tiny $b$};
\node [opacity=1] at (7.3,2) {\tiny $0$};
\node [opacity=1] at (2.5,.5) {\tiny $k+l$};
\end{tikzpicture}};
{\ar@{=>} (0,-33)*{}; (0,-27)*{}}; 
{\ar@{=>} (0,-13)*{}; (0,-7)*{}}; 
{\ar@{=>} (0,7)*{}; (0,13)*{}}; 
{\ar@{=>} (0,27)*{}; (0,33)*{}}; 
 \endxy
\; \longmapsfrom \;
(-1)^\zeta
\;\;
\xy
(0,0)*{
\begin{tikzpicture}[scale=.6]
\draw [ultra thick, blue] (0,-3) -- (0,0);
\draw [ultra thick, blue, ->] (0,0) .. controls (-.5,.1) and (-1,.2) .. (-1,3);
\draw [ultra thick, blue, ->] (0,0) .. controls (.5,.1) and (1,.2) .. (1,3);
\draw [ultra thick, black!30!green, directed=.55] (-2,.5) -- (-2,-.5);
\draw [ultra thick, black!30!green, directed=.8] (-2,-.5) .. controls (-1.9,-.6) and (-1.7,-.7) .. (-1.7,-1) .. controls (-1.7,-1.25) and (-1.3, -1.25) .. (-1.3,-1)  -- (-1.3,-.5);
\draw [ultra thick, black!30!green, directed=.55] (-1.3,-.5) .. controls (-1.4,-.4) and (-1.6,-.3) .. (-1.6,0) -- (-1.6,.75) .. controls (-1.6,1) and (-1.8,1) .. (-1.8,.75) .. controls (-1.8,.65) and (-1.9,.55) .. (-2,.5);
\draw [ultra thick, black!30!green, directed=.68] (-1.3,-.5) .. controls (-1.2,-.4) and (-1,-.3) .. (-1,0) .. controls (-1,.5) and (-.5,.25) .. (-.5,.75) .. controls (-.5,1.25) and (-1,1) .. (-1,1.5) .. controls (-1,2) and (-2,2) .. (-2,1) -- (-2,.5);
\draw [ultra thick, black!30!green, directed=.52] (-2,-.5) .. controls (-2.1,-.6) and (-2.3,-.7) .. (-2.3,-1) .. controls (-2.3,-2) and (1.5,-3) .. (1.5,-2.5) -- (1.5,2.5) .. controls (1.5,3) and (-2.3,2.5) .. (-2.3,1.5) -- (-2.3,1) .. controls (-2.3,.7) and (-2.1,.6) .. (-2,.5);
\node [blue] at (1,3.3) {\tiny $k$};
\node [blue] at (-1,3.3) {\tiny $l$};
\node [black!30!green] at (1.8,0) {\tiny $b$};
\node [black!30!green] at (-.2,.75) {\tiny $k$};
\node [black!30!green] at (-1.3,.25) {\tiny $l$};
\node at (2.5,2.7) {\tiny $[0,b,a]$};
\end{tikzpicture}
};
\endxy
\; = \;
\xy
(0,0)*{
\begin{tikzpicture}[scale=.6]
\draw [ultra thick, blue] (0,-3) -- (0,0);
\draw [ultra thick, blue, ->] (0,0) .. controls (-.5,.1) and (-1,.2) .. (-1,3);
\draw [ultra thick, blue, ->] (0,0) .. controls (.5,.1) and (1,.2) .. (1,3);
\node [blue] at (1,3.3) {\tiny $k$};
\node [blue] at (-1,3.3) {\tiny $l$};
\end{tikzpicture}
};
\endxy
\; \longmapsto \;
  \xy
 (0,19)*{
 \begin{tikzpicture} [scale=.35]
\draw [very thick, ->] (7,0) -- (0,0);
\draw [very thick, ->] (7,1) -- (0,1);
\draw [dotted] (7,2) -- (0,2);
\draw [very thick, directed= .6] (4,0) -- (4,1);
\draw [very thick, directed= .6] (3,0) -- (3,1);
\node [opacity=1] at (7.3,0) {\tiny $a$};
\node [opacity=1] at (7.3,1) {\tiny $b$};
\node [opacity=1] at (7.3,2) {\tiny $0$};
\node [opacity=1] at (2.5,.5) {\tiny $l$};
\node [opacity=1] at (4.5,.5) {\tiny $k$};
\end{tikzpicture}};
 (0,-19)*{
 \begin{tikzpicture} [scale=.35]
\draw [very thick, ->] (7,0) -- (0,0);
\draw [very thick, ->] (7,1) -- (0,1);
\draw [dotted] (7,2) -- (0,2);
\draw [very thick, directed= .6] (3.5,0) -- (3.5,1);
\node [opacity=1] at (7.3,0) {\tiny $a$};
\node [opacity=1] at (7.3,1) {\tiny $b$};
\node [opacity=1] at (7.3,2) {\tiny $0$};
\node [opacity=1] at (2.5,.5) {\tiny $k+l$};
\end{tikzpicture}};
{\ar@{=>} (0,-7)*{}; (0,7)*{}}?(.5)*\dir{}+(3,0)*{\scriptstyle{}}; 
 \endxy
\end{equation}
where $\zeta=\frac{(k+l)(k+l-1)+b(b-1)}{2}+b(k+l)$, allows us to assume that each point singularity is the image of 
a unique splitter/merger, and by equation \eqref{SplitCornerChange} we can assume that such a point singularity in the image 
of an $i$-colored merger/splitter lies on the $i^{th}$ (as opposed to $(i+1)^{th}$) slice, i.e. the slice closer to the back of the ladder foam.
We hence assume that all point singularities in $F_1$ and $F_2$ take the above form. 

We now observe that we can move point singularities, singular edges, and facets between slices. 
This follows from equation \eqref{PushSplitter} for point singularities, equations \eqref{PushSplitter}, 
\eqref{RungBend1}, \eqref{RungBend2}, and \eqref{PushCap} for singular seams, and the isotopies in Lemma \ref{easyisotopy} for facets. Additionally, 
equations \eqref{VertexRotate1}, \eqref{VertexRotate2}, and \eqref{CapSplit} allow seams to be moved independently, and equation \eqref{CapToMerge} 
can be used to pull a seam passing between slices onto one slice.

We now construct the ladder foams $F'_i$ for $i=1,2$. To begin, consider the `back sheet' of the foam $F_i$, i.e. our ladder foams lie in 
$[0,1]\times [-\infty,1] \times [0,1]$ and the back is the portion of the foam visible from the point $(1/2,2,1/2)$ 
(the slices are the planes $[0,1]\times \{2-i\} \times [0,1]$ for $i=1,\ldots,m'$). 
We now pull the entire back sheet of $F_i$, with the exception of a neighborhood of the (fixed) boundary, 
onto the first slice. Since $F_1$ and $F_2$ have identical boundaries, we can also ensure that the portions of the ladder foams' back sheets 
which don't lie on the first slice are identical. 

It then follows that the seams on the first slice differ only up to planar isotopy, so the portions of the foams lying in-between the first and second slice 
(which are determined by these seams) are isotopic. Using equations \eqref{VertexRotate1}, \eqref{VertexRotate2}, and \eqref{PushCap} we can 
then also ensure that the neighborhoods of the intersection of the in-between facets with the second slice are isotopic (outside a neighborhood of the foam boundaries). 
We next use these relations, together with \eqref{CapToMerge}, to pull the rest of each foam away from the second slice, 
guaranteeing that the second slice of each foam consists only of neighborhoods of the intersection of the facets in-between the first and second slices.
The only singular seams on the second slice are in small neighborhoods of where a singular seam enters the second slice from the first slice via a splitter/merger, 
before passing to the third slice. It again follows that the portions of the foams lying in-between the second and third slices are isotopic, and that neighborhoods of 
the intersections of these facets with the third slice are isotopic.

We next consider the complements of the portions of $F_1$ and $F_2$ we've considered so far (i.e. the portions lying on or after the third slice, 
and not in a neighborhood of the intersection of the third slice with the facets in-between the second and third slice). 
We can then perform the analog of the first step, pulling as much of the back of the remainder of each foam (outside a neighborhood of the foam boundary) 
onto the third slice as possible, without changing any part of the foams considered thus far. 
Note that if a portion of $F_1$ can be pulled onto the third slice, then the corresponding portion of $F_2$ can be pulled there as well, and that it's 
always possible to pull either a point singularity or (the remainder of) an entire seam or facet onto this slice. 

We iterate the above procedure, ensuring that the foams pass through even-numbered slices with isotopic neighborhoods, pushing the rest of the foam away, and then 
pulling as much of the remainder of the foam as possible onto the next odd-numbered slice. The procedure ensures that the in-between foam facets for $F_1$ and $F_2$ are 
always related by isotopy. This process eventually terminates; after dealing with each odd slice, there are fewer point singularities, singular seams, and/or 
facets in the remaining portion of the foam, since at each step we can at the very least pull a point singularity or (the remainder of) an entire seam or facet onto a slice. 
Since the ladder foams $F'_1$ and $F'_2$ remaining at the end of the procedure have isotopic portions on each slice and isotopic in-between facets, their preimage 
diagrams in $\Ucat_Q(\glm)$ are equal using biadjointness, $Q$-cyclicity, dot-sliding, and Reidemeister-like moves involving strands with distant colors (i.e. $|i-j| \geq 2$), 
which verifies the result.

We finally extend the above argument to the case of decorated foams. The main difficulty lies in the fact that the portions of a ladder foam lying on a slice cannot carry decorations. 
To remedy this, we'll formally view a decorated foam facet lying on a slice via the following shorthand:
\begin{equation}\label{slicedot}
\xy
(0,0)*{
\begin{tikzpicture} [scale=.5, fill opacity=0.2]
\path[fill=red] (-2,-2) rectangle (2,2);
\draw[very thick, directed=.55] (2,2) to (-2,2);
\draw[very thick, directed=.55] (2,-2) to (-2,-2);
\draw[very thick] (2,2) to (2,-2);
\draw[very thick] (-2,2) to (-2,-2);
\node[opacity=1] at (0,0) {$\bullet^{f}$};
\node[red,opacity=1] at (1.5,1.5) {$a$};
\end{tikzpicture}};
\endxy
\quad = \quad
\xy
(0,0)*{
\begin{tikzpicture} [scale=.5, fill opacity=0.2]
\draw (-.5,-1) to [out=90,in=180] (0,-.5) to [out=0,in=90] (.5,-1) to [out=270,in=0] (0,-1.5)
	to [out=180,in=270] (-.5,-1);
\draw (-1.25,0) to [out=90,in=180] (0,1.5) to [out=0,in=90] (1.25,0) to [out=270,in=0] (0,-2)
	to [out=180,in=270] (-1.25,0);
\path[fill=red] (-2.5,-3) rectangle (2.5,2.5);
\draw[very thick, directed=.55] (2.5,-3) to (-2.5,-3);
\draw[very thick] (2.5,2.5) to (2.5,-3);
\draw[very thick] (-2.5,2.5) to (-2.5,-3);
\draw[dashed] (2.5,0) to (1.25,0) to (.5,-1) to (-.5,-1) to (-1.25,0) to (-2.5,0);
\draw[very thick, directed=.55] (2.5,2.5) to (-2.5,2.5);
\node[opacity=1] at (.75,0) {$\bullet^{f}$};
\node[red,opacity=1] at (1.75,1.75) {$a$};
\end{tikzpicture}};
\endxy
\quad = \quad
\Phi \left(
(-1)^{\frac{a(a-1)}{2}}
\xy
(0,0)*{
\begin{tikzpicture} [scale=.5]
\draw [red, ultra thick, rdirected=.75] (1,0) to [out=90,in=0] (0,1) to [out=180,in=90] (-1,0) to [out=270,in=180] (0,-1) to [out=0,in=270] (1,0);
\node[opacity=1] at (1.15,0) {$\bullet^{f}$};
\node at (1.5,1.25) {\tiny $[a,0]$};
\end{tikzpicture}};
\endxy
\right) .
\end{equation}
We now repeat the above process for decorated foams, using equation \eqref{slicedot} when pulling a decorated facet onto a slice; since we only pull 
neighborhoods of the in-between facets onto even-numbered slices, there is always room for these decorations. The proof of the proposition in the decorated 
case is then the same as above, since the isotopies used are all compatible with decorations. The compatibility follows since the isotopy
\[
\xy
(0,0)*{
\begin{tikzpicture} [scale=.5, fill opacity=0.2]
\path[fill=red] (3,2) to (1.5,2) to (.5,3) to (-.5,3) to (-1.5,4) to (-3,4) to (-3,0) to (-1.5,0) to (-.5,-1) to (.5,-1) to
	(1.5,-2) to (3,-2);
\draw[very thick, directed=.35, directed=.7] (3,2) to (1.5,2) to (.5,3) to (-.5,3) to (-1.5,4) to (-3,4);
\draw[very thick] (3,-2) to (3,2);
\draw (1.5,-2) to (1.5,2);
\draw (.5,-1) to (.5,3);
\draw (-.5,-1) to (-.5,3);
\draw (-1.5,0) to (-1.5,4);
\draw[very thick] (-3,0) to (-3,4);
\draw[very thick, directed=.35, directed=.7] (3,-2) to (1.5,-2) to (.5,-1) to (-.5,-1) to (-1.5,0) to (-3,0);
\node[opacity=1] at (-1,1.5) {$\bullet^{f}$};
\end{tikzpicture}};
\endxy
\quad = \quad
\xy
(0,0)*{
\begin{tikzpicture} [scale=.5, fill opacity=0.2]
\path[fill=red] (3,2) to (1.5,2) to (.5,3) to (-.5,3) to (-1.5,4) to (-3,4) to (-3,0) to (-1.5,0) to (-.5,-1) to (.5,-1) to
	(1.5,-2) to (3,-2);
\draw[very thick, directed=.35, directed=.7] (3,2) to (1.5,2) to (.5,3) to (-.5,3) to (-1.5,4) to (-3,4);
\draw[very thick] (3,-2) to (3,2);
\draw (1.5,-2) to (1.5,2);
\draw (.5,-1) to (.5,3);
\draw (-.5,-1) to (-.5,3);
\draw (-1.5,0) to (-1.5,4);
\draw[very thick] (-3,0) to (-3,4);
\draw[very thick, directed=.35, directed=.7] (3,-2) to (1.5,-2) to (.5,-1) to (-.5,-1) to (-1.5,0) to (-3,0);
\node[opacity=1] at (1,.5) {$\bullet^{f}$};
\end{tikzpicture}};
\endxy
\]
is a consequence of the relation
$\xy
(0,0)*{
\begin{tikzpicture} [scale=.5]
\draw[ultra thick, black!30!green, ->] (-1,-1) to (-1,1);
\draw[ultra thick, blue, ->] (0,-1) to (0,1);
\node at (-.8,0) {$\bullet^{f}$};
\node[black!30!green] at (-1,-1.25) {\tiny $a$};
\node[blue] at (0,-1.25) {\tiny $a$};
\node at (1.25,.5) {\tiny$[0,0,a]$};
\end{tikzpicture}};
\endxy
= \;\;
\xy
(0,0)*{
\begin{tikzpicture} [scale=.5]
\draw[ultra thick, black!30!green, ->] (-1,-1) to (-1,1);
\draw[ultra thick, blue, ->] (0,-1) to (0,1);
\node at (.2,0) {$\bullet^{f}$};
\node[black!30!green] at (-1,-1.25) {\tiny $a$};
\node[blue] at (0,-1.25) {\tiny $a$};
\node at (1.5,.5) {\tiny$[0,0,a]$};
\end{tikzpicture}};
\endxy$
in $\Ucatc_Q(\glm)^{\geq 0}$, which in turn follows from equation \eqref{eq_r2_ij-gen-thick}.
\end{proof}

Propositions \ref{prop:fullness} and \ref{prop:foam_isotopy} indicate that the $2$-category $\foam{n}{m}$ is most naturally related to 
$\Ucatc_Q(\glm)^{(n)}$ in the limit $m \to \infty$. Since Lemma \ref{easyisotopy} implies that the objects $[a_1,\ldots,a_i,a_{i+1},0,\ldots,a_m]$ 
and $[a_1,\ldots,a_i,0,a_{i+1},\ldots,a_m]$ of $\foam{n}{m}$ are canonically isomorphic, we can identify them and hence consider an $\sln$ foam 
category in which the objects are sequences only of the elements $\{1,\ldots,n\}$, denoted by $\foam{n}{}$.

\begin{prop}\label{prop:direct_lim}
Let $_N\Ucatc_Q(\glm)^{(n)}$ denote the full $2$-subcategory of $\Ucatc_Q(\glm)^{(n)}$ given by objects \linebreak $[a_1,\ldots,a_m]$ with $\sum a_i = N$, then
the $2$-category $\foam{n}{}$ is equivalent to the direct limit 
\[_N\Ucatc_Q(\glnn{\infty})^{(n)} = \displaystyle \lim_{\longrightarrow} \left( \cdots \longrightarrow \textstyle _N\Ucatc_Q(\glm)^{(n)} \longrightarrow
\textstyle_N\Ucatc_Q(\glnn{m+1})^{(n)} \longrightarrow \cdots \right ).\]
\end{prop}
Objects in this direct limit $2$-category are infinite sequences $[a_1,a_2,\ldots]$ with $a_i \in \{0,\ldots,n\}$ so that $a_i=0$ for all but finitely many $i$ and $\sum a_i = N$, 
with $1$- and $2$-morphisms given as in $\Ucatc_Q(\slnn{\infty})$.
\begin{proof}
Proposition \ref{prop:foam_isotopy} indicates that isotopy relations in $\foam{n}{}$ follow as consequences of relations in $_N\Ucatc_Q(\glnn{\infty})^{(n)}$. 
The result then follows from Remark \ref{rem:rel}.
\end{proof}

\begin{rem}
Equation \eqref{digon_web_iso}, Relation \eqref{MVrel}, Lemma \ref{lemKLRc}, and Remark \ref{rem:SquareFormulas} imply that $\displaystyle \bigoplus_{N\geq 0} \foam{n}{}$ categorifies the category of left-directed (enhanced) webs. It follows from \cite[Theorem 5.3.1]{CKM} that (after truncating the $n$-labeled edges to tags) 
this gives the full subcategory of $\cat{Rep}(U_q(\mathfrak{sl}_n))$ generated (as a monoidal category) by the fundamental representations (i.e. without their duals).
\end{rem}

\begin{rem}
One can also consider the $2$-category $\foam{\infty}{}$ of $\slnn{\infty}$ foams, where we haven't imposed the condition that $n$-labeled 
facets are zero. The results of this section then imply that $\foam{\infty}{}$ is equivalent to $_N\Ucatc_Q(\glnn{\infty})^{\geq0}$, which 
is (a full $2$-subcategory of) the limit of categorified $q$-Schur algebras.
\end{rem}

%
\section{Link homology}\label{sec:link_hom}
%

In this section, we define our tangle (and as a special case link) invariant, and study its properties. 
Before doing so, we first check that the $2$-categories $\foam{n}{}$ are sufficiently non-degenerate. Specifically, 
we must show that the endomorphism ring of the identity morphism of the object $[n,\ldots,n]$ in $\foam[ns]{n}{}$ is non-zero. 
To this end, consider the $2$-representations of $\Ucat_Q(\slm)$ on the homotopy $2$-category of $\sln$ matrix factorizations 
$n\mathbf{hmf}$ given in \cite{MY}. These $2$-representations $\Gamma_n$ induce $2$-representations of $\Ucatc_Q(\glm)^{(n)}$ and 
satisfy the commutative diagram:
\[
\xymatrix{
\Ucatc_Q(\glm)^{(n)} \ar[rr] \ar[dr]^-{\Gamma_n} & & \Ucatc_Q(\glnn{m+1})^{(n)} \ar[dl]_-{\Gamma_n} \\
& n\mathbf{hmf} &}
\]
hence by Proposition \ref{prop:direct_lim} there exists an induced $2$-representation $\foam{n}{} \stackrel{\Psi}{\longrightarrow} n\mathbf{hmf}$. Indeed, 
one can explicitly specify\footnote{Note that the authors of that paper use the convention that $t_{i,i+1}=1$, so one 
must first apply an automorphism of $\Ucat_Q(\slm)$ to relate to our construction.} 
this $2$-representation using the definitions from \cite{MY} to assign matrix factorzations to $\sln$ webs and morphisms 
of matrix factorizations to generating foams. If $\1_{\mathbf{a}}$ denotes the identity web of the object $\mathbf{a}=[a_1,\ldots,a_k]$ in 
$\foam{n}{}$ then it follows from \cite[Corollary 7.1]{Wu} that $\dim \Hom_{n\mathbf{hmf}}(\Psi(\1_{\mathbf{a}}), \Psi(\1_{\mathbf{a}})) \neq 0$, 
which in turn implies that $\dim \Hom_{\foam{n}{}}(\1_{\mathbf{a}}, \1_{\mathbf{a}}) \neq 0$. As we'll see below, this suffices to define 
our link homology theory.

Equation $(14)$ from \cite{KhR} implies that under the $2$-representation $\Psi$, the image of an $n$-dotted 
$1$-labeled facet is zero. Since we use this $2$-representation to deduce the non-degeneracy
of $\foam{n}{}$, we'll impose this 
additional foam relation when constructing our link invariant, denoting the resulting foam $2$-category by $\foam{n}{}^\bullet$.
This relation holds for the (undeformed) constructions of $\slnn{2}$ and $\slnn{3}$ foam $2$-categories and a similar relation is 
satisfied for other existing constructions of $\sln$ link homology. 
Using equation \eqref{blisterrel}, this implies that an $a$-labeled sheet decorated by a Schur function 
$\pi_{\alpha}$ with any part larger than $n-a$ equals zero, i.e. such sheets carry an action of $H^*(Gr_a(\C^n))$.

\begin{rem}
Since the algebra of (non-degree preserving) endomorphisms of the matrix factorization assigned to 
a $a$-labeled edge has (graded) dimension $q^{a(n-a)} {n \brack a}$, 
we in fact see that 
\[\displaystyle \bigoplus_{k \in \Z} \Hom \left(
\xy
(0,0)*{
\begin{tikzpicture}
\node[rotate=-45] at (0,0){$\xy
(0,0)*{
\begin{tikzpicture} [scale=.2]
	\draw[very thick, directed=.55] (2,-2) -- (-2,-2);
\end{tikzpicture}
};
\endxy$};
\node at (.25,.25) {\small$a$};
\end{tikzpicture}
};
\endxy,
q^k
\xy
(0,0)*{
\begin{tikzpicture}
\node[rotate=-45] at (0,0){$\xy
(0,0)*{
\begin{tikzpicture} [scale=.2]
	\draw[very thick, directed=.55] (2,-2) -- (-2,-2);
\end{tikzpicture}
};
\endxy$};
\node at (.25,.25) {\small$a$};
\end{tikzpicture}
};
\endxy
\right) \cong H^*(Gr_a(\C^n))
\]
in $\foam{n}{}^\bullet$.
\end{rem}

Before proceeding to define our tangle invariant, we'll first discuss the subtleties of the identity $2$-morphisms of an $n$-labeled sheet. 
They play a crucial role in this theory, and give a higher representation theoretic explanation for the additional 
foam relation mentioned above.

\subsection{Some details concerning $n$-labeled facets}  \label{subsec:detailsNfacets} 

We next observe a representation-theoretic reason for considering the additional relation in $\foam{n}{}^\bullet$.
A priori, our link invariant will take values in the homotopy category of complexes over $\foam[ns]{n}{}$, for $s$ sufficiently large, 
in which we view the sequence $[n,\ldots,n]$ as corresponding to the `lowest weight' object $[0,\ldots,0,n,\ldots,n]$ in $\foam[ns]{n}{m}$.
In most $2$-representations of $\Ucat_Q(\slm)$, the (graded) vector space of (non-degree preserving) endomorphisms of the 
identity $1$-morphism of the highest weight object is $1$-dimensional in degree zero, and zero in other degrees. This is equivalent 
to setting any $n$-labeled facet with a non-identity decoration equal to zero.

In fact, imposing this relation on $n$-labeled facets implies the aforementioned relation (that $\bullet^n = 0$) for any $1$-labeled facet which 
interacts with an $n$-labeled facet. More precisely, when $a=n$, $b=0$, and $c=1$, equation \eqref{R2KLR1rel} implies that
\[
\sum_{d = 0}^n
(-1)^{d} \quad
\xy
(0,0)*{
\begin{tikzpicture} [scale=.6,fill opacity=0.2]
	\path[fill=yellow] (2.25,3) to (.75,3) to (.75,0) to (2.25,0);
	\path[fill=red] (.75,3) to [out=225,in=0] (-.5,2.5) to (-.5,-.5) to [out=0,in=225] (.75,0);
	\path[fill=red] (.75,3) to [out=135,in=0] (-1,3.5) to (-1,.5) to [out=0,in=135] (.75,0);	
	\draw [double] (2.25,0) to (.75,0);
	\draw [very thick,directed=.55] (.75,0) to [out=135,in=0] (-1,.5);
	\draw [very thick,directed=.55] (.75,0) to [out=225,in=0] (-.5,-.5);
	\draw[very thick, red, directed=.55] (.75,0) to (.75,3);
	\draw [very thick] (2.25,3) to (2.25,0);
	\draw [very thick] (-1,3.5) to (-1,.5);
	\draw [very thick] (-.5,2.5) to (-.5,-.5);
	\draw [double] (2.25,3) to (.75,3);
	\draw [very thick,directed=.55] (.75,3) to [out=135,in=0] (-1,3.5);
	\draw [very thick,directed=.55] (.75,3) to [out=225,in=0] (-.5,2.5);
	\node[opacity=1] at (-.375,3) {\small$\bullet^d$};
	\node[opacity=1] at (1.6,1.6) {\small$e_{n-d}$};
	\node[red, opacity=1] at (-.75,3.25) {\tiny{$1$}};
	\node[red, opacity=1] at (0,2.25) {\tiny{$_{n-1}$}};		
\end{tikzpicture}
};
\endxy
\quad = 0
\]
where here and for the duration we color $n$-labeled facets yellow.
The web isomorphism:
\[
\xy
(0,0)*{\begin{tikzpicture} [scale=.5]
 \draw [very thick, directed=.55] (3,1) -- (-1.75,1);
 \draw[double] (3,0) to (-1.75,0);
 \node at (3.5,0) {\small$n$};
 \node at (3.5,1) {\small$1$}; 
\end{tikzpicture}};
\endxy
\;\; \cong \;\;
\xy
(0,0)*{\begin{tikzpicture} [scale=.5]
 \draw [double] (3,0) -- (1.75,0);
 \draw[very thick, directed=.55] (3,1) to (1.25,1);
 \draw [double] (1.25,1) -- (0,1);
\draw [very thick, directed=.55] (1.75,0) -- (-.5,0);
\draw[very thick, directed=.55] (0,1) to (-1.75,1);
\draw[double] (-.5,0) to (-1.75,0);
 \draw [very thick, directed=.55] (1.75,0) -- (1.25,1);
 \draw[very thick, directed=.55] (0,1) to (-.5,0); 
 \node at (3.5,1) {\small$1$};
 \node at (3.5,0) {\small$n$}; 
\end{tikzpicture}};
\endxy
\]
from equation \eqref{FE_to_cb} gives the relation for $1$-labeled facets.

\begin{rem}
Since $n$-labeled edges correspond to the trivial representation of $\sln$, imposing the relation that their endomorphisms are 
$1$-dimensional is natural.
\end{rem}

We'll next see another way in which $n$-labeled facets can be viewed as trivial, by noting that decorations can `slide' through these facets. 
Composing equation \eqref{nH3rel} with a decorated cup generator \eqref{ccgens} shows that
\begin{equation*}
\xy
(0,0)*{
\begin{tikzpicture} [scale=.6,fill opacity=0.2]
	\path[fill=yellow] (-2.25,3) to (-.75,3) to (-.75,0) to (-2.25,0);
	\path[fill=red] (-.75,3) to [out=45,in=180] (.5,3.5) to (.5,.5) to [out=180,in=45] (-.75,0);
	\path[fill=red] (-.75,3) to [out=315,in=180] (1,2.5) to (1,-.5) to [out=180,in=315] (-.75,0);	
	\draw [double] (-2.25,0) to (-.75,0);
	\draw [very thick,rdirected=.55] (-.75,0) to [out=315,in=180] (1,-.5);
	\draw [very thick,rdirected=.55] (-.75,0) to [out=45,in=180] (.5,.5);
	\draw[very thick, red, rdirected=.55] (-.75,0) to (-.75,3);
	\draw [very thick] (-2.25,3) to (-2.25,0);
	\draw [very thick] (1,2.5) to (1,-.5);
	\draw [very thick] (.5,3.5) to (.5,.5);
	\draw [double] (-2.25,3) to (-.75,3);
	\draw [very thick,rdirected=.55] (-.75,3) to [out=315,in=180] (1,2.5);
	\draw [very thick,rdirected=.55] (-.75,3) to [out=45,in=180] (.5,3.5);
	\node[opacity=1] at (.5,0) {$h_d$};
\end{tikzpicture}
};
\endxy 
\quad = \quad (-1)^d \quad
\xy
(0,0)*{
\begin{tikzpicture} [scale=.6,fill opacity=0.2]
	\path[fill=yellow] (-2.25,3) to (-.75,3) to (-.75,0) to (-2.25,0);
	\path[fill=red] (-.75,3) to [out=45,in=180] (.5,3.5) to (.5,.5) to [out=180,in=45] (-.75,0);
	\path[fill=red] (-.75,3) to [out=315,in=180] (1,2.5) to (1,-.5) to [out=180,in=315] (-.75,0);	
	\draw [double] (-2.25,0) to (-.75,0);
	\draw [very thick,rdirected=.55] (-.75,0) to [out=315,in=180] (1,-.5);
	\draw [very thick,rdirected=.55] (-.75,0) to [out=45,in=180] (.5,.5);
	\draw[very thick, red, rdirected=.55] (-.75,0) to (-.75,3);
	\draw [very thick] (-2.25,3) to (-2.25,0);
	\draw [very thick] (1,2.5) to (1,-.5);
	\draw [very thick] (.5,3.5) to (.5,.5);
	\draw [double] (-2.25,3) to (-.75,3);
	\draw [very thick,rdirected=.55] (-.75,3) to [out=315,in=180] (1,2.5);
	\draw [very thick,rdirected=.55] (-.75,3) to [out=45,in=180] (.5,3.5);
	\node[opacity=1] at (0,3) {$e_d$};
\end{tikzpicture}
};
\endxy 
\end{equation*}
using equation \eqref{blisterrel} (and the Pieri formula \cite{Macdonald}).

One of the consequences of (and motivations for) the skew Howe approach to $\sln$ foams is to give an easy answer to the problem from 
\cite{MSV} of evaluating closed foams. 
In the setting of enhanced foams, this problem is reformulated as follows: do the foam relations suffice to express any foam whose top and 
bottom boundary consist entirely of 
$n$-labeled edges as a multiple of the identity foam. Forgetting the $n$-labeled facets, the former give closed foams, 
so we refer to them as \emph{enhanced closed foams}.

This result follows easily from the following, which describes the vector space of morphisms between grading shifts of 
identity tangles.
\begin{prop}\label{prop:id_end}
Let $\1_\mathbf{a}$ be the identity web of $\mathbf{a}=[a_1,\ldots,a_m]$ in $\foam{n}{}$, then the 
graded vector space $\displaystyle \bigoplus_{k \in \Z} \Hom(\1_\mathbf{a},q^k \1_\mathbf{a})$ is 
spanned by decorated (vertical) sheets, i.e. by foams which are equal to the identity after removing 
decorations. In particular, $\Hom(\1_\mathbf{a}, \1_\mathbf{a})$ is $1$-dimensional and 
$\Hom(\1_\mathbf{a},q^k \1_\mathbf{a})=0$ for $k<0$.
\end{prop}
\begin{proof}
Let $F: \1_\mathbf{a} \to q^k \1_\mathbf{a}$, then by Proposition \ref{prop:fullness}, there exists $m$ so that 
$F= \Phi_n^{m,N}(\cal{D})$ for some $\cal{D}: \onenn{\mathbf{a}'} \to  \onenn{\mathbf{a}'} \{k\}$, where 
$\mathbf{a}' = [a_1,\ldots,a_m,0,\ldots,0]$. If follows from \cite[Proposition 3.6]{KL3} that $\cal{D}$ is equal to 
a sum of products of dotted bubbles. The proof of Theorem \ref{thm:2functor} (or equations \eqref{R2KLR2rel} 
and \eqref{blisterrel}) shows that the image of a dotted bubble, which is a dotted `tube' between two vertical 
sheets, is equal to a linear combination of decorations on those sheets.
\end{proof}

\begin{cor} \label{cor:evaluationFoams}
Any enhanced closed foam in $\foam{n}{}^\bullet$ is a multiple of the trivial foam (consisting of parallel $n$-labeled sheets).
\end{cor}

\begin{rem}
If we do not impose the triviality of the $2$-endomorphism spaces of $n$-labeled strands, then enhanced closed foams would evaluate to 
multiples of decorated $n$-labeled sheets; for example, this gives rise to the deformed versions of $\slnn{2}$ and $\slnn{3}$ link homology. 
Such deformed versions of $\sln$ link homology can be defined in our setting; however, one must know that certain dotted bubbles in the
lowest weight space of $\Ucatc_Q(\glm)^{(n)}$ are non-zero to prove non-degeneracy.
\end{rem}

\begin{rem} 
Proposition \ref{prop:id_end} even provides us with an algorithmic way to evaluate any enhanced closed foam. Given such a foam, we can use Proposition \ref{prop:fullness} to find a preimage
in $\Ucatc_Q(\slm)$, i.e. a closed thick calculus diagram. We can then express the thick calculus elements in terms of thin diagrams, and a 
direct adaptation of \cite[Proposition 8.2]{Lau1} to the $\slm$ case (as outlined in \cite[Proposition 3.6]{KL3}) gives an 
algorithmic process for reducing the diagram to bubbles, from which the value of the foam can be immediately read off.
\end{rem}

\subsection{The $\sln$ tangle invariant} \label{subsec:tangleInv}

Let $\tau$ be a tangle with each component colored by a fundamental representation of $\sln$, i.e. a number between $1$ and $n-1$. 
If $c_1,\ldots,c_r$ are the colorings of the $r$ right-endpoints
\footnote{As with webs, we view tangles as mapping from their right endpoints to 
their left endpoints.} of $\tau$, the tangle invariant will take values in the homotopy category 
of complexes of the $2$-category\footnote{In this section, we'll work with the more general $2$-category $\foam{n}{}$ (without the 
extra relation on $1$-labeled facets), since all of our results actually hold there.}
$\foam{n}{}$ where $N=\sum_{i=1}^{r} \cal{O}(c_i) + ns$ for $s$ sufficiently large and
$\cal{O}(c_i)=c_i$ if the tangle component is directed out from the $i^{th}$ right endpoint and equals $n-c_i$ otherwise. We now define the 
tangle invariant $\llbracket \tau \rrbracket$, first locally on generating tangles, and then explain how to piece together composites of 
generating tangles. This construction is essentially the same as that from \cite{LQR1}, generalized to the $\sln$ foam $2$-categories.

Define a \emph{canonical sequence} to be an object of $\foam{n}{}$ of the form $[a_1,\ldots,a_k,n,\ldots,n]$ where 
$a_i \neq n$. We associate to each generating tangle a web mapping between two canonical sequences via the following rules. 
We map the identity tangle to the identity web of the object $[\cal{O}(c_1),\ldots,\cal{O}(c_r),n,\ldots,n]$ and 
for cap and cup tangles we set
\begin{align*}
\left \llbracket
\xy
(0,0)*{\begin{tikzpicture} [scale=.5]
\draw[very thick, directed=.99] (0,1) to [out=180,in=90] (-1.25,.5) to [out=270,in=180] (0,0);
\node at (-1.5,1) {\small$a$};
\end{tikzpicture}};
\endxy 
\;\; \right \rrbracket
\quad = \quad
\xy
(0,0)*{\begin{tikzpicture} [scale=.5]
 \draw [double] (-1.25,.5) to (-2.25,.5);
 \draw [very thick, directed=.55] (0,0) to [out=180,in=300] (-1.25,.5);
\draw [very thick, directed=.55] (0,1) to [out=180,in=60] (-1.25,.5);
 \node at (.875,0) {\small$n-a$};
 \node at (.5,1) {\small$a$};
 \node at (-2.75,.5) {\small$n$}; 
\end{tikzpicture}};
\endxy
\quad &, \quad
\left \llbracket \;\;
\xy
(0,0)*{\begin{tikzpicture} [scale=.5]
\draw[very thick, directed=.99] (0,0) to [out=0,in=270] (1.25,.5) to [out=90,in=0] (0,1);
\node at (1.5,1) {\small$a$};
\end{tikzpicture}};
\endxy
\right \rrbracket
\quad = \quad
\xy
(0,0)*{\begin{tikzpicture} [scale=.5]
 \draw [double] (2.25,.5) -- (1.25,.5);
 \draw [very thick, directed=.55] (1.25,.5) to [out=240,in=0] (0,0);
\draw [very thick, directed=.55] (1.25,.5) to [out=120,in=0] (0,1);
 \node at (-.875,0) {\small$n-a$};
 \node at (-.5,1) {\small$a$};
 \node at (2.75,.5) {\small$n$}; 
\end{tikzpicture}};
\endxy \\
\left \llbracket
\xy
(0,0)*{\begin{tikzpicture} [scale=.5]
\draw[very thick, rdirected=.05] (0,1) to [out=180,in=90] (-1.25,.5) to [out=270,in=180] (0,0);
\node at (-1.5,1) {\small$a$};
\end{tikzpicture}};
\endxy
\;\; \right \rrbracket
\quad = \quad
\xy
(0,0)*{\begin{tikzpicture} [scale=.5]
 \draw [double] (-1.25,.5) to (-2.25,.5);
 \draw [very thick, directed=.55] (0,0) to [out=180,in=300] (-1.25,.5);
\draw [very thick, directed=.55] (0,1) to [out=180,in=60] (-1.25,.5);
 \node at (.875,1) {\small$n-a$};
 \node at (.5,0) {\small$a$};
 \node at (-2.75,.5) {\small$n$}; 
\end{tikzpicture}};
\endxy
\quad &, \quad
\left \llbracket \;\;
\xy
(0,0)*{\begin{tikzpicture} [scale=.5]
\draw[very thick, rdirected=.05] (0,0) to [out=0,in=270] (1.25,.5) to [out=90,in=0] (0,1);
\node at (1.5,1) {\small$a$};
\end{tikzpicture}};
\endxy
\right \rrbracket
\quad = \quad
\xy
(0,0)*{\begin{tikzpicture} [scale=.5]
 \draw [double] (2.25,.5) -- (1.25,.5);
 \draw [very thick, directed=.55] (1.25,.5) to [out=240,in=0] (0,0);
\draw [very thick, directed=.55] (1.25,.5) to [out=120,in=0] (0,1);
 \node at (-.875,1) {\small$n-a$};
 \node at (-.5,0) {\small$a$};
 \node at (2.75,.5) {\small$n$}; 
\end{tikzpicture}};
\endxy \quad .
\end{align*}

To define the invariant on crossings, we first introduce the following shifted\footnote{Given a cochain complex $C^{\bullet}$, we follow the 
convention that $C[k]^\bullet = C^{\bullet+k}$.} versions of the Rickard 
complexes \eqref{Rickardp} and \eqref{Rickardn} in $\UcatD_Q(\glnn{2})$:
\begin{align}
\tilde{\cal{T}} \onenn{[a,b]}&= \cal{T}\onenn{(a-b)} [\min(a,b)]\{-\min(a,b)\} \label{Rickard_resc_+}\\
\tilde{\cal{T}}^{-1} \onenn{[a,b]}&= \cal{T}\onenn{(a-b)} [-\min(a,b)]\{\min(a,b)\} \label{Rickard_resc_-}
\end{align}
and analogous complexes $\tilde{\cal{T}_i}$ and $\tilde{\cal{T}_i}^{-1}$ in $\Ucatc_Q(\glm)$.
These shifts are motivated by the work of the first author \cite{Queff_aff}, in which it is shown at the decategorified level that 
analogous shifts give an invariant of tangled webs, in particular the relation
\begin{equation}\label{tangle_fork}
\xy
(0,0)*{\begin{tikzpicture} [scale=.5]
 \draw [very thick, ->] (-1.25,.5) to (-3,.5);
 \draw [very thick, directed=.55] (.75,0) to [out=180,in=300] (-1.25,.5);
\draw [very thick, directed=.55] (.75,1) to [out=180,in=60] (-1.25,.5);
\draw[very thick] (.75,2) to [out=180,in=60] (-1.7,.6);
\draw[very thick, ->] (-1.8,.4) to [out=240,in=0] (-3,-.5);
 \node at (1,0) {\small$c$};
 \node at (1,1) {\small$b$};
  \node at (1,2) {\small$a$};
 \node at (-3.875,.5) {\small$b+c$}; 
\end{tikzpicture}};
\endxy
\quad \sim \quad
\xy
(0,0)*{\begin{tikzpicture} [scale=.5]
 \draw [very thick, directed=.55] (-1.25,.5) to (-3,.5);
 \draw [very thick, directed=.55] (.75,0) to [out=180,in=300] (-1.25,.5);
\draw [very thick, directed=.7] (.75,1) to [out=180,in=60] (-1.25,.5);
\draw[very thick] (.75,2) to [out=180,in=60] (-.15,1.2);
\draw[very thick] (-.3,.9) to [out=240,in=60] (-.65,.2);
\draw[very thick, ->] (-.8,-.1) to [out=240,in=0] (-3,-.5);
 \node at (1,0) {\small$c$};
 \node at (1,1) {\small$b$};
 \node at (1,2) {\small$a$};
 \node at (-3.875,.5) {\small$b+c$}; 
\end{tikzpicture}};
\endxy
\end{equation}
and its variations are satisfied. This fact follows from a result of Cautis \cite[Lemma 5.2]{Cautis}, building 
on his previous work with Kamnitzer \cite{CK}. In the above notation, this result is the following.
\begin{lem} \label{lem:webMove}
The complexes $\tilde{\cal{T}}_i$ braid and for $|i-j|=1$ satisfy the relations:
\begin{equation}\label{VertexRel}
\cal{E}_j \tilde{\cal{T}}_i \tilde{\cal{T}}_j \onea \simeq \tilde{\cal{T}}_i \tilde{\cal{T}}_j\cal{E}_i \onea 
\quad,\quad 
\cal{F}_j\tilde{\cal{T}}_i \tilde{\cal{T}}_j \onea \simeq \tilde{\cal{T}}_i \tilde{\cal{T}}_j \cal{F}_i \onea. 
\end{equation}
in any integrable $2$-representation.
\end{lem}
\begin{proof}
The first part is obvious: the same crossings appear in both parts of the braiding relation, hence the same rescalings occur. 
The second part is a case-by-case adaptation of \cite[Lemma 5.2]{Cautis}.
\end{proof}

We now assign the following complexes to positive and negative left-directed crossings:
\begin{align*}
\left \llbracket \;\;
\xy
(0,0)*{
\begin{tikzpicture} [scale=.5]
\draw [very thick, directed=.99] (2,0) -- (0,1);
\draw [very thick] (2,1) -- (1.2, .6);
\draw [very thick, directed=.99] (.8,.4) -- (0,0);
\node at (2.3,1) {\small $a$};
\node at (2.3,0) {\small $b$};
\end{tikzpicture}};
\endxy
\right \rrbracket
\quad = \quad
\Phi_n(\tilde{\cal{T}}\mathbf{1}_{[a,b]})
\qquad , \qquad
\left \llbracket \;\;
\xy
(0,0)*{
\begin{tikzpicture} [decoration={markings, mark=at position 1 with {\arrow{>}}; },scale=.5]
\draw [very thick, postaction={decorate}] (2,1) -- (0,0);
\draw [very thick] (2,0) -- (1.2, .4);
\draw [very thick, postaction={decorate}] (.8,.6) -- (0,1);
\node at (2.3,1) {\small $a$};
\node at (2.3,0) {\small $b$};
\end{tikzpicture}};
\endxy
\right \rrbracket
\quad = \quad
\Phi_n(\tilde{\cal{T}}^{-1} \mathbf{1}_{[a,b]})
\end{align*}
for example, this gives
\begin{equation}\label{1,1_pos_cross}
\left \llbracket \;\;
\xy
(0,0)*{
\begin{tikzpicture} [scale=.5]
\draw [very thick, directed=.99] (2,0) -- (0,1);
\draw [very thick] (2,1) -- (1.2, .6);
\draw [very thick, directed=.99] (.8,.4) -- (0,0);
\node at (2.3,1) {\small $1$};
\node at (2.3,0) {\small $1$};
\end{tikzpicture}};
\endxy
\right \rrbracket
\quad = \quad
\left(
\xymatrix{
q^{-1} \;\;
\xy
(0,0)*{\begin{tikzpicture} [scale=.5]
\draw[very thick, directed=.45] (2.5,1) to [out=185,in=355] (0,1);
\draw[very thick, directed=.45] (2.5,0) to [out=175,in=5] (0,0);
 \node at (-.25,1) {\tiny $1$};
 \node at (-.25,0) {\tiny $1$};
\end{tikzpicture}};
\endxy
\;
\ar[rrr]^-{
\xy
(0,0)*{
\begin{tikzpicture} [scale=.4,fill opacity=0.2]
	\path [fill=red] (4.25,2) to (4.25,-.5) to [out=170,in=10] (-.5,-.5) to (-.5,2) to
		[out=0,in=225] (.75,2.5) to [out=270,in=180] (1.625,1.25) to [out=0,in=270] 
			(2.5,2.5) to [out=315,in=180] (4.25,2);
	\path [fill=red] (3.75,3) to (3.75,.5) to [out=190,in=350] (-1,.5) to (-1,3) to [out=0,in=135]
		(.75,2.5) to [out=270,in=180] (1.625,1.25) to [out=0,in=270] 
			(2.5,2.5) to [out=45,in=180] (3.75,3);
	\path[fill=blue] (2.5,2.5) to [out=270,in=0] (1.625,1.25) to [out=180,in=270] (.75,2.5);
	\draw [very thick,directed=.55] (4.25,-.5) to [out=170,in=10] (-.5,-.5);
	\draw [very thick, directed=.55] (3.75,.5) to [out=190,in=350] (-1,.5);
	\draw [very thick, red, directed=.75] (2.5,2.5) to [out=270,in=0] (1.625,1.25);
	\draw [very thick, red] (1.625,1.25) to [out=180,in=270] (.75,2.5);
	\draw [very thick] (3.75,3) to (3.75,.5);
	\draw [very thick] (4.25,2) to (4.25,-.5);
	\draw [very thick] (-1,3) to (-1,.5);
	\draw [very thick] (-.5,2) to (-.5,-.5);
	\draw [very thick,directed=.55] (2.5,2.5) to (.75,2.5);
	\draw [very thick,directed=.55] (.75,2.5) to [out=135,in=0] (-1,3);
	\draw [very thick,directed=.55] (.75,2.5) to [out=225,in=0] (-.5,2);
	\draw [very thick,directed=.55] (3.75,3) to [out=180,in=45] (2.5,2.5);
	\draw [very thick,directed=.55] (4.25,2) to [out=180,in=315] (2.5,2.5);
\end{tikzpicture}};
\endxy
}
& & &
\;\;
\uwave{
\xy
(0,9)*{};
(0,-7)*{};
(0,0)*{
\xy
(0,0)*{\begin{tikzpicture} [scale=.5]
\draw[very thick, directed=.45] (3.5,1) to [out=180,in=60] (2.25,.5);
\draw[very thick, directed=.45] (3.5,0) to [out=180,in=300] (2.25,.5);
 \draw [very thick, directed=.55] (2.25,.5) -- (1.25,.5);
 \draw [very thick, directed=.65] (1.25,.5) to [out=240,in=0] (0,0);
\draw [very thick, directed=.65] (1.25,.5) to [out=120,in=0] (0,1);
 \node at (-.25,1) {\tiny $1$};
 \node at (-.25,0) {\tiny $1$};
\node at (3.75,1) {\tiny $1$};
\node at (3.75,0) {\tiny $1$};
\end{tikzpicture}};
\endxy}
\endxy}
} \right)
\end{equation}
\begin{equation}\label{1,1_neg_cross}
\left \llbracket \;\;
\xy
(0,0)*{
\begin{tikzpicture} [scale=.5]
\draw [very thick, directed=.99] (2,1) to (0,0);
\draw [very thick] (2,0) -- (1.2, .4);
\draw [very thick, directed=.99] (.8,.6) to (0,1);
\node at (2.3,1) {\small $1$};
\node at (2.3,0) {\small $1$};
\end{tikzpicture}};
\endxy
\right \rrbracket
\quad = \quad
\left(
\xymatrix{
\uwave{
\xy
(0,9)*{};
(0,-7)*{};
(0,0)*{
\xy
(0,0)*{\begin{tikzpicture} [scale=.5]
\draw[very thick, directed=.45] (3.5,1) to [out=180,in=60] (2.25,.5);
\draw[very thick, directed=.45] (3.5,0) to [out=180,in=300] (2.25,.5);
 \draw [very thick, directed=.55] (2.25,.5) -- (1.25,.5);
 \draw [very thick, directed=.65] (1.25,.5) to [out=240,in=0] (0,0);
\draw [very thick, directed=.65] (1.25,.5) to [out=120,in=0] (0,1);
 \node at (-.25,1) {\tiny $1$};
 \node at (-.25,0) {\tiny $1$};
\node at (3.75,1) {\tiny $1$};
\node at (3.75,0) {\tiny $1$};
\end{tikzpicture}};
\endxy}
\endxy}
\ar[rrr]^-{
\xy
(0,0)*{
\begin{tikzpicture} [scale=.4,fill opacity=0.2]
	\path [fill=red] (4.25,-.5) to (4.25,2) to [out=170,in=10] (-.5,2) to (-.5,-.5) to 
		[out=0,in=225] (.75,0) to [out=90,in=180] (1.625,1.25) to [out=0,in=90] 
			(2.5,0) to [out=315,in=180] (4.25,-.5);
	\path [fill=red] (3.75,.5) to (3.75,3) to [out=190,in=350] (-1,3) to (-1,.5) to 
		[out=0,in=135] (.75,0) to [out=90,in=180] (1.625,1.25) to [out=0,in=90] 
			(2.5,0) to [out=45,in=180] (3.75,.5);
	\path[fill=blue] (.75,0) to [out=90,in=180] (1.625,1.25) to [out=0,in=90] (2.5,0);
	\draw [very thick,directed=.55] (2.5,0) to (.75,0);
	\draw [very thick,directed=.55] (.75,0) to [out=135,in=0] (-1,.5);
	\draw [very thick,directed=.55] (.75,0) to [out=225,in=0] (-.5,-.5);
	\draw [very thick,directed=.55] (3.75,.5) to [out=180,in=45] (2.5,0);
	\draw [very thick,directed=.55] (4.25,-.5) to [out=180,in=315] (2.5,0);
	\draw [very thick, red, directed=.75] (.75,0) to [out=90,in=180] (1.625,1.25);
	\draw [very thick, red] (1.625,1.25) to [out=0,in=90] (2.5,0);
	\draw [very thick] (3.75,3) to (3.75,.5);
	\draw [very thick] (4.25,2) to (4.25,-.5);
	\draw [very thick] (-1,3) to (-1,.5);
	\draw [very thick] (-.5,2) to (-.5,-.5);
	\draw [very thick,directed=.55] (4.25,2) to [out=170,in=10] (-.5,2);
	\draw [very thick, directed=.55] (3.75,3) to [out=190,in=350] (-1,3);
\end{tikzpicture}
};
\endxy
}
& & &
q \;\;
\xy
(0,0)*{\begin{tikzpicture} [scale=.5]
\draw[very thick, directed=.45] (2.5,1) to [out=185,in=355] (0,1);
\draw[very thick, directed=.45] (2.5,0) to [out=175,in=5] (0,0);
 \node at (-.25,1) {\tiny $1$};
 \node at (-.25,0) {\tiny $1$};
\end{tikzpicture}};
\endxy
} \right)
\end{equation}
which are the relevant complexes for (uncolored) link homology. Note that we need not 
specify the images of any other crossings (i.e. those oriented up, down, or rightward
\footnote{We show below that the two possible ways of defining a rightward crossing 
are equivalent}), since they 
may be constructed by composing the left-directed crossings with cap and cup tangles.

Notice, however, that the above procedure will not necessarily assign (a complex of) webs mapping 
between canonical sequences to these other crossings, or to tensor products of the above tangle generators, 
e.g. we would naively map
\[
\xymatrix{
\xy
(0,0)*{\begin{tikzpicture} [scale=.5]
\draw[very thick, directed=.99] (0,2) to (-3,2);
\draw[very thick, directed=.99] (0,1) to [out=180,in=90] (-1.25,.5) to [out=270,in=180] (0,0);
\draw[very thick, directed=.99] (-3,-1) to (0,-1);
\node at (-1.625,1) {\small$a_2$};
\node at (-3.5,2) {\small$a_1$};
\node at (-3.5,-1) {\small$a_3$};
\end{tikzpicture}};
\endxy
\quad
\ar@{|~>}[rr]
& &
\quad
\xy
(0,0)*{\begin{tikzpicture} [scale=.5]
 \draw [double] (-1.25,.5) to (-2.5,.5);
 \draw [very thick, directed=.55] (.25,0) to [out=180,in=300] (-1.25,.5);
\draw [very thick, directed=.55] (.25,1) to [out=180,in=60] (-1.25,.5);
\draw[very thick, directed=.99] (.25,2) to (-2.5,2);
\draw[very thick, directed=.99] (.25,-1) to (-2.5,-1);
 \node at (1.25,0) {\small$n-a_2$};
 \node at (.75,1) {\small$a_2$};
 \node at (-3,.5) {\small$n$}; 
 \node at (.75,2) {\small$a_1$};
\node at (1.25,-1) {\small$n-a_3$};
\end{tikzpicture}};
\endxy
}
\]
and the codomain of the latter is not a canonical sequence. 
We thus compose with the natural equivalence between a sequence and its 
corresponding canonical sequence, given by (composites of) appropriate \emph{unshifted}
\footnote{Here, we use the unshifted Rickard complexes, rather than
crossings between $n$- and $a$-labeled strands. Of course, these equivalences are equal, 
up to degree shifts, to both the positive and negative crossings defined above; however,
there is no natural choice of positive or negative crossing to use for this equivalence. 
These webs can be viewed as the `average' of the positive and negative crossing.}
Rickard complexes:
\begin{equation}\label{n_equiv}
\Phi_n(\cal{T}^{\pm1}\mathbf{1}_{[n,a]}) = 
\xy
(0,0)*{\begin{tikzpicture} [scale=.5]
 \draw [double] (3,1) -- (1.75,1);
 \draw[very thick, directed=.55] (3,0) to (1.25,0);
 \draw [double] (1.25,0) -- (0,0);
\draw [very thick, directed=.55] (1.75,1) -- (0,1);
 \draw [very thick, directed=.55] (1.75,1) -- (1.25,0);
 \node at (3.5,0) {\small$a$};
 \node at (3.5,1) {\small$n$}; 
\end{tikzpicture}};
\endxy
\quad , \quad
\Phi_n(\cal{T}^{\pm1} \mathbf{1}_{[a,n]}) = 
\xy
(0,0)*{\begin{tikzpicture} [scale=.5]
 \draw [double] (1.25,0) -- (0,0);
\draw [very thick, directed=.55] (1.25,1) -- (-.5,1);
\draw[very thick, directed=.55] (0,0) to (-1.75,0);
\draw[double] (-.5,1) to (-1.75,1);
 \draw[very thick, directed=.55] (0,0) to (-.5,1); 
 \node at (1.75,0) {\small$n$};
 \node at (1.75,1) {\small$a$}; 
\end{tikzpicture}};
\endxy
\end{equation}
so, for example, the invariant of the upward crossing
$
\left \llbracket \;\;
\xy
(0,0)*{
\begin{tikzpicture} [scale=.5]
\draw [very thick, rdirected=.1] (2,1) -- (0,0);
\draw [very thick] (2,0) -- (1.2, .4);
\draw [very thick, directed=.99] (.8,.6) -- (0,1);
\node at (2.3,1) {\small $1$};
\node at (2.3,0) {\small $1$};
\end{tikzpicture}};
\endxy
\right \rrbracket
$
equals the (horizontal) composition of the complex
\[
\xymatrix{
\xy
(0,0)*{\begin{tikzpicture} [scale=.5]
 \draw [double] (-2,1) to (-3,1);
 \draw [very thick, directed=.55] (3,.5) to (-.75,.5) to [out=180,in=300] (-2,1);
\draw [very thick, directed=.55] (3,1.5) to (-.75,1.5) to  [out=180,in=60] (-2,1);
 \draw [double] (3,-1) -- (2,-1);
 \draw [very thick, directed=.55] (2,-1) to [out=240,in=0] (.75,-1.5) to (-3,-1.5);
\draw [very thick, directed=.55] (2,-1) to [out=120,in=0] (.75,-.5) to (-3,-.5);
\node at (3.5,.5) {\small $1$};
\node at (3.875,1.5) {\small $n-1$};
\node at (3.5,-1) {\small $n$};
\node at (-3.5,-.5) {\small $1$};
\node at (-3.875,-1.5) {\small $n-1$};
\node at (-3.5,1) {\small $n$};
\end{tikzpicture}};
\endxy
\ar[rrr]^-{
\xy
(0,0)*{
\begin{tikzpicture} [scale=.4,fill opacity=0.2]
	\path [fill=red] (4.25,2) to (4.25,-.5) to [out=170,in=10] (-.5,-.5) to (-.5,2) to
		[out=0,in=225] (.75,2.5) to [out=270,in=180] (1.625,1.25) to [out=0,in=270] 
			(2.5,2.5) to [out=315,in=180] (4.25,2);
	\path [fill=red] (3.75,3) to (3.75,.5) to [out=190,in=350] (-1,.5) to (-1,3) to [out=0,in=135]
		(.75,2.5) to [out=270,in=180] (1.625,1.25) to [out=0,in=270] 
			(2.5,2.5) to [out=45,in=180] (3.75,3);
	\path[fill=blue] (2.5,2.5) to [out=270,in=0] (1.625,1.25) to [out=180,in=270] (.75,2.5);
	\draw [very thick,directed=.55] (4.25,-.5) to [out=170,in=10] (-.5,-.5);
	\draw [very thick, directed=.55] (3.75,.5) to [out=190,in=350] (-1,.5);
	\draw [very thick, red, directed=.75] (2.5,2.5) to [out=270,in=0] (1.625,1.25);
	\draw [very thick, red] (1.625,1.25) to [out=180,in=270] (.75,2.5);
	\draw [very thick] (3.75,3) to (3.75,.5);
	\draw [very thick] (4.25,2) to (4.25,-.5);
	\draw [very thick] (-1,3) to (-1,.5);
	\draw [very thick] (-.5,2) to (-.5,-.5);
	\draw [very thick,directed=.55] (2.5,2.5) to (.75,2.5);
	\draw [very thick,directed=.55] (.75,2.5) to [out=135,in=0] (-1,3);
	\draw [very thick,directed=.55] (.75,2.5) to [out=225,in=0] (-.5,2);
	\draw [very thick,directed=.55] (3.75,3) to [out=180,in=45] (2.5,2.5);
	\draw [very thick,directed=.55] (4.25,2) to [out=180,in=315] (2.5,2.5);
\end{tikzpicture}};
\endxy
}
&&& \; q \;\;
\xy
(0,0)*{\begin{tikzpicture} [scale=.5]
 \draw [double] (-2.5,1) to (-3.5,1);
 \draw [very thick, directed=.55] (-.5,0) to [out=120,in=300] (-2.5,1);
\draw [very thick, directed=.55] (3.5,1.5) to (-1.25,1.5) to [out=180,in=60] (-2.5,1);
\draw[very thick, directed=.55] (3.5,.5) to [out=180,in=60] (.5,0);
\draw [very thick, directed=.55] (.5,0) -- (-.5,0);
\draw[very thick, directed=.55] (-.5,0) to [out=240,in=0] (-3.5,-.5);
 \draw [double] (3.5,-1) -- (2.5,-1);
 \draw [very thick, directed=.55] (2.5,-1) to [out=240,in=0] (1.25,-1.5) to (-3.5,-1.5);
\draw [very thick, directed=.55] (2.5,-1) to [out=120,in=300] (.5,0);
\node at (4,.5) {\small $1$};
\node at (4.375,1.5) {\small $n-1$};
\node at (4,-1) {\small $n$};
\node at (-4,-.5) {\small $1$};
\node at (-4.375,-1.5) {\small $n-1$};
\node at (-4,1) {\small $n$};
\end{tikzpicture}};
\endxy
}
\]
with the web
\[
\xy
(0,0)*{\begin{tikzpicture} [scale=.5]
\draw [double] (3,1) -- (1.75,1);
\draw[very thick, directed=.55] (3,0) to (1.25,0);
\draw[very thick, directed=.55] (3,-1) to (-.5,-1);
\draw [very thick, directed=.55] (1.75,1) -- (1.25,0);
\draw [double] (1.25,0) -- (0,0);
\draw[very thick, directed=.55] (0,0) to (-1.75,0);
\draw[very thick, directed=.55] (0,0) to (-.5,-1);
\draw[double] (-.5,-1) to (-1.75,-1); 
\draw [very thick, directed=.55] (1.75,1) -- (-1.75,1);
\node at (3.875,-1) {\small$n-1$};
\node at (3.5,1) {\small$n$}; 
\node at (3.5,0) {\small$1$}; 
\node at (-2.625,0) {\small$n-1$};
\node at (-2.25,-1) {\small$n$}; 
\node at (-2.25,1) {\small$1$}; 
\end{tikzpicture}};
\endxy
\]
and similarly for other crossings. In fact, one can check that, up to isomorphism and grading shifts, 
non-left directed crossings are given by the Rickard complexes for the crossings in which 
any right-directed $a$-labeled strands are replaced by left-directed $(n-a)$-labeled strands.

The above formulae assign to any tangle $\tau$ a complex $\llbracket \tau \rrbracket$ of webs mapping 
between canonical sequences, by extending the local assignment for crossings to diagrams containing many 
crossing planar algebraically, i.e. in a manner similar to taking tensor products of complexes of abelian groups.

\begin{thm}\label{thm:tangle_inv}
Let $\tau$ be a framed, oriented tangle with right endpoints colored $c_1,\ldots,c_r$, then the complex
$\llbracket \tau \rrbracket$, viewed as an element of the homotopy category of complexes 
of $\foam[\sum_i \cal{O}(c_i)+ns]{n}{}$ for $s$ sufficiently large, is an invariant of $\tau$. 
\end{thm}

\begin{rem}
To be precise, the natural setting for the invariant is the ``direct limit" of the $2$-categories
$\displaystyle \lim_{\stackrel{\longrightarrow}{s}} \foam[ns+\sum \cal{O}(c_i)]{n}{}$ in which 
we take the tensor product of all $0$-, $1$-, and $2$-morphisms with infinitely many $n$-labels/edges/sheets.
\end{rem}

\begin{proof}[Proof (of Theorem \ref{thm:tangle_inv}).]
It suffices to check that applying a (framed) oriented tangle Reidemeister move to $\tau$ does not change $\llbracket \tau \rrbracket$ 
up to homotopy equivalence. Following \cite{Kassel}, it suffices to check homotopy invariance under the moves:

\begin{equation}\label{R0}
\xy
(0,0)*{

};
\endxy
\end{equation}
and to observe that the invariant behaves appropriately with regard to the change of framing moves:
\begin{equation} \label{R1}
\xy
(0,0)*{
%
};
\endxy
\end{equation}
where strands can carry any labels.

Relations \eqref{R2++} and \eqref{R3} are braidlike Reidemeister moves, and follow since the Rickard complexes satisfy the braid relations 
in any integrable $2$-representation \cite{CK}. The first tangle isotopy in \eqref{R0} follows from the web isomorphism
\begin{equation}\label{n_iso}
\xy
(0,0)*{
%
};
\endxy
\end{equation}
which is a consequence of equation \eqref{FE_to_cb}. A similar computation shows the other isotopy.

Throughout the proofs of the other moves, which we now sketch, we will also make use of Lemma \ref{lem:webMove}.
We first verify moves \eqref{R2+-} and \eqref{R2-+}; we'll only verify the first move in \eqref{R2+-} since the other cases follow similarly. 
The LHS corresponds to the following tangled web:
\[
\xy
(0,0)*{
%
};
\endxy
\]
We then use equation \eqref{EF_to_cb} to remove the square in the middle and perform a Reidemeister 2 move to obtain the web:
\[
\xy
(0,0)*{
%
};
\endxy \quad .
\]

To deduce \eqref{TMove2} and \eqref{R1}, we first note the relation:
\begin{equation}\label{fork_crossing}
\left \llbracket
\xy
(0,0)*{%
};
\endxy}}
\right)
\end{align*}
using the isomorphisms given by \eqref{MVrel} and \eqref{digon_web_iso}. The matrix giving $d$ in the latter 
equation has zeros below the subdiagonal since there are no $1$-morphisms 
$
\xy
(0,0)*{%
};
\endxy
$
for $l<0$, and a computation shows that the terms on the subdiagonal equal $\pm id$. 
Using a Gaussian elimination homotopy equivalence to cancel terms
(see e.g. \cite[Lemma 4.2]{BN3}) then shows that this complex is homotopy equivalent to 
one with only $q^a \;\;
\xy
(0,0)*{%
};
\endxy
$ in homological degree zero. The general form of equation \eqref{fork_crossing} follows from the $b=1$ case 
inductively, and formulas for the negative crossing (and with the web on the right of the crossing) follow similarly.
Using this, together with Lemma \ref{lem:webMove} and equation \eqref{n_iso} we compute that the change of 
framing relation \eqref{R1} is 
\begin{align*}
\left \llbracket
\xy
(0,0)*{
%
};
\endxy
\right) [-a]
\end{align*}
which follows using a similar argument to \cite[Lemma 14.7]{Wu}, together with a grading-based argument to 
cancel `multiples' on each side of the homotopy equivalence (since our categories are not Krull-Schmidt).

Finally, relation \eqref{TMove2} follows from those derived thus far.
Up to shifts coming from replacing the equivalences in \eqref{n_equiv} with crossings, 
the left hand side of the first relation in \eqref{TMove2} is 
\[
\xy
(0,0)*{
%
};
\endxy .
\]
Applying equation \eqref{fork_crossing} and its analog, noticing that the degree shifts compensate, 
this is the same as:
\[
\xy
(0,0)*{
%
};
\endxy .
\]
Incorporating the shift to replace the right-most $n$-crossing with the equivalence in \eqref{n_equiv}, 
and removing the square using equation \eqref{FE_to_cb} gives:
\[
q^{-k}\;
\xy
(0,0)*{
%
};
\endxy
\;[k]
\]
We can now perform the same operation at the top of the picture: pull the top back strand down behind the $n$-labeled strand, 
(this gives a $q^{k}$ factor and a homological shift $[-k]$, since it creates a negative $n$-crossing, which cancels the previous $q^{-k}$ factor and the homological shift), 
remove the created square, and delete both top and bottom curls, noticing again that the coefficients compensate. We are now left with:
\[
\xy
(0,0)*{
%
};
\endxy
\]
which we identify with the other side of relation \eqref{TMove2}, again up to shifts required to replace crossings with the
equivalences from \eqref{n_equiv}. It's easy to see that these shifts exactly cancel the shifts from the first step above.
\end{proof}

\begin{rem}\label{rem:ben}
The invariant $\llbracket \tau \rrbracket$ should be viewed as the `universal' (undeformed) $\sln$ link invariant, akin to the 
construction given in \cite{BN2} for $\slnn{2}$ link homology. As in that setting, one benefit of the above construction is that 
the invariant can be computed locally, first simplifying the complexes assigned to sub-tangles before glueing them together.

There are two additional benefits to our combinatorial presentation of $\sln$ link homology, over existing constructions. First, 
recall from \cite[Sections 9 and 10]{KhR} that the main obstruction to the (proper) functoriality of Khovanov-Rozansky homology 
is the fact that the map induced by a saddle cobordism is only well-defined up to $\pm1$. The analogs of a saddle in our $2$-category are 
the generators in equation \eqref{uzgens}, with $b=n-a$. That these foams lack the problematic rotational symmetry suggests that our 
link homology theory may be properly functorial; indeed, this is the case in the $\slnn{2}$ \cite{Blan} and $\slnn{3}$ case \cite{Clark}. 
However, as the current paper is already quite long, we'll reserve this consideration for future work.

Second, note that all of our foam relations involve integer coefficients. We can thus consider the complex $\llbracket \tau \rrbracket$
in an integral version of $\foam{n}{}$, where $2$-morphisms are given by (matrices of) $\Z$-linear combinations of foams. 
In the uncolored case, a straightforward check, akin to the proofs in \cite[Theorem 7.1]{MSV}, shows that all of 
the homotopy equivalences in Theorem \ref{thm:tangle_inv} hold over the integers.
The construction in the following section then produces an abelian group valued link homology theory.
\end{rem}

In fact, this integrality result holds true in the colored case as well:

\begin{prop}
The complex $\llbracket \tau \rrbracket$ considered in the homotopy category of $\foam{n}{}$ defined over the integers is a tangle invariant.
\end{prop}
\begin{proof}[Proof (sketch)]
It is enough to prove the Reidemeister $2$ move and the web-move \eqref{tangle_fork}, as all others can be derived from these.
The proof for these two moves is achieved simultaneously (over $\Z$ and for any value of $n$) by iteration.

One can first check the homotopy equivalence:
\[
\left \llbracket
\xy
(0,0)*{
\begin{tikzpicture}[scale=.6]
\draw [very thick] (3,1) -- (2,1);
\draw [very thick, ->] (2,1) .. controls (1,2) .. (0,2);
\draw [very thick, ->] (2,1) .. controls (1,0) .. (0,0);
\draw [very thick] (3,2) .. controls (2.35,2) .. (1.65,1.55);
\draw [very thick] (1.45,1.35) .. controls (1,.8) and (.9,.5) .. (.8,.15);
\draw [very thick, ->] (.75,-.1) .. controls (.6,-.6) and  (.4,-1) .. (0,-1);
\node at (3.6,1) {\tiny $k+1$};
\node at (3.3,2) {\tiny $p$};
\node at (-.3,2) {\tiny $k$};
\node at (-.3,0) {\tiny $1$};
\end{tikzpicture}
};
\endxy
\right \rrbracket
\quad
\simeq
\quad
\left \llbracket
\xy
(0,0)*{
\begin{tikzpicture}[scale=.6]
\draw [very thick] (3,1) -- (2,1);
\draw [very thick, ->] (2,1) .. controls (1,2) .. (0,2);
\draw [very thick, ->] (2,1) .. controls (1,0) .. (0,0);
\draw [very thick] (3,2) .. controls (2.6,2) .. (2.38,1.1);
\draw [very thick, ->] (2.35,.9) .. controls (2.2,0) and (1.5,-1) .. (0,-1);
\node at (3.6,1) {\tiny $k+1$};
\node at (3.3,2) {\tiny $p$};
\node at (-.3,2) {\tiny $k$};
\node at (-.3,0) {\tiny $1$};
\end{tikzpicture}
};
\endxy
\right \rrbracket
\]
via an explicit computation. The complex assigned to the left-hand side can be expanded, and, 
through a repeated use of Gaussian elimination homotopy equivalences and 
Remark \ref{rem:SquareFormulas} (or more precisely Stosic's formula from \cite{KLMS}), proven to be homotopy equivalent to the 
right-hand side.

An inductive argument then shows that the (general, colored) Reidemeister 2 move follows from cases which can be checked directly. 
Finally, one uses both the reduced case of the web-move and the Reidemeister 2 move to prove the general version of \eqref{tangle_fork}.
\end{proof}

\subsection{$\sln$ link homology} 

We now explain how to recover $\sln$ link homology from the tangle invariant described above;  to compare with the Khovanov-Rozansky 
link invariant, we now set $\Bbbk = \Q$.
Given a framed, oriented tangle $\tau$, consider the complex
$\llbracket \tau \rrbracket$ in the homotopy category of $\foam{n}{}^\bullet$, i.e. in the homotopy category of the 
$\Hom$-category whose objects are webs with boundary dictated by $\tau$'s and whose morphisms are foams between such webs.

Applying any functor from this category to the category of (graded) vector spaces produces a complex of 
vector spaces, and we can take its homology. If our tangle is a link $\cal{L}$, the complex lies in 
the category of webs whose boundary consist only of $n$-labeled points (with an equal number in the domain 
and codomain), and there is a preferred representable functor to graded vector spaces, given by
$
\displaystyle
\bigoplus_{k \in \Z}
\Hom
\left(
q^k \;\;
\xy
(0,0)*{
\begin{tikzpicture}[scale=.4]
\draw[double] (2,1) to (0,1);
\draw[double] (2,0) to (0,0);
\node[scale=.75] at (1,.75) {$\vdots$};
\node at (2.375,1) {\small$_n$};
\node at (2.375,0) {\small$_n$};
\end{tikzpicture}
};
\endxy
, -
\right)
$. To precisely match the conventions from \cite{KhR} and \cite{Wu}, we must also first pass to the 
complex $\llbracket \cal{L} \rrbracket^{\vee}$, in which we've reversed the quantum and homological 
gradings (and hence also reversed the direction of the morphisms) in $\llbracket \cal{L} \rrbracket$.

\begin{thm}
Let $\mathcal{L}$ be a colored, framed link, then the homology of the complex 
\[
\bigoplus_{k \in \Z}
\Hom
\left(
q^k \;\;
\xy
(0,0)*{
\begin{tikzpicture}[scale=.4]
\draw[double] (2,1) to (0,1);
\draw[double] (2,0) to (0,0);
\node[scale=.75] at (1,.75) {$\vdots$};
\node at (2.375,1) {\small$_n$};
\node at (2.375,0) {\small$_n$};
\end{tikzpicture}
};
\endxy
, \llbracket \cal{L} \rrbracket^{\vee}
\right)
\]
of graded vector spaces is isomorphic (up to shifts) to the colored Khovanov-Rozansky homology of $\cal{L}$.
\end{thm}
\begin{proof}
The crucial observation is the one implicit in \cite[Section 10]{Cautis} and explicitly described in 
\cite[Section 4.2]{LQR1}, that any two link homology theories coming from a ``skew Howe" 
$2$-functor\footnote{We're implicitly using here that for the codomain of such a $2$-functor, any $1$-morphism between 
the highest/lowest weight object and itself is isomorphic to a direct sum of grading shifts of the identity $1$-morphism.
This follows by the work of Cautis-Kamnitzer-Morrison \cite[Theorem 4.4.1]{CKM}.}
must be isomorphic, provided their highest (or lowest) weight objects have identity $1$-morphism 
$\1_{\lambda}$ so that $\Hom(\1_{\lambda},q^{l} \1_\lambda)$ is one dimensional when $l=0$ and 
is zero otherwise. This holds in $\foam[ns]{n}{}^\bullet$, and also for the construction of link homology 
resulting from the $2$-functor $\Gamma_n = \Psi_n \circ \Phi_n$, with braidings given by the Rickard 
complexes.

It suffices to show that the latter link homology is isomorphic to colored Khovanov-Rozansky homology 
(see \cite{KhR} for the original construction and \cite{Wu} and \cite{Yon} for the colored case). 
To this end, assume that $\cal{L}$ is the closure of a braid $\sigma$:
\[
\cal{L} = 
\xy
(0,0)*{
\begin{tikzpicture}[scale=.5,rotate=270]
\draw[very thick] (-.5,-2) rectangle (2.5,2);
\node at (1,0) {\Large $\sigma$};
\draw [very thick,rdirected=.5] (2,2) .. controls (2,3) and (3,3) .. (3,2) -- (3,-2) .. controls (3,-3) and (2,-3) .. (2,-2);
\draw [very thick,rdirected=.5] (0,2) .. controls (0,4) and (5,4) .. (5,2) -- (5,-2) .. controls (5,-4) and (0,-4) .. (0,-2);
\node at (3.875,0) {$\vdots$};
\end{tikzpicture}
};
\endxy
\]
using the Markov theorem. The complex of matrix factorizations assigned to the (unclosed) braid $\sigma$ 
in colored Khovanov-Rozansky homology is easily observed to be isomorphic (up to shifts in homological and quantum degree 
coming from the renormalization in \cite[Definition 12.16]{Wu}) 
to the complex of matrix factorizations $\Psi_n(\llbracket \sigma \rrbracket^{\vee})$. 

It thus suffices to show that the vector spaces assigned to the terms in these complexes by closing 
the braid in both the Khovanov-Rozanksy and skew Howe setting are 
isomorphic\footnote{Here we ignore the additional $\Z/2\Z$ grading on these vector spaces. The isomorphisms we describe do not 
preserve this degree.}. If 
$W$ is a resolution of $\sigma$, i.e. a web appearing in $\llbracket \sigma \rrbracket^{\vee}$, then the colored Khovanov-Rozansky 
construction assigns the vector space
\[
\bigoplus_{k \in \Z}
\Hom_{n\mathbf{hmf}}
\left(
q^k \varnothing,
\xy
(0,0)*{
\begin{tikzpicture}[scale=.3,rotate=270]
\draw[very thick] (-.5,-2) rectangle (2.5,2);
\node at (1,0) {\Large $W$};
\draw [very thick,rdirected=.6] (2,2) .. controls (2,3) and (3,3) .. (3,2) -- (3,-2) .. controls (3,-3) and (2,-3) .. (2,-2);
\draw [very thick, rdirected=.55] (0,2) .. controls (0,4) and (5,4) .. (5,2) -- (5,-2) .. controls (5,-4) and (0,-4) .. (0,-2);
\node at (3.5,0) {$_{_{\vdots}}$};
\end{tikzpicture}
};
\endxy
\right)
\]
where we graphically identify a planar graph\footnote{The latter isn't a left-directed web, i.e. is not a $1$-morphism in our foam categories.} 
and its associated matrix factorization in this formula. By contrast, the skew Howe 
approach assigns the vector space
\[
\bigoplus_{k \in \Z}
\Hom_{n\mathbf{hmf}}
\left(
q^k
\xy
(0,0)*{
\begin{tikzpicture}[scale=.4]
\draw[double] (2,1) to (0,1);
\draw[double] (2,0) to (0,0);
\node[scale=.75] at (1,.75) {$\vdots$};
\node at (2.375,1) {\small$_n$};
\node at (2.375,0) {\small$_n$};
\end{tikzpicture}
};
\endxy,
\xy
(0,0)*{
\begin{tikzpicture}[scale=.3,rotate=270]
\draw[very thick] (-.5,-2) rectangle (2.5,2);
\node at (1,0) {\Large $W$};
\draw [very thick,directed=.6] (2,2) .. controls (2,3) and (3,3) .. (3,2) -- (3,-2) .. controls (3,-3) and (2,-3) .. (2,-2);
\draw [very thick, directed=.55] (0,2) .. controls (0,4) and (5,4) .. (5,2) -- (5,-2) .. controls (5,-4) and (0,-4) .. (0,-2);
\node at (2.25,4) {$_{_{\vdots}}$};
\node at (2.25,-4) {$_{_{\vdots}}$};
\node at (3.5,0) {$_{_{\vdots}}$};
\draw[double] (2.5,2.75) to (2.5,3.5);
\draw[double] (2.5,-2.75) to (2.5,-3.5);
\draw[double] (1,3.25) to (1,5);
\draw[double] (4,3.25) to (4,5);
\draw[double] (1,-3.25) to (1,-5);
\draw[double] (4,-3.25) to (4,-5);
\end{tikzpicture}
};
\endxy
\right)
\]
which (up to shifts) is isomorphic to that above, using duality in the category of matrix factorizations and
the fact that the matrix factorizations corresponding to the planar graphs:
\[
\xy
(0,0)*{
\begin{tikzpicture}[scale=.4]
\draw[very thick, directed=.55] (0,-2) to (0,2);
\node at (.75,0) {\small $a$};
\end{tikzpicture}
};
\endxy
\quad \text{ and } \quad
\xy
(0,0)*{
\begin{tikzpicture}[scale=.4]
\draw[very thick, directed=.55] (0,-2) to (0,-1);
\draw[very thick, directed=.55] (0,1) to [out=315,in=45] (0,-1);
\draw[double] (0,1) to [out=225,in=135] (0,-1);
\draw[very thick, directed=.55] (0,1) to (0,2);
\node at (1.5,0) {\small $_{n-a}$};
\end{tikzpicture}
};
\endxy
\]
are isomorphic in $n\mathbf{hmf}$, and equation \eqref{n_iso}.
\end{proof}

\subsection{Comparison with previous work}\label{subsec:compMSV} 

\subsubsection{Work of Mackaay-Stosic-Vaz}
In \cite{MSV}, Mackaay, Stosic and Vaz give a foam-based construction of (uncolored) Khovanov-Rozansky link homology. 
Their framework utilizes foams containing only facets colored $1$, $2$ and $3$, which is enough for proving all Reidemeister moves in the uncolored case; 
however, they don't give a presentation by generators and relations, instead using the Kapustin-Li formula to deduce a set of useful relations on foams.
Since both constructions allow for a foam description of uncolored Khovanov-Rozansky link homology, one might hope that these constructions agree. A direct comparison between these foam categories is uneasy, 
mostly because the foams they use are not restricted to be leftward-directed, that is, are not necessarily in the image of our $2$-functor. 
We can compensate for this issue using $n$-labeled facets, giving an in-spirit identification of a rightward $k$-labeled strand with a leftward $(n-k)$-labeled one.

One then needs a process for canonically building a ladder foam out of any (MSV) foam, categorifying the process we use in Section \ref{subsec:tangleInv}. 
Such a process is easy to present, up to sign, but obtaining coherent signs, and moreover proving that the final result does not depend on the initial presentation of the foam, is more complicated.
Up to this putative process, we can at least give a correspondence between the Mackaay-Stosic-Vaz relations and their analogs in our $2$-category, as follows:
\begin{itemize}
\item the \emph{dot conversion} relations agree with ours: this has been discussed in Section \ref{subsec:detailsNfacets};
\item the \emph{dot migration} relations are special cases of equation \eqref{dotmigrel} involving only $1$, $2$ and $3$-labeled facets;
\item the \emph{cutting neck} relations are the first ones that are not of ladder type, but can easily be obtained by enhancing with an $n$-labeled sheet as a consequence of Relation \eqref{R2KLR1rel}, with $a+b=n$ and $b=c$;
\item the \emph{sphere relations} are consequences of our blister relations \eqref{blisterrel};
\item the two \emph{$\Theta$-foam} relations of \cite[Proposition 6.2]{MSV} reduce to sphere relations once one applies our KLR relation \eqref{R2KLR2rel};
\item our Matveev-Piergalini relation \ref{MVrel} is simply an extension of the one from \cite{MSV} (which exhausts their Proposition 6.2);
\item the other \emph{$\Theta$-foam} relation from \cite[Lemma 6.3]{MSV} can be recovered by applying our blister relation to compress the $1$- and $2$-labeled facets together, 
and then we are back to the \emph{sphere relation} for a $3$-labeled sphere;
\item the first \emph{digon removal} relation of \cite[Proposition 6.9]{MSV} is of ladder type (and the convention in \cite{MSV} is to use the left-hand rule for ordering along seams), 
but the sign differs from our equivalent relation \eqref{nH3rel}: a signed version of this identification would need an extra rescaling; 
\item Relation $(DR_{3_1})$ can be obtained from relation \eqref{nH3rel};
\item  the other \emph{digon removal} relations are not of ladder types and more easily follow from the isomorphisms underlying Remark \ref{rem:SquareFormulas};
\item the \emph{first square removal} relation is also not of ladder type and requires the use of $n$-labeled facets to again be identified with the isomorphisms underlying Remark \ref{rem:SquareFormulas}.
\end{itemize}

\subsubsection{Work of Cautis}
In \cite{Cautis}, Cautis shows how to construct a link invariant from any categorical $2$-representation of $\slm$ categorifying the skew Howe representation 
$\textstyle \bigwedge_q^N(\C_q^n \otimes \C_q^m)$. He then goes on to show that certain $2$-categories arising from the study of the affine Grassmannian and 
Nakajima quiver varieties (in type $A$) give examples of such categorical $2$-representations. As the definition of a categorical $2$-representation is weaker than 
that given in Definition \ref{def:2rep}, he goes on to conjecture that these $2$-categories also give such $2$-representations. Any such ``skew Howe" $2$-functor 
will factor through $\foam{n}{m}$, hence $\sln$ foams should act on these $2$-categories as well.

\subsubsection{Work of Chatav and Webster}
In \cite{Chatav}, Chatav shows that Bar-Natan's $\slnn{2}$ foam $2$-category gives a graphical description of $2$-intertwiners between Webster's categorifications of 
$\slnn{2}$ representations \cite{Webster}. One can hope to extend Chatav's work to show that our foams give $2$-intertwiners between tensor products of 
$\sln$ representations; however, this is far from a straightforward generalization. In particular, one must define the functors corresponding to the first two
generating webs in equation \eqref{webgen} and make sense of the role of $n$-labeled edges in Webster's construction.
Nevertheless, since our work can be viewed as a categorification of the category of left-directed $\sln$ webs, the relation to Webster's work seems likely.



%

%
\end{document}